\documentclass[11pt]{amsart}

\usepackage[foot]{amsaddr}
\usepackage[english]{babel}
\usepackage[utf8]{inputenc}
\usepackage[T1]{fontenc}
\usepackage[textwidth=17.5cm,textheight=22.275cm,
	height=24.275cm,width=18cm]{geometry}
\usepackage{amsmath,amssymb,mathtools,mathrsfs}
\usepackage[dvipsnames,svgnames,table,x11names]{xcolor}
\usepackage{newtxtext}
\usepackage{newtxmath}
\usepackage[bb=boondox,cal=boondoxo,scr=boondoxo]{mathalfa}
\mathtoolsset{showonlyrefs}
\numberwithin{equation}{section}
\NewEnvironmentCopy{standardbmatrix}{bmatrix}
\usepackage[renew-dots,renew-matrix]{nicematrix}
\usepackage{enumitem}
\usepackage{microtype}
\usepackage{tensor}
\usepackage{pgfplots}
\pgfplotsset{compat=1.18}
\usepackage{tikz}
\usetikzlibrary{calc,shadows,shapes.callouts,shapes.geometric,shapes.misc,positioning,patterns,decorations.pathmorphing,decorations.markings,decorations.fractals,decorations.pathreplacing,shadings,fadings}
\usepackage{longtable}
\usepackage{fancybox}
\usepackage{tcolorbox}
\usepackage{rotating,multirow,pdflscape,float}
\usepackage{arydshln,mathdots}
\usepackage{bigints}
\usepackage[hyperfootnotes]{hyperref}
\usepackage{orcidlink}
\usepackage{csquotes}

\usepackage[
backend=biber,
style=numeric,
sorting=nyt,
giveninits=true,
maxnames=99,
doi=true,
url=false,
isbn=false
]{biblatex}

\DeclareFieldFormat[article]{title}{\mkbibemph{#1}}

\DeclareFieldFormat[article]{journaltitle}{#1}

\DeclareFieldFormat[article]{volume}{\textbf{#1}}

\DeclareFieldFormat{pages}{#1}

\renewbibmacro{in:}{%
	\ifentrytype{article}
	{}
	{\printtext{\bibstring{in}\intitlepunct}}%
}

\renewbibmacro*{journal+issuetitle}{%
	\usebibmacro{journal}%
	\setunit*{\addspace}%
	\printfield{volume}%
	\iffieldundef{number}
	{}
	{\setunit{\addcomma\space}%
	 \printtext{no.\addnbspace}%
	 \printfield{number}}%
	\setunit{\addspace}%
	\printtext[parens]{\printfield{year}}%
	\setunit{\addspace}%
	\printfield{pages}%
	\newunit
}

\renewbibmacro*{note+pages}{%
	\ifentrytype{article}
	{}
	{\printfield{note}%
		\setunit{\bibpagespunct}%
		\printfield{pages}}%
}
\hypersetup{
	colorlinks=true,
	linkcolor=NavyBlue,
	urlcolor=RoyalPurple,
	citecolor=OliveGreen,
	pdftitle={Positive Bidiagonal Factorizations:
		Combinatorial Structure, Spectral Theory, and Matrix Continued Fractions},
	pdfauthor={Manuel Ma\~nas}
}
\allowdisplaybreaks
\newcommand{\dx}{\,\mathrm{d}x}

\newcommand{\e}{\mathrm{e}}

\newcommand{\one}{\mathbf{1}}
\newcommand{\pFq}[5]{\;{}_{#1}F_{#2}\left(\begin{matrix}#3\\#4\end{matrix};#5\right)}

\newcommand{\Nzero}{\mathbb N_0}

\renewcommand{\top}{\mathsf T}

\theoremstyle{plain}
\newtheorem{theorem}{Theorem}[section]
\newtheorem{proposition}[theorem]{Proposition}
\newtheorem{lemma}[theorem]{Lemma}
\newtheorem{corollary}[theorem]{Corollary}

\newtheoremstyle{definitionstyle}
{6pt}{12pt}{\normalfont}{}{\bfseries}{.}{0.5em}{}
\theoremstyle{definitionstyle}
\newtheorem{definition}[theorem]{Definition}
\newtheorem{example}[theorem]{Example}

\newtheoremstyle{remarkstyle}
{6pt}{12pt}{\normalfont}{}{\itshape}{.}{0.5em}{}
\theoremstyle{remarkstyle}
\newtheorem{remark}[theorem]{Remark}

\begin{document}

\title[Positive bidiagonal factorizations]
{Positive Bidiagonal Factorizations:\\
Combinatorial Structure, Spectral Theory, and Matrix Continued Fractions}
\author[M. Ma\~nas]{Manuel Ma\~nas\,\orcidlink{0000-0003-3764-5737}}
\address{Department of Theoretical Physics, Faculty of Physical Sciences,
	Complutense University of Madrid, 28040 Madrid, Spain}
\email{manuel.manas@ucm.es}

\subjclass[2020]{Primary 40A15, 41A21; Secondary 05A05, 68R15, 37K10,
	44A60, 65D15, 30B70, 65F30, 15A23, 15B48, 47A10, 42C05}
\keywords{positive bidiagonal factorization; Christoffel word; planar network;
	cyclic Darboux transformation; matrix continued fraction; mixed-type multiple
	orthogonality; integrable lattice; B\"acklund transformation; tau function;
	Stieltjes moment sequence; Weyl matrix; Pad\'e approximation}

\date{August 19, 2026}\enlargethispage{2cm}

\begin{abstract}
	Presenting a semi-infinite banded matrix through a prescribed positive
	bidiagonal factorization (PBF)
	\[
		T=L_1\cdots L_pU_q\cdots U_1
	\]
	reveals spectral and approximation-theoretic structure that is not visible
	from the entries of \(T\) alone.  This structure is developed for arbitrary
	\(p\) and \(q\), without assuming boundedness except for analytic resolvents.

	First, every cyclic permutation of the prescribed factors again admits a
	normalized PBF with all lower factors preceding all upper factors.  Thus the
	cyclic Darboux orbit remains within the scope of bounded or unbounded
	mixed-type Favard theory.  With suitable one-sided maximal-degree data at
	consecutive cyclic products, one factor transfer determines a degree-one
	matrix Christoffel transformation whose determinant
	is a nonzero multiple of \(x\); a complete circuit gives multiplication by
	\(x\), up to a constant triangular change of the terminal normalization.
	The fixed simultaneous row and column moment shifts form a finite rectangle
	when the two boundary chains are weakly normal, a condition that is automatic
	when the scalar moment matrix is totally positive.

	Second, there is an explicit interval of cyclic products for which the
	directed positive adjacent refactorizations produce a distinguished
	factorization of minimum height span.  In general this factorization is
	interlaced and occurs before the terminal lower--upper PBF in the sequence
	of exchanges.  Combinatorially, the minimizing factor pattern is a power of
	a lower Christoffel word, and its prefix-dominance property determines the
	maximal family of cyclic products that can reach a minimum-height positive
	refactorization.  The associated
	semi-infinite matrices have only two
	nonzero block displacements and \(\gcd(p,q)\) coordinates at each height;
	in the coprime case only the \(q\)-th superdiagonal and the \(p\)-th
	subdiagonal remain.  Their projected resolvents recover the corresponding
	cyclic products, and their finite nonzero spectra lie on a
	\((p+q)\)-ray star.  A planar-network expansion proves coefficientwise
	nonnegativity of the radial Hankel minors and hence yields positive radial
	Stieltjes moment sequences.  In the compactly supported radial case, in
	which the representing measure is determinate, the Fourier components admit
	finite complex representing measures supported on the corresponding compact
	star precisely when,
	entry by entry and in every nonzero mode, the radial measure has no atom at
	the origin and satisfies the negative fractional-moment condition.
	In the coprime case the two-diagonal recurrences also carry a local
	arithmetic-support Lax hierarchy.  A radial deformation of the star moments
	gives an exact determinant formula for every evolving factor coefficient.
	Their \((p+q)\)-th powers synchronize the Toda flows of all cyclic products
	of the refactorized factors, while the Christoffel transfers give compatible
	discrete B\"acklund transformations.  On the real first-time axis, the
	Pi\~neiro specialization yields global positive coefficientwise solutions
	expressed by determinants of Kummer functions;
	the intersection of the Jacobi-like positive-factorization chamber with
	the requisite AT region gives analogous global positive coefficientwise
	solutions expressed by \({}_qF_q\) determinants.

	Third, retaining the prescribed order and one coordinate level after each
	factor produces a sparse matrix.  Grouping consecutive height levels makes
	it block tridiagonal and yields a factor-resolved matrix continued fraction.
	Its finite height-block truncations have an explicit Pad\'e-type order of
	contact and denominator factorization, admit backward evaluation, and carry
	a computable a posteriori error bound.  For bounded nonnegative factors these
	convergents are monotone, their denominators are nonsingular \(M\)-matrices,
	and their limit is the least nonnegative solution; stochastic normalization
	gives convergence in the unit disk.  A second grouping gives a
	\(q\)-state branched recursion that reduces, for \(q=1\), to the branched
	Stieltjes fraction of ordinary multiple orthogonality.  In the bounded setting,
	under the Favard moment identification, a rectangular projection is the Weyl
	matrix of the mixed-type system, so all of its Stieltjes transforms are
	approximated simultaneously.  For unbounded factorizations the coefficientwise
	formal approximation and finite algorithm remain valid.  Analytic
	finite-section convergence, when available, is to the projected resolvent of
	a specified closed realization under the core and uniform-stability
	hypotheses, and agrees with a measure-defined Weyl matrix only under an
	additional resolvent--Cauchy-transform identification.  Pi\~neiro and
	Jacobi-like systems provide explicit
	spectral, approximation, and integrable applications.
\end{abstract}
\maketitle

\tableofcontents

\section{Introduction}
\label{sec:introduction}

Continued fractions provide a fundamental link between orthogonal
polynomials, spectral theory, and rational approximation.  For a Jacobi
matrix, the three-term recurrence produces a Stieltjes continued fraction for the
transform of the orthogonality measure.  Its convergents are Pad\'e
approximants, their denominators are the orthogonal polynomials, and their
poles are the eigenvalues of finite Jacobi truncations.  A continued fraction
therefore converts local recurrence coefficients into information about the
measure, the spectrum, polynomial zeros, and rational approximation.

The step-line recurrence in ordinary multiple orthogonality is lower
Hessenberg rather than tridiagonal.  Its \(p\)-branched continued fractions
encode the additional subdiagonals and simultaneously approximate a vector
of Stieltjes transforms.  The production/output formalism of
Deutsch--Ferrari--Rinaldi gives a path interpretation of lower-Hessenberg
recurrences \cite{DeutschFerrariRinaldi2005}; its relation with multiple
orthogonality and branched Jacobi fractions is developed in
\cite{PetreolleSokalZhu2023,Sokal2024}.  When the recurrence has \(p\) lower
bidiagonal factors and one upper factor, resolving those factors replaces
direct \(p\)-\L{}ukasiewicz paths by weighted \(p\)-Dyck paths and yields the
\(p\)-branched Stieltjes fraction
\cite{Lima2023,BranquinhoDiazFoulquieLimaManas2026JAT,
Flajolet1980,Viennot1983,BanderierFlajolet2002}.

A bidiagonal factorization, however, carries more than continued-fraction
data.  Its ordered factors define a chain of Darboux--Christoffel
transformations, while their positivity controls moment and spectral
properties that are not visible from the entries of their product alone.
The spectral and approximation-theoretic roles of the factors are developed
here simultaneously.

Mixed-type multiple orthogonality has several lower and upper recurrence
diagonals, and its Weyl object is a rectangular matrix.  This leads to the
three questions addressed here.  Consider a semi-infinite recurrence matrix
with a prescribed bidiagonal factorization
\[
 T=L_1\cdots L_pU_q\cdots U_1,
\]
where the \(L_a\) are lower bidiagonal and the \(U_b\) are upper
bidiagonal.  When their nontrivial entries are positive, this is called a
positive bidiagonal factorization, or PBF.  The questions are:
Does positive bidiagonal factorization survive when the factors are moved cyclically?  Can the
factor data be reorganized into a sparse recurrence carrying its own moment
and spectral structure?  And can the prescribed factor order be converted
into a continued fraction for the resolvent that retains the individual
factor coefficients, rather than only their sums of products in \(T\)?  All
three questions have affirmative answers for arbitrary \(p\) and \(q\).  The
main results are grouped accordingly.

The decisive organizing mechanism is combinatorial.  Words in the two factor
types are interpreted as rational lattice paths.  Christoffel words identify
the positive refactorizations of minimum height span, while prefix dominance
determines which cyclic products can reach them
\cite{BerstelLauveReutenauerSaliola2009}.  A second path construction, based
on planar networks and nonintersecting paths, proves the Hankel positivity
underlying the radial Stieltjes measures
\cite{Lindstrom1973,GesselViennot1985}.  Thus combinatorics first constructs
the sparse recurrences and then supplies their moment positivity.

\begin{enumerate}[label=\textnormal{(\arabic*)},leftmargin=*]
	
	\item \emph{Cyclic permutation preserves the positive factorization and
		transports the spectral data.}
	A local refactorization of an upper--lower pair shows that every cyclic
	permutation of the prescribed bidiagonal factors again admits a PBF
	(Lemma~\ref{lem:positive-adjacent-refactorization} and
	Theorem~\ref{thm:cyclic-PBF-closure}).  Thus one PBF generates a complete
	cyclic family of positive banded matrices.  This is an algebraic result and
	does not require boundedness.  Under the corresponding bounded or unbounded
	Favard hypotheses, every member of this family has positive mixed-type
	spectral data
	\cite{BranquinhoFoulquieManas2023Spectral,
		BranquinhoFoulquieManas2024BTP,
		BranquinhoFoulquieManas2026Unbounded}.
	
	When suitable one-sided maximal-degree data are fixed at two consecutive
	members of the cyclic family, moving one bidiagonal factor induces a
	matrix Christoffel transformation of degree one.  Its determinant is a
	nonzero multiple of \(x\), and one complete circuit through all factors
	multiplies the matrix of measures by \(x\), up to a constant triangular
	change of initial conditions
	(Theorem~\ref{thm:PBF-adapted-cyclic-Christoffel} and
	Corollary~\ref{cor:PBF-adapted-complete-circuit}).  Fixing instead the
	standard row and column Christoffel transformations produces a finite
	two-sided family of transformed moment systems.  Its existence requires
	weak normality of the two boundary chains; this condition is automatic
	when the scalar step-line moment matrix is totally positive
	(Theorem~\ref{thm:balanced-Christoffel-grid-path}).
	
	\item \emph{Selected cyclic products yield sparse recurrences with
		spectral data on a star.}
	For an explicitly determined interval of cyclic products, a finite sequence
	of exchanges reaches a distinguished positive pattern of minimum height
	span, generally before the terminal lower--upper PBF is reached
	(Theorem~\ref{thm:positive-balanced-refactorization}).  This
	factorization is retained because it produces, directly from its bidiagonal
	coefficients, a
	sparse semi-infinite matrix \(\mathscr H_{(a)}\) with only two nonzero
	block diagonals.  If \(d=\gcd(p,q)\), each height contains \(d\)
	coordinates and the two block displacements are
	\[
	+\frac{q}{d}
	\qquad\text{and}\qquad
	-\frac{p}{d}.
	\]
	When \(p\) and \(q\) are coprime, \(\mathscr H_{(a)}\) is scalar and its
	only nonzero entries are
	\[
	(\mathscr H_{(a)})_{h,h+q}=1,
	\qquad
	(\mathscr H_{(a)})_{h,h-p}=c_h^{(a)}>0.
	\]
	
	The matrices \(\mathscr H_{(a)}\) are not merely auxiliary
	linearizations.  They are recurrence matrices of explicitly constructed
	moment-functional systems; their radial moments are Stieltjes moment
	sequences, and the nonzero spectra of their finite truncations lie on a
	\((p+q)\)-ray star.  For fixed compactly supported radial representing
	measures, a separate entrywise criterion determines when the corresponding
	Fourier components admit finite complex representing measures supported on
	the corresponding compact star.
	Furthermore, a projected resolvent of \(\mathscr H_{(a)}\) recovers the
	resolvent of the corresponding cyclic product \(T_{(a)}\)
	(Corollary~\ref{cor:admissible-cyclic-star-CF-family}).  The sparse
	recurrence therefore retains the spectral information of \(T_{(a)}\)
	while exposing an additional star structure.
	When \(p\) and \(q\) are coprime, it also defines an arithmetic-support
	difference operator.  Section~\ref{sec:integrable-hierarchies-PBF} proves
	that this operator carries an infinite hierarchy of local commuting Lax
	flows.  The star moment matrix reconstructs the evolving coefficients by a
	four-minor \(\tau\)-quotient, so the PBF supplies an exact solution rather
	than only initial data.  The \((p+q)\)-th power evolves all cyclic products
	of the refactorized factors by synchronized Toda equations, and the
	Christoffel transfers commute with these flows as discrete B\"acklund
	transformations.  On the real first-time axis, the Pi\~neiro case gives
	global positive coefficientwise solutions expressed by determinants of
	Kummer functions;
	the intersection of the Jacobi-like positive-factorization chamber with the
	requisite AT range gives analogous global positive coefficientwise
	hypergeometric determinant solutions.
	
\item \emph{The original bidiagonal factorization yields a matrix continued
	fraction for the Weyl matrix.}
The construction in the preceding item replaces the original bidiagonal
factors by positively refactorized ones.  A second construction starts
directly from
\(T=L_1\cdots L_pU_q\cdots U_1\) and leaves these factors unchanged.
It resolves one application of \(T\) into its \(p+q\) successive
bidiagonal steps, so that every coefficient of every factor remains
visible.  Introducing one intermediate coordinate after each factor,
ordering these coordinates by an integer height, and grouping consecutive
heights turns the resulting sparse matrix into a block-tridiagonal matrix.
Its Schur complements give the matrix continued fraction of
Theorem~\ref{thm:matrix-CF-height-blocks}.  A different grouping of the
same coordinates gives the \(q\)-state recursion of
Theorem~\ref{thm:q-channel-matrix-BCF}; when \(q=1\), this recursion
reduces to the classical \(p\)-branched Stieltjes fraction.

Finite height-block truncations give rational approximants for which the
Pad\'e-type order of contact, denominator factorization, backward evaluation
algorithm, and a posteriori error bound are obtained explicitly.  For bounded
nonnegative factors, these convergents are entrywise
monotone, their Schur-complement denominators are nonsingular
\(M\)-matrices, and their limit is the least nonnegative solution of the
	continued-fraction equation.  The coefficientwise approximation and the
	finite algorithm remain valid for unbounded factorizations.
	Operator-resolvent convergence then requires a specified closed realization,
	the stated core hypothesis, and uniform finite-section stability.

	In the bounded setting, under the moment identification, a rectangular
	projection of the resolvent of this sparse matrix recovers the Weyl matrix
	of the mixed-type system.  The finite continued fractions therefore
	approximate simultaneously all the Stieltjes transforms in the underlying
	matrix of measures.  In the unbounded setting, convergence to the projected
	resolvent and identification of that resolvent with the Cauchy transform of
	a chosen Favard measure are separate assertions.  The two associated block
	output matrices recover the coefficient arrays of both mixed-type
	polynomial families
(Theorems~\ref{thm:mixed-block-output-matrices} and
\ref{thm:general-Weyl-CF}).
\end{enumerate}

We next place these results in the context of the state of the art.  The
relevant literature concerns sparse recurrences and spectral problems on
star-like sets, vector and matrix continued fractions for banded operators,
matrix-analytic methods for block recursions, and mixed-type multiple
orthogonality.  We review these four directions in turn and indicate precisely
where the present constructions differ from their closest precedents.

Two-diagonal lower-Hessenberg matrices, corresponding here to \(q=1\), occur
in multiple orthogonality on star-like sets
\cite{DelvauxLopezGarcia2013,Lima2026Symmetric}; rotational spectral patterns
for finite matrices with two separated diagonals were studied in
\cite{McMillen2008}.  The family \(\mathscr H_{(a)}\) constructed above is
their two-sided mixed-type extension.  For a Christoffel factorization, the general
factorization theorem of \cite{BranquinhoFoulquieManas2026Factorization}
identifies its cyclic products with successive Christoffel shifts.  The
present positive refactorization theorem proves, for bounded and unbounded
factors alike, that the complete cyclic orbit remains in the PBF class.  The
spectral consequences use two
distinct parts of the mixed-type theory.  The Favard theorem, established
first for bounded matrices and subsequently extended to unbounded matrices
under a coherence hypothesis for the finite PBFs, supplies positive measures
\cite{BranquinhoFoulquieManas2023Spectral,
BranquinhoFoulquieManas2026Unbounded}.  In the bounded case these measures
are compactly supported.  In the genuinely unbounded case, the
	coherence hypothesis for the shifted finite PBFs in the general theorem is
	automatic here: every leading principal submatrix inherits the leading
	submatrices of the same semi-infinite factors, with zero shift.  Hence the
	boundary matrices occurring in the coherence condition are independent of
	the truncation size once their fixed leading blocks are present.
Separately, the BTP criterion gives strong normality when
\(\nu^{-\top}\) is upper \(\Delta\mathrm{TP}\) and \(\xi^{-1}\) is lower
\(\Delta\mathrm{TP}\)
\cite{BranquinhoFoulquieManas2024BTP}.  The cited paper uses the word
\emph{normality} for this stronger maximal-degree property; it is called
\emph{strong normality} here.  For general bidiagonal factorizations the
cyclic resolvent identities remain valid.  Positivity is the additional
hypothesis used to obtain positive spectral data.

Vector and matrix continued fractions for banded operators have a substantial
literature.  Relevant precedents include the vector continued fractions and
Favard theory of Kaliaguine \cite{Kaliaguine1995}, the genetic-sum
representation of Aptekarev--Kaliaguine--Van Iseghem
\cite{AptekarevKaliaguineVanIseghem2000}, the matrix continued fractions of
Sorokin--Van Iseghem \cite{SorokinVanIseghem1999}, matrix Schur algorithms
\cite{ChenHu1998}, and noncommutative convergence questions
\cite{RaissouliKacha2000}.  L\'opez-Garc\'\i a and Prokhorov related the
Jacobi--Perron algorithm for Hessenberg resolvents to weighted
\L{}ukasiewicz paths and also treated a partial \(p\)-Dyck subfamily
\cite{LopezGarciaProkhorov2024}.

The closest operator-theoretic precedent is
\cite{KimLopezGarciaProkhorov2025}, where \(q\times p\) generating matrices
for two-sided paths and resolvent blocks of banded difference operators are
represented by matrix continued fractions.  There the path weights are the
entries of the banded operator and the rational approximants use its
principal truncations.  Here each recurrence step is first resolved into the
prescribed bidiagonal steps, and the convergents truncate the height-ordered
cyclic matrix rather than \(T\).  Consequently the present fraction keeps
the factor coefficients, the cyclic Darboux products, and factor positivity
visible.  The two approaches share a Schur-complement mechanism but start
from different data and retain different structure.

After the height grouping, the resulting block-tridiagonal matrices also
resemble the transition matrices of quasi-birth-and-death processes and the
matrices treated by matrix-analytic methods
\cite{Neuts1981,LatoucheRamaswami1999}.  In that theory, the block transition
probabilities are given data and the central minimal nonnegative object is
the \(G\)-matrix of the Markov process.  Here the blocks are constructed from
a prescribed bidiagonal factorization, retain its individual coefficients,
and the least nonnegative solution is a factor-resolved tail resolvent.  Thus
the parallel extends to the matrix recursions and their minimal nonnegative
solutions, whereas the input data and the interpreted object are different.

The mixed-type framework has antecedents in the matrix-functional setting of
Sorokin--Van Iseghem \cite{SorokinVanIseghem1997}.  Daems and Kuijlaars later
introduced the terminology of mixed-type multiple orthogonality and related
its Christoffel--Darboux kernel to nonintersecting Brownian motions
\cite{DaemsKuijlaars2007}.  The Gauss--Borel formulation of the two
mixed-type polynomial families and their multi-component 2D Toda dynamics
was developed in \cite{AlvarezFidalgoManas2011}.  Under the Favard hypotheses used below, a
rectangular projection of \((I-tT)^{-1}\) is generated by a \(q\times p\)
matrix of measures.  The continued fraction therefore gives simultaneous
rational approximants to all entries of its Weyl matrix.  Two normalized
block output matrices recover the coefficient arrays of the right and left
mixed-type polynomial systems from powers of \(T\) and \(T^{\top}\).  The
cyclic matrix records the intermediate factor positions and serves a
different purpose.

The Pi\~neiro and Jacobi-like systems illustrate what the general
construction adds once a PBF is available.  Their factor and
special-function formulas are imported as inputs from
\cite{BranquinhoFoulquieManas2026Pineiro,ManasMarkov2026,
Manas2026Hypergeometric}; the Pi\~neiro factors and their positivity chamber
also have the self-contained alternative derivation in
\cite[Section~7 and Appendix~A]{ManasMarkov2026}.  The
general construction is independent of those formulas.  The Pi\~neiro
section contains one complete rational example, including its
matrix measure, banded product, bidiagonal factors, finite fraction,
stochastic normalization, and Weyl approximation.  The Jacobi-like section
gives the general beta-convolution measure, its hypergeometric Weyl matrix,
and their rotational lift to the star.  In that family the upper
Christoffel transforms can leave the initial special-function family,
although cyclic closure preserves their PBFs.

The paper is organized as follows.
Section~\ref{sec:positive-refactorization-sparse-matrices} proves positive adjacent refactorization,
shows that every cyclic product has a terminal lower--upper PBF, and
identifies those cyclic products for which a finite initial sequence of the
same exchanges reaches a minimum-height factorization.  It defines directly
the sparse matrix for the prescribed order and the family
\(\mathscr H_{(a)}\) obtained from these minimum-height factorizations.
It then introduces the auxiliary cyclic matrix and records its Darboux
resolvent blocks.  Section~\ref{sec:height-ordering-block-CF} proves the two
direct placement rules, decomposes the polynomial solutions of
\(\mathscr H_{(a)}\) into
cyclic-product recurrences, and derives the block continued fraction together
with its approximation and convergence properties.
Section~\ref{sec:branched-q-state-CF} gives the branched and
\(q\)-state recursions for the original two-run order.
Section~\ref{sec:general-Weyl-bridge}
identifies the rectangular Weyl matrix inside the cyclic resolvent and
develops its finite rational approximation in the bounded and unbounded
settings.  Section~\ref{sec:balanced-star-Christoffel} constructs the
full star moment systems associated with \(\mathscr H_{(a)}\), proves
coefficientwise nonnegativity of the radial Hankel minors and positive
Stieltjes representations,
establishes scalar weak normality in the coprime case and a uniform block
Gauss--Borel realization in general, proves finite ray localization, and
identifies the finite rectangle of simultaneous row and column moment shifts.
The condensed Pi\~neiro and Jacobi-like applications occupy
Sections~\ref{sec:case-study-mixed-Pineiro} and~\ref{sec:Jacobi-like-CF}.
Section~\ref{sec:integrable-hierarchies-PBF} identifies the integrable
hierarchies selected by the coprime sparse recurrences, reconstructs their
solutions from the star moments, proves compatibility with the discrete
Christoffel transformations, and gives explicit hypergeometric determinant
solutions for the Pi\~neiro and Jacobi-like specializations.

\section{Positive refactorization and the two sparse matrices}
\label{sec:positive-refactorization-sparse-matrices}

This section constructs the two sparse matrices used in the rest of the
article.  The first retains the prescribed factors and all their coefficients.
The second, denoted by \(\mathscr H_{(a)}\), is obtained from a
minimum-height factorization of each selected cyclic product
\(T_{(a)}\), \(a\in\mathcal A_{\mathrm{cyc}}\).  The distinction between
this factorization and the terminal lower--upper PBF of the same product is
made explicit in Subsection~\ref{subsec:balanced-Christoffel-word}.

The section is organized as follows. Subsection~\ref{subsec:positive-adjacent-refactorization} proves that a
positive upper--lower pair can be interchanged while preserving normalized
bidiagonal form and positivity.  Subsection~\ref{subsec:balanced-Christoffel-word}
first identifies the selected cyclic products and then uses these exchanges
to construct their minimum-height factorizations.
Subsection~\ref{subsec:direct-matrices-from-PBF} then defines the two sparse
matrices directly.  Only after these constructions are complete does
Subsection~\ref{subsec:cyclic-linearization} introduce the cyclic matrix.
That matrix will be used to prove the identities connecting the sparse matrices with \(T\),
the cyclic products in the displayed range, and their resolvents.  The two
sparse matrices themselves are assembled directly from their bidiagonal
coefficients.

Before beginning the construction, we fix two terminological conventions
used throughout the article.  Our terminology for total positivity follows
Fallat and Johnson
\cite{FallatJohnson2011}.  A real matrix is \emph{totally nonnegative}
(TN) if all its minors are nonnegative, and \emph{totally positive} (TP) if
all its minors are positive.  For classical background, see also
\cite{Karlin1968,Pinkus2010}.  Pinkus uses TP and STP for the classes called
TN and TP here, respectively; statements cited from that monograph are
expressed in our TN--TP terminology.  A nontrivial triangular matrix cannot
be TP in this sense.  Following Fallat and Johnson, we write
\(\Delta\mathrm{TP}\) for a triangular TN matrix whose minors are positive
whenever they are not forced to vanish solely by triangularity.  Throughout
the article, \emph{triangular totally positive} means
\(\Delta\mathrm{TP}\).

\begin{definition}[Weak and strong normality along a step-line]
	\label{def:weak-strong-step-line-normality}
	A position of a prescribed mixed-type step-line is called
	\emph{weakly normal} if each of its right and left mixed-type
	orthogonality problems has a one-dimensional solution space;
	equivalently, each polynomial vector exists and is unique up to
	multiplication by a nonzero constant.  It is called
	\emph{strongly normal} if, in addition, every component whose prescribed
	polynomial space is nontrivial attains its maximal allowed degree.  A
	step-line is weakly, respectively strongly, normal through order \(N\)
	when this property holds at every one of its first \(N\) positions.
\end{definition}

\subsection{Positive adjacent refactorization}
\label{subsec:positive-adjacent-refactorization}

Let \(p,q\in\mathbb N\), and put
\[
M\coloneq p+q.
\]
The pair \((p,q)\) remains fixed throughout the general construction.  Once
defined, \(\mathscr C_{p,q}\) and \(\mathscr H_{p,q}\) are abbreviated to
\(\mathscr C\) and \(\mathscr H\); their subscripts are restored for special
values of \((p,q)\) and for the explicit families.  All matrices below act
algebraically on
\(\mathcal E\coloneq\mathbb C^{\mathbb N_0}\), the vector space of all complex
sequences \(x=(x_n)_{n\ge0}\), without a support condition.  For each
\(n\in\mathbb N_0\), let \(e_n\) denote the sequence having a \(1\) in
coordinate \(n\) and \(0\) in every other coordinate.  Write
\(\one=(1,1,\ldots)\in\mathcal E\) for the constant sequence of ones and,
for \(r\in\mathbb N\),
\(\one_r=(1,\ldots,1)\in\mathbb C^r\).  Every row of a
bidiagonal factor has at most two nonzero entries, and every finite product of
the factors is banded.  Hence each coordinate of the image of any
\(x\in\mathcal E\) is a finite sum, so no convergence condition is needed for
this algebraic action.  Formal resolvents are understood coefficient by
coefficient.  In the analytic statements, the matrices are restricted to the
indicated \(\ell^s\) spaces and are explicitly assumed to be bounded there.
For \(1\le s\le\infty\), the notation \(\|\cdot\|_s\) denotes the
\(\ell^s\)-norm for vectors and the corresponding subordinate norm for
matrices and operators; the relevant coordinate spaces are determined by the
domain and codomain displayed in each statement.
For a semi-infinite matrix \(X=[X_{i,j}]_{i,j\in\mathbb N_0}\) and
\(N\in\mathbb N_0\), write
\[
	X^{[N]}
	\coloneq
	\bigl[X_{i,j}\bigr]_{0\le i,j\le N}
\]
for its leading principal truncation, which has order \(N+1\).
Let \(L_1,\ldots,L_p\) be lower bidiagonal matrices and
\(U_1,\ldots,U_q\) upper bidiagonal matrices acting on \(\mathcal E\), with
entries
\begin{align*}
	(L_a)_{n,n}&=1,
	&
	(L_a)_{n+1,n}&=\ell_n^{(a)},
	&&
	a\in\{1,\ldots,p\},
	\quad
	n\in\mathbb N_0,
	\\
	(U_b)_{n,n}&=u_n^{(b)},
	&
	(U_b)_{n,n+1}&=1,
	&&
	b\in\{1,\ldots,q\},
	\quad
	n\in\mathbb N_0.
\end{align*}
All other entries vanish.  Put
\(A\coloneq L_1\cdots L_p\) and
\(B\coloneq U_q\cdots U_1\).  Consider the banded matrix
\begin{equation}
	\label{eq:general-two-sided-product}
	T
	\coloneq
	AB
	=
	L_1\cdots L_p U_q\cdots U_1.
\end{equation}
Enumerate the bidiagonal factors in their order of appearance in this
product by
\[
	F_a\coloneq L_a,
	\qquad 1\le a\le p,
	\qquad
	F_{p+b}\coloneq U_{q-b+1},
	\qquad 1\le b\le q.
\]
Then \(T=F_1F_2\cdots F_M\).
For \(a\in\{1,\ldots,M\}\), define the cyclic product starting with
\(F_a\) by
\[
	T_{(a)}
	\coloneq
	F_aF_{a+1}\cdots F_MF_1\cdots F_{a-1}.
\]
Thus \(T_{(1)}=T\).  The product \(T_{(a)}\) is obtained by moving
\(F_1,\ldots,F_{a-1}\) from the beginning of the product to the end,
without changing their relative order.  Such a change of the initial factor
is called a cyclic permutation of the factors.
The matrix \(T\) has at most \(p\) strict subdiagonals and \(q\) strict
superdiagonals.
The displayed factorization is called nonnegative when every
\(\ell_n^{(a)}\) and \(u_n^{(b)}\) is nonnegative, and positive when all of
them are strictly positive.  The fixed diagonal entries of the lower factors
and the fixed superdiagonal entries of the upper factors are the normalizing
entries of the factorization.

The entries of \(T\) combine contributions from all \(p+q\) factors.  The
two sparse matrices introduced below instead record one factor at a time and
retain the factor from which each coefficient comes.

Finite bidiagonal factorizations are a basic tool in the classical theory of
totally nonnegative matrices; see, among others,
\cite{GantmacherKrein,Pinkus2010,FallatJohnson2011}; the related finite
\(qd\) and \(LR\) transformations are reviewed in
\cite{GutknechtParlett2011}.  Those results concern
factorization and parametrization of a fixed finite totally nonnegative
matrix.  The result below instead gives an explicit normalized
semi-infinite refactorization \(UL=\widetilde L\widetilde U\) that preserves 
strict positivity coefficient by coefficient.  Its iteration
will imply that every cyclic permutation of an arbitrary prescribed PBF
again admits a PBF.

\begin{lemma}[Positive adjacent refactorization]
	\label{lem:positive-adjacent-refactorization}
	Let \(U\) and \(L\) be normalized upper and lower bidiagonal matrices,
	respectively, with
	\[
	U_{n,n}=u_n>0,
	\qquad U_{n,n+1}=1,
	\qquad L_{n,n}=1,
	\qquad L_{n+1,n}=\ell_n>0.
	\]
	Then there is a unique normalized positive pair
	\(\widetilde L,\widetilde U\) such that
	\begin{equation}
		UL=\widetilde L\widetilde U.
	\end{equation}
	Its nontrivial entries are determined recursively by
	\begin{align}
		\label{eq:positive-adjacent-refactorization-coefficients}
		\widetilde u_0
		&=u_0+\ell_0,
		&
		\widetilde\ell_n
		&=\frac{u_{n+1}\ell_n}{\widetilde u_n},
		&
		\widetilde u_{n+1}
		&=u_{n+1}+\ell_{n+1}-\widetilde\ell_n,
		\qquad n\in\mathbb N_0.
	\end{align}
	More precisely,
	\begin{equation}
		\label{eq:positive-adjacent-refactorization-inequalities}
		\widetilde u_n>\ell_n,
		\qquad
		0<\widetilde\ell_n<u_{n+1},
		\qquad n\in\mathbb N_0.
	\end{equation}
	If the two original coefficient sequences are bounded, then so are the
	refactorized sequences.
\end{lemma}

\begin{proof}
	The product \(UL\) has first superdiagonal equal to \(1\), diagonal
	entries \(u_n+\ell_n\), and first subdiagonal entries
	\(u_{n+1}\ell_n\).  On the other hand,
	\(\widetilde L\widetilde U\) has first superdiagonal equal to \(1\),
	diagonal entries
	\(\widetilde u_n+\widetilde\ell_{n-1}\), with
	\(\widetilde\ell_{-1}\coloneq0\), and first subdiagonal entries
	\(\widetilde\ell_n\widetilde u_n\).  Equality of the zeroth diagonal
	entry, followed successively by equality of the \(n\)-th subdiagonal and
	the \((n+1)\)-st diagonal entries, gives exactly
	\eqref{eq:positive-adjacent-refactorization-coefficients}.  This also
	proves uniqueness.

	It remains to prove that no denominator vanishes and that every new entry
	is positive.  The initial value satisfies
	\(\widetilde u_0=u_0+\ell_0>\ell_0\).  If
	\(\widetilde u_n>\ell_n\), then
	\[
	0<\widetilde\ell_n
	=\frac{u_{n+1}\ell_n}{\widetilde u_n}
	<u_{n+1},
	\]
	and consequently
	\[
	\widetilde u_{n+1}
	=u_{n+1}+\ell_{n+1}-\widetilde\ell_n
	>\ell_{n+1}.
	\]
	Induction proves \eqref{eq:positive-adjacent-refactorization-inequalities}
	and the strict positivity of the two refactorized factors.

	Finally, if \(u_n,\ell_n\le\Gamma\) for all \(n\), then
	\(\widetilde\ell_n<u_{n+1}\le\Gamma\), while
	\(\widetilde u_0\le2\Gamma\) and
	\(\widetilde u_{n+1}<u_{n+1}+\ell_{n+1}\le2\Gamma\).
	Thus both new coefficient sequences are bounded.
\end{proof}

The preceding local exchange can be iterated to restore the lower--upper
order after any cyclic permutation of the factors.

\begin{theorem}[Cyclic closure of positive bidiagonal factorizations]
	\label{thm:cyclic-PBF-closure}
	Suppose that the factorization
	\eqref{eq:general-two-sided-product} is positive.  Every cyclic product
	\(T_{(a)}\), with \(a\in\{1,\ldots,M\}\), admits a normalized positive
	bidiagonal factorization with \(p\) lower
	factors followed by \(q\) upper factors.  The resulting lower--upper
	factors are independent of the sequence of adjacent refactorizations.
	If the coefficients of the
	original factors are bounded, the coefficients in every one of these
	cyclic factorizations are bounded.
\end{theorem}

\begin{proof}
	Fix a cyclic product and read its factors from left to right.  Whenever an
	upper factor occurs immediately before a lower factor, apply
		Lemma~\ref{lem:positive-adjacent-refactorization} to replace that pair by
	a normalized positive lower--upper pair with the same product.  Each
	interchange decreases by one the number of upper--lower inversions in the
	factor order.  After finitely many interchanges no such inversion remains,
	so all \(p\) lower factors precede all \(q\) upper factors.  Every
	intermediate product equals the original cyclic product and every factor
	remains positive.

	We prove independence of the exchange sequence by induction on the number
	of upper--lower inversions.  There is nothing to prove when this number is
	zero.  If two exchange sequences begin with the same exchange, the induction
	hypothesis applies after that common first step.  If they begin with
	different exchanges, the two occurrences of \(UL\) are disjoint, because
	two such occurrences cannot overlap.  The corresponding refactorizations
	act on disjoint factor pairs and therefore commute.  Applying both exchanges
	in either order gives the same factor list, now with fewer inversions, so the
	induction hypothesis shows that the two sequences have the same terminal
	factors.

	Boundedness is preserved at each interchange by the
	last assertion of Lemma~\ref{lem:positive-adjacent-refactorization}; only
	finitely many interchanges are required.
\end{proof}

\begin{remark}[Boundary correction for separately truncated factors]
	\label{rem:adjacent-refactorization-finite-boundary}
	The refactorization identity concerns semi-infinite matrices.  If the
	factors are truncated separately before multiplication, the last diagonal
	entry acquires a boundary correction.  Indeed, if \(P_N\) denotes the coordinate projection onto
	the indices \(0,\ldots,N\), then in general
	\(P_NULP_N
	\neq
	(P_NUP_N)(P_NLP_N)\),
	because the product on the right omits the contribution
	\(U_{N,N+1}L_{N+1,N}=\ell_N\) to its last diagonal entry.
	When all lower factors precede all upper
	factors, every leading principal submatrix does factor as the product of
	the corresponding leading submatrices, as proved in
	Lemma~\ref{lem:global-PBF-compatible-truncations}.
\end{remark}

\subsection{Selected cyclic products and their minimum-height factorizations}
\label{subsec:balanced-Christoffel-word}

Each matrix \(T_{(a)}\) is obtained by choosing a different first factor in
the prescribed product and reading the factors cyclically from that point.
Starting from this cyclic order, each adjacent exchange
\(UL=\widetilde L\widetilde U\) preserves the product and positivity but
replaces the two factors involved.  Continuing these exchanges until no upper
factor precedes a lower factor gives the terminal lower--upper PBF of
\(T_{(a)}\) for every \(a\), by Theorem~\ref{thm:cyclic-PBF-closure}.

For the indices selected below, a finite initial sequence of exchanges reaches
a pattern of minimum height span.  We call a factorization with such a
factor-type pattern a \emph{minimum-height factorization}.  We retain the
factorization obtained at that stage.  Its factors are generally interlaced;
they coincide with the
terminal PBF only when this pattern is already lower--upper.  Theorem~
\ref{thm:positive-balanced-refactorization} proves that the resulting factors
are independent of the exchange sequence.  They still have product
\(T_{(a)}\), but now construct the two-block-diagonal matrix
\(\mathscr H_{(a)}\), whose projected resolvent recovers that of \(T_{(a)}\)
and whose powers define the star moment system of Section~
\ref{sec:balanced-star-Christoffel}.  We now determine the admissible indices.

Put
\[
	d\coloneq\gcd(p,q),
	\qquad
	P\coloneq\frac{p}{d},
	\qquad
	Q\coloneq\frac{q}{d},
	\qquad
	R\coloneq P+Q=\frac{M}{d}.
\]

\begin{definition}[Factor-position heights]
	\label{def:factor-position-heights}
	Let \(\sigma=(\sigma_1,\ldots,\sigma_M)\) be an ordered sequence
	containing \(p\) symbols \(L\) and \(q\) symbols \(U\).  Set
	\[
		\eta_0\coloneq0,
		\qquad
		\eta_k
		\coloneq
		\eta_{k-1}
		+
		\begin{cases}
			Q,&\sigma_k=L,\\
			-P,&\sigma_k=U,
		\end{cases}
		\qquad
		k\in\{1,\ldots,M\}.
	\]
	The integers \(\eta_0,\ldots,\eta_M\) are the
	\emph{factor-position heights} of \(\sigma\).  Since \(pQ=qP\), one has
	\(\eta_M=\eta_0=0\).  The \emph{height span} of \(\sigma\) is
	\[
		\max_{0\le k\le M}\eta_k
		-
		\min_{0\le k\le M}\eta_k.
	\]
\end{definition}

\begin{definition}[Selected cyclic products]
	\label{def:admissible-cyclic-products}
	Define
	\begin{equation}
	\label{eq:admissible-cyclic-index-set}
	\begin{gathered}
		c_-\coloneq\left\lceil\frac{P}{Q}\right\rceil,
		\qquad
		c_+\coloneq\left\lceil\frac{Q}{P}\right\rceil,\\
		a_-\coloneq p-c_-+1,
		\qquad
		a_+\coloneq p+1+\min\{q-1,c_+\},
		\qquad
		\mathcal A_{\mathrm{cyc}}
		\coloneq\{a_-,a_-+1,\ldots,a_+\}.
	\end{gathered}
	\end{equation}
	The matrices \(T_{(a)}\) with
	\(a\in\mathcal A_{\mathrm{cyc}}\) are called the selected cyclic
	products.
\end{definition}

\begin{lemma}[Factor pattern of minimum height span]
	\label{lem:minimum-height-factor-pattern}
	For \(k\in\{0,\ldots,M\}\), let \(r_k\) be the remainder of \(kQ\)
	upon division by \(R\), and define
	\(\sigma_{\mathrm{min}}
	=(\sigma_{\mathrm{min},1},\ldots,\sigma_{\mathrm{min},M})\) by
	\[
		\sigma_{\mathrm{min},k+1}
		=
		\begin{cases}
			L,&r_k<P,\\
			U,&r_k\ge P,
		\end{cases}
		\qquad k\in\{0,\ldots,M-1\}.
	\]
	The sequence \(\sigma_{\mathrm{min}}\) contains \(p\) symbols \(L\) and
	\(q\) symbols \(U\).  Its factor-position heights satisfy
	\(\eta_k=r_k\), and its height span is \(R-1\), the minimum possible.
	Moreover, it begins with \(c_-\) symbols \(L\) and ends with
	\(c_+\) symbols \(U\).

	There are exactly \(R\) distinct factor-type sequences with height span
	\(R-1\), namely the distinct cyclic rotations of
	\(\sigma_{\mathrm{min}}\).  Among them,
	\(\sigma_{\mathrm{min}}\) is the unique sequence whose heights all belong
	to \(\{0,\ldots,R-1\}\).  If \(N_L(\rho;k)\) denotes the number of
	letters \(L\) among the first \(k\) positions of a minimizing sequence
	\(\rho\), then
	\[
		N_L(\sigma_{\mathrm{min}};k)
		\ge N_L(\rho;k),
		\qquad k\in\{0,\ldots,M\}.
	\]
	The first \(R\) symbols of \(\sigma_{\mathrm{min}}\) form the lower
	Christoffel word with \(P\) letters \(L\) and \(Q\) letters \(U\), with
	\(L\) horizontal and \(U\) vertical; hence
	\(\sigma_{\mathrm{min}}\) is the \(d\)-th power of that primitive word
	\cite{BerstelLauveReutenauerSaliola2009}.
\end{lemma}

\begin{proof}
	Let \(\eta_0,\ldots,\eta_M\) be the factor-position heights of
	\(\sigma_{\mathrm{min}}\).  We prove by induction that
	\(\eta_k=r_k\).  This is clear for \(k=0\).  If \(r_k<P\), then
	position \(k+1\) is lower and \(\eta_{k+1}=r_k+Q<R\).  If
	\(r_k\ge P\), then that position is upper and
	\(\eta_{k+1}=r_k-P=r_k+Q-R\).  In either case
	\(\eta_{k+1}=r_{k+1}\).

	Since \(\gcd(Q,R)=1\), the residues
	\(r_0,\ldots,r_{R-1}\) are a permutation of
	\(0,\ldots,R-1\).  Hence one block of \(R\) positions contains exactly
	\(P\) lower and \(Q\) upper positions.  Repetition over the \(d\) residue
	cycles gives exactly \(p=dP\) lower and \(q=dQ\) upper positions.
	The sequence begins with \(c_-=\lceil P/Q\rceil\) lower positions.
	Reading the residue cycle backwards from zero shows that it ends with
	\(c_+=\lceil Q/P\rceil\) upper positions.

	The factor-position heights are the residues \(r_k\), so their span is
	\(R-1\).  To prove minimality, observe that the two increments \(+Q\) and
	\(-P\) are congruent to \(Q\) modulo \(R=P+Q\).  Therefore the first
	\(R\) factor-position heights of any sequence represent all residue
	classes modulo \(R\).  Any \(R\) integers with distinct residues modulo
	\(R\) have span at least \(R-1\).

	Let \(\rho\) be any sequence whose height span is \(R-1\), and denote its
	heights by \(\theta_k\).  The integers
	\(\theta_0,\ldots,\theta_{R-1}\) have all the residues modulo \(R\) and
	lie in an interval of length \(R-1\).  Since \(\theta_0=0\), that interval
	is
	\[
		\{-j,-j+1,\ldots,R-1-j\}
	\]
	for a unique \(j\in\{0,\ldots,R-1\}\).  For every \(k\), the height
	\(\theta_k\) is therefore the unique representative of \(r_k\) in this
	interval.  Thus the increments, and hence the letters of \(\rho\), are
	uniquely determined by \(j\).  Choosing the unique position \(s\) with
	\(r_s=j\) shows that \(\rho\) is the cyclic rotation of
	\(\sigma_{\mathrm{min}}\) beginning at position \(s+1\).  Conversely,
	every cyclic rotation has the same height span.  The primitive block has
	length \(R\), because \(P\) and \(Q\) are coprime; hence there are exactly
	\(R\) distinct rotations.  The value \(j=0\) gives
	\(\sigma_{\mathrm{min}}\), and it is the unique rotation whose heights are
	all nonnegative.

	Finally, both possible height increments are congruent to \(Q\) modulo
	\(R\), so \(\theta_k\equiv r_k\pmod R\).  The representative of \(r_k\)
	in the interval above is either \(r_k\) or \(r_k-R\); consequently
	\(\theta_k\le r_k\).  For any factor-type sequence,
	\[
		N_L(\rho;k)=\frac{\theta_k+Pk}{R}.
	\]
	Since the corresponding height of \(\sigma_{\mathrm{min}}\) is \(r_k\),
	the asserted prefix dominance follows.  The residue construction of the
	first block is precisely the lower Christoffel construction, and the
	periodicity of \(r_k\) gives its \(d\)-fold repetition.
\end{proof}

\begin{figure}
	\centering
	\begin{tikzpicture}[
		x=0.76cm,
		y=0.46cm,
		>=stealth,
		heightnode/.style={circle,draw=black,fill=white,minimum size=5.2mm,
			inner sep=0pt,font=\scriptsize},
		lowerstep/.style={->,thick,NavyBlue},
		upperstep/.style={->,thick,BrickRed},
		steplabel/.style={midway,fill=white,inner sep=1pt,font=\scriptsize},
		panel/.style={font=\small},
		input/.style={font=\scriptsize,align=left},
		annotation/.style={font=\scriptsize}
	]
		\begin{scope}[yshift=6.4cm]
			\node[panel,anchor=west] at (0,4.75)
				{\textup{(a)} \((p,q)=(3,2)\)};
			\node[input,anchor=west] at (0,2)
				{\(T_{(2)}:\ LLUUL\qquad T_{(3)}:\ LUULL\)\\[2pt]
				 \(T_{(4)}:\ UULLL\qquad T_{(5)}:\ ULLLU\)};
			\draw[->,thick] (7.1,2) -- node[above,font=\scriptsize,align=center]
				{positive adjacent\\refactorizations} (9.2,2);
			\node[annotation,anchor=south] at (12,4.45)
				{\(\sigma_{\mathrm{min}}=LLULU\)};
			\fill[NavyBlue!3] (9.65,-0.35) rectangle (14.35,4.35);
			\draw[densely dashed,Gray] (9.65,0) -- (14.35,0);
			\draw[densely dashed,Gray] (9.65,4) -- (14.35,4);
			\node[heightnode] (a0) at (10,0) {\(0\)};
			\node[heightnode] (a1) at (10.8,2) {\(2\)};
			\node[heightnode] (a2) at (11.6,4) {\(4\)};
			\node[heightnode] (a3) at (12.4,1) {\(1\)};
			\node[heightnode] (a4) at (13.2,3) {\(3\)};
			\node[heightnode] (a5) at (14,0) {\(0\)};
			\draw[lowerstep] (a0) -- node[steplabel,above,sloped] {\(L\)} (a1);
			\draw[lowerstep] (a1) -- node[steplabel,above,sloped] {\(L\)} (a2);
			\draw[upperstep] (a2) -- node[steplabel,above,sloped] {\(U\)} (a3);
			\draw[lowerstep] (a3) -- node[steplabel,above,sloped] {\(L\)} (a4);
			\draw[upperstep] (a4) -- node[steplabel,above,sloped] {\(U\)} (a5);
		\end{scope}

		\begin{scope}[yshift=3.2cm]
			\node[panel,anchor=west] at (0,4.75)
				{\textup{(b)} \((p,q)=(2,3)\)};
			\node[input,anchor=west] at (0,2)
				{\(T_{(2)}:\ LUUUL\qquad T_{(3)}:\ UUULL\)\\[2pt]
				 \(T_{(4)}:\ UULLU\qquad T_{(5)}:\ ULLUU\)};
			\draw[->,thick] (7.1,2) -- node[above,font=\scriptsize,align=center]
				{positive adjacent\\refactorizations} (9.2,2);
			\node[annotation,anchor=south] at (12,4.45)
				{\(\sigma_{\mathrm{min}}=LULUU\)};
			\fill[NavyBlue!3] (9.65,-0.35) rectangle (14.35,4.35);
			\draw[densely dashed,Gray] (9.65,0) -- (14.35,0);
			\draw[densely dashed,Gray] (9.65,4) -- (14.35,4);
			\node[heightnode] (b0) at (10,0) {\(0\)};
			\node[heightnode] (b1) at (10.8,3) {\(3\)};
			\node[heightnode] (b2) at (11.6,1) {\(1\)};
			\node[heightnode] (b3) at (12.4,4) {\(4\)};
			\node[heightnode] (b4) at (13.2,2) {\(2\)};
			\node[heightnode] (b5) at (14,0) {\(0\)};
			\draw[lowerstep] (b0) -- node[steplabel,above,sloped] {\(L\)} (b1);
			\draw[upperstep] (b1) -- node[steplabel,above,sloped] {\(U\)} (b2);
			\draw[lowerstep] (b2) -- node[steplabel,above,sloped] {\(L\)} (b3);
			\draw[upperstep] (b3) -- node[steplabel,above,sloped] {\(U\)} (b4);
			\draw[upperstep] (b4) -- node[steplabel,above,sloped] {\(U\)} (b5);
		\end{scope}

		\begin{scope}
			\node[panel,anchor=west] at (0,4.1)
				{\textup{(c)} \((p,q)=(4,2)\), \(d=2\)};
			\node[input,anchor=west] at (0,2)
				{\(T_{(3)}:\ LLUULL\qquad T_{(4)}:\ LUULLL\)\\[2pt]
				 \(T_{(5)}:\ UULLLL\qquad T_{(6)}:\ ULLLLU\)};
			\draw[->,thick] (7.1,2) -- node[above,font=\scriptsize,align=center]
				{positive adjacent\\refactorizations} (9.2,2);
			\node[annotation,anchor=south] at (12,2.45)
				{\(\sigma_{\mathrm{min}}=LLULLU\)};
			\fill[NavyBlue!3] (9.65,-0.35) rectangle (14.35,2.35);
			\draw[densely dashed,Gray] (9.65,0) -- (14.35,0);
			\draw[densely dashed,Gray] (9.65,2) -- (14.35,2);
			\draw[densely dashed,Gray] (12,-0.45) -- (12,2.45);
			\node[heightnode] (c0) at (9.8,0) {\(0\)};
			\node[heightnode] (c1) at (10.55,1) {\(1\)};
			\node[heightnode] (c2) at (11.3,2) {\(2\)};
			\node[heightnode] (c3) at (12,0) {\(0\)};
			\node[heightnode] (c4) at (12.75,1) {\(1\)};
			\node[heightnode] (c5) at (13.5,2) {\(2\)};
			\node[heightnode] (c6) at (14.2,0) {\(0\)};
			\draw[lowerstep] (c0) -- node[midway,above=3pt,sloped,
				inner sep=0pt,font=\scriptsize,text=NavyBlue] {\(L\)} (c1);
			\draw[lowerstep] (c1) -- node[midway,above=3pt,sloped,
				inner sep=0pt,font=\scriptsize,text=NavyBlue] {\(L\)} (c2);
			\draw[upperstep] (c2) -- node[midway,above=3pt,sloped,
				inner sep=0pt,font=\scriptsize,text=BrickRed] {\(U\)} (c3);
			\draw[lowerstep] (c3) -- node[midway,above=3pt,sloped,
				inner sep=0pt,font=\scriptsize,text=NavyBlue] {\(L\)} (c4);
			\draw[lowerstep] (c4) -- node[midway,above=3pt,sloped,
				inner sep=0pt,font=\scriptsize,text=NavyBlue] {\(L\)} (c5);
			\draw[upperstep] (c5) -- node[midway,above=3pt,sloped,
				inner sep=0pt,font=\scriptsize,text=BrickRed] {\(U\)} (c6);
		\end{scope}
	\end{tikzpicture}
	\caption{Selected cyclic products and their minimum-height
		refactorizations.  The words on the left are the factor-type orders in
		the products \(T_{(a)}\); the path on the right belongs to the resulting
		minimum-height factorization.  Blue \(L\)-steps have
		increment \(+Q\), and red \(U\)-steps have increment \(-P\).  Panel
		\textup{(c)} displays the repeated residue cycle when \(d>1\).}
	\label{fig:minimum-height-factor-patterns}
\end{figure}
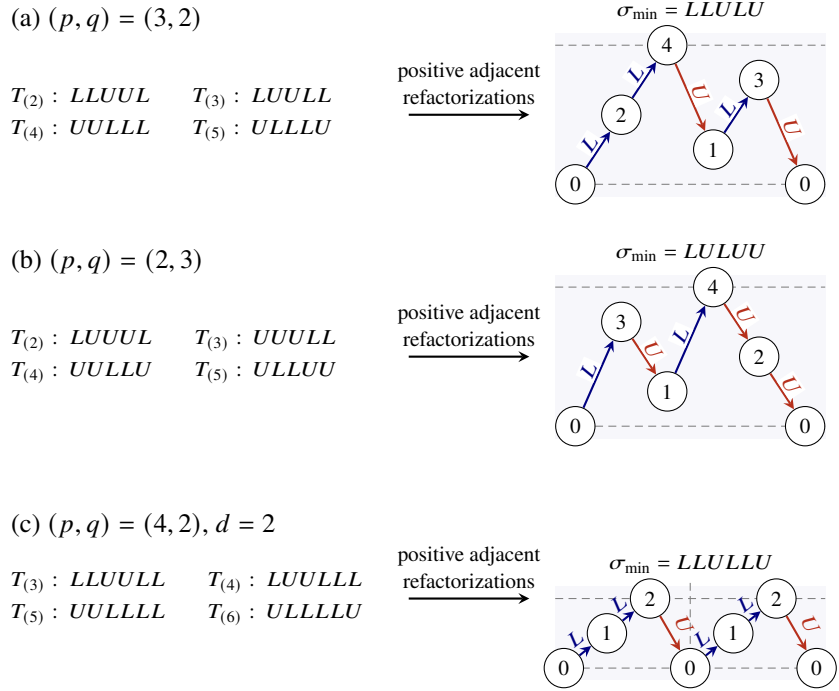

Figure~\ref{fig:minimum-height-factor-patterns} shows the cyclic factor-type
orders of the selected products and the common minimum-height pattern reached
after a suitable finite sequence of adjacent refactorizations.  We retain the
factorization obtained at this stage in order to construct
\(\mathscr H_{(a)}\).  Further adjacent refactorizations are possible and
eventually produce the terminal lower--upper PBF.

\begin{theorem}[Positive refactorization into the minimum-height pattern]
	\label{thm:positive-balanced-refactorization}
	For \(a\in\{1,\ldots,M\}\), the following conditions are equivalent:
	\begin{enumerate}[label=\textnormal{(\roman*)}]
		\item \(a\in\mathcal A_{\mathrm{cyc}}\);
		\item the cyclically ordered factors of \(T_{(a)}\) can be transformed,
		by finitely many directed adjacent refactorizations
		\(UL=\widetilde L\widetilde U\), into normalized positive bidiagonal
		factors whose lower--upper type sequence is
		\(\sigma_{\mathrm{min}}\).
	\end{enumerate}
	When these conditions hold, the resulting factors
	\(\widetilde F^{(a)}_1,\ldots,\widetilde F^{(a)}_M\) are independent of
	the sequence of adjacent refactorizations and satisfy
	\begin{equation}
		\label{eq:selected-cyclic-star-factorization}
		T_{(a)}
		=
		\widetilde F^{(a)}_1\cdots\widetilde F^{(a)}_M.
	\end{equation}
	If the original factor coefficients are bounded, then the coefficients of
	the resulting factors are bounded.
	Moreover, this choice is maximal among the minimum-height patterns: if the
	factor-type sequence of a cyclic product can be transformed by directed
	exchanges into any sequence of height span \(R-1\), then it can be
	transformed into \(\sigma_{\mathrm{min}}\).
\end{theorem}

\begin{proof}
	By Lemma~\ref{lem:minimum-height-factor-pattern},
	\(\sigma_{\mathrm{min}}\) begins with \(c_-\) lower positions and ends
	with \(c_+\) upper positions.  We determine which cyclic factor orders can
	be transformed into this pattern by directed exchanges.

	We first record the reachability criterion used below.  For a factor-type
	sequence \(\rho\), let \(N_L(\rho;k)\) be the number of lower positions in
	its first \(k\) entries.  If \(\widehat\rho\) is obtained from \(\rho\) by one
	directed exchange \(UL\mapsto LU\), then
	\(N_L(\widehat\rho;k)\ge N_L(\rho;k)\) for every \(k\).  Hence a
	necessary condition for transforming \(\rho\) into a target sequence
	\(\tau\) is
	\[
		N_L(\tau;k)\ge N_L(\rho;k),
		\qquad k\in\{0,\ldots,M\}.
	\]
	It is also sufficient.  Indeed, let
	\(i_1<\cdots<i_p\) and \(j_1<\cdots<j_p\) be the positions of the lower
	factors in \(\rho\) and \(\tau\), respectively.  The prefix inequalities
	are equivalent to \(j_m\le i_m\) for every \(m\).  Proceed by induction on
	\(m\).  Before stage \(m\), the first \(m-1\) lower factors occupy the
	positions \(j_1,\ldots,j_{m-1}\), while the \(m\)-th lower factor still
	occupies position \(i_m\): the preceding moves cross only upper factors
	lying to its left and do not change the positions of the later lower
	factors.  Since all earlier target positions are strictly smaller than
	\(j_m\), every position from \(j_m\) through \(i_m-1\) is upper.  Successive
	directed exchanges \(UL\mapsto LU\) therefore move the \(m\)-th lower factor
	from \(i_m\) to \(j_m\) without disturbing the factors already placed.
	After the \(p\) stages the target sequence is reached.

	Consider first \(a\in\{1,\ldots,p\}\).  The factor types in \(T_{(a)}\)
	occur as
	\[
		L^{p-a+1}U^qL^{a-1}.
	\]
	Put \(x=p-a+1\).  If \(x>c_-\), the first \(c_-+1\) positions of the initial
	sequence are lower, whereas \(\sigma_{\mathrm{min}}\) has only \(c_-\) lower positions
	in that initial segment; reachability therefore fails.  Conversely, assume
	\(x\le c_-\).  Up to position \(x\), both sequences are lower.  From position
	\(x+1\) through position \(x+q\), the initial arrangement has exactly \(x\)
	lower positions, while
	\(N_L(\sigma_{\mathrm{min}};k)\ge\min\{k,c_-\}\ge x\).  After position
	\(x+q\), any arrangement containing only \(q\)
	upper symbols has at least \(k-q\) lower symbols in its first \(k\)
	positions, which is exactly the number in the initial arrangement.  Hence the
	prefix condition is equivalent to \(x\le c_-\), or \(a\ge a_-\).
	
	Now let \(a=p+b\), with \(b\in\{1,\ldots,q\}\).  The cyclic product begins
	with the type sequence
	\[
		U^{q-b+1}L^pU^{b-1}.
	\]
	Put \(z=b-1\), the length of its final upper group.  If \(z>c_+\), the
	suffix of length \(z\) of \(\sigma_{\mathrm{min}}\) contains a lower symbol, whereas the
	corresponding suffix of the initial arrangement contains only upper symbols.  At
	position \(M-z\), the initial sequence has already accumulated all \(p\)
	lower symbols but \(\sigma_{\mathrm{min}}\) has not, so the prefix condition fails.
	Conversely, let \(z\le c_+\).  Before the initial lower group begins, its
	prefix contains no lower symbol.  During that lower group, the final \(z\)
	positions of \(\sigma_{\mathrm{min}}\) are upper, so its first \(k\) positions contain at
	most \(q-z\) upper symbols and therefore at least as many lower symbols as
	the initial arrangement.  After the initial lower group ends, the omitted suffix
	of \(\sigma_{\mathrm{min}}\) lies inside its final upper group, so its prefix
	contains all \(p\) lower symbols.  Thus the prefix condition is equivalent
	to \(z=b-1\le c_+\).  Since \(b\le q\), this is equivalent to
	\(a\le a_+\).  The two cases give exactly
	\(\mathcal A_{\mathrm{cyc}}\).
	
	For such an index \(a\), choose any sequence of adjacent exchanges taking
	the initial type sequence to \(\sigma_{\mathrm{min}}\), and apply
	Lemma~\ref{lem:positive-adjacent-refactorization} at every exchange.  The
	product remains \(T_{(a)}\), and all factors remain normalized and
	positive.  The final assertion of
	Lemma~\ref{lem:positive-adjacent-refactorization} preserves uniform
	boundedness of the coefficient sequences at each exchange.  Hence, if the
	original sequences are bounded, the finitely many exchanges produce
	bounded factor coefficients.

	Every directed exchange decreases by one the number of inversions
	consisting of an upper position followed by a lower position.  Therefore
	every exchange sequence with the prescribed initial and final type
	sequences has the same length and is automatically reduced.
	
	It remains to prove independence of the chosen exchange sequence.  Give
	temporary labels to the factors in order to follow them through the
	calculation.  Every upper--lower pair whose relative order changes is
	exchanged exactly once, and no other pair is exchanged.  The exchanges
	involving a fixed factor have a forced order.  Thus every reduced exchange
	sequence is a linear extension of the same finite partially ordered set of
	required crossings.  Two crossings that share an upper factor are ordered
	by the unchanged relative order of the lower factors that it crosses; the
	same argument applies to two crossings that share a lower factor.  Hence
	two incomparable crossings involve disjoint factor pairs.  Any two linear
	extensions differ by successive interchanges of adjacent incomparable
	crossings, and the corresponding bidiagonal refactorizations commute.
	Consequently all reduced exchange sequences give the same final matrices.

	Finally, let \(\rho\) be any minimum-height target reachable from the type
	sequence of a cyclic product.  By
	Lemma~\ref{lem:minimum-height-factor-pattern},
	\[
		N_L(\sigma_{\mathrm{min}};k)
		\ge N_L(\rho;k)
	\]
	for every \(k\).  The reachability criterion and the inequalities from the
	initial sequence to \(\rho\) therefore show that the same initial sequence
	can be transformed into \(\sigma_{\mathrm{min}}\).  This proves the final
	maximality assertion.
\end{proof}

For the remainder of the construction, fix
\(a\in\mathcal A_{\mathrm{cyc}}\).
Use the minimum-height factorization
\eqref{eq:selected-cyclic-star-factorization}.  The matrix built
from these factors will be denoted by \(\mathscr H_{(a)}\).  The index
\((a)\) records the selected cyclic product \(T_{(a)}\).  Its entries are taken from
the new factors \(\widetilde F^{(a)}_1,\ldots,\widetilde F^{(a)}_M\), not
from the original factors with their labels permuted.  They coincide with the
terminal lower--upper PBF factors only when the minimum-height pattern itself
is lower--upper.
In particular, \(p+1\in\mathcal A_{\mathrm{cyc}}\) and
\(T_{(p+1)}=BA\). 
The next subsection constructs \(\mathscr H_{(a)}\) from these factors, and
Theorem~\ref{thm:balanced-two-extreme-diagonals} proves that it has exactly
two nonzero block diagonals.

\subsection{Direct construction of the two sparse matrices}
\label{subsec:direct-matrices-from-PBF}

We now define the two matrices directly from the bidiagonal coefficients;
the cyclic matrix will enter later in the resolvent identities.

\begin{definition}[The two sparse matrices]
\label{def:direct-sparse-matrices}
Retain the prescribed order
\[
F_1,\ldots,F_M
=
L_1,\ldots,L_p,U_q,\ldots,U_1
\]
and define
\[
	\omega_j=(j-1)Q,
	\qquad j\in\{1,\ldots,p\},
	\qquad
	\omega_{p+b}=(q-b+1)P,
	\qquad b\in\{1,\ldots,q\}.
\]
For \(h\in\mathbb N_0\), let
\[
	I_h
	\coloneq
	\left\{j\in\{1,\ldots,M\}:
	h-\omega_j\in R\mathbb N_0\right\},
	\qquad
	n_j(h)\coloneq\frac{h-\omega_j}{R}
	\quad(j\in I_h).
\]
The rows and columns of the first sparse matrix are indexed by the pairs
\((h,j)\) with \(j\in I_h\), ordered first by \(h\) and then by \(j\).
For any factor-position index \(k\in\{1,\ldots,M\}\), write
\(k^+=k+1\) when \(k<M\) and \(k^+=1\) when \(k=M\).  Define
\(\mathscr H\) by
\begin{equation}
	\label{eq:direct-prescribed-matrix-definition}
	\mathscr H_{(h,j),(h+Q,j^+)}=1,
\end{equation}
and by the following coefficient entries:
\begin{align}
	\mathscr H_{(h,j),(h-P,j^+)}
	&=\ell_{n_j(h)-1}^{(j)},
	&&j\in\{1,\ldots,p\},\quad n_j(h)\ge1,
	\label{eq:direct-prescribed-lower-entry}
	\\
	\mathscr H_{(h,p+b),(h-P,(p+b)^+)}
	&=u_{n_{p+b}(h)}^{(q-b+1)},
	&&b\in\{1,\ldots,q\}.
	\label{eq:direct-prescribed-upper-entry}
\end{align}
All entries not specified in
\eqref{eq:direct-prescribed-matrix-definition}--
\eqref{eq:direct-prescribed-upper-entry} are zero.  A term with a negative
height is absent.  This defines \(\mathscr H\) directly from the original
factor coefficients, including its finite initial part.  When \(d=1\), each
nonempty \(I_h\) contains one factor position; when \(d>1\), the same formula
gives finite blocks indexed by the elements of \(I_h\).

For the second matrix, use the selected index \(a\) fixed after
Theorem~\ref{thm:positive-balanced-refactorization} and write
\[
	T_{(a)}
	=\widetilde F^{(a)}_1\cdots\widetilde F^{(a)}_M.
\]
For \(k\in\{0,\ldots,R-1\}\), let \(r_k\) be the remainder of \(kQ\)
upon division by \(R\).  For \(j\in\{1,\ldots,M\}\), let
\[
	\widetilde\omega_j
	\coloneq r_{(j-1)\bmod R},
\]
and put
\[
	I_h^{(a)}
	\coloneq
	\left\{j\in\{1,\ldots,M\}:
	h-\widetilde\omega_j\in R\mathbb N_0\right\},
	\qquad
	\widetilde n_j(h)
	\coloneq\frac{h-\widetilde\omega_j}{R}.
\]
The rows and columns of \(\mathscr H_{(a)}\) are indexed by
\((h,j)\), with \(j\in I_h^{(a)}\), in the same order as above.
For every such pair, place the normalizing entry of
\(\widetilde F^{(a)}_j\), which is equal to \(1\), at
\[
	(\mathscr H_{(a)})_{(h,j),(h+Q,j^+)}=1.
\]
If \(\widetilde F^{(a)}_j\) is lower bidiagonal and
\(\widetilde n_j(h)\ge1\), place its subdiagonal entry at
\[
	(\mathscr H_{(a)})_{(h,j),(h-P,j^+)}
	=
	(\widetilde F^{(a)}_j)_{\widetilde n_j(h),\widetilde n_j(h)-1}.
\]
If \(\widetilde F^{(a)}_j\) is upper bidiagonal, place its diagonal entry at
\[
	(\mathscr H_{(a)})_{(h,j),(h-P,j^+)}
	=
	(\widetilde F^{(a)}_j)_{\widetilde n_j(h),\widetilde n_j(h)}.
\]
All other entries are zero.  Here every
\(I_h^{(a)}\) has exactly \(d\) elements.  Thus the definition gives
two block diagonals from the first row; when \(d=1\), it gives directly
\[
	(\mathscr H_{(a)})_{h,h+q}=1,
	\qquad
	(\mathscr H_{(a)})_{h,h-p}=c_h^{(a)}>0
	\quad(h\ge p).
\]
\end{definition}

Section~\ref{sec:height-ordering-block-CF} proves these placement rules in
detail and derives the block continued fraction.  It extends the
factor-position heights of
Definition~\ref{def:factor-position-heights} to the coordinate labels and
uses them to verify the two target indices and to relate the directly
constructed matrices to the products \(T\) and \(T_{(a)}\).

\subsection{Auxiliary cyclic realization}
\label{subsec:cyclic-linearization}

The two sparse matrices have already been defined from the factor
coefficients.  We now introduce an auxiliary matrix that will prove their
relation with the original product and its cyclic permutations.  Enlarge the
space to the direct sum
\(\mathcal E^{\oplus M}=\mathcal E\oplus\cdots\oplus\mathcal E\), with
\(M=p+q\) ordered summands.  The \(a\)-th copy of \(\mathcal E\) means the
\(a\)-th summand in this direct sum.

Recall the ordered factors \(F_1,\ldots,F_M\) introduced in
Subsection~\ref{subsec:positive-adjacent-refactorization}, so that
\(T=F_1F_2\cdots F_M\).
Introduce the cyclic matrix
\begin{equation*}
	\mathscr C_{p,q}
	\coloneq
	\begin{bNiceMatrix}
		0 & F_1 & 0 & \cdots & 0\\
		0 & 0 & F_2 & \ddots & \vdots\\
		\vdots & \ddots & \ddots & \ddots & 0\\
		0 & \cdots & 0 & 0 & F_{M-1}\\
		F_M & 0 & \cdots & 0 & 0
	\end{bNiceMatrix}.
\end{equation*}
For \(a\in\{1,\ldots,M\}\), let
\(\Pi_a:\mathcal E^{\oplus M}\to\mathcal E\) be the coordinate projection,
so that \(\Pi_a(x^{(1)},\ldots,x^{(M)})=x^{(a)}\).  Its transpose
\(\Pi_a^{\top}\) embeds \(\mathcal E\) into the \(a\)-th summand.  Thus
\(\Pi_1\) selects the first block row of the displayed cyclic matrix and
\(\Pi_1^{\top}\) selects its first block column.

\begin{lemma}[The cyclic matrix and the factor products]
	\label{lem:general-cyclic-linearization}
	One complete circuit through the factor positions gives
	\[
	\mathscr C^M
	=
	\operatorname{diag}
	\bigl(T_{(1)},\ldots,T_{(M)}\bigr).
	\]
	In particular,
	\[
	\Pi_1\mathscr C^{M}\Pi_1^{\top}=T.
	\]
	Consequently, as an identity of formal power series in \(\tau\),
	\begin{equation*}
		\Pi_1(I-\tau\mathscr C)^{-1}\Pi_1^{\top}
		=
		(I-\tau^M T)^{-1}.
	\end{equation*}
\end{lemma}

\begin{proof}
	The nonzero blocks of \(\mathscr C\) are
	\[
	(\mathscr C)_{a,a+1}=F_a,
	\qquad
	a\in\{1,\ldots,M-1\},
	\]
	and
	\[
	(\mathscr C)_{M,1}=F_M.
	\]
	Reading the block indices cyclically, induction on \(s\) gives
	\[
	(\mathscr C^{s})_{a,a+s}
	=
	F_aF_{a+1}\cdots F_{a+s-1}.
	\]
	For \(s\in\{1,\ldots,M-1\}\), this is the only nonzero block in row
	\(a\), and it is not on the main block diagonal.  Therefore
	\[
	\Pi_1\mathscr C^{s}\Pi_1^{\top}=0,
	\qquad
	s\in\{1,\ldots,M-1\}.
	\]

	After \(M\) steps the block position returns to its starting point, and
	\[
	(\mathscr C^{M})_{a,a}
	=
	F_aF_{a+1}\cdots F_MF_1\cdots F_{a-1}.
	\]
	These are precisely \(T_{(1)},\ldots,T_{(M)}\), which proves the first
	assertion and its \((1,1)\)-block identity.  Repeating complete circuits
	shows, for \(n\ge0\) and \(s\in\{0,\ldots,M-1\}\), that
	\[
	\Pi_1\mathscr C^{Mn+s}\Pi_1^{\top}
	=
	\begin{cases}
		T^n,&s=0,\\
		0,&s\ne0.
	\end{cases}
	\]
	Substitution in the formal Neumann series now gives
	\[
		\Pi_1(I-\tau\mathscr C)^{-1}\Pi_1^{\top}
=
		\sum_{m\ge0}
		\tau^m
		\Pi_1\mathscr C^{m}\Pi_1^{\top}
=
		\sum_{n\ge0}
		\tau^{Mn}T^n
=
		(I-\tau^M T)^{-1}.
	\]
\end{proof}

\subsection{Darboux interpretation and off-diagonal resolvent blocks}
\label{subsec:cyclic-Darboux-transforms}

The cyclic products \(T_{(a)}\) introduced above also form the Darboux orbit
of the prescribed factorization: passing from \(T_{(a)}\) to
\(T_{(a+1)}\) moves the first factor to the right end.  By
Lemma~\ref{lem:general-cyclic-linearization}, these products occur as the
diagonal blocks of \(\mathscr C^M\).

When the prescribed factors are Christoffel connection matrices, the
matrices \(T_{(a)}\) are the recurrence matrices of the corresponding
transformed systems
\cite{BranquinhoFoulquieManas2026Factorization}.  Their precise spectral
realization, including the bounded and unbounded Favard cases, is developed
in Section~\ref{sec:balanced-star-Christoffel}.

The diagonal blocks describe complete circuits through the factors.  The
off-diagonal blocks of the cyclic resolvent describe the intervening partial
products.  We now make this statement precise.  For
\(a,b\in\{1,\ldots,M\}\), let \(s(a,b)\in\{0,\ldots,M-1\}\) be the unique
cyclic distance determined by
\[
b\equiv a+s(a,b)\pmod M,
\]
where factor positions are represented by \(1,\ldots,M\).  Define
\[
F_{a\to b}
\coloneq
F_aF_{a+1}\cdots F_{a+s(a,b)-1},
\]
with cyclic indices; when \(a=b\), this is the empty product and hence the
identity.  The next proposition records both the diagonal and off-diagonal
resolvent blocks in terms of these objects.

	\begin{proposition}[Resolvent blocks for the cyclic Darboux transforms]
		\label{prop:Darboux-shifted-resolvent-blocks}
		For \(a,b\in\{1,\ldots,M\}\), the partial cyclic product defined above
		satisfies
	\begin{equation}
		\label{eq:Darboux-partial-intertwining}
		T_{(a)}F_{a\to b}=F_{a\to b}T_{(b)}.
	\end{equation}
	Moreover, as an identity of formal power series in \(\tau\),
	\begin{equation}
		\label{eq:Darboux-shifted-resolvent-right}
		\Pi_a(I-\tau\mathscr C)^{-1}\Pi_b^{\top}
		=
		\tau^{s(a,b)}F_{a\to b}(I-\tau^M T_{(b)})^{-1}.
	\end{equation}
	Equivalently,
	\begin{equation}
		\label{eq:Darboux-shifted-resolvent-left}
		\Pi_a(I-\tau\mathscr C)^{-1}\Pi_b^{\top}
		=
		\tau^{s(a,b)}(I-\tau^M T_{(a)})^{-1}F_{a\to b}.
	\end{equation}
	In particular,
	\begin{equation}
		\label{eq:Darboux-diagonal-resolvent-block}
		\Pi_a(I-\tau\mathscr C)^{-1}\Pi_a^{\top}
		=
		(I-\tau^M T_{(a)})^{-1}.
	\end{equation}
\end{proposition}

\begin{proof}
	The identity \eqref{eq:Darboux-partial-intertwining} follows by expanding both
	sides. Indeed, the left-hand side is
	\[
	F_a\cdots F_{a+s(a,b)-1}
	F_{a+s(a,b)}\cdots F_{a-1}
	F_a\cdots F_{a+s(a,b)-1},
	\]
	whereas the right-hand side places \(F_{a\to b}\) in front of the cyclic
	product starting at \(b\). Since
	\(b\equiv a+s(a,b)\pmod M\), the two products coincide.

	For the resolvent identity, expand the resolvent into its Neumann series.
	A term contributing to the \((a,b)\)-block must have length
	\(s(a,b)+Mn\), with \(n\in\mathbb N_0\). The first \(s(a,b)\)
	elementary steps contribute \(F_{a\to b}\). Each subsequent full cycle,
	starting and ending at the \(b\)-th copy, contributes one power of
	\(T_{(b)}\). Therefore
	\[
	\Pi_a(I-\tau\mathscr C)^{-1}\Pi_b^{\top}
	=
	\sum_{n\ge0}\tau^{s(a,b)+Mn}F_{a\to b}(T_{(b)})^n,
	\]
	which gives \eqref{eq:Darboux-shifted-resolvent-right}. The equivalent form
	\eqref{eq:Darboux-shifted-resolvent-left} follows from the intertwining
	relation \eqref{eq:Darboux-partial-intertwining}. The diagonal case
	\(a=b\) gives \eqref{eq:Darboux-diagonal-resolvent-block}.
\end{proof}

\section{Proof of the sparse forms and block continued fractions}
\label{sec:height-ordering-block-CF}

The matrices \(\mathscr H\) and \(\mathscr H_{(a)}\) were defined
directly from the bidiagonal coefficients in
Subsection~\ref{subsec:direct-matrices-from-PBF}.  This section proves those
definitions have the asserted two-diagonal form and relates them to the
auxiliary cyclic matrices.  A coordinate is labeled by the factor being
applied and its row index.  The integer height assigned to that label shows
that every nonzero entry has displacement \(+Q\) or \(-P\).

After the two sparse forms have been proved, consecutive heights are grouped
into finite blocks.  This makes \(\mathscr H\) block tridiagonal, and ordinary
finite-dimensional Schur complementation produces the continued fraction for
its resolvent.

The first three subsections verify the direct placement rules for the
prescribed order and for \(\mathscr H_{(a)}\), and describe the polynomial
recurrences contained in the latter matrix.  The remaining subsections
construct the continued fraction from the prescribed factor order, then
establish its finite approximation, convergence, positivity, backward
evaluation, and residual estimates.  Notation without the cyclic subscript
\((a)\) refers to the prescribed order
\(L_1\cdots L_pU_q\cdots U_1\) and hence to the original factor coefficients.

\subsection{Verification of the direct coefficient placement}
\label{subsec:integer-height-rational-steps}

We verify the direct construction of
Subsection~\ref{subsec:direct-matrices-from-PBF}.  To retain the information
carried by the individual bidiagonal factors, one application of their
product is resolved into \(M=p+q\) elementary stages.  A coordinate is
therefore specified by the factor position \(a\in\{1,\ldots,M\}\) and the
row index \(n\in\mathbb N_0\); the pair \((a,n)\) is the coordinate label
used in the proof.

Recall the factor-position heights
\(\eta_0,\ldots,\eta_M\) from
Definition~\ref{def:factor-position-heights}.  Extend them to the coordinate
labels by defining
\[
	h(a,n)\coloneq Rn+\eta_{a-1}.
\]
Thus the factor position determines the residue of the height modulo \(R\),
whereas \(n\) determines its multiple of \(R\).  This choice makes every
nonzero entry of an individual bidiagonal factor produce one of the two
height displacements \(+Q\) and \(-P\).  For the prescribed two-run order,
\(\eta_{a-1}=\omega_a\); for a selected refactorization,
\(\eta_{a-1}=r_{a-1}\).

The height is integer-valued; the word
``rational'' below refers to the coprime pair of increments \(+Q\) and \(-P\),
or equivalently to their rational slope.  The type of a transition is determined by its
position in the bidiagonal factor, not by its numerical value.  The
\emph{normalizing entries} are the fixed diagonal entries of the lower
factors and the fixed superdiagonal entries of the upper factors; by
the factor normalization specified above, these entries are identically
equal to \(1\).  They will be represented as upward steps.  The
\emph{coefficient entries} are the subdiagonal entries
\(\ell_n^{(a)}\) of the lower factors and the diagonal entries
\(u_n^{(b)}\) of the upper factors; they will be represented as downward
steps carrying those coefficients.  In particular, if a coefficient
\(\ell_n^{(a)}\) or \(u_n^{(b)}\) happens to equal \(1\), it is still a
coefficient entry and still gives a downward step.

Minors of weighted path matrices admit a combinatorial reading through the
Lindstr{\"o}m--Gessel--Viennot lemma \cite{Lindstrom1973,GesselViennot1985};
the determinant factorizations below are instead obtained from Schur
complements, which produce directly the ordered pivot blocks required by the
backward recursion, and all path arguments here are one-particle expansions.

Throughout the path interpretation the convention is the same as in the
Hessenberg case of \cite{BranquinhoDiazFoulquieLimaManas2026JAT}: a matrix
entry is read in row-to-column orientation. Thus an entry \(A_{i,j}\) is
regarded as an edge from row state \(i\) to column state \(j\). This is a
path-reading convention only; all algebraic identities for the cyclic
linearization and for resolvents remain in the usual matrix notation. With
this convention, a subdiagonal entry of a lower bidiagonal factor gives a
downward edge, and a superdiagonal entry of an upper bidiagonal factor gives
an upward edge.

Recall from Subsection~\ref{subsec:balanced-Christoffel-word} that
\[
	d\coloneq \gcd(p,q),
	\qquad
	P\coloneq \frac{p}{d},
	\qquad
	Q\coloneq \frac{q}{d},
	\qquad
	R\coloneq P+Q=\frac{M}{d}.
\]
The division by \(d=\gcd(p,q)\) removes the common factor of the two
unscaled height increments
\[
+q,
\qquad
-p.
\]
Indeed,
\[
q=dQ,
\qquad
p=dP.
\]
Thus every edge in the reduced height graph has increment
\[
+Q,
\qquad
-P.
\]
When \(d>1\), the integer \(d\) appears as the multiplicity of factor labels
over the same reduced height.  Thus \(P\) and \(Q\)
determine the two reduced height steps, while \(d\) measures how many
factor labels may project to the same height residue.

For example, if \((p,q)=(4,2)\), then \(d=2\), \(P=2\), and \(Q=1\). The
unreduced steps \(+2\) and \(-4\) become
\[
+1,
\qquad
-2,
\]
but with two labels over each reduced residue class.

Use the ordered factors \(F_1,\ldots,F_M\) defined in
Subsection~\ref{subsec:positive-adjacent-refactorization}.
A coordinate state of \(\mathcal E^{\oplus M}\) is written
\[
e_{(a,n)}\coloneq\Pi_a^{\top}e_n,
\qquad
a\in\{1,\ldots,M\},
\quad
n\in\mathbb N_0.
\]
Thus the label \((a,n)\) records the summand associated with the factor
\(F_a\) and the coordinate \(n\) inside that summand.
The complete set of labeled coordinate states is denoted by
\[
\mathscr V_{p,q}
\coloneq
\{(a,n):a\in\{1,\ldots,M\},\ n\in\mathbb N_0\}.
\]
For any set \(\mathscr S\) of such labels, the notation
\(\operatorname{span}\mathscr S\) means
\[
\operatorname{span}\mathscr S
\coloneq
\operatorname{span}\{e_{(a,n)}:(a,n)\in\mathscr S\}.
\]

For the prescribed two-run order, write the state height as
\[
	h(a,n)=Rn+\omega_a,
\]
where \(\omega_a=h(a,0)\) is the height of the factor position \(a\).
The numbers \(\omega_a\) are determined by requiring the normalizing entries
to give rises of height \(+Q\).

For the lower factors, the normalizing entries are the diagonal entries of
\(L_1,\ldots,L_p\). They do not change \(n\). Hence the height shifts must
increase by \(Q\):
\[
\omega_1=0,
\qquad
\omega_a=(a-1)Q,
\qquad
a\in\{1,\ldots,p+1\}.
\]
For the upper factors, the normalizing entries are the superdiagonal entries
of \(U_q,\ldots,U_1\). In the row-to-column path convention, such an
entry gives an edge from index \(n\) to index \(n+1\). Hence the term
\(Rn\) increases by \(R=P+Q\), and the height shift must decrease by
\(P\) in order to make the total increment equal to \(+Q\). Thus, for
\(b\in\{1,\ldots,q\}\),
\begin{equation}
	\label{eq:general-height-shifts}
	\omega_{p+b}
	=
	(q-b+1)P.
\end{equation}
For the common endpoint \(a=p+1\), the two descriptions give the same value:
from the lower side \(\omega_{p+1}=pQ\), while from the upper side
\(\omega_{p+1}=qP\). These coincide because \(pQ=qP=pq/d\).

The pairs
\[
(h(a,n);a)
\]
are height states with an additional factor label. The label \(a\) is kept
because it identifies the bidiagonal factor containing the coefficient of the
transition. If the height alone determines \(a\), the label can be suppressed.
If not, it must be retained.

\begin{proposition}[The two possible height increments]
	\label{prop:general-height-steps}
	After relabeling the cyclic linearization by the height states
	\((h(a,n);a)\), every nonzero elementary transition has height increment
	either \(+Q\) or \(-P\). More precisely:
	\begin{enumerate}[label=\textnormal{(\roman*)}]
		\item for \(a\in\{1,\ldots,p\}\), the fixed diagonal entry of \(L_a\)
		gives a transition of height \(+Q\), while the subdiagonal entry
		\(\ell_n^{(a)}\) gives a coefficient transition of height \(-P\);
		\item for \(b\in\{1,\ldots,q\}\), the factor in cyclic position
		\(p+b\) is \(U_{q-b+1}\); its fixed superdiagonal entry gives a
		transition of height \(+Q\), while its diagonal entry
		\(u_n^{(q-b+1)}\) gives a coefficient transition of height \(-P\).
	\end{enumerate}
	The classification is positional: statements (i) and (ii) remain unchanged
	when one of the coefficients \(\ell_n^{(a)}\) or \(u_n^{(b)}\) has the
	numerical value \(1\).
\end{proposition}

\begin{proof}
	For the lower factors, the diagonal entry of \(L_a\) gives an edge from
	index \(n\) at \(F_a\) to index \(n\) at \(F_{a+1}\). Therefore \(n\) is unchanged,
	and the height increment is
	\[
	\omega_{a+1}-\omega_a=Q.
	\]
	The subdiagonal entry of \(L_a\) gives an edge from row index \(n\) to
	column index \(n-1\). Thus the term
	\(Rn\) decreases by \(R\), while the height shift increases by \(Q\). The
	total height increment is
	\[
	Q-R=Q-(P+Q)=-P.
	\]

	For the upper factors, the fixed superdiagonal entry gives an edge from row
	index \(n\) to column index \(n+1\). The term \(Rn\) therefore increases by \(R=P+Q\), while the
	height shift decreases by \(P\). The total height increment is
	\[
	R-P=Q.
	\]
	The diagonal entry \(u_n^{(b)}\) preserves \(e_n\). Hence only the height
	shift changes, and it changes by \(-P\). Therefore the diagonal entry gives
	a coefficient transition of height \(-P\).
\end{proof}

Proposition~\ref{prop:general-height-steps} determines the two possible height increments.
The next question is whether the height alone identifies the factor that
contains an edge weight.  This is true in the coprime case, apart from
finitely many missing heights near the origin, but factor labels remain necessary when
\(\gcd(p,q)>1\).

\begin{proposition}[Projection to the reduced height graph]
	\label{prop:general-coprime-projection-boundary}
	Assume first that \(d=1\). Then \(P=p\), \(Q=q\), and \(R=M=p+q\). The
	height shifts
	\[
	\omega_1,\ldots,\omega_M
	\]
	are pairwise distinct modulo \(M\). Hence, whenever a height belongs to the
	image of the map \((a,n)\mapsto h(a,n)\), that height determines the factor
	label \(a\). In this case the labeled height graph projects to a scalar
	rational path graph with steps
	\[
	+q,
	\qquad
	-p,
	\]
	apart from finitely many missing heights near the origin.

	If \(d>1\), different factor labels may give the same height modulo \(R\).
	Then the height alone does not determine the current factor, and the labels
	must be retained.
\end{proposition}

\begin{proof}
	When \(d=1\), formula~\eqref{eq:general-height-shifts} gives, for every
	factor position \(a\in\{1,\ldots,M\}\),
	\[
	\omega_a\equiv(a-1)q\pmod M.
	\]
	For \(1\le a\le p+1\), this follows from
	\(\omega_a=(a-1)q\).  For an upper-factor position \(a=p+b\),
	formula~\eqref{eq:general-height-shifts} gives
	\[
	(p+b-1)q-\omega_{p+b}
	=(p+b-1)q-(q-b+1)p
	=(b-1)M,
	\]
	which proves the same congruence.  Moreover,
	\(\gcd(q,M)=\gcd(q,p+q)=\gcd(p,q)=1\).  Multiplication by \(q\)
	therefore permutes the residue classes modulo \(M\), so the \(M\) shifts
	\(\omega_1,\ldots,\omega_M\) are pairwise distinct modulo \(M\). Thus two states with
	the same height must have the same factor label and the same index \(n\).

	However, the restriction \(n\ge0\) may remove some small heights from the
	image. For example, height \(1\) is missing when
	\((p,q)=(2,3)\). Away from this finite
	set of missing initial heights, the projected height graph is the usual rational-step graph.

	If \(d>1\), there are \(M=dR\) factor labels but only \(R\) residue classes
	modulo \(R\). Thus different labels may have the same height residue, and
	the label is necessary to know which bidiagonal coefficient is attached to
	the transition.
\end{proof}

The preceding congruence constructs the prescribed-order matrix directly
from the factor coefficients.

\begin{corollary}[Direct placement of the prescribed factor coefficients]
	\label{cor:direct-prescribed-coefficient-placement}
	Assume that \(d=1\), so that \(M=p+q\).  For every height \(h\), let
	\(k(h)\in\{0,\ldots,M-1\}\) be the unique integer satisfying
	\[
		k(h)q\equiv h\pmod M,
	\]
	and put
	\[
		a(h)\coloneq k(h)+1,
		\qquad
		n(h)\coloneq\frac{h-\omega_{a(h)}}{M}.
	\]
	The height \(h\) occurs in the prescribed ordering precisely when
	\(n(h)\ge0\).  At every such height there is an entry equal to \(1\)
	joining \(h\) to \(h+q\).  When the indicated coefficient exists, the
	entry joining \(h\) to \(h-p\) is
	\begin{equation*}
		c_h^{\mathrm{pre}}
		\coloneq
		\begin{cases}
			\ell_{n(h)-1}^{(a(h))},
			& a(h)\in\{1,\ldots,p\},\quad n(h)\ge1,\\[2pt]
			u_{n(h)}^{(q-b+1)},
			& a(h)=p+b,\quad b\in\{1,\ldots,q\}.
		\end{cases}
	\end{equation*}
	Consequently, after omitting the finitely many heights for which
	\(n(h)<0\), the two nonzero diagonals of \(\mathscr H\) are obtained by
	reading the original factor coefficients in the periodic order determined
	by \(k(h)q\equiv h\pmod M\).  In particular, every height
	\(h\ge p q\) occurs.
\end{corollary}

\begin{proof}
	Proposition~\ref{prop:general-coprime-projection-boundary} gives
	\(\omega_a\equiv(a-1)q\pmod M\), so the congruence selects the unique
	factor position whose states can have height \(h\).  Solving
	\(h=Mn+\omega_a\) gives \(n=n(h)\).  A lower factor has its coefficient
	entry in row \(n\), column \(n-1\), with value
	\(\ell_{n-1}^{(a)}\); an upper factor in position \(p+b\) has its
	coefficient entry on row and column \(n\), with value
	\(u_n^{(q-b+1)}\).  Proposition~\ref{prop:general-height-steps} gives the
	two target heights.  Finally,
	\(\max_a\omega_a=pq\), so \(n(h)\ge0\) for every selected factor when
	\(h\ge pq\).
\end{proof}

For example, during one period the factors selected by increasing height
are
\[
\begin{array}{c|ccccc}
	(p,q) & h\equiv0 & h\equiv1 & h\equiv2 & h\equiv3 & h\equiv4\\
	\hline
	(3,2) & L_1 & U_2 & L_2 & U_1 & L_3\\
	(2,3) & L_1 & U_3 & U_1 & L_2 & U_2
\end{array}
\qquad (\bmod\ 5).
\]
The formula for \(n(h)\) automatically removes an occurrence that would
require a negative row index; for instance, height \(1\) is absent in the
initial part of both examples.  After this finite initial part, the same
five-factor order repeats with the next coefficient of each factor.

When \(d>1\), the same rule is block-valued.  For a height residue
\(\rho\in\{0,\ldots,R-1\}\), the congruence
\[
	(a-1)Q\equiv\rho\pmod R
\]
selects exactly \(d\) factor positions.  For each selected position one
solves \(h=Rn+\omega_a\), discards negative values of \(n\), and places the
corresponding bidiagonal coefficient in the block from height \(h\) to
height \(h-P\); the normalizing entries fill the block from \(h\) to
\(h+Q\).  Thus the block matrix is also obtained directly from the periodic
ordering of the original coefficients.  The coordinate-state notation
keeps track of the \(d\) rows and columns within each block.

For the continued fraction it is useful to group several consecutive heights
into one block.  Before doing so, the exact-height ordering reveals a more
rigid recurrence structure that will also be used in the spectral
interpretation.  For \(h\in\mathbb N_0\), put
\[
\mathscr F_h
\coloneq
\{(a,n):h(a,n)=h\},
\qquad
\mathcal E_h
\coloneq
\operatorname{span}\mathscr F_h.
\]
Empty spaces \(\mathcal E_h\) are allowed in the finite initial part of the
decomposition.

\begin{definition}[Height ordering used in the proof]
	\label{def:height-ordered-cyclic-operator}
	Recall that \(e_{(a,n)}=\Pi_a^{\top}e_n\) denotes the coordinate vector
	corresponding to coordinate \(n\) in factor copy \(a\).  Enumerate these
	vectors as
	\[
	e_{(a_0,n_0)},e_{(a_1,n_1)},e_{(a_2,n_2)},\ldots
	\]
	by requiring
	\[
	h(a_j,n_j)<h(a_k,n_k)
	\quad\Longrightarrow\quad
	j<k,
	\]
	and, when \(h(a_j,n_j)=h(a_k,n_k)\), by placing the smaller factor
	label first.
	
	Let \(\mathscr P_h\) be the permutation matrix determined by
	\[
	\mathscr P_h e_{(a_j,n_j)}=e_j,
	\qquad j\in\mathbb N_0.
	\]
\end{definition}

\begin{proposition}[Agreement with the auxiliary cyclic matrix]
	\label{prop:direct-cyclic-agreement}
	The matrix \(\mathscr H_{p,q}\) defined directly in
	Subsection~\ref{subsec:direct-matrices-from-PBF} satisfies
	\[
	\mathscr H_{p,q}
	=
	\mathscr P_h\mathscr C_{p,q}\mathscr P_h^{\top}.
	\]
\end{proposition}

\begin{proof}
	Fix a coordinate label \((a,n)\), and write
	\[
		h=Rn+\omega_a.
	\]
	The nonzero entries in the corresponding row of \(\mathscr C_{p,q}\)
	are precisely the nonzero entries in row \(n\) of the factor \(F_a\), and
	their column labels belong to the next factor position \(a^+\).

	Suppose first that \(a\in\{1,\ldots,p\}\), so that \(F_a=L_a\).
	The diagonal entry of \(L_a\) is equal to \(1\).  It joins \((a,n)\) to
	\((a^+,n)\), whose height is \(h+Q\).  If \(n\ge1\), the other nonzero
	entry in row \(n\) is
	\[
		(L_a)_{n,n-1}=\ell_{n-1}^{(a)}.
	\]
	It joins \((a,n)\) to \((a^+,n-1)\), whose height is \(h-P\).
	Since \(n=n_a(h)\), these are exactly the entries in
	\eqref{eq:direct-prescribed-matrix-definition} and
	\eqref{eq:direct-prescribed-lower-entry}.

	Now suppose that \(a=p+b\), so that \(F_a=U_{q-b+1}\).  The fixed
	superdiagonal entry of this factor is equal to \(1\).  It joins \((a,n)\)
	to \((a^+,n+1)\), whose height is again \(h+Q\).  Its diagonal entry is
	\[
		(U_{q-b+1})_{n,n}=u_n^{(q-b+1)}.
	\]
	It joins \((a,n)\) to \((a^+,n)\), whose height is \(h-P\).  Since
	\(n=n_{p+b}(h)\), these are exactly the entries in
	\eqref{eq:direct-prescribed-matrix-definition} and
	\eqref{eq:direct-prescribed-upper-entry}.

	There are no further entries because every factor is bidiagonal.  The
	permutation \(\mathscr P_h\) merely places the labels \((a,n)\) in the
	order used in the direct definition of \(\mathscr H_{p,q}\).  Therefore
	\(\mathscr P_h\mathscr C_{p,q}\mathscr P_h^{\top}\) and
	\(\mathscr H_{p,q}\) agree entry by entry.
\end{proof}

\begin{theorem}[Two-diagonal form of the prescribed matrix]
	\label{thm:exact-height-double-band}
	Group together the coordinate labels having the same height.  Relative to
	the resulting sequence of coordinate
	spaces
	\[
	\mathcal E_0,\mathcal E_1,\mathcal E_2,\ldots,
	\]
	the directly defined matrix \(\mathscr H\) has only two nonzero block
	diagonals:
	\begin{equation}
		\label{eq:exact-height-double-band}
		[\mathscr H]_{\mathcal E_h\times\mathcal E_k}=0
		\qquad\textnormal{unless}\qquad
		k-h\in\{Q,-P\}.
	\end{equation}
	The block with \(k=h+Q\) consists of normalizing entries, whereas the
	block with \(k=h-P\) consists of the individual coefficient entries of
	the bidiagonal factors.  Moreover,
	\[
	0\le\dim\mathcal E_h\le d,
	\]
	and \(\dim\mathcal E_h=d\) whenever
	\(h\ge\max_{1\le a\le M}\omega_a=pQ\).

	If \(d=1\), every sufficiently large exact height carries one state, and
	the scalar tail of \(\mathscr H\) is a two-diagonal band matrix: its only
	nonzero diagonals are the \(Q\)-th superdiagonal and the \(P\)-th
	subdiagonal.  In particular, for \(q=1\) this is the two-diagonal
	lower-Hessenberg matrix associated with the usual star-like
	\(p\)-orthogonality.  For arbitrary \(p\) and \(q\),
	\eqref{eq:exact-height-double-band} is its two-sided block analogue.
\end{theorem}

\begin{proof}
	Proposition~\ref{prop:general-height-steps} states that every nonzero
	transition has increment \(+Q\) or \(-P\), and it also identifies which
	type of factor entry produces each increment.  This proves
	\eqref{eq:exact-height-double-band} and the assertion about its two block
	diagonals.

	For \(a\le p+1\), the definition gives
	\(\omega_a=(a-1)Q\).  If \(a=p+b\), with
	\(b\in\{1,\ldots,q\}\), then the upper-factor formula for the shift gives
	\[
	(p+b-1)Q-\omega_{p+b}=(b-1)(P+Q)=(b-1)R.
	\]
	Consequently,
	\(\omega_a\equiv(a-1)Q\pmod R\) for every factor label \(a\).  Since
	\(\gcd(Q,R)=1\) and \(M=dR\), every residue class modulo \(R\) occurs
	for exactly \(d\) factor labels.  For fixed \(h\), each such label
	determines at most one index,
	\[
	n=\frac{h-\omega_a}{R}.
	\]
	Hence \(\dim\mathcal E_h\le d\).  If
	\(h\ge\max_a\omega_a\), all the resulting indices are nonnegative, so
	all \(d\) labels occur and \(\dim\mathcal E_h=d\).  When \(d=1\), the
	exact-height blocks are therefore scalar after the finite initial part,
	and the two increments become the stated scalar diagonal offsets.
\end{proof}

The rotational symmetry of this two-diagonal matrix is already present
before height ordering.  Let \(\zeta\) be an \(M\)-th root of unity and set
\[
\mathscr D_\zeta
\coloneq
\operatorname{diag}
\bigl(I,\zeta I,\ldots,\zeta^{M-1}I\bigr)
\]
on the \(M\) factor copies.
Write
\[
\widehat{\mathscr D}_\zeta
\coloneq
\mathscr P_h\mathscr D_\zeta\mathscr P_h^{\top}
\]
for the same diagonal transformation in the ordering by height.

\begin{proposition}[Rotational symmetry and cyclic products]
	\label{prop:cyclic-star-symmetry}
	Let \(\zeta\) be an \(M\)-th root of unity.  Then
	\begin{equation}
		\label{eq:cyclic-star-similarity}
		\mathscr D_\zeta^{-1}\mathscr C\mathscr D_\zeta
		=
		\zeta\mathscr C,
		\qquad
		\widehat{\mathscr D}_\zeta^{-1}\mathscr H
		\widehat{\mathscr D}_\zeta
		=
		\zeta\mathscr H,
	\end{equation}
	Moreover,
	\begin{equation}
		\label{eq:height-power-cyclic-products}
		\mathscr H^M
		=
		\mathscr P_h
		\operatorname{diag}
		\bigl(T_{(1)},\ldots,T_{(M)}\bigr)
		\mathscr P_h^{\top}.
	\end{equation}
	
	Consequently, every finite principal submatrix of \(\mathscr H\) has
	spectrum invariant under multiplication by \(\zeta\).  In particular,
	this holds for the leading principal submatrices obtained by retaining
	the first finitely many height blocks.  Equation
	\eqref{eq:height-power-cyclic-products} expresses \(\mathscr H^M\), up
	to the permutation \(\mathscr P_h\), as the direct sum of the cyclic
	products \(T_{(1)},\ldots,T_{(M)}\).
\end{proposition}

\begin{proof}
	In the cyclic matrix, the block in row \(a\) and column \(a+1\) is
	\(F_a\), with the factor-copy indices read modulo \(M\).  Conjugation by
	\(\mathscr D_\zeta\) multiplies each of these blocks by \(\zeta\).  For
	the block joining the last copy to the first, this follows from
	\(\zeta^M=1\).  Hence
	\[
	\mathscr D_\zeta^{-1}\mathscr C\mathscr D_\zeta
	=
	\zeta\mathscr C.
	\]
	Conjugating this identity by \(\mathscr P_h\) gives the second identity
	in \eqref{eq:cyclic-star-similarity}.
	
	Lemma~\ref{lem:general-cyclic-linearization} gives
	\[
	\mathscr C^M
	=
	\operatorname{diag}
	\bigl(T_{(1)},\ldots,T_{(M)}\bigr).
	\]
	Conjugation by \(\mathscr P_h\) proves
	\eqref{eq:height-power-cyclic-products}.
	
	Finally, \(\widehat{\mathscr D}_\zeta\) is diagonal in the coordinate
	ordering by height.  If \(S\) is any finite set of coordinate
	indices, restriction of the second identity in
	\eqref{eq:cyclic-star-similarity} gives
	\[
	\left(\widehat{\mathscr D}_\zeta|_S\right)^{-1}
	\mathscr H_{S\times S}
	\left(\widehat{\mathscr D}_\zeta|_S\right)
	=
	\zeta\,\mathscr H_{S\times S}.
	\]
	Thus \(\mathscr H_{S\times S}\) and
	\(\zeta\mathscr H_{S\times S}\) are similar, which proves the asserted
	spectral invariance.
\end{proof}

This symmetry has a discrete Fourier interpretation.  Put
\[
	\zeta_M\coloneq\exp\left(\frac{2\pi\mathrm i}{M}\right).
\]
The diagonal matrices \(\mathscr D_{\zeta_M^j}\),
\(j\in\{0,\ldots,M-1\}\), implement the characters of the cyclic group
\(\mathbb Z/M\mathbb Z\) on the factor-copy index.  The discrete Fourier
transform on that index converts this diagonal character action into cyclic
translation.  Accordingly,
\[
	\frac1M\sum_{j=0}^{M-1}\zeta_M^{j(n-s)}
	=
	\begin{cases}
		1,&n\equiv s\pmod M,\\
		0,&n\not\equiv s\pmod M,
	\end{cases}
\]
is the Fourier projection onto the congruence class \(s\).  For the
prescribed-order Weyl matrix, the zero Fourier mode produces the rotational
lift in Theorem~\ref{thm:cyclic-star-Weyl}.  The full star moment system of
\(\mathscr H_{(a)}\) assembles all congruence classes.  Its component
corresponding to \(T_{(a)}\) is the transformed system
\(\mathrm d\Psi^{(a)}\); for \(a=p+1\), this is
\(x\,\mathrm d\Psi\).  The remaining classes give the Fourier components
studied in Proposition~\ref{prop:balanced-character-measure-lift}.  Thus the rotational
spectral symmetry, the decomposition by powers modulo \(M\), and the
Fourier components of the measures on the star are different forms of the same discrete
Fourier decomposition.

For positive scalar double-band matrices, rotationally arranged spectra are
also obtained in \cite{McMillen2008}.  In the Hessenberg case \(q=1\), the
recurrence
\[
\lambda P_h(\lambda)
=P_{h+1}(\lambda)+c_hP_{h-p}(\lambda)
\]
and its orthogonality on a star-like set are classical; see
\cite{DelvauxLopezGarcia2013,Lima2026Symmetric}.
Theorem~\ref{thm:exact-height-double-band} shows that a mixed-type bidiagonal
factorization produces the two-sided version, with \(Q\) forward levels and
\(P\) backward levels and with all intervening diagonals equal to zero.  The
continued fraction below is obtained by grouping this exact-height matrix
into neighboring finite blocks, which provides a Schur-complement evaluation
of its resolvent.

\begin{remark}[Terminology and rational Dyck paths]
	\label{rem:mixed-rational-path-terminology}
	The term ``rational Dyck path'' is avoided because it already denotes
	other path families in rational Catalan combinatorics
	\cite{GorskyMazinVazirani2017,DaiFuQiu2026,
	BarcucciBerniniBilottaPinzani2025}.  Here the paths are simply weighted
	paths with increments \(+Q\) and \(-P\), and the factor label is retained
	to identify the bidiagonal factor supplying each coefficient.  When
	\(q=1\), these increments become \(+1\) and \(-p\), and the construction
	reduces to the standard weighted \(p\)-Dyck path model for the
	\(p\)-branched Stieltjes fraction
	\cite{Flajolet1980,Viennot1983,BanderierFlajolet2002}.
\end{remark}

\subsection{Verification of the two-block-diagonal form of
\texorpdfstring{\(\mathscr H_{(a)}\)}{H(a)}}
\label{subsec:balanced-exact-height}

We now verify the direct construction of
\(\mathscr H_{(a)}\).  Let
\[
	T_{(a)}
	=\widetilde F^{(a)}_1\cdots\widetilde F^{(a)}_M
\]
be the factorization supplied by
Theorem~\ref{thm:positive-balanced-refactorization}, and let
\(\mathscr C_{(a)}\) be the cyclic matrix having
\(\widetilde F^{(a)}_k\) in block position \((k,k+1)\), with the block
indices read cyclically.  In the refactorized factor arrangement, read the
remainders \(r_k\) periodically modulo \(R\).
The height reached before the factor in position \(j\) is then
\(r_{j-1}\).  Accordingly, the height of its coordinate state is
\[
	h^{(a)}(j,n)\coloneq Rn+r_{j-1}.
\]
This is the definition \(h(j,n)=Rn+h_{j-1}\) from
Subsection~\ref{subsec:integer-height-rational-steps}, specialized to the
minimum-height factorization.

For \(h\ge0\), let \(\mathcal E_h^{(a)}\) be the span of the coordinate
labels having height \(h\) in this order.  These are precisely the rows
and columns indexed by \(I_h^{(a)}\) in the direct definition.
Ordering the auxiliary cyclic matrix by these spaces gives the directly
defined matrix \(\mathscr H_{(a)}\), by the same entry-by-entry
verification as in Proposition~\ref{prop:direct-cyclic-agreement}.
Recall that a monomial matrix is
a square matrix with exactly one nonzero entry in each row and each column;
it is positive when all those entries are positive.

\begin{theorem}[The two block diagonals of \(\mathscr H_{(a)}\)]
	\label{thm:balanced-two-extreme-diagonals}
	For every \(h\in\mathbb N_0\), there are exactly \(d\) coordinate states
	of height \(h\); equivalently,
	\[
	\dim\mathcal E_h^{(a)}=d.
	\]
	
	With the row-to-column convention fixed above, the only possible
	transitions from height \(h\) are the following:
	\begin{enumerate}[label=\textnormal{(\roman*)}]
		\item Every state of height \(h\) has one transition of weight \(1\)
		to height \(h+Q\).  These transitions give a bijection between the
		\(d\) states of heights \(h\) and \(h+Q\).  Consequently, the block
		\[
		\left[
		\mathscr H_{(a)}
		\right]_{
			\mathcal E_h^{(a)}
			\times
			\mathcal E_{h+Q}^{(a)}
		}
		\]
		has exactly one nonzero entry in every row and every column, and all
		these entries are equal to \(1\).
		
		\item If \(h\ge P\), every state of height \(h\) has one transition
		with a strictly positive factor coefficient to height \(h-P\).
		These transitions give a bijection between the \(d\) states of
		heights \(h\) and \(h-P\).  Consequently, the block
		\[
		\left[
		\mathscr H_{(a)}
		\right]_{
			\mathcal E_h^{(a)}
			\times
			\mathcal E_{h-P}^{(a)}
		}
		\]
		has exactly one strictly positive entry in every row and every
		column.  If \(h<P\), no such transition exists.
		
		\item All other blocks are zero.
	\end{enumerate}
\end{theorem}

\begin{proof}
	Write the height uniquely as
	\[
	h=Rn+\rho,
	\qquad
	n\in\mathbb N_0,
	\qquad
	\rho\in\{0,\ldots,R-1\}.
	\]
	In this factor arrangement, every residue
	\(\rho\in\{0,\ldots,R-1\}\) occurs exactly \(d\) times among the
	\(M=dR\) factor positions.  Therefore there are exactly \(d\) positions
	\(k\) for which
	\[
	r_{k-1}=\rho.
	\]
	For each of them, the coordinate state \((k,n)\) has height
	\[
	Rn+r_{k-1}=Rn+\rho=h.
	\]
	Conversely, every coordinate state of height \(h\) must be obtained in
	this way.  Hence
	\(\dim\mathcal E_h^{(a)}=d\).
	
	We now examine the two nonzero entries in the row of each bidiagonal
	factor.  The index \(k+1\) is read cyclically when \(k=M\).
	
	If \(\widetilde F_k^{(a)}\) is lower bidiagonal, then
	\(r_{k-1}<P\) and
	\[
	r_k=r_{k-1}+Q.
	\]
	Its fixed diagonal entry, which is equal to \(1\), gives
	\[
	(k,n)\longmapsto(k+1,n).
	\]
	The height of the target is
	\[
	Rn+r_k
	=
	Rn+r_{k-1}+Q
	=
	h+Q.
	\]
	Its strictly positive subdiagonal coefficient gives, when \(n\ge1\),
	\[
	(k,n)\longmapsto(k+1,n-1),
	\]
	whose target has height
	\[
	R(n-1)+r_k
	=
	h+Q-R
	=
	h-P.
	\]
	
	If \(\widetilde F_k^{(a)}\) is upper bidiagonal, then
	\(r_{k-1}\ge P\) and
	\[
	r_k=r_{k-1}+Q-R=r_{k-1}-P.
	\]
	Its fixed superdiagonal entry, which is equal to \(1\), gives
	\[
	(k,n)\longmapsto(k+1,n+1),
	\]
	and the target has height
	\[
	R(n+1)+r_k
	=
	h+Q.
	\]
	Its strictly positive diagonal coefficient gives
	\[
	(k,n)\longmapsto(k+1,n),
	\]
	whose target has height
	\[
	Rn+r_k=h-P.
	\]
	
	Thus every state has exactly one weight-\(1\) transition to height
	\(h+Q\).  A transition to height \(h-P\) exists for every state exactly
	when \(h\ge P\): at a lower-factor position this is equivalent to
	\(n\ge1\), while every state belonging to an upper-factor position
	already has height at least \(P\).
	
	Finally, advancing the factor position \(k\) cyclically is bijective.
	Hence the transitions just described match the \(d\) source states
	bijectively with the \(d\) target states.  Since a bidiagonal factor has
	no other nonzero entries, no other height transitions occur.
\end{proof}

\begin{corollary}[Coprime double-band matrix]
	\label{cor:balanced-coprime-double-band}
	If \(d=1\), the height determines the factor position and
	\[
		(\mathscr H_{(a)})_{h,h+q}=1
		\quad(h\ge0),
		\qquad
		(\mathscr H_{(a)})_{h,h-p}=c_h^{(a)}>0
		\quad(h\ge p),
	\]
	with every other entry equal to zero.
\end{corollary}

\begin{proof}
	When \(d=1\), one has \(P=p\), \(Q=q\), and
	\(\dim\mathcal E_h^{(a)}=1\) for every \(h\ge0\).  Hence the
	permutation block from \(\mathcal E_h^{(a)}\) to
	\(\mathcal E_{h+Q}^{(a)}\) is the scalar \(1\), while the
	positive monomial block from \(\mathcal E_h^{(a)}\) to
	\(\mathcal E_{h-P}^{(a)}\) is a positive scalar, denoted by
	\(c_h^{(a)}\).  The asserted entries and the vanishing of all the others follow
	from Theorem~\ref{thm:balanced-two-extreme-diagonals}.
\end{proof}

\begin{corollary}[Two-block-diagonal recurrences and continued fractions along the constructed cyclic family]
	\label{cor:admissible-cyclic-star-CF-family}
	For every \(a\in\mathcal A_{\mathrm{cyc}}\), let
	\(\mathscr H_{(a)}\) be the height-ordered matrix obtained
	from the refactorization of \(T_{(a)}\) in
	Theorem~\ref{thm:positive-balanced-refactorization}.  Then
	\(\mathscr H_{(a)}\) has only the two block displacements
	\(+Q\) and \(-P\), as described in
	Theorem~\ref{thm:balanced-two-extreme-diagonals}.

	Let \(\mathsf R_a\) denote the coordinate projection onto all coordinates
	belonging to the first factor position in this refactorization.  Then,
	as a formal power-series identity,
	\begin{equation}
		\label{eq:admissible-cyclic-star-resolvent}
		\mathsf R_a
		(I-\tau\mathscr H_{(a)})^{-1}
		\mathsf R_a^{\top}
		=
		(I-\tau^M T_{(a)})^{-1}.
	\end{equation}
	If the matrices involved define bounded operators, the same identity
	holds on the common operator-resolvent domain.

	Consequently, grouping consecutive heights of
	\(\mathscr H_{(a)}\) into finite blocks and applying Schur
	complementation gives a factor-resolved height-block continued fraction
	whose \(\mathsf R_a\)-projection recovers the resolvent of
	\(T_{(a)}\) through \eqref{eq:admissible-cyclic-star-resolvent}.
	Whenever \(T_{(a)}\) is supplied with mixed-type Favard data, the
	corresponding rectangular projection gives simultaneous rational
	approximants to its Weyl matrix.
\end{corollary}

\begin{proof}
	Theorem~\ref{thm:balanced-two-extreme-diagonals} gives the two stated
	block displacements.  Before the height permutation, the matrix
	\(\mathscr H_{(a)}\) is the cyclic linearization of the minimum-height
	factorization of \(T_{(a)}\).  Applying
	Lemma~\ref{lem:general-cyclic-linearization} to those factors and then
	conjugating by the height permutation gives
	\eqref{eq:admissible-cyclic-star-resolvent}.  The block continued fraction
	now follows from Theorem~\ref{thm:matrix-CF-height-blocks}, and the Weyl
	statement follows from the rectangular projection in
	Proposition~\ref{prop:general-rectangular-projection}, applied to the same
	cyclic product \(T_{(a)}\) with its Favard normalization matrices.
\end{proof}

The entries \(c_h^{(a)}\) in the preceding corollary can be read directly from
the refactorized bidiagonal matrices.

\begin{corollary}[Direct placement in \(\mathscr H_{(a)}\)]
	\label{cor:direct-balanced-coefficient-placement}
	Write
	\(h=Rn+\rho\), with
	\(\rho\in\{0,\ldots,R-1\}\), and let
	\(k(\rho)\in\{0,\ldots,R-1\}\) be determined by
	\[
		k(\rho)Q\equiv\rho\pmod R.
	\]
	The factor positions associated with height \(h\) are
	\[
		a_j(h)\coloneq k(\rho)+1+jR,
		\qquad j\in\{0,\ldots,d-1\}.
	\]
	For each of these positions, the coefficient placed in the block from
	height \(h\) to height \(h-P\) is
	\begin{equation}
		\label{eq:direct-balanced-coefficient-placement}
		\begin{cases}
			(\widetilde F^{(a)}_{a_j(h)})_{n,n-1},
			& \widetilde F^{(a)}_{a_j(h)}\text{ is lower bidiagonal},\\[2pt]
			(\widetilde F^{(a)}_{a_j(h)})_{n,n},
			& \widetilde F^{(a)}_{a_j(h)}\text{ is upper bidiagonal}.
		\end{cases}
	\end{equation}
	The first case occurs only when \(n\ge1\); together with the second case,
	it supplies all \(d\) entries of the positive monomial block whenever
	\(h\ge P\).  The block from height \(h\) to height \(h+Q\) is the
	permutation block formed by the normalizing entries.  Hence
	\(\mathscr H_{(a)}\) is obtained directly from the entries of
	\(\widetilde F^{(a)}_1,\ldots,\widetilde F^{(a)}_M\) by this periodic
	modular rule.
\end{corollary}

\begin{proof}
	The congruence selects exactly the positions for which
	\(r_{a-1}=\rho\).  At each such position the coordinate index is the
	quotient \(n\) in \(h=Rn+\rho\).  The two alternatives in
	\eqref{eq:direct-balanced-coefficient-placement} are precisely the
	coefficient entries of a lower or upper bidiagonal factor.  Their target
	height is \(h-P\), while the normalizing entry has target height \(h+Q\),
	by Theorem~\ref{thm:balanced-two-extreme-diagonals}.
\end{proof}

In particular, when \(d=1\), one reads a single factor at each height and
the scalar coefficient in Corollary~\ref{cor:balanced-coprime-double-band}
is
\[
	c_h^{(a)}=
	\begin{cases}
		(\widetilde F^{(a)}_{a(h)})_{n,n-1},
		& \widetilde F^{(a)}_{a(h)}\text{ is lower bidiagonal},\\[2pt]
		(\widetilde F^{(a)}_{a(h)})_{n,n},
		& \widetilde F^{(a)}_{a(h)}\text{ is upper bidiagonal},
	\end{cases}
	\qquad
	h=(p+q)n+\rho,
\]
where \(a(h)=k(\rho)+1\).  Definition~\ref{def:direct-sparse-matrices}
determines the factor position
\(a(h)\) directly from the residue of \(h\).  Thus, once the adjacent
refactorization has produced the factors in the star arrangement, the
entries of \(\mathscr H_{(a)}\) are read from those factors by the displayed
modular rule.

The adjacent exchanges preserve \(T_{(a)}\), although they generally change
its individual bidiagonal factors.  For \(a\in\{1,\ldots,M\}\), put
\[
	C_a\coloneq F_1\cdots F_{a-1},
	\qquad
	D_a\coloneq F_a\cdots F_M.
\]
Then
\[
	T=C_aD_a,
	\qquad
	T_{(a)}=D_aC_a.
\]
All these products are banded and are therefore well defined entry by entry.

\begin{theorem}[Resolvent transport along the cyclic family]
	\label{prop:balanced-ABBA-resolvent-dressing}
	For every \(a\in\{1,\ldots,M\}\), the following identities hold as formal
	power series in \(t\):
	\begin{align}
		(I-tT)^{-1}
		&=I+tC_a(I-tT_{(a)})^{-1}D_a,
		\label{eq:balanced-ABBA-resolvent-dressing}\\
		(I-tT_{(a)})^{-1}
		&=I+tD_a(I-tT)^{-1}C_a,
		\label{eq:cyclic-resolvent-reverse-dressing}\\
		(I-tT)^{-1}C_a
		&=C_a(I-tT_{(a)})^{-1},
		\label{eq:cyclic-resolvent-intertwining-right}\\
		D_a(I-tT)^{-1}
		&=(I-tT_{(a)})^{-1}D_a.
		\label{eq:cyclic-resolvent-intertwining-left}
	\end{align}
	
	For \(a\in\mathcal A_{\mathrm{cyc}}\), let \(\mathsf R_a\) be the
	coordinate projection in
	\eqref{eq:admissible-cyclic-star-resolvent}.  With \(t=\tau^M\),
	\begin{equation}
		\label{eq:original-resolvent-through-admissible-star-matrix}
		(I-\tau^M T)^{-1}
		=
		I+\tau^M C_a\mathsf R_a
		(I-\tau\mathscr H_{(a)})^{-1}
		\mathsf R_a^{\top}D_a.
	\end{equation}
	Consequently,
	\begin{align}
		\label{eq:original-rectangular-resolvent-through-admissible-star-matrix}
		&\mathsf E_{[q]}(I-\tau^M T)^{-1}\mathsf E_{[p]}^{\top}
	=
		\mathsf E_{[q]}\mathsf E_{[p]}^{\top}
		+\tau^M\mathsf E_{[q]}C_a\mathsf R_a
		(I-\tau\mathscr H_{(a)})^{-1}
		\mathsf R_a^{\top}D_a\mathsf E_{[p]}^{\top}.
	\end{align}
	Thus every selected cyclic product gives a factor-resolved continued
	fraction for the same rectangular resolvent of \(T\).
	
	If the displayed matrices define bounded operators on the chosen sequence
	space, the same formulas hold as operator identities wherever the
	corresponding inverses exist.
\end{theorem}

\begin{proof}
	For \(n\ge1\),
	\[
		(C_aD_a)^n=C_a(D_aC_a)^{n-1}D_a,
	\]
	which gives \eqref{eq:balanced-ABBA-resolvent-dressing} after summing the
	formal Neumann series.  Interchanging \(C_a\) and \(D_a\) gives
	\eqref{eq:cyclic-resolvent-reverse-dressing}.  The coefficient identities
	\[
		(C_aD_a)^nC_a=C_a(D_aC_a)^n,
		\qquad
		D_a(C_aD_a)^n=(D_aC_a)^nD_a
	\]
	prove the two intertwining formulas.  Substitution of
	\eqref{eq:admissible-cyclic-star-resolvent} into
	\eqref{eq:balanced-ABBA-resolvent-dressing}, followed by the rectangular
	coordinate projections, proves
	\eqref{eq:original-resolvent-through-admissible-star-matrix} and
	\eqref{eq:original-rectangular-resolvent-through-admissible-star-matrix}.
	For bounded operators the same identities follow by multiplying the
	corresponding inverse operators.
\end{proof}

The resolvents of consecutive cyclic products are connected by one factor.
If
\[
	K_a\coloneq F_{a+1}\cdots F_MF_1\cdots F_{a-1},
\]
with indices read modulo \(M\), then
\[
	T_{(a)}=F_aK_a,
	\qquad
	T_{(a+1)}=K_aF_a,
\]
and therefore
\begin{equation}
	\label{eq:consecutive-cyclic-resolvent-intertwining}
	(I-tT_{(a)})^{-1}F_a
	=
	F_a(I-tT_{(a+1)})^{-1}.
\end{equation}
Here \(a+1\) is read cyclically.  This is the one-factor Darboux relation
joining consecutive members of the cyclic family.

\subsection{Polynomial recurrences contained in
\texorpdfstring{\(\mathscr H_{(a)}\)}{H(a)}}
\label{subsec:balanced-Fourier-Darboux-decimation}

The matrix \(\mathscr H_{(a)}\) is linked to the original factorization at
the level of recurrence equations.  One complete circuit through its factor
copies recovers a cyclic product of the refactorized factors.  The next results make
this statement precise for right and left polynomial solutions of those
equations.  The discussion in this subsection is algebraic; orthogonality and
representing measures enter in Section~\ref{sec:balanced-star-Christoffel}.

For \(k\in\{1,\ldots,M\}\), set
\[
	\widetilde T_{(a;k)}
	\coloneq
	\widetilde F_k^{(a)}\widetilde F_{k+1}^{(a)}
	\cdots\widetilde F_M^{(a)}
	\widetilde F_1^{(a)}\cdots\widetilde F_{k-1}^{(a)}.
\]
These are the cyclic products of the refactorization of
\(T_{(a)}\), and \(\widetilde T_{(a;1)}=T_{(a)}\).
Let a right vector on the \(M\) factor copies be written as
\(
	\boldsymbol X=(\boldsymbol X_1^{\top},\ldots,
	\boldsymbol X_M^{\top})^{\top}
\), and write a left vector as
\(
	\boldsymbol Y=(\boldsymbol Y_1,\ldots,\boldsymbol Y_M)
\).

\begin{proposition}[Polynomial solutions along a complete factor circuit]
	\label{prop:balanced-factor-copy-eigenvalue-equations}
	Write
	\[
	\boldsymbol X(\lambda)
	=
	\bigl(
	\boldsymbol X_1(\lambda)^{\top},
	\ldots,
	\boldsymbol X_M(\lambda)^{\top}
	\bigr)^{\top},
	\qquad
	\boldsymbol Y(\lambda)
	=
	\bigl(
	\boldsymbol Y_1(\lambda),
	\ldots,
	\boldsymbol Y_M(\lambda)
	\bigr),
	\]
	and adopt the cyclic convention
	\[
	\boldsymbol X_{M+1}=\boldsymbol X_1,
	\qquad
	\boldsymbol Y_{M+1}=\boldsymbol Y_1.
	\]
	Then
	\[
	\mathscr C_{(a)}\boldsymbol X(\lambda)
	=
	\lambda\boldsymbol X(\lambda)
	\]
	if and only if
	\begin{equation}
		\label{eq:balanced-right-factor-copy-equations}
		\widetilde F_k^{(a)}\boldsymbol X_{k+1}(\lambda)
		=
		\lambda\boldsymbol X_k(\lambda),
		\qquad k\in\{1,\ldots,M\}.
	\end{equation}
	In this case, each component satisfies
	\begin{equation}
		\label{eq:balanced-right-cyclic-eigenvalue-equation}
		\widetilde T_{(a;k)}\boldsymbol X_k(\lambda)
		=
		\lambda^M\boldsymbol X_k(\lambda).
	\end{equation}
	
	Similarly,
	\[
	\boldsymbol Y(\lambda)\mathscr C_{(a)}
	=
	\lambda\boldsymbol Y(\lambda)
	\]
	if and only if
	\begin{equation}
		\label{eq:balanced-left-factor-copy-equations}
		\boldsymbol Y_k(\lambda)\widetilde F_k^{(a)}
		=
		\lambda\boldsymbol Y_{k+1}(\lambda),
		\qquad k\in\{1,\ldots,M\}.
	\end{equation}
	In this case,
	\begin{equation}
		\label{eq:balanced-left-cyclic-eigenvalue-equation}
		\boldsymbol Y_k(\lambda)\widetilde T_{(a;k)}
		=
		\lambda^M\boldsymbol Y_k(\lambda).
	\end{equation}
\end{proposition}

\begin{proof}
	By the block definition of \(\mathscr C_{(a)}\), its \(k\)-th block
	row contains only the factor \(\widetilde F_k^{(a)}\), in block column
	\(k+1\), with the indices read cyclically.  Therefore the \(k\)-th block
	component of
	\(\mathscr C_{(a)}\boldsymbol X(\lambda)
	=\lambda\boldsymbol X(\lambda)\)
	is precisely
	\[
	\widetilde F_k^{(a)}\boldsymbol X_{k+1}(\lambda)
	=
	\lambda\boldsymbol X_k(\lambda).
	\]
	This proves the first equivalence.
	
	Starting with the equation of index \(k\) and substituting successively
	the equations for the following \(M-1\) factor positions gives
	\[
	\widetilde F_k^{(a)}
	\widetilde F_{k+1}^{(a)}
	\cdots
	\widetilde F_{k+M-1}^{(a)}
	\boldsymbol X_k(\lambda)
	=
	\lambda^M\boldsymbol X_k(\lambda),
	\]
	where all factor indices are read modulo \(M\).  The product on the
	left is \(\widetilde T_{(a;k)}\), which proves
	\eqref{eq:balanced-right-cyclic-eigenvalue-equation}.
	
	The left eigenvalue equation is treated in the same way.  Its \(k\)-th
	block component is
	\[
	\boldsymbol Y_k(\lambda)\widetilde F_k^{(a)}
	=
	\lambda\boldsymbol Y_{k+1}(\lambda).
	\]
	Composing these identities once around the complete factor circuit
	and using the cyclic convention gives
	\[
	\boldsymbol Y_k(\lambda)
	\widetilde T_{(a;k)}
	=
	\lambda^M\boldsymbol Y_k(\lambda),
	\]
	as claimed.
\end{proof}

The rotational symmetry separates the powers of \(\lambda\) according to
their residues modulo \(M\).  To state this precisely, choose ordered right
boundary states
\[
(\alpha_r,n_r^{\mathrm R}),
\qquad r\in\{1,\ldots,q\},
\]
and ordered left boundary states
\[
(\beta_s,n_s^{\mathrm L}),
\qquad s\in\{1,\ldots,p\},
\]
where the first component is the factor-copy label and the second is the
coordinate within that copy.  Normalize the corresponding polynomial
solutions by
\[
(\boldsymbol X_{\alpha_r,j}(\lambda))_{n_r^{\mathrm R}}
=\delta_{rj},
\qquad
(\boldsymbol Y_{i,\beta_s}(\lambda))_{n_s^{\mathrm L}}
=\delta_{is}.
\]
For an integer \(\ell\), let \([\ell]_M\in\{0,\ldots,M-1\}\) denote its residue
modulo \(M\).

\begin{theorem}[Fourier decomposition into cyclic-product recurrences]
	\label{thm:balanced-Fourier-Darboux-decimation}
	Assume that the normalized right and left polynomial solutions described
	above exist and are unique.  Their blocks on factor copy \(k\) have the
	forms
	\begin{align}
		\label{eq:balanced-right-Fourier-decimation}
		\boldsymbol X_{k,j}(\lambda)
		&=
		\lambda^{[k-\alpha_j]_M}
		\boldsymbol V_{k,j}(\lambda^M),\\
		\label{eq:balanced-left-Fourier-decimation}
		\boldsymbol Y_{i,k}(\lambda)
		&=
		\lambda^{[\beta_i-k]_M}
		\boldsymbol W_{i,k}(\lambda^M).
	\end{align}
	With \(x=\lambda^M\), the radial polynomial sequences satisfy
	\begin{equation}
		\label{eq:balanced-radial-cyclic-recurrences}
		\widetilde T_{(a;k)}\boldsymbol V_{k,j}(x)
		=
		x\boldsymbol V_{k,j}(x),
		\qquad
		\boldsymbol W_{i,k}(x)\widetilde T_{(a;k)}
		=
		x\boldsymbol W_{i,k}(x).
	\end{equation}
	Moreover, the individual refactorized factors connect consecutive radial
	systems; the copy index \(k+1\) is read cyclically:
	\begin{equation}
		\label{eq:balanced-right-radial-connections}
		\widetilde F_k^{(a)}\boldsymbol V_{k+1,j}
		=
		\begin{cases}
			x\boldsymbol V_{k,j},
			&k\equiv\alpha_j-1\pmod M,\\
			\boldsymbol V_{k,j},&\textnormal{otherwise},
		\end{cases}
	\end{equation}
	and
	\begin{equation}
		\label{eq:balanced-left-radial-connections}
		\boldsymbol W_{i,k}\widetilde F_k^{(a)}
		=
		\begin{cases}
			x\boldsymbol W_{i,k+1},
			&k\equiv\beta_i\pmod M,\\
			\boldsymbol W_{i,k+1},&\textnormal{otherwise}.
		\end{cases}
	\end{equation}
	Thus the polynomial solution of the cyclic recurrence of the refactorized factors contains,
	as its Fourier components, polynomial solutions of all the cyclic products
	\(\widetilde T_{(a;1)},\ldots,\widetilde T_{(a;M)}\).
\end{theorem}

\begin{proof}
	Let
	\[
		\mathscr D_\zeta
		=
		\operatorname{diag}
		(I,\zeta I,\ldots,\zeta^{M-1}I),
		\qquad \zeta^M=1.
	\]
	The proof of Proposition~\ref{prop:cyclic-star-symmetry}, applied to the
	cyclic matrix of the refactorized factors, gives
	\[
		\mathscr D_\zeta^{-1}\mathscr C_{(a)}
		\mathscr D_\zeta
		=
		\zeta\mathscr C_{(a)}.
	\]
	Hence \(\mathscr D_\zeta\boldsymbol X(\lambda)\) is a right polynomial
	solution at \(\zeta\lambda\).  Restoring the identity normalization at
	the right boundary and using uniqueness gives
	\[
		\boldsymbol X(\zeta\lambda)
		=
		\mathscr D_\zeta\boldsymbol X(\lambda)
		\operatorname{diag}
		(\zeta^{1-\alpha_1},\ldots,\zeta^{1-\alpha_q}).
	\]
	Thus every power occurring in \(\boldsymbol X_{k,j}\) is congruent to
	\(k-\alpha_j\) modulo \(M\), which proves
	\eqref{eq:balanced-right-Fourier-decimation}.  The left relation
	\[
		\boldsymbol Y(\zeta\lambda)
		=
		\operatorname{diag}
		(\zeta^{\beta_1-1},\ldots,\zeta^{\beta_p-1})
		\boldsymbol Y(\lambda)\mathscr D_\zeta^{-1}
	\]
	proves \eqref{eq:balanced-left-Fourier-decimation}.

	Substitution in
	Proposition~\ref{prop:balanced-factor-copy-eigenvalue-equations} proves
	\eqref{eq:balanced-radial-cyclic-recurrences}.  In the right factor
	equation, the residues \([k-\alpha_j]_M+1\) and
	\([k+1-\alpha_j]_M\) agree as integers unless the second residue wraps
	from \(M\) to \(0\).  In that case their powers differ by \(M\), leaving
	the factor \(\lambda^M=x\).  This is precisely
	\eqref{eq:balanced-right-radial-connections}.  The same calculation in
	the left factor equation proves
	\eqref{eq:balanced-left-radial-connections}.
\end{proof}

The refactorization is a factorization of
\(T_{(a)}=D_aC_a\), whereas the original matrix is
\(T=C_aD_a\).  The corresponding polynomial solutions are related without
any spectral hypothesis.

Suppose that polynomial vectors \(\boldsymbol P_T(x)\) and
\(\boldsymbol Q_T(x)\) satisfy
	\(T\boldsymbol P_T(x)=x\boldsymbol P_T(x)\),
	\(\boldsymbol Q_T(x)T=x\boldsymbol Q_T(x)\).
Define
\begin{align*}
	\boldsymbol P^{(a)}(x)
	&\coloneq
	D_a\boldsymbol P_T(x),&
	\boldsymbol Q^{(a)}(x)
	&\coloneq
	\boldsymbol Q_T(x)C_a.
\end{align*}
These products are well defined coordinatewise because \(C_a\) and \(D_a\)
are banded.

\begin{proposition}[Polynomial transport to a cyclic product in the constructed range]
	\label{prop:balanced-ABBA-polynomial-transport}
	The polynomial vectors defined above satisfy
	\(	T_{(a)}\boldsymbol P^{(a)}
		=x\boldsymbol P^{(a)}\),
		\(\boldsymbol Q^{(a)}T_{(a)}
		=x\boldsymbol Q^{(a)}\).
\end{proposition}

\begin{proof}
	The right recurrence \(C_aD_a\boldsymbol P_T=x\boldsymbol P_T\) gives
	\[
		D_aC_a(D_a\boldsymbol P_T)
		=
		D_a(C_aD_a)\boldsymbol P_T
		=
		xD_a\boldsymbol P_T.
	\]
	This proves the right identity.  Similarly,
	\[
		(\boldsymbol Q_TC_a)D_aC_a
		=
		\boldsymbol Q_T(C_aD_a)C_a
		=
		x\boldsymbol Q_TC_a,
	\]
	which proves the left identity.
\end{proof}

Consequently, \(D_a\boldsymbol P_T\) and
\(\boldsymbol Q_TC_a\) are polynomial solutions for \(T_{(a)}\).  If their
right and left initial-condition blocks are constant and nonsingular,
uniqueness of the recurrence problems identifies them, up to constant
nonsingular changes of initial conditions, with the radial system at the
factor copy for which \(\widetilde T_{(a;1)}=T_{(a)}\).  Without this
additional condition no constant-gauge identification is asserted.  The
other radial systems are obtained by the explicit factor connections
\eqref{eq:balanced-right-radial-connections} and
\eqref{eq:balanced-left-radial-connections}.  Thus the polynomial solutions
of the recurrence \(\mathscr H_{(a)}\) combine, through their residue classes modulo
\(M\), the recurrence solutions of the cyclic products of the
refactorization.  These identities give the recurrence solutions
algebraically.  Under the Favard--Christoffel hypotheses, the radial systems
with the prescribed admissible initial conditions are the mixed-type
orthogonal polynomials for the transform associated with \(T_{(a)}\).  For
\(a=p+1\), this is the familiar transform \(x\,\mathrm d\Psi\).  What remains, and is treated
in Section~\ref{sec:balanced-star-Christoffel}, is the moment problem that
assembles all cyclic components into the full recurrence \(\mathscr H_{(a)}\) on the
star.

For reference, Table~\ref{tab:height-CF-notation} collects the notation used
for the state heights, the two blockings, and the original and reordered
coordinate projections.
\begin{table}[H]
\centering
\caption{Principal notation for the height and boundary-state constructions.}
\label{tab:height-CF-notation}
\small
\renewcommand{\arraystretch}{1.25}
\begin{tabular}{@{}p{0.28\textwidth}p{0.66\textwidth}@{}}
\hline
Symbol & Meaning \\
\hline
\(d,P,Q,R,w\) &
\(d=\gcd(p,q)\), reduced steps \(+Q\) and \(-P\),
\(R=P+Q\), and height-block width \(w=\max(P,Q)\). \\
\(h(j,n)\), \(h^{(a)}(j,n)\) &
Heights in the prescribed and minimum-height factor orders. \\
\(\boldsymbol v\) &
Positive eigenvector used only for the diagonal normalization. \\
\(\mathscr F_h\), \(\mathcal E_h\),
\(\mathcal E_h^{(a)}\) &
Labeled states of height \(h\) and their coordinate spaces in the prescribed
and minimum-height orders. \\
\(\Pi_a\), \(\mathsf E_{[r]}\) &
Projection onto the \(a\)-th copy of the original sequence space; projection
onto its first \(r\) coordinates. \\
\(\mathsf Q_r^{(a)}\) &
The same \(r\)-coordinate projection after height reordering. \\
\(\mathsf J_j\) &
Zero-filling injection of a block column into the coordinates of the height
block \(\mathscr V_j\). \\
\(\mathscr V_m\), \(\mathscr V_{[m,N]}\) &
Height block of width \(w\); union of the height blocks
\(\mathscr V_m,\ldots,\mathscr V_N\). \\
\(\mathscr H^{\langle m\rangle}\),
\(\mathscr H^{[m,N]}\) &
Semi-infinite tail beginning at \(\mathscr V_m\); finite principal matrix on
blocks \(m,\ldots,N\). \\
\(\mathsf A_m,\mathsf B_m,\mathsf C_m\), \(\Phi_m\) &
Diagonal and rectangular neighboring-block couplings for the width-\(w\)
blocking; its leading tail-resolvent block. \\
\(\mathcal U_m\), \(E_m\), \(\mathbf G_m\) &
Group of \(q\) boundary states over \(Q\) heights; its coordinate isometry;
its \(q\times q\) tail resolvent. \\
\(\mathbf A_m,\mathbf B_m,\mathbf D_m^{(\eta)}\) &
The \(q\times q\) internal, upward, and depth-\(\eta\) return matrices in the
boundary-state recursion. \\
\hline
\end{tabular}
\end{table}

\subsection{The matrix continued fraction obtained from height blocks}
\label{subsec:matrix-CF-height-blocks}

Before constructing the continued fraction, the correspondence is summarized.
One application of the original recurrence matrix \(T\) is decomposed into
\(M=p+q\) elementary factor steps. Each elementary step is one edge of the
labeled height graph. The fixed normalizing entries are the upward edges,
whereas the coefficient entries are the downward edges and carry the
corresponding bidiagonal coefficients. Thus:
\begin{itemize}
\item a bidiagonal factor corresponds to an elementary step;
\item a coefficient entry corresponds to a downward edge carrying that
coefficient as its weight;
\item one complete cycle
\[
L_1\cdots L_pU_q\cdots U_1
\]
corresponds to one complete sequence of \(M=p+q\) elementary steps.
\end{itemize}
The continued fraction is obtained by applying Schur complements to finite
sections of the resolvent matrix of this height graph. Since the jumps are
generally \(+Q\) and \(-P\), the graph is first grouped into finite height
blocks; this is what turns the reordered cyclic linearization into a
block-tridiagonal matrix.  The infinite tail identity is then obtained by
coefficientwise stabilization of these finite calculations.

The rational-step height graph has elementary steps \(+Q\) and \(-P\).  Unless
\(P=Q=1\), this graph is not tridiagonal with respect to individual heights:
one elementary transition may jump over several height levels.  Thus, if the
basis were ordered by single heights, the corresponding matrix would not be
tridiagonal.

To recover a tridiagonal structure, group consecutive heights into blocks
large enough so that one elementary transition can move at most to a
neighboring block.  Put
\begin{equation*}
	w\coloneq \max(P,Q).
\end{equation*}
For \(m\in\mathbb N_0\), define the \(m\)-th height block by
\begin{equation*}
	\mathscr V_m
	\coloneq
	\bigl\{(a,n):
	mw\le h(a,n)<(m+1)w\bigr\}.
\end{equation*}
Here
\[
|\mathscr V_m|
\coloneq
\#\bigl\{(a,n):mw\le h(a,n)<(m+1)w\bigr\}
\]
is the number of labeled states in the block, and hence the dimension of
its coordinate space.  In particular, \(I_{|\mathscr V_m|}\) denotes the
identity matrix of that order.  The congruence
\[
\omega_a\equiv(a-1)Q\pmod R,
\qquad a\in\{1,\ldots,M\},
\]
follows from \(\omega_a=(a-1)Q\) at the lower-factor positions.  For an
upper-factor position
\(a=p+b\), formula~\eqref{eq:general-height-shifts} and \(pQ=qP\) give
\[
(p+b-1)Q-\omega_{p+b}=(b-1)(P+Q)=(b-1)R,
\]
which proves the same congruence.
Because \(\gcd(Q,R)=1\), each residue class modulo \(R\) occurs for
exactly \(d\) factor labels.  Hence every height supports at most \(d\)
labeled states.  Since a block contains \(w\) consecutive heights,
\[
0\le |\mathscr V_m|\le dw=\max\{p,q\}.
\]
In fact, no block is empty.  Starting from the state \((1,0)\) and following
only normalizing entries produces a state at every height in
\(Q\mathbb N_0\).  Every interval of \(w\ge Q\) consecutive heights contains
a multiple of \(Q\), and hence contains at least one labeled state.  Therefore
\[
1\le |\mathscr V_m|\le\max\{p,q\},
\qquad m\in\mathbb N_0.
\]
Moreover, if
\(mw\ge\max_{1\le a\le M}\omega_a\), then all \(d\) labels belonging to
each height in the block occur with a nonnegative sequence index, and
therefore
\[
|\mathscr V_m|=\max\{p,q\}.
\]
Since \(\max_a\omega_a=pQ\) and \(w\ge Q\), this equality holds, in
particular, for every \(m\ge p\).
Thus only finitely many initial blocks can have smaller dimension.  Grouping
the basis of the cyclic linearization in this way makes its matrix block
tridiagonal, so the standard finite-dimensional Schur-complement recursion
applies to every finite section.

The consecutive-height blocks retain the exact-height order of
Definition~\ref{def:height-ordered-cyclic-operator}.  In block notation this
order is
\[
\mathscr V_0,\mathscr V_1,\mathscr V_2,\ldots,
\]
with increasing exact height inside each block and with the factor label
retained whenever a height carries more than one state.  This displays the
already defined matrix \(\mathscr H_{p,q}\) in finite neighboring blocks.

\begin{example}[The height pattern for \texorpdfstring{$(p,q)=(2,3)$}{(p,q)=(2,3)}]
	\label{ex:height-blocks-CF-23}
	For \((p,q)=(2,3)\), one has
	\[
	d=1,
	\qquad
	P=2,
	\qquad
	Q=3,
	\qquad
	R=5,
	\qquad
	w=3.
	\]
	The first labeled states, ordered by height and then grouped into intervals
	of width \(w\), are listed in Table~\ref{tab:height-blocks-CF-23}.
	\begin{table}[htbp]
	\centering
	\caption{The first height blocks for \((p,q)=(2,3)\).  The columns give
		the block number, its height interval, and the corresponding labeled
		state set.}
	\label{tab:height-blocks-CF-23}
	\(
	\begin{NiceArray}{c|c|c}[cell-space-limits=3pt]
		m&\textnormal{height interval}&\mathscr V_m\\
		\hline
		0&\{0,1,2\}&\{(1,0),(5,0)\}\\
		1&\{3,4,5\}&\{(2,0),(4,0),(1,1)\}\\
		2&\{6,7,8\}&\{(3,0),(5,1),(2,1)\}\\
		3&\{9,10,11\}&\{(4,1),(1,2),(3,1)\}
	\end{NiceArray}
	\)
	\end{table}
	The state at height \(1\) is absent because of
	the restriction \(n\ge0\).  The normalizing edges of height \(+3\) and the
	coefficient edges of height \(-2\) connect only the same or adjacent blocks.  Therefore
	Table~\ref{tab:height-blocks-CF-23} already displays the block-tridiagonal
	structure.  The general finite-dimensional Schur-complement recursion below,
	followed by coefficientwise stabilization, therefore applies to this pair
	\((p,q)\).
\end{example}

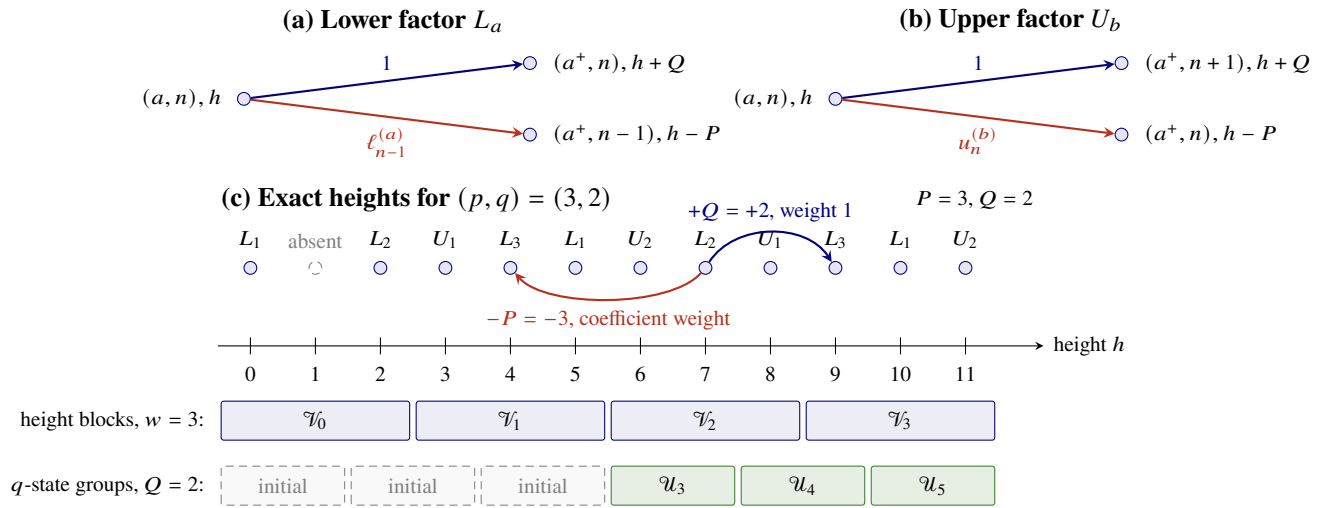
\begin{figure}[H]
	\centering
	\begin{tikzpicture}[
		x=0.86cm,
		y=0.86cm,
		>=stealth,
		state/.style={circle,draw=NavyBlue,fill=NavyBlue!10,inner sep=1.7pt},
		missing/.style={circle,draw=Gray,densely dashed,inner sep=1.7pt},
		smalllab/.style={font=\scriptsize},
		edge/.style={->,thick}
	]
		\node[font=\small\bfseries] at (-4.8,3.8) {(a) Lower factor \(L_a\)};
		\node[state] (ls) at (-7.1,2.6) {};
		\node[left=2pt of ls,smalllab] {\((a,n),h\)};
		\node[state] (lu) at (-2.7,3.15) {};
		\node[right=2pt of lu,smalllab] {\((a^+,n),h+Q\)};
		\node[state] (ld) at (-2.7,2.05) {};
		\node[right=2pt of ld,smalllab] {\((a^+,n-1),h-P\)};
		\draw[edge,NavyBlue] (ls) -- node[above,smalllab] {\(1\)} (lu);
		\draw[edge,BrickRed] (ls) -- node[below,smalllab] {\(\ell_{n-1}^{(a)}\)} (ld);

		\node[font=\small\bfseries] at (4.7,3.8) {(b) Upper factor \(U_b\)};
		\node[state] (us) at (2.0,2.6) {};
		\node[left=2pt of us,smalllab] {\((a,n),h\)};
		\node[state] (uu) at (6.4,3.15) {};
		\node[right=2pt of uu,smalllab] {\((a^+,n+1),h+Q\)};
		\node[state] (ud) at (6.4,2.05) {};
		\node[right=2pt of ud,smalllab] {\((a^+,n),h-P\)};
		\draw[edge,NavyBlue] (us) -- node[above,smalllab] {\(1\)} (uu);
		\draw[edge,BrickRed] (us) -- node[below,smalllab] {\(u_n^{(b)}\)} (ud);

		\node[font=\small\bfseries,anchor=west] at (-7.6,1.05)
			{(c) Exact heights for \((p,q)=(3,2)\)};
		\node[smalllab,anchor=east] at (5.2,1.05) {\(P=3\), \(Q=2\)};
		\draw[->] (-7.5,-1.2) -- (5.2,-1.2) node[right,smalllab] {height \(h\)};
		\foreach \h/\lab in {
			0/L_1,1/\varnothing,2/L_2,3/U_1,4/L_3,5/L_1,
			6/U_2,7/L_2,8/U_1,9/L_3,10/L_1,11/U_2}
		{
			\pgfmathsetmacro{\xx}{-7+\h}
			\draw (\xx,-1.33) -- (\xx,-1.07);
			\node[smalllab,below] at (\xx,-1.36) {\(\h\)};
			\ifnum\h=1
				\node[missing] (h\h) at (\xx,0) {};
				\node[smalllab,above=1pt of h\h,text=Gray] {absent};
			\else
				\node[state] (h\h) at (\xx,0) {};
				\node[smalllab,above=1pt of h\h] {\(\lab\)};
			\fi
		}
		\draw[edge,NavyBlue]
			(h7) .. controls (0.45,0.70) and (1.55,0.70) ..
			node[above,smalllab] {\(+Q=+2\), weight \(1\)} (h9);
		\draw[edge,BrickRed]
			(h7) .. controls (-0.55,-0.62) and (-2.45,-0.62) ..
			node[below,smalllab] {\(-P=-3\), coefficient weight} (h4);

		\node[smalllab,anchor=east] at (-7.55,-2.35)
			{height blocks, \(w=3\):};
		\foreach \m/\xa/\xb in {0/-7.45/-4.55,1/-4.45/-1.55,2/-1.45/1.45,3/1.55/4.45}
		{
			\draw[rounded corners=1pt,fill=NavyBlue!7,draw=NavyBlue]
				(\xa,-2.65) rectangle (\xb,-2.05);
			\node[smalllab] at ({(\xa+\xb)/2},-2.35) {\(\mathscr V_{\m}\)};
		}

		\node[smalllab,anchor=east] at (-7.55,-3.35)
			{\(q\)-state groups, \(Q=2\):};
		\foreach \m/\xa/\xb in {0/-7.45/-5.55,1/-5.45/-3.55,2/-3.45/-1.55}
		{
			\draw[rounded corners=1pt,densely dashed,fill=Gray!5,draw=Gray]
				(\xa,-3.65) rectangle (\xb,-3.05);
			\node[smalllab,text=Gray] at ({(\xa+\xb)/2},-3.35) {initial};
		}
		\foreach \m/\xa/\xb in {3/-1.45/0.45,4/0.55/2.45,5/2.55/4.45}
		{
			\draw[rounded corners=1pt,fill=OliveGreen!9,draw=OliveGreen]
				(\xa,-3.65) rectangle (\xb,-3.05);
			\node[smalllab] at ({(\xa+\xb)/2},-3.35) {\(\mathcal U_{\m}\)};
		}
	\end{tikzpicture}
	\caption{States, heights, and the two groupings.  Panels (a) and (b)
		show that every nonzero bidiagonal entry produces one of the two height
		displacements \(+Q\) and \(-P\).  Panel (c) shows the beginning of the
		exact-height order for \((p,q)=(3,2)\).  Blocks of width
		\(w=\max(P,Q)\) give the block-tridiagonal matrix used in the Schur
		continued fraction; groups of width \(Q\), complete from \(m=p\)
		onward, give the \(q\)-state recursion.}
	\label{fig:height-states-and-groupings}
\end{figure}

Figure~\ref{fig:height-states-and-groupings} shows the same construction for
\((p,q)=(3,2)\), where the two groupings used in the article have different
widths.  The lower row of the figure anticipates the groups of \(Q\)
consecutive heights used in the \(q\)-state recursion of
Subsection~\ref{subsec:q-channel-branched-refinement}.

\begin{lemma}[Block-tridiagonal height form]
	\label{lem:general-rational-block-tridiagonal}
	Relative to the decomposition
	\[
	\operatorname{span}\mathscr V_{p,q}
	=
	\bigoplus_{m\ge0}\operatorname{span}\mathscr V_m,
	\]
	the matrix \(\mathscr H\) has block-tridiagonal form
	\begin{equation*}
		\mathscr H
		=
		\begin{bNiceMatrix}[margin=2pt,cell-space-limits=2pt]
			\mathsf A_0&\mathsf B_0&0&\cdots\\
			\mathsf C_0&\mathsf A_1&\mathsf B_1&\ddots\\
			0&\mathsf C_1&\mathsf A_2&\ddots\\
			\vdots&\ddots&\ddots&\ddots[shorten-end=10pt]
		\end{bNiceMatrix}.
	\end{equation*}
	Every nonzero entry of the blocks is either \(1\), or one of the individual
	bidiagonal coefficients \(\ell_n^{(a)}\), \(u_n^{(b)}\).
\end{lemma}

\begin{proof}
	Every elementary transition changes height by either \(+Q\) or \(-P\).
	Since
	\[
	P\le w,
	\qquad
	Q\le w,
	\]
	a transition starting from a state in \(\mathscr V_m\) can end only in
	\(\mathscr V_{m-1}\), \(\mathscr V_m\), or \(\mathscr V_{m+1}\).  Hence, in
	the block ordering, only the diagonal, upper, and lower block diagonals can
	be nonzero.  The assertion about the entries follows from
	Proposition~\ref{prop:general-height-steps}, because the height reordering
	does not change the actual coefficients of the elementary transitions.
\end{proof}

The sizes of the three block types will be used repeatedly.  With the
row-to-column convention fixed above,
\begin{equation*}
\mathsf A_m\in
\mathbb C^{|\mathscr V_m|\times|\mathscr V_m|},
\qquad
\mathsf B_m\in
\mathbb C^{|\mathscr V_m|\times|\mathscr V_{m+1}|},
\qquad
\mathsf C_m\in
\mathbb C^{|\mathscr V_{m+1}|\times|\mathscr V_m|}.
\end{equation*}
Thus the coupling blocks \(\mathsf B_m\) and \(\mathsf C_m\) are generally
rectangular, although every Schur complement below is taken in a square
finite-dimensional matrix.

The variable \(\tau\) marks one elementary factor step of the cyclic
linearization.  Thus a path of length \(\ell\) in the enlarged graph contributes
a factor \(\tau^\ell\).  Since one complete multiplication by the original
matrix \(T\) corresponds to \(M=p+q\) elementary factor steps, the variable
for the original recurrence matrix will later be
\[
t=\tau^M.
\]
Equivalently, in the Weyl or Stieltjes transform variable \(z\), use
\[
z^{-1}=\tau^M.
\]
Thus \(\tau\) counts individual factor steps.  After projection back to the
first copy of the sequence space, the variables for the original operator
are \(t=\tau^M\) and \(z=t^{-1}\).
Here \(\tau\) is only the elementary-step variable; it is unrelated to the
Christoffel determinants or the integrable tau functions introduced later.

\begin{definition}[Coordinate submatrices, compressions, and finite sections]
	\label{def:principal-compression}
	Let \(K\) be a matrix written in the labeled-state basis.  For two sets
	of states \(\mathscr A\) and \(\mathscr B\), finite or infinite, write
	\[
	[K]_{\mathscr A\times\mathscr B}
	\]
	for the matrix obtained by retaining the rows indexed by \(\mathscr A\)
	and the columns indexed by \(\mathscr B\).  If
	\(\mathscr A=\mathscr B=\mathscr S\), the principal compression of
	\(K\) to \(\mathscr S\) is the principal submatrix
	\([K]_{\mathscr S\times\mathscr S}\).  Thus a compression simply deletes
	all rows and columns outside \(\mathscr S\); the coordinate subspace
	spanned by \(\mathscr S\) may therefore fail to be invariant under \(K\).

	For the height blocks introduced above, set
	\[
	\mathscr V_{\ge m}\coloneq\bigcup_{k\ge m}\mathscr V_k,
	\qquad
	\mathscr V_{[m,N]}\coloneq\bigcup_{k=m}^{N}\mathscr V_k
	\quad (m\le N),
	\]
	and use the notation
	\[
	\mathscr H^{\langle m\rangle}
	\coloneq
	[\mathscr H]_{\mathscr V_{\ge m}\times\mathscr V_{\ge m}},
	\qquad
	\mathscr H^{[m,N]}
	\coloneq
	[\mathscr H]_{\mathscr V_{[m,N]}\times\mathscr V_{[m,N]}}.
	\]
	The first matrix is semi-infinite and acts, coordinatewise, on sequences
	indexed by \(\mathscr V_{\ge m}\).  The second is called the
	\emph{finite section} from block \(m\) through block \(N\): it is simply
	the finite principal matrix obtained by truncating the tail after block
	\(N\).  It acts
	on the finite-dimensional coordinate space
	\[
	\mathbb C^{\mathscr V_{[m,N]}}
	=
	\operatorname{span}
	\{e_{(a,n)}:(a,n)\in\mathscr V_{[m,N]}\}
	=
	\bigoplus_{k=m}^{N}\operatorname{span}\mathscr V_k,
	\qquad
	\dim\mathbb C^{\mathscr V_{[m,N]}}
	=
	\sum_{k=m}^{N}|\mathscr V_k|.
	\]
	Consequently, angle brackets always indicate an infinite tail, whereas
	square brackets with two endpoints indicate the finite principal submatrix
	that retains the stated range of blocks.  In both
	cases the endpoints refer to height-block indices, not to the original
	sequence indices.  In particular,
	\(\mathscr H^{[0,N]}\) retains the blocks \(0,\ldots,N\).  This notation
	will be used consistently below.  If the full matrix carries a
	family label, the compression is taken after fixing the family, as in
	\((\mathscr H_{p,q}^{\mathrm P})^{\langle m\rangle}\) or
	\((\mathscr H_{p,q}^{\mathrm J})^{[m,N]}\).
\end{definition}

Until boundedness is imposed in
Subsection~\ref{subsec:finite-convergents}, every infinite resolvent in this
subsection is understood through its formal Neumann series.  This is
well-defined coefficient by coefficient because each state has only finitely
many outgoing edges, so only finitely many paths of a fixed length begin at a
given state.
Define the initial block of its resolvent by
\begin{equation*}
	\Phi_m(\tau)
	\coloneq
	\left[
	(I-\tau\mathscr H^{\langle m\rangle})^{-1}
	\right]_{\mathscr V_m\times\mathscr V_m}.
\end{equation*}
This matrix records all paths which start in the block \(\mathscr V_m\),
remain in the tail
\[
\mathscr V_m\cup\mathscr V_{m+1}\cup\cdots,
\]
and return to \(\mathscr V_m\).  The entries of \(\Phi_m(\tau)\) are therefore
generating functions for returns between labeled states of the same height
block.

For \(N\ge m\), define also the resolvent block of the finite section by
\begin{equation*}
	\Phi_m^{[N]}(\tau)
	\coloneq
	\left[
	\left(I-\tau\mathscr H^{[m,N]}\right)^{-1}
	\right]_{\mathscr V_m\times\mathscr V_m}.
\end{equation*}
The matrix inverted here is square and finite-dimensional.  These finite
resolvent blocks will be used to derive the recursion for \(\Phi_m\), and all
Schur complements below are finite-dimensional.

\begin{theorem}[Matrix continued fraction from height blocks]
	\label{thm:matrix-CF-height-blocks}
	For every \(m\in\mathbb N_0\), the resolvent blocks
	\(\Phi_m\) satisfy
	\begin{equation*}
		\Phi_m(\tau)
		=
		\left(
		I_{|\mathscr V_m|}
		-
		\tau\mathsf A_m
		-
		\tau^2\mathsf B_m\Phi_{m+1}(\tau)\mathsf C_m
		\right)^{-1}.
	\end{equation*}
	Consequently, \(\Phi_0(\tau)\) is given by an iterated matrix continued
	fraction whose coefficients are the individual bidiagonal entries of
	\[
	L_1,\ldots,L_p,U_q,\ldots,U_1.
	\]
\end{theorem}

\begin{proof}
	Fix \(N>m\) and separate the first height block of the finite section
	from its remaining blocks.  With respect to the finite-dimensional
	decomposition
	\[
	\operatorname{span}\mathscr V_m
	\oplus
	\bigoplus_{k=m+1}^{N}\operatorname{span}\mathscr V_k,
	\]
	the square matrix \(I-\tau\mathscr H^{[m,N]}\) has the block
	form
	\[
	I-\tau\mathscr H^{[m,N]}
	=
	\begin{bNiceMatrix}
		I_{|\mathscr V_m|}-\tau\mathsf A_m
		&
		-\tau\mathsf B_m
		\\
		-\tau\mathsf C_m
		&
		I-\tau\mathscr H^{[m+1,N]}
	\end{bNiceMatrix}.
	\]
	Here \(\mathsf A_m\) describes transitions staying inside
	\(\mathscr V_m\), \(\mathsf B_m\) describes transitions from
	\(\mathscr V_m\) to \(\mathscr V_{m+1}\), and \(\mathsf C_m\) describes
	transitions from \(\mathscr V_{m+1}\) back to \(\mathscr V_m\).

	All four blocks have finite dimensions.  The lower-right block has constant
	term \(I\), so it is invertible as a finite matrix over
	\(\mathbb C[[\tau]]\).  Applying the ordinary Schur-complement formula gives
	\[
	\Phi_m^{[N]}(\tau)
	=
	\left(
	I_{|\mathscr V_m|}
	-
	\tau\mathsf A_m
	-
	\tau^2\mathsf B_m
	\left[
	(I-\tau\mathscr H^{[m+1,N]})^{-1}
	\right]_{\mathscr V_{m+1}\times\mathscr V_{m+1}}
	\mathsf C_m
	\right)^{-1}.
	\]
	By the definition of the finite resolvent block on the next tail, this is
	\[
	\Phi_m^{[N]}(\tau)
	=
	\left(
	I_{|\mathscr V_m|}
	-
	\tau\mathsf A_m
	-
	\tau^2\mathsf B_m\Phi_{m+1}^{[N]}(\tau)\mathsf C_m
	\right)^{-1}.
	\]
	This calculation involves only the finite square matrix
	\(I-\tau\mathscr H^{[m,N]}\).

	It remains to pass to the infinite tail.  The coefficient of \(\tau^j\)
	in \(\Phi_m\) is the sum of all length-\(j\) paths that begin and end in
	\(\mathscr V_m\) and stay in the tail.  Since the height-block matrix is
	block tridiagonal, a path that reaches the first deleted block
	\(\mathscr V_{N+1}\) and returns to \(\mathscr V_m\) has length at least
	\(2(N-m+1)\).  Consequently,
	\[
	[\tau^j]\Phi_m^{[N]}(\tau)
	=
	[\tau^j]\Phi_m(\tau)
	\qquad\textnormal{for }j<2(N-m+1).
	\]
	The analogous statement holds for the tail beginning at
	\(\mathscr V_{m+1}\).  Hence both finite resolvent blocks in the preceding
	recursion converge coefficientwise to their corresponding formal tail
	resolvent blocks as \(N\to\infty\).  Finally, the coefficient of degree
	\(j\) in the inverse of a matrix power series with constant term \(I\)
	depends only on its coefficients through degree \(j\).  Taking the
	coefficientwise limit in the finite Schur-complement identity therefore
	gives
	\[
	\Phi_m(\tau)
	=
	\left(
	I_{|\mathscr V_m|}
	-
	\tau\mathsf A_m
	-
	\tau^2\mathsf B_m\Phi_{m+1}(\tau)\mathsf C_m
	\right)^{-1}.
	\]

	The factor \(\tau^2\) has a simple meaning.  To leave the block
	\(\mathscr V_m\), enter the next tail, and return to \(\mathscr V_m\), a
	path must use one elementary step from \(\mathscr V_m\) to
	\(\mathscr V_{m+1}\), represented by \(\tau\mathsf B_m\), and one elementary
	step back from \(\mathscr V_{m+1}\) to \(\mathscr V_m\), represented by
	\(\tau\mathsf C_m\).  The excursions inside the next tail are summed by
	\(\Phi_{m+1}(\tau)\).  Thus the product
	\[
	\tau\mathsf B_m\,\Phi_{m+1}(\tau)\,\tau\mathsf C_m
	=
	\tau^2\mathsf B_m\Phi_{m+1}(\tau)\mathsf C_m
	\]
	is the contribution of all paths that leave \(\mathscr V_m\), make an
	arbitrary excursion in the next tail, and return to \(\mathscr V_m\).
\end{proof}

\begin{remark}[The variables \(\tau\), \(t\), and \(z\)]
	\label{rem:tau-t-z-variables}
	The matrix continued fraction above is written in the elementary-step
	variable \(\tau\).  Each application of \(\mathscr H\) corresponds
	to one bidiagonal factor.

	The original recurrence matrix \(T\), however, is recovered only after one
	full cycle of \(M\) elementary steps.  Hence, after projection onto the first
	copy of the cyclic linearization, only powers \(\tau^{Mn}\) survive.  This
	is why the moment variable for the original matrix is
	\(t=\tau^M\).
	In the Weyl matrix one uses the resolvent
	\[
	(zI-T)^{-1}
	=
	\frac{1}{z}(I-z^{-1}T)^{-1}.
	\]
	Thus the corresponding relation with the elementary-step variable is
	\(z^{-1}=\tau^M\).
	The matrix continued fraction is therefore first constructed in the variable
	\(\tau\), and its rectangular projection onto the first copy becomes a function of
	\(t=\tau^M\), or equivalently of \(z\) through \(z^{-1}=\tau^M\).
\end{remark}

\begin{remark}[The case \texorpdfstring{$q=1$}{q=1} and the scalar tails]
	\label{rem:Phi-G-q1}
	When \(q=1\), the labeled height graph becomes the ordinary weighted
	\(p\)-Dyck path graph associated with the bidiagonal factors.  The direct
	production graph of \(T\) is instead the contracted \(p\)-\L{}ukasiewicz
	graph.  The height-block resolvents \(\Phi_m\)
	remain matrix-valued when the blocks of width \(p\) are retained.  What
	becomes scalar is the finer resolvent obtained by cutting the graph at one
	height at a time.  These single-height tails, denoted by \(\widehat G_k\)
	below, satisfy the \(p\)-branched Stieltjes recursion; the precise reduction
	is proved in Subsection~\ref{subsec:q1-reduction-branched}.
\end{remark}

\begin{remark}[Relation with branched continued fractions]
	\label{rem:mixed-analogue-branched}
	The construction above extends the branched continued fractions of
	the Hessenberg case \(q=1\) to arbitrary \((p,q)\)
	\cite{PetreolleSokalZhu2023,Sokal2024,Lima2023,
	BranquinhoDiazFoulquieLimaManas2026JAT}.
	For \(q=1\), the recurrence matrix is the scalar lower-Hessenberg
	production matrix.  Its direct paths are \(p\)-\L{}ukasiewicz paths and give
	the branched Jacobi fraction; the additional bidiagonal factorization
	unfolds them into weighted \(p\)-Dyck paths and gives the branched
	Stieltjes fraction used here.  For general \((p,q)\), the same recurrence
	matrix \(T\), grouped into blocks of size \(q\), is lower block Hessenberg.
	Section~\ref{sec:general-Weyl-bridge} introduces its block output array and
	proves the corresponding mixed-type production/output identities, giving a
	block extension of the scalar formalism.

	The production matrix is \(T\), whereas the sparse matrix \(\mathscr H\)
	records the individual bidiagonal factors before they are multiplied.
	The reduced steps \(+Q\) and \(-P\), together with the labels attached to
	the cyclic factor positions, prevent a scalar first-return recursion in
	general.  The powers of \(\mathscr H\) still enumerate the labeled paths,
	but the first-return generating functions are matrices indexed by boundary
	states.  When \(T\) comes from a
	Gauss--Borel factorization of a mixed-type moment matrix,
	Theorem~\ref{thm:mixed-block-output-matrices} makes this correspondence
	precise: two block output matrices recover the inverses of the complete
	coefficient arrays of the two mixed-type polynomial vectors.
\end{remark}

\subsection{Finite convergents and denominator factorization}
\label{subsec:finite-convergents}

The Schur-complement recursion gives the convergent
\(\Phi_m^{[N]}\) when the height graph is cut after block \(N\).  This
subsection determines its denominator factorization and the first coefficient
that can be lost at the truncation boundary.  The following two subsections
use these algebraic results to prove analytic convergence and positivity, and
then to formulate the Pad\'e-type approximation, backward evaluation, and
residual error estimates.  Only the block-tridiagonal structure is used, so
the results apply to every factorization considered above.

For \(0\le m\le N\), the finite coordinate space is
\[
\mathbb C^{\mathscr V_{[m,N]}}
=
\bigoplus_{k=m}^{N}\operatorname{span}\mathscr V_k.
\]
Thus,
by Definition~\ref{def:principal-compression},
\(\mathscr H^{[m,N]}\) is obtained from \(\mathscr H\) by retaining exactly
the rows and columns in those blocks.  The convergent for the tail beginning
at block \(m\), truncated after block \(N\), was defined above by
\begin{equation*}
	\Phi_m^{[N]}(\tau)
	\coloneq
	\left[
	\left(I-\tau\mathscr H^{[m,N]}\right)^{-1}
	\right]_{\mathscr V_m\times\mathscr V_m}.
\end{equation*}
The subscript \(m\) records the first retained block, and the bracketed
superscript \([N]\) records the last one.  The matrix \(\Phi_m^{[N]}\) is the
\(\mathscr V_m\times\mathscr V_m\) block of the resolvent of the finite
matrix \(\mathscr H^{[m,N]}\).
Schur complementation from the last block upwards gives
\begin{align}
	\label{eq:finite-tail-terminal}
	\Phi_N^{[N]}(\tau)
	&=
	\left(I-\tau\mathsf A_N\right)^{-1},
	\\
	\label{eq:finite-tail-backward-recursion}
	\Phi_m^{[N]}(\tau)
	&=
	\left(
	I-\tau\mathsf A_m
	-\tau^2\mathsf B_m
	\Phi_{m+1}^{[N]}(\tau)
	\mathsf C_m
	\right)^{-1},
	\qquad 0\le m<N.
\end{align}
Thus \(\Phi_0^{[N]}\) is obtained by evaluating the matrix continued
fraction from the terminal value \eqref{eq:finite-tail-terminal} back to the
initial height block.

The matrices inverted during this backward evaluation are the successive
Schur-complement pivots.  Recording them separately identifies a common
denominator of the convergent and relates its possible poles to the spectrum
of the finite height section.  For fixed \(0\le m\le N\), whenever the
required inverses exist, define these pivots backwards by
\begin{align}
	\label{eq:finite-CF-denominator-terminal}
	\mathsf D_N^{[N]}(\tau)
	&\coloneq I-\tau\mathsf A_N,
	\\
	\label{eq:finite-CF-denominator-recursion}
	\mathsf D_k^{[N]}(\tau)
	&\coloneq
	I-\tau\mathsf A_k
	-\tau^2\mathsf B_k
	\bigl(\mathsf D_{k+1}^{[N]}(\tau)\bigr)^{-1}
	\mathsf C_k,
	\qquad m\le k<N.
\end{align}
The superscript \([N]\) records that the recursion starts at the terminal
block \(N\).

\begin{proposition}[Schur-complement denominators and determinant factorization]
	\label{prop:finite-CF-determinants}
	Fix \(0\le m\le N\), and assume that the denominator matrices defined
	above exist.  Then
	\(\Phi_k^{[N]}=(\mathsf D_k^{[N]})^{-1}\), and
	\begin{equation}
		\label{eq:finite-CF-determinant-factorization}
		\det\left(I-\tau\mathscr H^{[m,N]}\right)
		=
		\prod_{k=m}^{N}\det\mathsf D_k^{[N]}(\tau).
	\end{equation}
	Consequently, every entry of \(\Phi_m^{[N]}\) is rational in \(\tau\);
	its poles are contained among the reciprocals of the nonzero eigenvalues
	of \(\mathscr H^{[m,N]}\).
\end{proposition}

\begin{proof}
	Eliminate the height blocks in the order
	\(\mathscr V_N,\mathscr V_{N-1},\ldots,\mathscr V_m\).
	The last diagonal block of
	\(I-\tau\mathscr H^{[m,N]}\) is
	\(I-\tau\mathsf A_N=\mathsf D_N^{[N]}\).  After this block has been
	eliminated, its contribution to the preceding diagonal block is
	\[
		(-\tau\mathsf B_{N-1})
		(\mathsf D_N^{[N]})^{-1}
		(-\tau\mathsf C_{N-1})
		=
		\tau^2\mathsf B_{N-1}
		(\mathsf D_N^{[N]})^{-1}\mathsf C_{N-1}.
	\]
	The new diagonal block at level \(N-1\) is therefore
	\(\mathsf D_{N-1}^{[N]}\).  Eliminating successively the blocks
	\(\mathscr V_{N-1},\ldots,\mathscr V_{m+1}\) shows
	inductively that the pivot block at level \(k\) is exactly the matrix
	\(\mathsf D_k^{[N]}\) defined in
	\eqref{eq:finite-CF-denominator-recursion}.

	The leading block of the inverse of the remaining tail
	\(\mathscr V_k\oplus\cdots\oplus\mathscr V_N\) is the inverse of its
	first Schur complement.  Hence
	\(\Phi_k^{[N]}=(\mathsf D_k^{[N]})^{-1}\).  These block eliminations also give
	the block factorization
	\[
	I-\tau\mathscr H^{[m,N]}
	=
	\mathcal F_{m,N}^{\mathrm{up}}(\tau)
	\operatorname{diag}\left(
	\mathsf D_m^{[N]}(\tau),\ldots,\mathsf D_N^{[N]}(\tau)
	\right)
	\mathcal F_{m,N}^{\mathrm{low}}(\tau),
	\]
	where \(\mathcal F_{m,N}^{\mathrm{up}}\) and
	\(\mathcal F_{m,N}^{\mathrm{low}}\) are respectively
	upper and lower block unitriangular.  Their determinants equal \(1\).
	Taking determinants leaves only the determinants of the pivot blocks and
	therefore proves
	\eqref{eq:finite-CF-determinant-factorization}.  Cramer's rule applied to
	the finite matrix \(I-\tau\mathscr H^{[m,N]}\) now shows that every
	entry of the inverse, and in particular every entry of
	\(\Phi_m^{[N]}\), is rational in \(\tau\).  Its scalar denominator divides
	\(\det(I-\tau\mathscr H^{[m,N]})\).  The zeros of this determinant
	are the reciprocals of the nonzero eigenvalues of the finite height
	section.  An entry may have cancellations between its numerator and
	denominator, so the set of its poles is contained in, but need not equal,
	that reciprocal spectrum.
\end{proof}

The determinant factorization describes the rational function
\(\Phi_m^{[N]}\), but not how accurately it approximates the tail resolvent
\(\Phi_m\).  The path expansion identifies the coefficient at the earliest
degree at which a return through the omitted block can occur.  That
coefficient may vanish because the required block path need not lift to a
path between individual height states; the precise structural condition is
recorded after the proposition.

\begin{proposition}[Coefficient at the first possible truncation-error degree]
	\label{prop:finite-CF-leading-error}
	As a formal power series at the origin,
	\begin{equation}
		\label{eq:finite-CF-leading-error}
		\Phi_m(\tau)-\Phi_m^{[N]}(\tau)
		=
		\tau^{2(N-m+1)}
		\mathsf B_m\mathsf B_{m+1}\cdots\mathsf B_N
		\mathsf C_N\mathsf C_{N-1}\cdots\mathsf C_m
		+
		\mathrm{O}\left(\tau^{2(N-m+1)+1}\right),
		\qquad \tau\longrightarrow0.
	\end{equation}
	If the displayed block product is nonzero, then
	\(\Phi_m-\Phi_m^{[N]}\) has a zero of exact order \(2(N-m+1)\) at
	\(\tau=0\).
\end{proposition}

\begin{proof}
	The coefficient of \(\tau^k\) in \(\Phi_m\) is the
	\(\mathscr V_m\times\mathscr V_m\) block of
	\((\mathscr H^{\langle m\rangle})^k\).  It is the sum of the
	weights of all length-\(k\) paths that start and finish in
	\(\mathscr V_m\) and remain in the tail.  The corresponding coefficient
	of \(\Phi_m^{[N]}\) is the same sum restricted to paths that never leave
	\(\mathscr V_m\cup\cdots\cup\mathscr V_N\).

	It follows that a term in the difference must come from a path that
	visits \(\mathscr V_{N+1}\).  Since \(\mathscr H\) is block
	tridiagonal, the block index changes by at most one at each step.  The
	path must therefore make at least \(N-m+1\) upward block transitions to reach
	\(\mathscr V_{N+1}\), and at least \(N-m+1\) downward block
	transitions to return to \(\mathscr V_m\).  No deleted path can have length
	smaller than \(2(N-m+1)\).

	At length \(2(N-m+1)\), there is no time for a horizontal step or a reversal
	before reaching the omitted block \(\mathscr V_{N+1}\).  Hence there is only one possible
	sequence of block indices:
	\[
	m,m+1,\ldots,N,N+1,N,\ldots,m+1,m.
	\]
	With the convention that matrix products are read in the chronological
	order of the block transitions, the outward segment has weight
	\(\mathsf B_m\mathsf B_{m+1}\cdots\mathsf B_N\), and the returning
	segment has weight
	\(\mathsf C_N\mathsf C_{N-1}\cdots\mathsf C_m\).  Their product is
	therefore the coefficient displayed in
	\eqref{eq:finite-CF-leading-error}.  Every other deleted path has at
	least one further step and contributes only to powers
	\(\tau^{2(N-m+1)+1}\) and higher.  This proves the expansion.  If the displayed
	product is nonzero, the coefficient of \(\tau^{2(N-m+1)}\) cannot vanish, so
	the difference has a zero of exact order \(2(N-m+1)\) at \(\tau=0\).
\end{proof}

\begin{remark}[When the first possible omitted coefficient vanishes]
	\label{rem:leading-omitted-product-vanishing}
	Put \(s=N-m+1\).  The coefficient displayed in
	Proposition~\ref{prop:finite-CF-leading-error} follows the block-index path
	\[
	m,m+1,\ldots,m+s,m+s-1,\ldots,m.
	\]
	Its nonvanishing requires the following structural condition, independently
	of the values of the coefficient entries:
	\begin{equation*}
		s\,|P-Q|<w,
		\qquad w=\max(P,Q).
	\end{equation*}
	Indeed, if the product is nonzero, at least one labeled path follows the
	displayed sequence of blocks.  Write its initial height as \(mw+x\), with
	\(0\le x<w\).  If \(P\ge Q\), each of the \(s\) rises of size \(Q\)
	must cross one block boundary, which is possible exactly when
	\(x\ge s(P-Q)\).  If instead \(Q>P\), the rises cross automatically, whereas
	the \(s\) falls of size \(P\) cross successively exactly when
	\(x+s(Q-P)<w\).  The existence of such an integer \(x\) forces the
	displayed inequality.

	Conversely, suppose that all coefficient entries are strictly positive,
	that \(s|P-Q|<w\), and that
	\begin{equation*}
		(m+1)w>(p-1)Q.
	\end{equation*}
	It remains to lift a corresponding reduced-height path to the labeled graph.  If
	\(Q\ge P\), choose \(x=0\).  The required initial state exists because
	following normalizing entries from \((1,0)\) produces a state at every
	height in \(Q\mathbb N_0\).  If instead \(P>Q\), choose \(x=P-1\), so the
	initial height is \(h_0=(m+1)P-1\).  Every residue modulo \(R\) occurs
	among the shifts \(\omega_a\).  A matching shift different from \(pQ\)
	is at most \((p-1)Q\), and hence at most \(h_0\) by the additional
	hypothesis.  If the only matching shift is \(pQ\), then
	\(h_0\equiv pQ\pmod R\) and \(h_0\ge(p-1)Q\); since
	\(pQ-(p-1)Q=Q<R\), these relations also force \(h_0\ge pQ\).
	Thus the initial height again carries a labeled state.

	All outward rises use normalizing entries and therefore exist.  The sources
	of the returning falls lie in the blocks
	\(\mathscr V_{m+1},\ldots,\mathscr V_{m+s}\), so their heights are at
	least \((m+1)w>(p-1)Q\).  A lower-factor state with sequence index
	\(n=0\) has height \(\omega_a=(a-1)Q\le(p-1)Q\); hence no returning
	fall can be lost through a nonexistent lower-factor subdiagonal entry.
	Upper-factor falls use diagonal entries and also exist at \(n=0\).
	Consequently, the reduced-height path lifts to a labeled path, and strict
	positivity prevents cancellation in the block product.  This proves
	sufficiency.  Without the additional inequality, it remains possible for a
	returning fall to start at a lower-factor state with \(n=0\), and the
	height condition alone need not be sufficient.

	For both \((p,q)=(3,2)\) and \((p,q)=(2,3)\), one has \(w=3\) and
	\(|P-Q|=1\), so the reduced-height condition leaves only \(s=1,2\) as
	possible.  At \(m=0\), however, the additional inequality above fails.
	Direct inspection at \(m=N=0\) in both examples shows that every omitted
	return of length \(2\) or \(3\) has zero weight.  In block notation, the
	length-two pattern is \(0,1,0\).  The three possible length-three patterns
	that visit the omitted block are
	\[
	0,0,1,0,
	\qquad
	0,1,1,0,
	\qquad
	0,1,0,0.
	\]
	The first and third contain the vanishing length-two return as a
	subproduct, and direct multiplication shows that the middle pattern also
	has zero weight.  On the other hand, the length-four pattern
	\(0,1,2,1,0\) has nonzero weight.  Equivalently, the displayed block
	product is zero for \(s=1\) and nonzero for \(s=2\), but the explicit
	length-three check is needed when the retained section ends at \(N=0\).
	Thus, when \(m=N=0\),
	\(\Phi_0-\Phi_0^{[0]}\) has a zero of exact order \(4\) at \(\tau=0\) in
	both systems.  For \(s\ge3\),
	the unconditional necessary condition fails, so the displayed coefficient
	vanishes and the actual first omitted term occurs at a higher degree.
	Proposition~\ref{prop:finite-CF-leading-error} identifies the coefficient at
	the first \emph{possible} degree and gives exactness when this coefficient
	is nonzero.
\end{remark}

\subsection{Analytic convergence and positivity}
\label{subsec:analytic-convergence-positivity}

The preceding proposition is a statement about formal power series.  Analytic
convergence on the chosen sequence space follows when the height-ordered
operator is bounded, equivalently when the scalar entries in the given
bidiagonal factorization are uniformly bounded.  In the nonnegative case, a
positive eigenvector gives a second useful normalization in which the
infinity norm is known exactly.

Define the uniform size of the nontrivial factor entries by
\[
	\Gamma_{\mathrm{fac}}
	\coloneq
	\sup\Bigl\{
	|\ell_n^{(a)}|,\ |u_n^{(b)}|:
	n\ge0,\ 1\le a\le p,\ 1\le b\le q
	\Bigr\}
	\in[0,\infty].
\]

\begin{lemma}[Boundedness in terms of the bidiagonal coefficients]
	\label{lem:height-operator-norm-bound}
	The following conditions are equivalent:
	\begin{enumerate}[label=\textnormal{(\roman*)}]
		\item \(\Gamma_{\mathrm{fac}}<\infty\);
		\item \(\mathscr C\) is bounded on
		\(\ell^s(\{1,\ldots,M\}\times\mathbb N_0)\) for one, and hence for
		every, \(1\le s\le\infty\);
		\item \(\mathscr H\) is bounded on
		\(\ell^s(\mathscr V_{p,q})\) for one, and hence for every,
		\(1\le s\le\infty\).
	\end{enumerate}
	If the bidiagonal factorization is nonnegative, these conditions are also
	equivalent to boundedness of \(T\) on \(\ell^s(\mathbb N_0)\) for one,
	and hence for every, \(1\le s\le\infty\).
	When these conditions hold,
	\begin{equation*}
		\Gamma_{\mathrm{fac}}
		\le
		\|\mathscr H\|_s
		=
		\|\mathscr C\|_s
		\le
		1+\Gamma_{\mathrm{fac}},
		\qquad 1\le s\le\infty,
	\end{equation*}
	and the original recurrence matrix is bounded on \(\ell^s\), with
	\begin{equation*}
		\|T\|_s
		\le
		\|\mathscr H\|_s^{M}
		\le
		(1+\Gamma_{\mathrm{fac}})^{M}.
	\end{equation*}
\end{lemma}

\begin{proof}
	In the cyclic ordering, each row and each column contains at most two
	nonzero entries: one has modulus \(1\), and the other is one of the
	bidiagonal coefficients.  Hence
	\(\|\mathscr C\|_\infty\le1+\Gamma_{\mathrm{fac}}\) and
	\(\|\mathscr C\|_1\le1+\Gamma_{\mathrm{fac}}\).  Interpolation between these two
	bounds gives
	\(\|\mathscr C\|_s\le1+\Gamma_{\mathrm{fac}}\) for every
	\(1\le s\le\infty\).  Height ordering is a coordinate permutation and
	is therefore an isometry on each of these sequence spaces, so
	\(\|\mathscr H\|_s=\|\mathscr C\|_s\).

	Conversely, every number \(\ell_n^{(a)}\) and \(u_n^{(b)}\) occurs as an
	entry of \(\mathscr C\), and therefore also as an entry of
	\(\mathscr H\).  Apply the operator to the corresponding coordinate
	vector.  Since that vector has \(\ell^s\)-norm one, the modulus of the
	selected output coordinate does not exceed the operator norm.  Taking the
	supremum over all factor entries gives
	\(\Gamma_{\mathrm{fac}}\le\|\mathscr H\|_s\).  This proves the
	equivalence and both bounds for the linearization.

	Finally, Lemma~\ref{lem:general-cyclic-linearization} gives
	\(T=\Pi_1\mathscr C^M\Pi_1^{\top}\).  The coordinate projection and
	embedding have norm one, and hence
	\[
		\|T\|_s
		\le\|\mathscr C\|_s^M
		=\|\mathscr H\|_s^M.
	\]
	This proves in particular that boundedness of the factors implies
	boundedness of \(T\), without a sign assumption.

	Suppose conversely that the factorization is nonnegative and that \(T\)
	is bounded on one \(\ell^s\) space.  In the matrix product defining
	\(T\), choose the fixed normalizing superdiagonal entry in every upper factor, the
	entry \(\ell_n^{(a)}\) in \(L_a\), and the fixed normalizing diagonal entry in every
	other lower factor.  This path contributes \(\ell_n^{(a)}\) to the entry
	\(T_{n+1,n+q}\).  Since all other contributions are nonnegative,
	\[
		0\le\ell_n^{(a)}\le T_{n+1,n+q}\le\|T\|_s.
	\]
	Similarly, for \(n\ge q-b\), choose the diagonal entry
	\(u_n^{(b)}\) in \(U_b\), the fixed normalizing superdiagonal entry in every other
	upper factor, and the fixed normalizing diagonal entry in every lower factor.  This
	contributes \(u_n^{(b)}\) to
	\(T_{n-q+b,n+b-1}\), and hence
	\[
		0\le u_n^{(b)}\le T_{n-q+b,n+b-1}\le\|T\|_s.
	\]
	For each \(b\), the omitted indices \(0\le n<q-b\) form a finite set.
	Consequently all factor coefficients are uniformly bounded, so
	\(\Gamma_{\mathrm{fac}}<\infty\).  The equivalence already proved then
	gives boundedness on every \(\ell^s\) space.
\end{proof}

Thus the boundedness condition used below can be checked before the height
blocks are assembled: it is exactly
\(\Gamma_{\mathrm{fac}}<\infty\).  The converse statement involving only
the product \(T\) is false for signed factorizations, because cancellation in
the product can conceal unbounded individual factors.  For example, take
\(p=q=1\), put
\[
\ell_{2k}=k+1,
\qquad
\ell_{2k+1}=\frac1{k+1},
\qquad
u_0=-1,
\qquad
u_n=-\ell_{n-1}\quad(n\ge1).
\]
Then the coefficients of both factors are unbounded, whereas
\(T=LU\) has unit superdiagonal, diagonal
\((-1,0,0,\ldots)\), and subdiagonal entries
\(-1\), \(-1\), \(-2\), \(-1\), \(-3/2\), \(-1\),
\(-4/3\), \(-1\), \(\ldots\);
hence \(T\) is bounded on every \(\ell^s\).  The continued
fraction resolves those individual factors, so its analytic convergence
must retain control of them.  For a nonnegative bidiagonal factorization,
the last part of the lemma shows instead that boundedness may be checked on
the given product \(T\) itself.  The diagonal normalization below gives a
second direct route: after normalization the height-ordered matrix is
row-stochastic and hence automatically bounded.

\begin{definition}[Diagonal normalization associated with a positive eigenvector]
	\label{def:stochastic-height-normalization}
	Let \(\boldsymbol v\) be a strictly positive vector satisfying
	\(T\boldsymbol v=\boldsymbol v\).  Put
	\[
	\boldsymbol v^{[r]}\coloneq F_{r+1}\cdots F_M\boldsymbol v,
	\qquad 0\le r<M,
	\qquad
	\boldsymbol v^{[M]}\coloneq\boldsymbol v,
	\]
	so that \(\boldsymbol v^{[0]}=T\boldsymbol v=\boldsymbol v\).  Assume that
	every intermediate vector \(\boldsymbol v^{[r]}\) is strictly positive,
	and let
\(\mathcal D_r\coloneq\operatorname{diag}(\boldsymbol v^{[r]})\).  The normalized factors are
	\begin{equation*}
		\widehat F_r
		\coloneq
		\mathcal D_{r-1}^{-1}F_r\mathcal D_r,
		\qquad 1\le r\le M.
	\end{equation*}
	The associated normalized cyclic matrix is defined blockwise by
	\[
	(\widehat{\mathscr C})_{r,r+1}=\widehat F_r
	\quad (1\le r<M),
	\qquad
	(\widehat{\mathscr C})_{M,1}=\widehat F_M,
	\]
	with every other block equal to zero.  Its height-ordered form is
	\begin{equation*}
		\widehat{\mathscr H}
		\coloneq
		\mathscr P_h\widehat{\mathscr C}\mathscr P_h^{\top},
	\end{equation*}
	where \(\mathscr P_h\) is the height-ordering permutation introduced above.
	If instead \(T\boldsymbol v=\lambda\boldsymbol v\) with \(\lambda>0\),
	the same construction is applied to \(\lambda^{-1}T\), for example by
	replacing one factor by its multiple \(\lambda^{-1}\).  This auxiliary
	rescaling changes the fixed normalization of that one factor, but it
	does not change its bidiagonal sparsity or any cyclic resolvent identity.
	The height increments depend only on the positions of the two bidiagonal
	entries, so the height blocks are unchanged; the rescaled entry is simply
	retained as an edge weight.
	The resolvent of
	\(T\) is then recovered from that of \(\lambda^{-1}T\) by replacing its
	moment variable \(t\) with \(\lambda t\).  If \(t=\tau^M\), the corresponding
	cyclic variable therefore satisfies
	\(\widetilde\tau^M=\lambda\tau^M\).
\end{definition}

\begin{proposition}[Properties of the diagonal normalization]
	\label{prop:stochastic-height-gauge}
	Use the notation of
	Definition~\ref{def:stochastic-height-normalization}, and assume in addition
	that all factors \(F_1,\ldots,F_M\) are nonnegative.  Then every
	\(\widehat F_r\) is a nonnegative row-stochastic bidiagonal matrix.  The
	normalized matrices are related to the original ones by diagonal conjugation;
	this is an operator similarity whenever the diagonal maps and their inverses
	are bounded.  Moreover,
	\begin{equation*}
		\|\widehat{\mathscr H}\|_{\infty}=1,
		\qquad
		\|\widehat{\mathscr H}\|_{2}\le\sqrt2.
	\end{equation*}
\end{proposition}

\begin{proof}
	The strict positivity of the intermediate vectors makes every \(\mathcal D_r\)
	invertible.  Since diagonal multiplication by positive matrices preserves
	entrywise signs, every entry of \(\widehat F_r\) is nonnegative.

	By construction,
	\[
		F_r\boldsymbol v^{[r]}
		=
		F_rF_{r+1}\cdots F_M\boldsymbol v
		=
		\boldsymbol v^{[r-1]}.
	\]
	Consequently,
	\[
		\widehat F_r\one
		=
		\mathcal D_{r-1}^{-1}F_r\mathcal D_r\one
		=
		\mathcal D_{r-1}^{-1}F_r\boldsymbol v^{[r]}
		=
		\mathcal D_{r-1}^{-1}\boldsymbol v^{[r-1]}
		=
		\one.
	\]
	Hence every row of \(\widehat F_r\) sums to \(1\).  Diagonal
	multiplication cannot create new nonzero entries, so
	\(\widehat F_r\) is still bidiagonal.  This proves that it is a
	nonnegative row-stochastic bidiagonal factor.

	At the \(r\)-th position in the factor cycle, the diagonal change of coordinates is
	\(\mathcal D_r\).  If
	\(\mathscr D=\operatorname{diag}(\mathcal D_0,\ldots,\mathcal D_{M-1})\), checking one
	factor transition at a time gives
	\(\widehat{\mathscr C}
	=\mathscr D^{-1}\mathscr C\mathscr D\).
	The last-to-first transition uses \(\mathcal D_M=\mathcal D_0\), because
	\(\boldsymbol v^{[M]}=\boldsymbol v^{[0]}=\boldsymbol v\).  Height ordering applies the same permutation to
	the rows and columns and therefore carries this identity to the
	height-ordered matrices:
	\[
	\widehat{\mathscr H}
	=
	(\mathscr P_h\mathscr D\mathscr P_h^{\top})^{-1}
	\mathscr H
	(\mathscr P_h\mathscr D\mathscr P_h^{\top}).
	\]
	It is an operator similarity precisely when
	\(\mathscr D\) and \(\mathscr D^{-1}\) act boundedly on the chosen
	sequence space.

	Every row of the normalized cyclic operator sums to \(1\), and a
	permutation preserves row sums.  Hence
	\(\|\widehat{\mathscr H}\|_\infty=1\).
	Every column contains at most two nonzero entries.  Each of them lies in
	\([0,1]\), since it belongs to a stochastic row, so the column sum is at
	most \(2\).  Thus
	\(\|\widehat{\mathscr H}\|_1\le2\), and the standard interpolation
	estimate gives
	\[
	\|\widehat{\mathscr H}\|_2^2
	\le
	\|\widehat{\mathscr H}\|_1
	\|\widehat{\mathscr H}\|_\infty
	\le2.
	\]
\end{proof}

\begin{remark}[Finite resolvents under diagonal normalization]
	\label{rem:finite-diagonal-normalization-resolvents}
	Every fixed finite section is related by a finite-dimensional diagonal
	similarity.  Let
	\(\mathscr D_h\coloneq
	\mathscr P_h\mathscr D\mathscr P_h^{\top}\), and let
	\(\mathscr D_{h,[m,N]}\) be its restriction to the coordinates in
	\(\mathscr V_{[m,N]}\).  Then
	\begin{equation*}
	\widehat{\mathscr H}^{[m,N]}
	=
	\mathscr D_{h,[m,N]}^{-1}\mathscr H^{[m,N]}
	\mathscr D_{h,[m,N]},
	\end{equation*}
	and therefore
	\begin{equation*}
	\left(I-\tau\widehat{\mathscr H}^{[m,N]}\right)^{-1}
	=
	\mathscr D_{h,[m,N]}^{-1}
	\left(I-\tau\mathscr H^{[m,N]}\right)^{-1}
	\mathscr D_{h,[m,N]}.
	\end{equation*}
	Thus, for fixed input and output coordinates, the normalized and
	unnormalized finite resolvent entries differ by fixed positive diagonal
	factors, independent of \(N\) once those coordinates are retained.  The
	same relation holds coefficientwise for the formal infinite resolvents.
	Consequently, convergence of normalized finite entries transfers to the
	corresponding unnormalized entries after this fixed rescaling.  A bounded
	operator similarity on the whole sequence space additionally requires
	boundedness of \(\mathscr D_h\) and \(\mathscr D_h^{-1}\).  The Pi\~neiro
	experiment below is carried out directly with \(\widehat{\mathscr H}\).
\end{remark}

The preceding lemma translates the operator hypothesis into the scalar data
supplied by the factorization.  The theorem is stated with the actual operator
norm because that gives the sharper estimate; the directly computable bound
\(\|\mathscr H\|_s\le1+\Gamma_{\mathrm{fac}}\) may always be substituted.
The shortest-path argument below gives both the exact coefficient range and
an error bound given by a geometric series.

\begin{theorem}[Finite convergents and a norm error bound]
	\label{thm:finite-CF-convergence}
	Let \(0\le m\le N\).  The following identity holds for the formal power
	series at the origin:
	\begin{equation}
		\label{eq:finite-CF-coefficient-exactness}
		\Phi_m(\tau)-\Phi_m^{[N]}(\tau)
		=
		\mathrm{O}\left(\tau^{2(N-m+1)}\right)
		\qquad(\tau\longrightarrow0).
	\end{equation}
	Equivalently, the two power series have identical coefficients through
	degree \(2(N-m+1)-1\).

	Suppose in addition that \(\mathscr H\) defines a bounded operator on
	\(\ell^s(\mathscr V_{p,q})\), for some
	\(1\le s\le\infty\).  If
	\(\theta\coloneq|\tau|\|\mathscr H\|_s<1\), then
	\begin{equation}
		\label{eq:finite-CF-geometric-error}
		\left\|
		\Phi_m(\tau)-\Phi_m^{[N]}(\tau)
		\right\|_s
		\le
		\frac{2\theta^{2(N-m+1)}}{1-\theta}.
	\end{equation}
	In particular, for every fixed \(m\), the convergents
	\(\Phi_m^{[N]}\) converge locally uniformly to the tail resolvent
	\(\Phi_m\) in
	\(|\tau|<\|\mathscr H\|_s^{-1}\).
	Consequently, a convergence disc expressed only in the entries of the
	given factorization is
	\begin{equation*}
		|\tau|<\frac{1}{1+\Gamma_{\mathrm{fac}}}.
	\end{equation*}
	More precisely, if
	\(\eta\coloneq|\tau|(1+\Gamma_{\mathrm{fac}})<1\), then
	\begin{equation*}
		\left\|
		\Phi_m(\tau)-\Phi_m^{[N]}(\tau)
		\right\|_s
		\le
		\frac{2\eta^{2(N-m+1)}}{1-\eta}.
	\end{equation*}
\end{theorem}

\begin{proof}
	Consider first a coefficient of degree \(k\) in the two resolvent
	expansions.  It is the corresponding block of
	\((\mathscr H^{\langle m\rangle})^k\) or of
	\((\mathscr H^{[m,N]})^k\).  Since the operator is
	block tridiagonal, a product contributing to the
	\(\mathscr V_m\times\mathscr V_m\) block changes the block index by at
	most one at each multiplication.

	A product that is present in the infinite tail but absent from the finite
	compression must visit \(\mathscr V_{N+1}\).  Starting at
	\(\mathscr V_m\), reaching \(\mathscr V_{N+1}\) requires at least
	\(N-m+1\) outward block transitions.  Returning to
	\(\mathscr V_m\) requires at least the same number of inward transitions.
	Consequently, such a product has length at least
	\(2(N-m+1)\).  The coefficients of all lower degrees therefore agree,
	which proves \eqref{eq:finite-CF-coefficient-exactness}.

	For the norm estimate, expand both resolvents in Neumann series.  Principal
	compression does not increase the subordinate norm, so
	\(\|\mathscr H^{[m,N]}\|_s\le\|\mathscr H\|_s\).  The equality of
	the coefficients below degree \(2(N-m+1)\) gives
	\begin{equation*}
		\left\|
		\Phi_m(\tau)-\Phi_m^{[N]}(\tau)
		\right\|_s
		\le
		2\sum_{k=2(N-m+1)}^\infty
		|\tau|^k\|\mathscr H\|_s^k
=
		\frac{2\theta^{2(N-m+1)}}{1-\theta}.
	\end{equation*}
	The bound is uniform when \(\theta\) ranges over a compact subinterval of
	\([0,1)\), and local uniform convergence follows.  By
	Lemma~\ref{lem:height-operator-norm-bound},
	\(\theta\le\eta\).  Since
	\(x\mapsto2x^{2(N-m+1)}/(1-x)\) is increasing on \([0,1)\), replacing
	\(\theta\) by \(\eta\) gives the final estimate stated in the theorem.
\end{proof}

The coefficient bound provides a convergence disc directly from the input
factorization, but it need not exhaust the resolvent set.  Beyond that disc,
convergence follows if the inverses of the finite sections are uniformly
bounded.  The same stability principle also applies to a closed unbounded
realization of the height matrix.  In that case the finite sequences must be
a core, so that vectors in the operator domain and their images can both be
approximated by finite sequences.  The theorem below assumes uniform
stability of the explicit finite matrices assembled from the bidiagonal
factors.  This is the standard stability condition in the finite-section
method for infinite matrices; see \cite{Lindner2006}.

\begin{theorem}[Convergence under a uniform finite-section resolvent bound]
	\label{thm:finite-CF-stable-domain}
	Fix \(m\) and \(1\le s<\infty\), and put
	\(X_m\coloneq\ell^s(\mathscr V_{\ge m})\).  Let \(\mathcal H_m\) be a
	closed, densely defined realization on \(X_m\) of the matrix
	\(\mathscr H^{\langle m\rangle}\), and assume that
	\(c_{00}(\mathscr V_{\ge m})\) is a core for \(\mathcal H_m\).
	Let \(K\) be a compact subset of
	\begin{equation*}
	\Omega_m
	\coloneq
	\left\{\tau\in\mathbb C:
	I-\tau\mathcal H_m:D(\mathcal H_m)\longrightarrow X_m
	\textnormal{ is bijective with bounded inverse}\right\}.
	\end{equation*}
	Suppose that the finite-section inverses are uniformly bounded on \(K\):
	\begin{equation}
		\label{eq:finite-CF-uniform-stability}
		\sup_{\substack{N\ge m\\\tau\in K}}
		\left\|
		\left(I-\tau\mathscr H^{[m,N]}\right)^{-1}
		\right\|_s
		<\infty.
	\end{equation}
	Here the norm in \eqref{eq:finite-CF-uniform-stability} is the induced
	operator norm on \(\ell^s(\mathscr V_{[m,N]})\).  For \(\tau\in K\), put
	\begin{equation*}
	\Phi_m(\tau)
	\coloneq
	\left[(I-\tau\mathcal H_m)^{-1}\right]_{
	\mathscr V_m\times\mathscr V_m}.
	\end{equation*}
	Then
	\begin{equation}
		\label{eq:finite-CF-stable-domain-convergence}
		\sup_{\tau\in K}
		\left\|\Phi_m(\tau)-\Phi_m^{[N]}(\tau)\right\|_s
		\longrightarrow0.
	\end{equation}
	More generally, for any fixed finite state sets
	\(I,J\subset\mathscr V_{\ge m}\), one has, once \(N\) contains both sets,
	\begin{equation}
		\label{eq:finite-CF-stable-window-convergence}
		\sup_{\tau\in K}
		\left\|
		\left[\left(I-\tau\mathscr H^{[m,N]}\right)^{-1}\right]_{I\times J}
		-
		\left[\left(I-\tau\mathcal H_m\right)^{-1}\right]_{I\times J}
		\right\|_s
		\longrightarrow0.
	\end{equation}
	The norm in \eqref{eq:finite-CF-stable-domain-convergence} is the induced
	operator norm on the fixed space \(\ell^s(\mathscr V_m)\); in
	\eqref{eq:finite-CF-stable-window-convergence} it is the induced norm from
	\(\ell^s(J)\) to \(\ell^s(I)\).
	Thus \(\Phi_m^{[N]}\to\Phi_m\) locally uniformly on every compact part of
	\(\Omega_m\) where the finite-section inverses are uniformly bounded.
	In particular, the assertion applies to the bounded realization
	\(\mathcal H_m=\mathscr H^{\langle m\rangle}\) considered above.
\end{theorem}

\begin{proof}
	Let \(\mathcal P_N:X_m\to\ell^s(\mathscr V_{[m,N]})\) be the coordinate
	projection, and let \(\mathcal J_N\) extend a finite sequence by zero.
	Fix a finitely supported vector \(u\), and set
	\[
		x(\tau)
		\coloneq
		(I-\tau\mathcal H_m)^{-1}u,
		\qquad \tau\in K.
	\]
	The set \(\{x(\tau):\tau\in K\}\) is compact in the graph norm
	\(\|x\|_s+\|\mathcal H_mx\|_s\).  Indeed, the resolvent is continuous in
	operator norm on \(\Omega_m\), and, for \(\tau\ne0\),
	\[
	\mathcal H_m(I-\tau\mathcal H_m)^{-1}
	=
	\tau^{-1}\bigl((I-\tau\mathcal H_m)^{-1}-I\bigr).
	\]
	If \(0\in K\), the definition of \(\Omega_m\) forces
	\(D(\mathcal H_m)=X_m\); the closed graph theorem then makes
	\(\mathcal H_m\) bounded, and the same conclusion follows directly.

	Because the finite sequences form a core, a compact set in the graph norm
	can be approximated uniformly by finitely many finite sequences.  Thus,
	for every \(\varepsilon>0\), there are
	\(y_1,\ldots,y_r\in c_{00}(\mathscr V_{\ge m})\) such that, for every
	\(\tau\in K\), one of these vectors satisfies
	\begin{equation*}
	\|x(\tau)-y_j\|_s+
	\|\mathcal H_m(x(\tau)-y_j)\|_s<\varepsilon.
	\end{equation*}
	Choose \(N\) large enough that all the \(y_j\) are supported in
	\(\mathscr V_{[m,N]}\).  Since the matrix is banded, increasing \(N\), if
	necessary, also ensures
	\[
	\mathscr H^{[m,N]}y_j=\mathcal P_N\mathcal H_my_j
	\qquad(1\le j\le r).
	\]

	Let
	\[
		x_N(\tau)
		\coloneq
		\left(I-\tau\mathscr H^{[m,N]}\right)^{-1}\mathcal P_Nu.
	\]
	For a vector \(y_j\) chosen as above, the two resolvent equations give
	\begin{equation*}
	\left(I-\tau\mathscr H^{[m,N]}\right)
	\bigl(x_N(\tau)-y_j\bigr)=
	\mathcal P_N(I-\tau\mathcal H_m)
	\bigl(x(\tau)-y_j\bigr).
	\end{equation*}
	Put \(C_K\) equal to the supremum in
	\eqref{eq:finite-CF-uniform-stability} and
	\(R_K\coloneq\sup_{\tau\in K}|\tau|\).  It follows that
	\[
	\|x_N(\tau)-y_j\|_s
	\le C_K(1+R_K)\varepsilon.
	\]
	Consequently,
	\[
	\sup_{\tau\in K}
	\|\mathcal J_Nx_N(\tau)-x(\tau)\|_s
	\le \bigl(1+C_K(1+R_K)\bigr)\varepsilon.
	\]
	As \(\varepsilon\) is arbitrary, the zero-extended finite resolvent
	columns converge uniformly on \(K\) for every finitely supported input
	vector.  Take the inputs successively in the coordinate basis indexed by
	\(J\), and project the outputs onto the coordinates in \(I\).  Since both
	sets are finite, convergence of these columns is equivalent to convergence
	in the induced matrix norm.  This proves
	\eqref{eq:finite-CF-stable-window-convergence}; taking
	\(I=J=\mathscr V_m\) gives
	\eqref{eq:finite-CF-stable-domain-convergence}.
\end{proof}

Recall from Definition~\ref{def:principal-compression} that
\(\mathscr V_{\ge m}=\bigcup_{k\ge m}\mathscr V_k\).  The sets in this
union are pairwise disjoint because they correspond to
disjoint height intervals.
The argument takes place on \(\ell^2(\mathscr V_{\ge m})\); its orthonormal
basis is formed by the labeled states in the height blocks
\(\mathscr V_m,\mathscr V_{m+1},\ldots\).  On this space, the numerical range
gives a concrete region where the bounds in
Theorem~\ref{thm:finite-CF-stable-domain} hold uniformly.  For a bounded
operator \(A\) on this space, its numerical range is
\[
W(A)\coloneq
\{\langle Ax,x\rangle:\ x\in\ell^2(\mathscr V_{\ge m}),\ \|x\|_2=1\}.
\]

\begin{corollary}[A numerical-range convergence domain]
	\label{cor:finite-CF-numerical-range}
	Assume that \(\mathscr H^{\langle m\rangle}\) is bounded on
	\(\ell^2(\mathscr V_{\ge m})\).  Then, for fixed \(m\), the convergents
	\(\Phi_m^{[N]}(\tau)\) converge locally uniformly to \(\Phi_m(\tau)\) as
	\(N\to\infty\) on the set
	\[
	\tau=0
	\qquad\textnormal{or}\qquad
		\tau^{-1}\notin
		\overline{W\bigl(\mathscr H^{\langle m\rangle}\bigr)}.
		\]
\end{corollary}

\begin{proof}
	The numerical range of every principal compression is contained in the
	numerical range of the full tail.  If a compact set of values
	\(z=\tau^{-1}\) stays at distance \(\delta>0\) from its closure, then
	the full operator \(zI-\mathscr H^{\langle m\rangle}\) and every finite
	compression are invertible, by the standard numerical-range resolvent
	estimate, and
	\[
	\left\|(zI-\mathscr H^{[m,N]})^{-1}\right\|_2
	\le\frac1\delta
	\]
	uniformly in \(N\).  Since
	\(I-\tau\mathscr H=\tau(zI-\mathscr H)\), for \(\tau\ne0\), it follows that
	\[
	\left\|(I-\tau\mathscr H^{[m,N]})^{-1}\right\|_2
	\le \frac{1}{|\tau|\delta}.
	\]
	On a compact subset of the stated domain that does not contain \(0\),
	\(|\tau|\) is bounded away from zero, so this is the uniform bound
	required by Theorem~\ref{thm:finite-CF-stable-domain}.  A neighborhood of
	\(\tau=0\) is covered directly by
	Theorem~\ref{thm:finite-CF-convergence}.
\end{proof}

When the entries are nonnegative, comparison of path weights yields more than
norm convergence.  It enlarges the basic disc from the operator norm to the
spectral radius and gives monotonicity of the convergents together with the
\(M\)-matrix property of the denominator blocks.

\begin{definition}[\(Z\)-matrices and nonsingular \(M\)-matrices]
	\label{def:Z-and-M-matrices}
	A real square matrix is a \(Z\)-matrix if all its off-diagonal entries are
	nonpositive.  It is a nonsingular \(M\)-matrix if it is a \(Z\)-matrix,
	is invertible, and its inverse is entrywise nonnegative.
\end{definition}

This is the standard inverse-positivity characterization of a nonsingular
\(M\)-matrix; see \cite{BermanPlemmons1994}.

\begin{theorem}[Convergence and monotonicity for nonnegative coefficients when
\texorpdfstring{$|\tau|\rho(\mathscr H)<1$}{|tau| rho(H) < 1}]
	\label{thm:positive-CF-convergents}
	Assume that \(\mathscr H\) is a bounded operator on
	\(\ell^s(\mathscr V_{p,q})\), for some \(1\le s\le\infty\), and that
	all its entries are nonnegative.  Let
	\(\varrho\coloneq\rho(\mathscr H)\) be its spectral radius on
	\(\ell^s(\mathscr V_{p,q})\).  If
	\(|\tau|\varrho<1\), then the convergents \(\Phi_m^{[N]}\) converge locally
	uniformly to \(\Phi_m\).  More precisely, choose any
	\(\gamma>\varrho\) such that \(|\tau|\gamma<1\), and put
	\[
	C_\gamma\coloneq
	\sup_{k\ge0}\frac{\|\mathscr H^{\,k}\|_s}{\gamma^k}<\infty.
	\]
	Then
	\begin{equation}
		\label{eq:positive-CF-spectral-radius-error}
		\left\|\Phi_m(\tau)-\Phi_m^{[N]}(\tau)\right\|_s
		\le
		C_\gamma\frac{(|\tau|\gamma)^{2(N-m+1)}}{1-|\tau|\gamma}.
	\end{equation}
	For real \(0\le\tau<\varrho^{-1}\), every finite convergent is
	nonnegative and, for fixed \(m\),
	\begin{equation}
		\label{eq:positive-CF-monotone-convergence}
		0\le
		\Phi_m^{[N]}(\tau)
		\le
	\Phi_m^{[N+1]}(\tau)
	\le
	\Phi_m(\tau)
	\end{equation}
	entrywise.  In the same real range, the denominator matrices in the
	backward recursion, defined in
	\eqref{eq:finite-CF-denominator-terminal}--
	\eqref{eq:finite-CF-denominator-recursion}, are nonsingular
	\(M\)-matrices in the sense of
	Definition~\ref{def:Z-and-M-matrices}.
\end{theorem}

\begin{proof}
	First consider the constant \(C_\gamma\).  By the spectral-radius formula,
	\(\lim_{k\to\infty}\|\mathscr H^{\,k}\|_s^{1/k}=\varrho\).
	Since \(\gamma>\varrho\), one has
	\(\|\mathscr H^{\,k}\|_s\le \gamma^k\) for all sufficiently large \(k\),
	up to a fixed multiplicative constant.  Enlarging that constant to cover
	the finitely many remaining powers gives \(C_\gamma<\infty\).

	Next compare the infinite and finite coefficients.  Every entry of a
	power of a nonnegative matrix is a sum of nonnegative path weights.
	Every path contained in the finite section is also a path of the
	infinite tail, with unchanged weight.  Hence, after embedding the
	finite block in the infinite sequence space,
	\[
		0\le
		\bigl[(\mathscr H^{[m,N]})^k\bigr]_
		{\mathscr V_m\times\mathscr V_m}
		\le
		\bigl[(\mathscr H^{\langle m\rangle})^k\bigr]_
		{\mathscr V_m\times\mathscr V_m}
	\]
	entrywise.  For nonnegative matrices, entrywise domination implies
	domination of every induced \(\ell^s\)-norm: for an arbitrary vector
	\(x\), apply the two matrices to \(\lvert x\rvert\).  The norm of the
	right-hand block is therefore bounded by
	\(\|\mathscr H^{\,k}\|_s\).  By
	Proposition~\ref{prop:finite-CF-leading-error}, the difference of the
	two series has no term below degree \(2(N-m+1)\).  Therefore
	\[
	\left\|\Phi_m-\Phi_m^{[N]}\right\|_s
	\le
	\sum_{k=2(N-m+1)}^{\infty}
	|\tau|^k\|\mathscr H^{\,k}\|_s
	\le
	C_\gamma\sum_{k=2(N-m+1)}^{\infty}(|\tau|\gamma)^k,
	\]
	which proves \eqref{eq:positive-CF-spectral-radius-error}.

	For real \(\tau\ge0\), all terms of the finite and infinite Neumann series are
	nonnegative.  Passing from the cut at \(N\) to the cut at \(N+1\)
	retains every old path and adds the paths that visit the new terminal
	block.  Passing to the infinite tail adds all remaining paths.  This
	gives each of the entrywise inequalities in
	\eqref{eq:positive-CF-monotone-convergence}.  The estimate just proved
	also shows that the increasing limit is \(\Phi_m\), locally
	uniformly for \(|\tau|\varrho<1\).

	The matrix \(I-\tau\mathscr H^{[m,N]}\) is a \(Z\)-matrix in the sense of
	Definition~\ref{def:Z-and-M-matrices}, and its
	inverse is a convergent nonnegative Neumann series, because the path
	comparison above also gives
	\(\rho(\mathscr H^{[m,N]})\le\varrho\).  Indeed, after zero extension,
	every power of the compression is entrywise dominated by the corresponding
	power of the full nonnegative operator.  The preceding \(\ell^s\)-norm
	comparison and the spectral-radius formula then give the stated
	inequality.  Thus \(\tau\rho(\mathscr H^{[m,N]})<1\), and the finite
	resolvent is entrywise nonnegative.  This proves that
	\(I-\tau\mathscr H^{[m,N]}\) is a nonsingular \(M\)-matrix.

	For any \(k\in\{m,\ldots,N\}\), the matrix
	\(\mathsf D_k^{[N]}\) is the Schur complement obtained after eliminating
	the blocks \(k+1,\ldots,N\).  Its off-diagonal entries are nonpositive:
	the original diagonal block is a \(Z\)-matrix and the correction
	\(\tau^2\mathsf B_k(\mathsf D_{k+1}^{[N]})^{-1}\mathsf C_k\) is
	nonnegative.  Its inverse is
	\(\Phi_k^{[N]}\), the leading block of a nonnegative finite resolvent.
	Hence \(\mathsf D_k^{[N]}\) is a nonsingular \(M\)-matrix for every
	level \(k\), as asserted.
\end{proof}

\begin{corollary}[Convergence after row-stochastic normalization]
	\label{cor:stochastic-normalized-CF-convergence}
	Under the hypotheses of
	Proposition~\ref{prop:stochastic-height-gauge}, form the tail resolvents
	\(\widehat\Phi_m\) and their finite convergents
	\(\widehat\Phi_m^{[N]}\) from \(\widehat{\mathscr H}\).  Then, for
	\(\lvert\tau\rvert<1\),
	\begin{equation}
		\label{eq:stochastic-normalized-CF-error}
		\left\|
		\widehat\Phi_m(\tau)-\widehat\Phi_m^{[N]}(\tau)
		\right\|_{\infty}
		\le
		\frac{2\lvert\tau\rvert^{2(N-m+1)}}{1-\lvert\tau\rvert}.
	\end{equation}
	Consequently, for every fixed \(m\), the convergents tend locally
	uniformly to \(\widehat\Phi_m\) throughout the unit disc under the stated
	stochastic-normalization hypotheses.  For real
	\(0\le\tau<1\), the convergence is entrywise monotone and every
	denominator in the backward recursion is a nonsingular \(M\)-matrix.

	The same normalized operator is bounded on \(\ell^2\).  Hence the
	convergence also holds locally uniformly when
	\[
		\tau=0
		\qquad\textnormal{or}\qquad
		\tau^{-1}\notin
		\overline{W\bigl(\widehat{\mathscr H}^{\langle m\rangle}\bigr)}.
	\]
	In this second domain, the uniform bound on all finite-section inverses is
	a consequence of the numerical-range estimate rather than a separate
	hypothesis.
\end{corollary}

\begin{proof}
	Proposition~\ref{prop:stochastic-height-gauge} gives
	\(\|\widehat{\mathscr H}\|_\infty=1\).  Applying
	Theorem~\ref{thm:finite-CF-convergence} in \(\ell^\infty\) gives
	\eqref{eq:stochastic-normalized-CF-error}.  Moreover,
	\(\widehat{\mathscr H}\one=\one\), so its spectral radius on
	\(\ell^\infty\) is one.  The entrywise monotonicity and the
	\(M\)-matrix assertion now follow from
	Theorem~\ref{thm:positive-CF-convergents}.

	The same proposition gives
	\(\|\widehat{\mathscr H}\|_2\le\sqrt2\).  Thus
	Corollary~\ref{cor:finite-CF-numerical-range} applies and proves the last
	assertion, including the uniform finite-section resolvent bound used in
	its proof.
\end{proof}

\subsection{Pad\'e-type approximation, backward evaluation, and residual control}
\label{subsec:Pade-backward-residual}

The coefficient agreement proved above is now translated into approximation
language.  The term ``Pad\'e-type'' is used because the denominator is fixed
by the finite height section rather than chosen through a separate left or
right matrix-denominator normalization.

\begin{definition}[Pad\'e-type order]
	\label{def:matrix-Pade-type-order}
	Let \(\mathcal F(\tau)\) be a matrix formal power series.  A rational matrix
	\(\mathcal R(\tau)\) is called a Pad\'e-type approximant of order \(r\) at the
	origin if
	\(\mathcal F(\tau)-\mathcal R(\tau)=\mathrm{O}(\tau^{r+1})\) as \(\tau\to0\).  The qualifier ``type'' records that
	the denominator is prescribed by the finite height section; no left or
	right matrix-denominator normalization is imposed.  Similarly, if
	\(\mathcal F(z)=\sum_{n\ge0}\mathsf M_nz^{-n-1}\), a rational matrix \(\mathcal R(z)\)
	has moment order \(r\) at infinity when
	\(\mathcal F(z)-\mathcal R(z)=\mathrm{O}(z^{-r-1})\) as \(z\to\infty\); equivalently, its first \(r\) moments agree
	with those of \(\mathcal F\).
\end{definition}

\begin{corollary}[Pad\'e-type approximation by the convergents]
	\label{cor:finite-CF-Pade-type}
	The convergent \(\Phi_m^{[N]}\) is a Pad\'e-type approximant of order
	\(2(N-m+1)-1\) to \(\Phi_m\).  A common scalar denominator for its entries
	is \(\det(I-\tau\mathscr H^{[m,N]})\), and this denominator factors
	according to \eqref{eq:finite-CF-determinant-factorization}.
\end{corollary}

\begin{proof}
	The approximation order is
	Proposition~\ref{prop:finite-CF-leading-error}; the denominator statement follows from
	Proposition~\ref{prop:finite-CF-determinants}.
\end{proof}

The block-tridiagonal structure used to obtain the denominator also makes the
approximants inexpensive to evaluate: the full finite section is never
inverted at once.  Instead, the continued fraction requires one bounded-size
matrix inversion per height block.  Put
\(w_0\coloneq\max\{p,q\}\), the eventual order of every height block.

\begin{corollary}[Cost of backward evaluation]
	\label{cor:finite-CF-complexity}
	Every height block has cardinality at most \(\max\{p,q\}\), with equality
	after finitely many initial blocks.  Consequently,
	\(\Phi_0^{[N]}\) can be evaluated from
	\eqref{eq:finite-tail-terminal}--\eqref{eq:finite-tail-backward-recursion}
	using \(N+1\) inversions of matrices of order at most
	\(w_0\).  With standard dense matrix algorithms, each
	inversion requires \(\mathrm{O}(w_0^3)\) arithmetic operations and
	\(\mathrm{O}(w_0^2)\) storage as \(w_0\to\infty\).  Across the \(N+1\)
	height blocks, the complete backward evaluation requires
	\(\mathrm{O}((N+1)w_0^3)\) operations as \(N+w_0\to\infty\), and
	\(\mathrm{O}(w_0^2)\) working storage as \(w_0\to\infty\).  For fixed \(p,q\), the
	operation count is linear in \(N\).  These estimates assume that the
	bidiagonal coefficients are generated as needed.
\end{corollary}

\begin{proof}
	The cardinality statement was proved when the blocks \(\mathscr V_m\)
	were defined.  The operation and storage bounds now follow directly from
	the backward recursion.
\end{proof}

The estimates obtained from the Neumann series above depend on a norm bound
for \(\mathscr H\) and can be
conservative.  A practical stopping criterion should instead measure the
coupling that remains across the truncation boundary.  The finite resolvent
column needed for this residual is obtained recursively from level \(m\) to
level \(N\), after computing the convergents \(\Phi_k^{[N]}\) from level
\(N\) down to level \(m\).  For every \(j\), let \(\mathsf J_j\) be the
zero-filling injection of \(\mathbb C^{\mathscr V_j}\) into the height-ordered
sequence space: \(\mathsf J_jy\) equals \(y\) on
the coordinates in \(\mathscr V_j\) and is zero on every other height block.
When only the tail beginning at \(\mathscr V_m\) is considered, the same
symbol denotes its restriction to that tail.  For \(m\le k\le N\), define the corresponding
finite resolvent column by
\[
\mathsf X_k^{[m,N]}(\tau)
\coloneq
\left[
\left(I-\tau\mathscr H^{[m,N]}\right)^{-1}
\right]_{\mathscr V_k\times\mathscr V_m}.
\]

\begin{proposition}[A posteriori estimate from the truncation boundary]
	\label{prop:finite-CF-a-posteriori}
	Assume that the finite tail matrices
	\[
		I-\tau\mathscr H^{[k,N]},
		\qquad m\le k\le N,
	\]
	are invertible.  Then the finite resolvent column defined above can be recovered
	from the backward recursion without another large
	inversion:
	\begin{align}
		\label{eq:finite-CF-Green-column-initial}
		\mathsf X_m^{[m,N]}&=\Phi_m^{[N]},
		\\
		\label{eq:finite-CF-Green-column-forward}
		\mathsf X_{k+1}^{[m,N]}
		&=
		\tau\Phi_{k+1}^{[N]}\mathsf C_k
		\mathsf X_k^{[m,N]},
		\qquad m\le k<N.
	\end{align}
	If, in addition, \(I-\tau\mathscr H^{\langle m\rangle}\) is invertible on
	\(\ell^s(\mathscr V_{\ge m})\), for some \(1\le s\le\infty\), the exact
	error identity is
	\begin{equation}
		\label{eq:finite-CF-exact-residual-error}
		\Phi_m-\Phi_m^{[N]}
		=
		\tau
		\left[
		\left(I-\tau\mathscr H^{\langle m\rangle}\right)^{-1}
		\right]_{\mathscr V_m\times\mathscr V_{N+1}}
		\mathsf C_N\mathsf X_N^{[m,N]}.
	\end{equation}
	In particular, if
	\(\theta=|\tau|\|\mathscr H\|_s<1\), then
	\begin{equation}
		\label{eq:finite-CF-a-posteriori-bound}
		\left\|\Phi_m-\Phi_m^{[N]}\right\|_s
		\le
		\frac{|\tau|}{1-\theta}
		\left\|\mathsf C_N\mathsf X_N^{[m,N]}\right\|_s.
	\end{equation}
\end{proposition}

\begin{proof}
	The first block of the finite resolvent column is the leading block of the
	finite resolvent, and hence it is \(\Phi_m^{[N]}\).  This proves
	\eqref{eq:finite-CF-Green-column-initial}.

	After the blocks \(\mathscr V_{k+2},\ldots,\mathscr V_N\) have been
	eliminated, the block row at level \(k+1\) of the finite resolvent
	equation reads
	\[
		\mathsf D_{k+1}^{[N]}\mathsf X_{k+1}^{[m,N]}
		-\tau\mathsf C_k\mathsf X_k^{[m,N]}=0.
	\]
	Multiplication by
	\((\mathsf D_{k+1}^{[N]})^{-1}=\Phi_{k+1}^{[N]}\) gives
	\eqref{eq:finite-CF-Green-column-forward}.  Thus the tail resolvents are
	computed from level \(N\) down to level \(m\), and the resolvent column is
	then computed from level \(m\) up to level \(N\).

	Extend the finite resolvent column by zero on
	\(\mathscr V_{N+1},\mathscr V_{N+2},\ldots\), and call the extended
	column \(\widetilde{\mathsf X}^{[m,N]}\).  Its block equations coincide
	with the infinite resolvent equation through level \(N\).  At level
	\(N+1\), the infinite operator still contains the coupling from
	\(\mathscr V_N\), while the extended column has a zero
	\(\mathscr V_{N+1}\)-block.  Therefore
	\[
		(I-\tau\mathscr H^{\langle m\rangle})
		\widetilde{\mathsf X}^{[m,N]}
		=
		\mathsf J_m-
		\tau\mathsf J_{N+1}
		\mathsf C_N\mathsf X_N^{[m,N]},
	\]
		where \(\mathsf J_j\) is the zero-filling coordinate injection defined above.

	The corresponding column of the infinite resolvent is
	\(\mathsf X=(I-\tau\mathscr H^{\langle m\rangle})^{-1}\mathsf J_m\).
	Subtracting the preceding equation from its resolvent equation and
	applying the infinite inverse gives
	\[
		\mathsf X-\widetilde{\mathsf X}^{[m,N]}
		=
		\tau
		(I-\tau\mathscr H^{\langle m\rangle})^{-1}
		\mathsf J_{N+1}\mathsf C_N\mathsf X_N^{[m,N]}.
	\]
	The \(\mathscr V_m\)-block of this identity is
	\eqref{eq:finite-CF-exact-residual-error}.

	Finally, the Neumann-series estimate gives
	\[
	\left\|
	\left(I-\tau\mathscr H^{\langle m\rangle}\right)^{-1}
	\right\|_s
	\le\frac1{1-\theta},
	\]
	because the norm of the tail compression does not exceed the norm used in
	the hypothesis.  Apply this estimate to the exact identity; the
	zero-extension maps \(\mathsf J_j\) and the coordinate projections have
	norm \(1\).  This proves
	\eqref{eq:finite-CF-a-posteriori-bound}.  Every quantity on the
	right-hand side is produced by these two recursions and the known norm
	bound, so the estimate is a computable stopping criterion.
\end{proof}

Classical Pincherle theorems identify the value of a scalar continued fraction
through a distinguished solution of its associated recurrence.  In the scalar
three-term setting, a nonzero solution \((y_n)\) is called \emph{minimal}, or
\emph{recessive}, if \(y_n/z_n\to0\) for one, and hence every, linearly
independent solution \((z_n)\).  Matrix
extensions usually use invertible square matrices in the corresponding transfer
recurrence \cite{Ahlbrandt1996}.  The height blocks here may have different sizes and
	the coupling blocks may be singular.  Nonnegativity gives an alternative
	that does not require those blocks to be inverted: the continued fraction
	selects the least entrywise nonnegative solution of the recurrence below.
	For fixed \(m\) and every \(\tau\) for which the tail resolvent exists, write
	\[
	\mathsf X_k(\tau)
	\coloneq
	\left[
	\left(I-\tau\mathscr H^{\langle m\rangle}\right)^{-1}
	\right]_{\mathscr V_k\times\mathscr V_m},
	\qquad k\ge m.
	\]

	\begin{theorem}[Least nonnegative solution and a matrix Pincherle principle]
		\label{thm:matrix-Pincherle-principle}
		Assume that \(\mathscr H\) is a nonnegative bounded operator on
		\(\ell^s(\mathscr V_{p,q})\), for some \(1\le s\le\infty\), and let
		\(0\le\tau<\rho(\mathscr H)^{-1}\), where the spectral radius is taken
		on this space.  Fix \(m\).  The blocks defined above have dimensions
		\(
		\mathsf X_k(\tau)
		\in\mathbb C^{|\mathscr V_k|\times|\mathscr V_m|}
		\), so every product in the recurrence below has the stated rectangular
		dimensions.
	Then \((\mathsf X_k)_{k\ge m}\) is the least entrywise nonnegative
	solution of
	\begin{align}
		\label{eq:matrix-Pincherle-boundary}
		(I-\tau\mathsf A_m)\mathsf X_m
		-\tau\mathsf B_m\mathsf X_{m+1}
		&=I,
		\\
		\label{eq:matrix-Pincherle-recurrence}
		-\tau\mathsf C_{k-1}\mathsf X_{k-1}
		+(I-\tau\mathsf A_k)\mathsf X_k
		-\tau\mathsf B_k\mathsf X_{k+1}
		&=0,
		\qquad k>m.
	\end{align}
	Moreover,
	\begin{equation}
		\label{eq:matrix-Pincherle-tail-ratio}
		\mathsf X_m=\Phi_m,
		\qquad
		\mathsf X_{m+1}
		=
		\tau\Phi_{m+1}\mathsf C_m\Phi_m.
	\end{equation}
	For each fixed \(k\ge m\), the corresponding finite blocks with the terminal condition
	\(\mathsf X_{N+1}^{[m,N]}=0\) increase entrywise to
	\(\mathsf X_k\) as \(N\to\infty\).
\end{theorem}

\begin{proof}
	The block rows of
	\[
	\left(I-\tau\mathscr H^{\langle m\rangle}\right)
	\begin{bNiceMatrix}
	\mathsf X_m\\
	\mathsf X_{m+1}\\
	\Vdots
	\end{bNiceMatrix}
	=
	\begin{bNiceMatrix}
	I\\
	0\\
	\Vdots
	\end{bNiceMatrix}
	\]
	are exactly \eqref{eq:matrix-Pincherle-boundary} and
	\eqref{eq:matrix-Pincherle-recurrence}.  The first identity in
	\eqref{eq:matrix-Pincherle-tail-ratio} is the definition of \(\Phi_m\).
	For the second identity, split the tail after \(\mathscr V_m\).  The
	lower-left block of the inverse of the resulting \(2\times2\) block
	operator is
	\[
		\tau
		(I-\tau\mathscr H^{\langle m+1\rangle})^{-1}
		\mathsf C_m\Phi_m.
	\]
	Projecting its output onto the first block of the smaller tail replaces
	the resolvent by its leading block \(\Phi_{m+1}\).  Hence
	\(\mathsf X_{m+1}=\tau\Phi_{m+1}\mathsf C_m\Phi_m\).

	It remains to prove the least-solution assertion.  If
	\((\mathsf Y_k)_{k\ge m}\) is any other nonnegative solution, the
	boundary equation and the recurrence can be written jointly as
	\(\mathsf Y=\mathsf J_m+\tau\mathscr H^{\langle m\rangle}
	\mathsf Y\), where \(\mathsf J_m\) is the block column with \(I\) in its
	first position.  This identity and its finite iterations are legitimate
	coordinatewise without assuming that the entries of \(\mathsf Y\) form a
	bounded operator: each row of the banded matrix
	\(\mathscr H^{\langle m\rangle}\) contains only finitely many nonzero
	entries.  Iterating this identity \(r+1\) times gives
	\[
		\mathsf Y
		=
		\sum_{j=0}^{r}\tau^j
		\bigl(\mathscr H^{\langle m\rangle}\bigr)^j\mathsf J_m
		+
		\tau^{r+1}
		\bigl(\mathscr H^{\langle m\rangle}\bigr)^{r+1}
		\mathsf Y.
	\]
	The last term is nonnegative, so \(\mathsf Y\) dominates every partial
	Neumann sum.  After zero extension from the tail space to
	\(\ell^s(\mathscr V_{p,q})\), every power of
	\(\mathscr H^{\langle m\rangle}\) is entrywise dominated by the
	corresponding power of \(\mathscr H\).  Monotonicity of the induced
	\(\ell^s\)-norm under entrywise domination and the spectral-radius
	formula therefore give
	\(\rho(\mathscr H^{\langle m\rangle})\le\rho(\mathscr H)\).
	For \(\tau=0\) the assertion is immediate.  If \(\tau>0\), choose
	\(\gamma\) with
	\(\rho(\mathscr H)<\gamma<\tau^{-1}\).  The spectral-radius formula,
	applied as in the proof of
	Theorem~\ref{thm:positive-CF-convergents}, gives a constant
	\(C_\gamma\) such that
	\[
		\bigl\|(\mathscr H^{\langle m\rangle})^j\bigr\|_s
		\le C_\gamma\gamma^j,
		\qquad j\in\mathbb N_0.
	\]
	Therefore the Neumann series converges in operator norm on the tail
	space, and its first block column is \(\mathsf X\).
	Letting \(r\) tend to infinity yields
	\(\mathsf Y\ge\mathsf X\), so \(\mathsf X\) is the least nonnegative
	solution.

	Finally, the terminal condition
	\(\mathsf X_{N+1}^{[m,N]}=0\) removes exactly the paths that cross the
	boundary after \(\mathscr V_N\).  The finite resolvent column therefore sums
	the paths confined to the blocks \(m,\ldots,N\).  Increasing \(N\)
	retains every old path and adds new paths with nonnegative weights.
	Theorem~\ref{thm:positive-CF-convergents} identifies their increasing
	limit with the infinite resolvent column.  This is the matrix Pincherle
	principle appropriate here:
	the least nonnegative solution is selected by finite truncations satisfying
	\(\mathsf X_{N+1}^{[m,N]}=0\), without
	requiring the generally rectangular coupling blocks
	\(\mathsf B_k,\mathsf C_k\) to be invertible.
\end{proof}

\section{Branched and \texorpdfstring{$q$}{q}-state matrix continued fractions}
\label{sec:branched-q-state-CF}

This section relates the matrix continued fraction of
Theorem~\ref{thm:matrix-CF-height-blocks} to branched first-return
recursions.  Subsection~\ref{subsec:q1-reduction-branched} treats the
case \(q=1\): replacing the height blocks by single-height blocks recovers the scalar
\(p\)-branched Stieltjes fraction with the original bidiagonal entries as its
weights.  Restricting the \((p+1)\)-step contraction to the first factor copy
gives the direct \(p\)-\L{}ukasiewicz production graph of \(T\) and hence its
branched Jacobi fraction; the remaining factor copies give the analogous
graphs of the cyclic products.  The scalar Stieltjes fraction appears before
this contraction, at the level of the individual factors.
Subsection~\ref{subsec:q-channel-branched-refinement} keeps instead
\(q\) boundary states, formed from \(Q=q/\gcd(p,q)\) consecutive reduced
heights and the \(\gcd(p,q)\) factor labels over them.  Its main result is a
\(q\times q\) matrix first-return recursion valid for arbitrary \(p\) and
\(q\); the scalar fraction and the two-state fraction are derived below as
special cases.

\subsection{Reduction to the branched Stieltjes fraction}
\label{subsec:q1-reduction-branched}

When \(q=1\), the preceding construction contains the \(p\)-branched
Stieltjes continued fraction associated with the displayed bidiagonal
factorization of the \((p+2)\)-term Hessenberg recurrence.  A general
Hessenberg recurrence gives a branched Jacobi fraction; the Stieltjes form
used here comes from the additional factorization.
The matrix continued fraction
of Theorem~\ref{thm:matrix-CF-height-blocks} is constructed after grouping
heights into blocks of width \(w=\max(P,Q)\).  When \(q=1\), one has
\(w=p\).  If this block grouping is kept, the Schur complements in every
finite truncation are still matrix-valued.  The scalar branched fraction is obtained by using the extra
simplification specific to \(q=1\): the height itself
already determines the factor label, and the Schur-complement decomposition
may be refined from height blocks to single-height tails.

Assume \(q=1\). Then
\[
d=1,
\qquad
P=p,
\qquad
Q=1,
\qquad
R=p+1.
\]
The ordered factors are
\[
F_1=L_1,\ldots,F_p=L_p,
\qquad
F_{p+1}=U_1.
\]
The height shifts are
\[
\omega_1=0,
\quad
\omega_2=1,
\quad
\ldots,
\quad
\omega_{p+1}=p.
\]
Thus
\[
h(a,n)=(p+1)n+a-1,
\qquad
a\in\{1,\ldots,p+1\}.
\]
Hence the map
\[
(a,n)\longmapsto (p+1)n+a-1
\]
is a bijection from the labeled height states onto \(\mathbb N_0\).  In this
case the height serves as the label, since it determines both the factor
position and the sequence index.

Define the scalar weights
\begin{equation}
	\label{eq:general-q1-weights}
	\gamma_{(p+1)n}=u_n^{(1)},
	\qquad
	\gamma_{(p+1)n+a}=\ell_n^{(a)},
	\qquad
	a\in\{1,\ldots,p\},
	\quad
	n\in\mathbb N_0.
\end{equation}
Under the height identification, the normalizing entries give rises
\[
h\longmapsto h+1,
\]
whereas the coefficient entries give falls
\[
h+p\longmapsto h
\]
with weight \(\gamma_h\).  Indeed, for \(a\in\{1,\ldots,p\}\) and
	\(n\in\mathbb N_0\), the coefficient entry of \(L_a\) sends
\((a,n+1)\) to \((a+1,n)\).  The source and target heights are
\[
(p+1)(n+1)+a-1=(p+1)n+a+p
\]
and
\[
(p+1)n+a,
\]
respectively.  Thus this transition is the fall from height \(h+p\) to
height \(h\), with \(h=(p+1)n+a\), and its weight is
\(\ell_n^{(a)}=\gamma_h\).  Similarly, if \(h=(p+1)n\), the diagonal entry
of \(U_1\) sends \((p+1,n)\) to \((1,n)\), which is the fall from height
\(h+p\) to height \(h\) with weight \(u_n^{(1)}=\gamma_h\).  Thus the
	projected graph is exactly the factor-resolved weighted \(p\)-Dyck path
	graph, with
rises \(+1\) and falls \(-p\).

The relation between the scalar tails and the resolvents introduced above is
now made explicit.  After the identification of labeled states with heights,
write \(e_k\coloneq e_{(a,n)}\) for the unique basis state satisfying
\(h(a,n)=k\).  In the notation of
Definition~\ref{def:principal-compression}, define the single-height tail by
\[
\mathscr H_{p,1}^{\langle k\rangle}
\coloneq
[\mathscr H_{p,1}]_{
\{j:j\ge k\}\times\{j:j\ge k\}}.
\]
Thus it is obtained by deleting all rows and columns whose height is less
than \(k\).  Here the angle brackets again indicate an infinite tail; the
index is a single height because every height contains one state when
\(q=1\).
Define, with the inverse again interpreted as a formal Neumann series,
\begin{equation}
	\label{eq:q1-G-resolvent-definition-tau}
	\widehat G_k(\tau)
	\coloneq
	e_k^{\top}
	\left(I-\tau\mathscr H_{p,1}^{\langle k\rangle}\right)^{-1}
	e_k.
\end{equation}
Equivalently, \(\widehat G_k(\tau)\) is the diagonal resolvent entry of the
cyclic linearized operator restricted to the tail of heights not below \(k\).
It generates precisely the paths that start at height \(k\), end at height
\(k\), and never go below \(k\).

	Every closed path in this tail contains, for each coefficient fall, exactly \(p\)
	normalizing rises.  Hence its length in the elementary-step variable is a multiple of
\(p+1\).  Therefore there is a unique formal power series \(G_k(t)\) such that
\begin{equation*}
	\widehat G_k(\tau)=G_k(\tau^{p+1}).
\end{equation*}
Thus the standard branched-Stieltjes variable is
\begin{equation}
	\label{eq:q1-tau-t-relation}
	t=\tau^{p+1}.
\end{equation}
In this normalization, \(t\) marks one coefficient fall, equivalently one
complete circuit through the \(p+1\) factor positions.  On the first factor
copy this circuit represents one application of \(T\); on the other factor
copies it represents the corresponding cyclic products.

This definition also clarifies the relation with the block tails
\(\Phi_m(\tau)\) of Theorem~\ref{thm:matrix-CF-height-blocks}.  For \(q=1\),
the height blocks have width \(p\).  Under the bijection between labeled
states and heights,
\[
\mathscr V_m
\longleftrightarrow
\{mp,mp+1,\ldots,mp+p-1\},
\]
with both sides ordered by increasing height.  Consequently,
\begin{equation*}
	G_{mp}(\tau^{p+1})
	=
	\left(\Phi_m(\tau)\right)_{mp,mp}.
\end{equation*}
For a general height \(k\), however, \(G_k\) is not simply a diagonal entry of
the block tail beginning at \(mp\), because \(\Phi_m\) allows paths to move inside the whole
block \(\mathscr V_m\), possibly below \(k\).  The scalar branched fraction is
obtained by replacing the block tails \(\Phi_m\) with the single-height tail
resolvents \eqref{eq:q1-G-resolvent-definition-tau}.

\begin{theorem}[The case \texorpdfstring{$q=1$}{q=1} and the branched Stieltjes recursion]
	\label{thm:matrix-CF-reduces-to-branched-CF}
	Assume \(q=1\). Then the height model underlying the matrix continued
	fraction reduces, after regrouping the states into single-height tails, to
	the scalar factor-resolved \(p\)-Dyck path graph. The scalar tail functions satisfy the
	\(p\)-branched Stieltjes recursion
	\begin{equation*}
		G_k(t)
		=
		1+\gamma_k t\,G_k(t)G_{k+1}(t)\cdots G_{k+p}(t),
	\end{equation*}
	or, equivalently,
	\begin{equation}
		\label{eq:q1-branched-tail}
		G_k(t)
		=
		\frac{1}{
			1-\gamma_k t\,G_{k+1}(t)\cdots G_{k+p}(t)
		}.
	\end{equation}
	In particular, \(G_0(t)\) is the \(p\)-branched Stieltjes series with
	weights \(\gamma_0,\gamma_1,\ldots\), where the weights are given by
	\eqref{eq:general-q1-weights}.
\end{theorem}

\begin{proof}
	The height identification above shows that, when \(q=1\), the labeled
	height graph is the scalar graph with elementary transitions
	\[
	h\longmapsto h+1,
	\qquad
	h+p\longmapsto h.
	\]
	The first transition has weight \(1\), while the second has weight
	\(\gamma_h\).  Therefore this is precisely the factor-resolved weighted
	\(p\)-Dyck path graph.  Restricting the \((p+1)\)-step contraction to the
	first factor copy gives the \(p\)-\L{}ukasiewicz production graph of \(T\);
	restriction to any other factor copy gives the production graph of the
	corresponding cyclic product.

	A nonempty path counted by \(G_k(t)\) has a unique first return to height
	\(k\).  This return is a fall from height \(k+p\).  In the part of the
	path preceding that fall, mark the last rise from each of the heights
	\(k,k+1,\ldots,k+p-1\).  These marked rises occur successively:
	\[
	k
	\longrightarrow
	k+1
	\longrightarrow
	\cdots
	\longrightarrow
	k+p.
	\]
	The first fall returning to height \(k\) then goes from height \(k+p\) to
	height \(k\).  This fall has weight \(\gamma_k\) and contributes one power
	of \(t\).
	Because this is a first-return decomposition, the path has no excursion
	based at \(k\) before the first marked rise: any nonempty such excursion
	would already return to \(k\) earlier.

	Between two consecutive marked rises the path contains an arbitrary
	excursion based at the current height, but never below that height.  If it
	went below, it would later have to use another rise from the preceding
	height, contrary to the choice of the marked rise.  After the last marked
	rise, the path may make an arbitrary excursion based at \(k+p\) before the
	returning fall; this gives \(G_{k+p}(t)\) and is not a segment between two
	marked rises.  After the fall from \(k+p\) to \(k\), the remainder is an
	arbitrary excursion based at \(k\), giving \(G_k(t)\).  Hence the scalar
	factors which appear are
	\[
	G_{k+1}(t),
	\quad
	G_{k+2}(t),
	\quad
	\ldots,
	\quad
	G_{k+p}(t),
	\quad
	G_k(t).
	\]
	Since these are scalar series, their order is immaterial.  Therefore the
	first-return decomposition gives
	\[
	G_k(t)
	=
	1+\gamma_k t\,G_k(t)G_{k+1}(t)\cdots G_{k+p}(t),
	\]
	where the term \(1\) corresponds to the empty path. Solving for \(G_k(t)\)
	gives \eqref{eq:q1-branched-tail}.
\end{proof}

\subsection{A \texorpdfstring{$q$}{q}-state first-return recursion}
\label{subsec:q-channel-branched-refinement}

The reduction in Subsection~\ref{subsec:q1-reduction-branched} shows that the
case \(q=1\) permits a finer grouping than the general height-block matrix
continued fraction: after passing to single-height tails,
one obtains the scalar \(p\)-branched Stieltjes fraction.  The corresponding
grouping with \(q\) boundary states for arbitrary \(q\) is described next.  The outcome is a
matrix-valued first-return recursion, equivalently a matrix-valued branched
continued fraction with \(q\) boundary states, obtained from
\((I-\tau\mathscr H)^{-1}\).  When \(Q<P\), it groups \(Q\) consecutive
heights instead of the \(P\) consecutive heights used in
Theorem~\ref{thm:matrix-CF-height-blocks}; when \(Q\ge P\), the two
groupings coincide once every admissible factor label is present.  It is called the \(q\)-state
recursion because each tail begins with \(q\) boundary states.

Recall that
\[
        d=\gcd(p,q),
        \qquad
        P=\frac p d,
        \qquad
        Q=\frac q d,
        \qquad
        R=P+Q.
\]
The elementary height steps are \(+Q\) for the normalizing entries and
\(-P\) for the coefficient entries.  The height-block continued fraction of
Theorem~\ref{thm:matrix-CF-height-blocks} groups heights in intervals of width
\(\max(P,Q)\), in order to obtain a block-tridiagonal matrix.  Here the
heights are grouped in intervals of width \(Q\).  This produces \(q=dQ\) boundary states:
\(Q\) consecutive reduced heights, with \(d\) possible factor labels at each
height.

For small heights, some labeled states are absent because \(n\ge0\).
Choose an index \(m_*\) such that every admissible label is present for
\(m\ge m_*\).  From
\eqref{eq:general-height-shifts} and from the lower-side formula
\(\omega_a=(a-1)Q\), the largest height shift is
\[
        \max_{1\le a\le M}\omega_a=pQ=qP.
\]
Hence one may take
\[
        m_*=p.
\]
We use this uniform choice throughout.  In the case \(q=1\), the height map
is already a bijection onto
\(\mathbb N_0\), so one may and will take \(m_*=0\).
Indeed, if \(m\ge p\) and
\[
        h\in\{mQ,mQ+1,\ldots,(m+1)Q-1\},
\]
then \(h\ge pQ\).  For any factor label \(a\) satisfying
\(h\equiv\omega_a\pmod R\), the integer
\[
        n=\frac{h-\omega_a}{R}
\]
is nonnegative because \(h\ge\max_a\omega_a\).
Thus no factor label allowed by the residue class is lost at the boundary
\(n=0\).

It remains to count the labels over each reduced height.  Consecutive offsets
satisfy
\[
        \omega_{a+1}\equiv\omega_a+Q\pmod R,
\]
also across the upper-factor part, because \(-P\equiv Q\pmod R\).
Therefore
\[
        \omega_a\equiv(a-1)Q\pmod R.
\]
Since \(\gcd(Q,R)=1\) and \(M=dR\), the list of \(M\) offsets contains
every residue class modulo \(R\) exactly \(d\) times.  Consequently, for
every \(m\ge m_*\), each reduced height in the interval
\[
        mQ,
        mQ+1,
        \ldots,
        (m+1)Q-1
\]
carries exactly \(d\) factor labels.  From this point on set
\[
        \mathcal U_m
        \coloneq
        \operatorname{span}
        \{e_{(a,n)}:mQ\le h(a,n)<(m+1)Q\},
        \qquad
        m\ge m_*.
\]
Then
\[
        \mathcal U_m\simeq\mathbb C^Q\otimes\mathbb C^d,
        \qquad
        \dim \mathcal U_m=q.
\]
Write the reduced-height coordinate as
\[
        c\in\{0,\ldots,Q-1\}
\]
and the remaining coordinate as one of the \(d\) factor labels above that
height.  Thus a boundary state is the pair formed by this coordinate and the
corresponding factor label.  Let
\[
		E_m:\mathbb C^q\longrightarrow \ell^2(\mathscr V_{p,q})
	\]
	be the coordinate isometry whose range is \(\mathcal U_m\), with the
	columns ordered by the boundary-state basis just described.  Thus
	\[
	E_m^{\top}E_m=I_q,
	\qquad
	E_mE_m^{\top}
	\textnormal{ is the coordinate projection onto }\mathcal U_m.
	\]
	To distinguish
this \(Q\)-height grouping from the height blocks of width
\(w=\max\{P,Q\}\) used earlier, define
\[
\mathscr H^{\langle m\rangle_Q}
\coloneq
[\mathscr H]_{
\{(a,n):h(a,n)\ge mQ\}
\times
\{(a,n):h(a,n)\ge mQ\}}.
\]
Equivalently, this compression retains the coordinate spaces
\(\mathcal U_m,\mathcal U_{m+1},\ldots\).  The subscript \(Q\) records that
the cutoff is the height \(mQ\), rather than the left endpoint \(mw\) of
the earlier block \(\mathscr V_m\).
The \(q\)-state tail resolvent is
\begin{equation}
	\label{eq:q-channel-tail-resolvent}
	\mathbf G_m(\tau)
	\coloneq
	E_m^{\top}
	\left(I-\tau\mathscr H^{\langle m\rangle_Q}\right)^{-1}
	E_m
	\in\mathbb C[[\tau]]^{q\times q},
	\qquad
	m\ge m_*.
\end{equation}
Thus the scalar tails of the case \(q=1\) are replaced by matrix-valued tails
with \(q\) boundary states.

Write
\begin{equation}
        \label{eq:P-divided-by-Q}
        P=\kappa Q+\sigma,
        \qquad
        \kappa\in\mathbb N_0,
        \qquad
        \sigma\in\{0,\ldots,Q-1\}.
\end{equation}
A fall of size \(-P\) from reduced-height coordinate \(c\) in
\(\mathcal U_{m+\eta}\) lands at height
\[
        mQ+(\eta-\kappa)Q+c-\sigma.
\]
Therefore it can land in \(\mathcal U_m\) only for
\[
        \eta=\kappa
        \qquad\textnormal{or}\qquad
        \eta=\kappa+1.
\]
More precisely, a fall of depth \(\kappa\) uses the reduced-height coordinates
\(c\ge\sigma\), whereas a fall of depth \(\kappa+1\) uses the coordinates
\(c<\sigma\).  If \(\sigma=0\), only the depth \(\kappa\) case occurs.  If
\(\kappa=0\), the falls of depth zero are internal transitions inside
\(\mathcal U_m\).

Using the source-row and target-column convention fixed above, define
\begin{align}
        \label{eq:q-channel-A-B-D-definitions}
        \mathbf A_m
        &\coloneq
        E_m^{\top}\mathscr HE_m,
&
        \mathbf B_m
        &\coloneq
        E_m^{\top}\mathscr HE_{m+1},
&
	\mathbf D_m^{(\eta)}
        &\coloneq
        E_{m+\eta}^{\top}\mathscr HE_m,
        \qquad
        \eta\in\{\kappa,\kappa+1\}.
\end{align}
All three matrices are square of order \(q\):
\[
\mathbf A_m,\ \mathbf B_m,\ \mathbf D_m^{(\eta)}
\in\mathbb C^{q\times q}.
\]
Here \(\mathbf A_m\) contains the one-step transitions that remain inside
\(\mathcal U_m\).  The block \(\mathbf B_m\) contains the normalizing upward
transitions from \(\mathcal U_m\) to \(\mathcal U_{m+1}\).  The block
	\(\mathbf D_m^{(\eta)}\) contains the coefficient falls from
	\(\mathcal U_{m+\eta}\) back to \(\mathcal U_m\).  When the corresponding subset
	of boundary states is empty, this block is the zero matrix.  The nonzero
	entries of \(\mathbf B_m\) are equal to \(1\) for the canonical factor
	normalization fixed at the beginning of
	Subsection~\ref{subsec:integer-height-rational-steps}, whereas the nonzero entries of
	\(\mathbf D_m^{(\eta)}\) are precisely the bidiagonal coefficients
	\(\ell_n^{(a)}\) and \(u_n^{(b)}\) selected by the factor label.  After a
	diagonal or eigenvalue normalization, the same definition applies but the
	nonzero entries of \(\mathbf B_m\) need not equal \(1\).

The diagonal block \(\mathbf A_m\) can be nonzero only when \(\kappa=0\).
To see this, recall that the only transitions inside or below one of the
spaces \(\mathcal U_j\), each formed from \(Q\) consecutive heights,
which are not normalizing rises are the coefficient falls of size \(-P\).  If
\(\kappa\ge1\), then \(P\ge Q\), and every such fall starting in
\(\mathcal U_m\) lands in some \(\mathcal U_j\) with \(j<m\); it cannot remain in
\(\mathcal U_m\).  The normalizing transitions go from \(\mathcal U_m\) to
\(\mathcal U_{m+1}\).  Therefore
\[
        \mathbf A_m=0,
        \qquad
        \kappa\ge1.
\]
When \(\kappa=0\), the falls with reduced-height coordinate \(c\ge\sigma\)
have depth zero and are precisely the internal coefficient transitions recorded
by \(\mathbf A_m\).

The next identity is the matrix extension of the scalar branched Stieltjes
recursion.  For \(q=1\), all its blocks are scalars and the ordered product of
the tail matrices \(\mathbf G_j\) becomes the product of scalar tails in the
\(p\)-branched continued fraction, after the change of variable
\(t=\tau^{p+1}\).  For \(q>1\), the scalar tails are replaced by
\(q\times q\) resolvent blocks, and their order in the product is essential.
Put
\[
\mathcal I_+
\coloneq
\{\eta\in\{\kappa,\kappa+1\}:\eta\ge1\},
\]
the set of positive return depths.

\begin{theorem}[The \texorpdfstring{$q$}{q}-state matrix extension of the branched Stieltjes recursion]
        \label{thm:q-channel-matrix-BCF}
        For every \(m\ge m_*\), the \(q\)-state tails
        \eqref{eq:q-channel-tail-resolvent} satisfy the following identity of
        matrix-valued formal power series:
        \begin{equation}
        \label{eq:q-channel-matrix-BCF}
        \mathbf G_m(\tau)
        =
        \Bigl(
        I_q
        -\tau\mathbf A_m
        -\sum_{\eta\in\mathcal I_+}
        \tau^{\eta+1}
        \mathbf B_m\mathbf G_{m+1}(\tau)
        \mathbf B_{m+1}\mathbf G_{m+2}(\tau)
        \cdots
        \mathbf B_{m+\eta-1}\mathbf G_{m+\eta}(\tau)
        \mathbf D_m^{(\eta)}
        \Bigr)^{-1}.
        \end{equation}
        When \(\eta=1\), the alternating product in the summand stops after
        \(\mathbf B_m\mathbf G_{m+1}(\tau)\); there is no intermediate
        \(\mathbf B\)-factor.
        If \(\sigma=0\), then \(\mathbf D_m^{(\kappa+1)}=0\).  If
        \(\kappa=0\), the depth-zero falls are included in \(\mathbf A_m\), and
        this is why \(0\notin\mathcal I_+\).
\end{theorem}

\begin{proof}
        A path counted by \(\mathbf G_m\) starts in \(\mathcal U_m\), stays in
        the tail
        \[
                \mathcal U_m\oplus\mathcal U_{m+1}\oplus\mathcal U_{m+2}\oplus\cdots,
        \]
        and ends in \(\mathcal U_m\).  Decompose such paths into primitive
        returns, where a primitive return is a nonempty path whose first return
        to \(\mathcal U_m\) occurs at its last step.

        A primitive return has two possible types.  It may be a one-step
        transition that stays inside \(\mathcal U_m\), and the corresponding
        contribution is \(\tau\mathbf A_m\).  Otherwise, it must first leave
	        \(\mathcal U_m\) by the upward normalizing transition
        \[
                \mathcal U_m\longrightarrow\mathcal U_{m+1},
        \]
        which contributes \(\tau\mathbf B_m\).  There is no excursion before
        this first step: an internal one-step transition would already be a
        return to \(\mathcal U_m\), while a downward transition would leave the
        retained tail.  Either possibility is incompatible with a primitive
        return of the second type.

        Since the only upward block move is from \(\mathcal U_j\) to
        \(\mathcal U_{j+1}\), a primitive return whose final fall comes from
        \(\mathcal U_{m+\eta}\) must pass successively through
        \[
                \mathcal U_{m+1},
                \mathcal U_{m+2},
                \ldots,
                \mathcal U_{m+\eta}.
        \]
        The initial passage from \(\mathcal U_m\) is its last one before the
        final fall, because a primitive return cannot visit \(\mathcal U_m\)
        earlier.  For each \(r\in\{1,\ldots,\eta-1\}\), also mark the last upward passage
        from \(\mathcal U_{m+r}\) to \(\mathcal U_{m+r+1}\) before the final
        fall.  After the marked passage from the preceding block, the path
        cannot visit \(\mathcal U_{m+r-1}\) again; otherwise the marked
        passage from that block would not have been its last one.  This remains
	        true when a coefficient fall jumps across several groups.  If such a fall
        lands below \(\mathcal U_{m+r-1}\), then, because every upward move
        advances by exactly one group, the path must later cross
        \(\mathcal U_{m+r-1}\to\mathcal U_{m+r}\) in order to reach the final
        falling state, contradicting the choice of the last marked passage.
        The segment before the next marked passage is therefore an arbitrary return to
        \(\mathcal U_{m+r}\) that stays in the tail beginning there.  It is
        summed by \(\mathbf G_{m+r}(\tau)\).

        After the last marked rise into \(\mathcal U_{m+\eta}\), the segment
	        before the final coefficient fall is an arbitrary return based at
        \(\mathcal U_{m+\eta}\) that stays in its tail; it contributes
        \(\mathbf G_{m+\eta}(\tau)\).  Thus this final segment is accounted for
        separately rather than as a segment between two marked rises.  The
        final fall contributes \(\tau\mathbf D_m^{(\eta)}\).  This last-exit
        decomposition is unique.  Hence the contribution of all
        primitive returns whose final fall has depth \(\eta\ge1\) is
        \[
        \tau^{\eta+1}
        \mathbf B_m\mathbf G_{m+1}(\tau)
        \mathbf B_{m+1}\mathbf G_{m+2}(\tau)
        \cdots
        \mathbf B_{m+\eta-1}\mathbf G_{m+\eta}(\tau)
        \mathbf D_m^{(\eta)}.
        \]
	        The exponent \(\eta+1\) counts the \(\eta\) upward normalizing steps and the
	        final coefficient fall.

        By \eqref{eq:P-divided-by-Q}, a single fall can return from
        \(\mathcal U_{m+\eta}\) to \(\mathcal U_m\) only for
        \(\eta=\kappa\) or \(\eta=\kappa+1\).  Therefore the generating matrix of
        primitive nonempty returns is
        \[
        \tau\mathbf A_m
        +
        \sum_{\eta\in\mathcal I_+}
        \tau^{\eta+1}
        \mathbf B_m\mathbf G_{m+1}(\tau)
        \mathbf B_{m+1}\mathbf G_{m+2}(\tau)
        \cdots
        \mathbf B_{m+\eta-1}\mathbf G_{m+\eta}(\tau)
        \mathbf D_m^{(\eta)}.
        \]
        An arbitrary return to \(\mathcal U_m\) is a concatenation of primitive
        returns.  Summing the resulting noncommutative geometric series gives
        \eqref{eq:q-channel-matrix-BCF}.
\end{proof}

The next corollary verifies the reduction to the scalar \(p\)-branched
Stieltjes fraction when \(q=1\).  The case \(q=2\) then gives the first
matrix-valued extension.

\begin{corollary}[Recovery of the scalar branched fraction]
        \label{cor:q-channel-recovers-scalar-BCF}
        For \(q=1\), take \(m_*=0\).  Then the \(q\)-state recursion in
        Theorem~\ref{thm:q-channel-matrix-BCF} is the scalar
        \(p\)-branched Stieltjes fraction of
        Theorem~\ref{thm:matrix-CF-reduces-to-branched-CF}.
\end{corollary}

\begin{proof}
        If \(q=1\), then
        \[
                d=1,
                \qquad
                Q=1,
                \qquad
                P=p.
        \]
        Hence \(\kappa=p\), \(\sigma=0\), and each boundary-state block has
        dimension one.  Moreover \(\mathbf A_m=0\), \(\mathbf B_m=1\), and
        \(\mathbf D_m^{(p)}=\gamma_m\), where \(\gamma_m\) is the weight of the
        fall \(m+p\to m\).  Formula \eqref{eq:q-channel-matrix-BCF} becomes
        \[
                \widehat G_m(\tau)
                =
                \frac{1}{
                1-
                \gamma_m\tau^{p+1}
                \widehat G_{m+1}(\tau)\cdots \widehat G_{m+p}(\tau)
                }.
        \]
        By \eqref{eq:q1-G-resolvent-definition-tau}, these are precisely the
        single-height resolvents introduced above.  Substituting
        \(\widehat G_j(\tau)=G_j(t)\) and setting \(t=\tau^{p+1}\) gives
        exactly the scalar \(p\)-branched Stieltjes recursion.
\end{proof}

\begin{corollary}[The two-state case]
        \label{cor:q2-two-channel-from-general}
        For \(q=2\), Theorem~\ref{thm:q-channel-matrix-BCF} gives a
        \(2\times2\) matrix-valued branched continued fraction with two
        boundary states.  If \(p=2r+1\), then
        \[
                d=1,
                \qquad
                Q=2,
                \qquad
                P=2r+1,
        \]
        and the two return depths are \(r\) and \(r+1\).  If \(p=2r\), then
        \[
                d=2,
                \qquad
                Q=1,
                \qquad
                P=r,
        \]
        and the two boundary states are the two factor labels over each
        reduced height; there is a single return depth \(r\).
        In the endpoint case \(p=1\), which corresponds to \(r=0\) in the
        first alternative, the depth-zero falls are internal transitions and
        are included in \(\mathbf A_m\).  Thus only the depth-one term occurs
        in the sum in \eqref{eq:q-channel-matrix-BCF}.
\end{corollary}

\begin{proof}
        If \(p=2r+1\), then \(\gcd(p,2)=1\), so \(P=2r+1\) and
        \(Q=2\).  The Euclidean division
        \(P=rQ+1\) gives the two possible depths \(r\) and \(r+1\), and
        each group contains \(Q=2\) reduced heights with one label over each.
        When \(r=0\), the explanation following
        \eqref{eq:q-channel-A-B-D-definitions} places the depth-zero falls in
        \(\mathbf A_m\), so they do not appear in the sum over positive
        depths.
        If \(p=2r\), then \(\gcd(p,2)=2\), so \(P=r\) and \(Q=1\).
        Each group then contains one reduced height with two labels over it.
        Since the remainder in the division of \(P\) by \(Q\) is zero, only
        the depth \(r\) occurs.
\end{proof}

\begin{remark}[Relation between the \texorpdfstring{$q$}{q}-state and height-block recursions]
        \label{rem:q-channel-lives-inside-matrix-CF}
		The grouping into \(q\) boundary states and the grouping into height
		blocks give two first-return descriptions of resolvent blocks of
		\(\mathscr H\).  The former is proved by the path decomposition in
		Theorem~\ref{thm:q-channel-matrix-BCF}; the latter is obtained from the
		finite-dimensional Schur complements in
		Theorem~\ref{thm:matrix-CF-height-blocks}.
        Theorem~\ref{thm:matrix-CF-height-blocks}
        uses height blocks \(\mathscr V_m\) large enough to make the operator
        block tridiagonal for arbitrary \((p,q)\).  Theorem~\ref{thm:q-channel-matrix-BCF}
        uses instead the grouping
        \[
                \mathcal U_m\simeq\mathbb C^q,
        \]
        which keeps together the \(q\) boundary states associated with a
        group of \(Q\) consecutive heights.  If \(Q<P\), the two
        constructions group the height states differently.  If
        \(Q\ge P\), equivalently \(q\ge p\), then \(w=Q\), and, for every
        \(m\ge m_*\), the spaces \(\mathcal U_m\) and \(\mathscr V_m\)
        contain exactly the same labeled states.  With the same ordering of
        those states, their coefficient blocks and tail recursions coincide.
        Thus the boundary-state construction is genuinely different from the
        height-block construction only when \(q<p\); the scalar branched
        fraction for \(q=1<p\) is the basic instance.
\end{remark}

\begin{remark}[The \(q\)-state tails and the rectangular Weyl matrix]
	\label{rem:q-state-tail-versus-Weyl}
	The matrices \(\mathbf G_m(\tau)\) are principal \(q\times q\) resolvent
	blocks of tails beginning at \(m\ge m_*\), where every reduced height carries
	all \(d\) factor labels.  The rectangular moment series at the origin is a
	different projection: it observes the first factor copy through
	\(\mathsf Q_q\) and \(\mathsf Q_p\), and the selected states may lie in the
	finite initial part below \(m_*\).  Hence it is not, in general, one of the
	matrices \(\mathbf G_m\).

	Remark~\ref{rem:finite-window-extraction-Darboux} gives an explicit
	finite-window extraction using the width-\(w\) height-block tail
	\(\Phi_{N+1}\).  If the remote tail is instead evaluated through the
	\(q\)-state recursion, it must still be coupled to the initial coordinates
	by a finite Schur-complement calculation and, when \(P>Q\), through the
	finite boundary slab crossed by the downward transitions.
\end{remark}

\section{Rectangular resolvents and Weyl approximation}
\label{sec:general-Weyl-bridge}

This section answers two approximation-theoretic questions.  How can the
rectangular moment series and Weyl matrix of \(T\) be recovered from the
factor-resolved cyclic linearization?  And what do the finite matrix
continued fractions approximate, both when the PBF is bounded and when it is
unbounded?  The first question connects the cyclic resolvent with mixed-type
multiple orthogonality; the second determines the precise analytic or
asymptotic meaning of its finite truncations.

Subsection~\ref{subsec:rectangular-cyclic-projection} proves that rectangular
projections of the cyclic resolvent recover the moment series of \(T\) and of
all its cyclic Darboux transforms.  It also identifies the two block output
matrices associated with the mixed-type recurrence and clarifies their
relation with the scalar production-matrix formalism.

Subsection~\ref{subsec:general-Weyl-spectral-approximation} applies the
mixed-type spectral Favard theory developed in
\cite{BranquinhoFoulquieManas2023Spectral,
BranquinhoFoulquieManas2026Unbounded}.  In the bounded setting it identifies
the relevant rectangular resolvent block with the Stieltjes transform and
transfers the convergence and Pad\'e-type results to finite rectangular
approximants.  In the unbounded setting the Poincar\'e moment approximation
of the measure-defined Weyl matrix is separated from analytic finite-section
convergence: the latter targets the projected resolvent of a specified
closed realization and gives the measure-defined Weyl matrix only under an
additional spectral identification.

\subsection{Rectangular projections of the cyclic resolvent}
\label{subsec:rectangular-cyclic-projection}

After projection onto a fixed copy, only complete cycles of length \(M\)
contribute, and each complete cycle gives one power of the corresponding
Darboux transform of \(T\).

For \(r\in\mathbb N\), let
\[
\mathsf E_{[r]}\coloneq\begin{bmatrix}I_r&0_{r\times\infty}\end{bmatrix}
\]
be the projection onto the first \(r\) coordinates of the original sequence
space. The rectangular moment series associated with \(T\) is the
\(q\times p\) matrix of formal power series
\begin{equation*}
	\mathsf G_{p,q}(t)
	\coloneq
	\mathsf E_{[q]}(I-tT)^{-1}\mathsf E_{[p]}^{\top}
	=
	\sum_{n=0}^{\infty}
	\mathsf E_{[q]}T^n\mathsf E_{[p]}^{\top}t^n,
\end{equation*}
whose coefficients \(\mathsf E_{[q]}T^n\mathsf E_{[p]}^{\top}\) are the
rectangular moments entering the spectral Favard theorem for mixed-type
multiple orthogonality. More generally, for the cyclic Darboux transforms
\(T_{(a)}\), set
\begin{equation*}
	\mathsf G_{p,q}^{(a)}(t)
	\coloneq
	\mathsf E_{[q]}(I-tT_{(a)})^{-1}\mathsf E_{[p]}^{\top},
	\qquad
	a\in\{1,\ldots,M\}.
\end{equation*}

In the production/output formalism of
Deutsch--Ferrari--Rinaldi \cite{DeutschFerrariRinaldi2005}, a
unit-lower-Hessenberg matrix \(\Pi\) is called the production matrix.  Its
output matrix is \(\mathcal O(\Pi)\), with entries
\[
\mathcal O(\Pi)_{n,k}=(\Pi^n)_{0,k}.
\]
Thus \(\mathcal O(\Pi)_{n,k}\) is the total weight of the \(n\)-step paths from state \(0\)
to state \(k\).  Since \(\Pi\) is unit lower Hessenberg,
\(\mathcal O(\Pi)_{n,k}=0\) for \(k>n\) and
\(\mathcal O(\Pi)_{n,n}=1\); hence the output matrix \(\mathcal O(\Pi)\), not the production
matrix \(\Pi\), is lower unitriangular.  In the polynomial setting, the rows
of \(\mathcal O(\Pi)^{-1}\) contain the coefficients of the associated monic polynomial
sequence, and multiplication by the variable is represented by \(\Pi\).  Its
connection with multiple orthogonality and branched continued fractions is
developed in
\cite{PetreolleSokalZhu2023,Sokal2024}.

The same distinction applies in the mixed-type setting.  When \(q>1\), the
recurrence matrix \(T\) is not lower Hessenberg in scalar coordinates, but it
is lower block Hessenberg after consecutive coordinates are grouped into
blocks of size \(q\).  Similarly, \(T^{\top}\) is lower block Hessenberg for
blocks of size \(p\).  Applying the output construction to these two block
production matrices gives two output matrices.  They will be identified below
with \(\mathscr L^{-1}\) and \(\mathscr U^{-\top}\); the latter are output
matrices, not production matrices.

For a semi-infinite matrix \(S\), an integer \(r\ge1\), and a nonsingular
\(r\times r\) matrix \(C\), define its normalized block output matrix by
\begin{equation*}
	\mathcal O_r(S;C)
	\coloneq
	\begin{bmatrix}
		C\mathsf E_{[r]}\\
		C\mathsf E_{[r]}S\\
		C\mathsf E_{[r]}S^2\\
		\vdots
	\end{bmatrix}.
\end{equation*}
Its \(n\)-th block row, numbered from zero, is
\(C\mathsf E_{[r]}S^n\).  When \(r=1\) and \(C=1\), this is the usual
output matrix of a production matrix: its \((n,k)\)-entry is
\((S^n)_{0,k}\).  Deleting the first block row of
\(\mathcal O_r(S;C)\) gives \(\mathcal O_r(S;C)S\), which is the block
analogue of the defining row-production identity.
For later use, write
\[
\boldsymbol\chi_{[r]}(x)
\coloneq
\begin{bNiceMatrix}[margin=2pt]
	I_r&xI_r&x^2I_r&\Cdots[shorten-end=3pt]
\end{bNiceMatrix}^{\top},
\]
and let \(\mathsf\Lambda_{[r]}\) be the block shift characterized by
\(\mathsf\Lambda_{[r]}\boldsymbol\chi_{[r]}(x)
=x\boldsymbol\chi_{[r]}(x)\).

We use the Gauss--Borel formulation of mixed-type multiple orthogonality
developed in \cite{AlvarezFidalgoManas2011}.
Let \(\boldsymbol\mu\) be a \(q\times p\) matrix of measures whose
step-line moment matrix has the Gauss--Borel factorization
\[
\mathscr M
=
\mathscr L^{-1}\mathscr U^{-1},
\]
where \(\mathscr L\) is lower unitriangular and \(\mathscr U\) is upper
triangular.  Suppose that the associated recurrence matrix has the two
representations
\[
T
=
\mathscr L\mathsf\Lambda_{[q]}\mathscr L^{-1}
=
\mathscr U^{-1}\mathsf\Lambda_{[p]}^{\top}\mathscr U.
\]
The notation is chosen to agree with the initial-condition matrices introduced
in \cite[Section~3]{BranquinhoFoulquieManas2023Spectral}.  There, \(\xi\)
collects the initial values of the \(q\) right recursion-polynomial sequences,
whereas \(\nu\) collects those of the \(p\) left recursion-polynomial
sequences.  In the present Gauss--Borel realization these matrices are no
longer independent data: they are the leading blocks of the two triangular
coefficient matrices.  Indeed, since
\[
\boldsymbol B(x)=\mathscr L\boldsymbol\chi_{[q]}(x),
\qquad
\boldsymbol A(x)=\boldsymbol\chi_{[p]}(x)^{\top}\mathscr U,
\]
the matrix \(\boldsymbol B(x)\in
\mathbb C[x]^{\mathbb N_0\times q}\) is semi-infinite by \(q\), whereas
\(\boldsymbol A(x)\in
\mathbb C[x]^{p\times\mathbb N_0}\) is \(p\) by semi-infinite.  Thus
the corresponding initial-condition matrices are
\[
\xi
\coloneq
\mathsf E_{[q]}\mathscr L\mathsf E_{[q]}^{\top},
\qquad
\nu^{\top}
\coloneq
\mathsf E_{[p]}\mathscr U\mathsf E_{[p]}^{\top}.
\]
The transpose in the second identity reflects the row convention for the left
polynomial family.  Consequently, the factors \(\xi^{-1}\) and \(\nu^{-1}\)
in the normalized output matrices remove precisely these initial conditions.

\begin{theorem}[Output matrices of the mixed-type recurrence]
	\label{thm:mixed-block-output-matrices}
	Under the preceding hypotheses,
	\begin{equation}
		\label{eq:mixed-block-output-identities}
		\mathcal O_q(T;\xi^{-1})=\mathscr L^{-1},
		\qquad
		\mathcal O_p(T^{\top};\nu^{-1})=\mathscr U^{-\top}.
	\end{equation}
	Consequently,
	\begin{equation}
		\label{eq:mixed-two-output-factorization}
		\mathscr M
		=
		\mathcal O_q(T;\xi^{-1})
		\mathcal O_p(T^{\top};\nu^{-1})^{\top}.
	\end{equation}
	Moreover, the complete coefficient arrays of the two step-line mixed-type
	polynomial matrices \(\boldsymbol B\) and \(\boldsymbol A\) are respectively
	\[
	\mathscr L=\mathcal O_q(T;\xi^{-1})^{-1},
	\qquad
	\mathscr U^{\top}=\mathcal O_p(T^{\top};\nu^{-1})^{-1}.
	\]
\end{theorem}

\begin{proof}
	The first recurrence representation gives
	\[
	\mathsf\Lambda_{[q]}\mathscr L^{-1}
	=
	\mathscr L^{-1}T.
	\]
	The leading block of the inverse of a lower triangular matrix is the
	inverse of its leading block.  Hence
	\[
	\mathsf E_{[q]}\mathscr L^{-1}
	=
	\xi^{-1}\mathsf E_{[q]}.
	\]
	The \(n\)-th block row of \(\mathscr L^{-1}\) is therefore
	\begin{equation*}
		\mathsf E_{[q]}\mathsf\Lambda_{[q]}^n\mathscr L^{-1}=\mathsf E_{[q]}\mathscr L^{-1}T^n=\xi^{-1}\mathsf E_{[q]}T^n.
	\end{equation*}
	Stacking these block rows proves the first identity in
	\eqref{eq:mixed-block-output-identities}.

	Transposing the second recurrence representation gives
	\[
	\mathsf\Lambda_{[p]}\mathscr U^{-\top}
	=
	\mathscr U^{-\top}T^{\top}.
	\]
	Since the leading \(p\times p\) block of \(\mathscr U^{-\top}\) is
	\(\nu^{-1}\), triangularity gives
	\[
	\mathsf E_{[p]}\mathscr U^{-\top}
	=
	\nu^{-1}\mathsf E_{[p]}.
	\]
	Repeating the preceding block-row argument yields
	\[
	\mathsf E_{[p]}\mathsf\Lambda_{[p]}^n\mathscr U^{-\top}
	=
	\nu^{-1}\mathsf E_{[p]}(T^{\top})^n,
	\]
	which proves the second identity.  Equation
	\eqref{eq:mixed-two-output-factorization} is now the Gauss--Borel
	factorization itself.  Finally, \(\mathscr L\) and \(\mathscr U^{\top}\)
	are the coefficient matrices in the displayed definitions of
	\(\boldsymbol B\) and \(\boldsymbol A\), so the last two identities
	follow by inversion.
\end{proof}

\begin{proposition}[The recurrence and the initial data determine the moment functional]
	\label{prop:recurrence-initial-data-moment-functional}
	Fix a semi-infinite banded matrix \(T\) and nonsingular matrices
	\(\xi\in\mathbb C^{q\times q}\) and
	\(\nu\in\mathbb C^{p\times p}\).  There is a unique
	\(q\times p\) matrix \(\boldsymbol{\mathcal F}_{T;\xi,\nu}\) of linear
	functionals on \(\mathbb C[x]\) whose moments satisfy
	\begin{equation}
		\label{eq:recurrence-initial-data-moment-functional}
		\bigl\langle
		\boldsymbol{\mathcal F}_{T;\xi,\nu},x^n
		\bigr\rangle
		=
		\xi^{-1}\mathsf E_{[q]}T^n
		\mathsf E_{[p]}^{\top}\nu^{-\top},
		\qquad n\in\mathbb N_0.
	\end{equation}
	Thus the triple \((T,\xi,\nu)\) determines the moment functional.  If a
	matrix of measures represents the moments in
	\eqref{eq:recurrence-initial-data-moment-functional}, it represents this
	unique functional.  In an indeterminate moment problem, several matrices of
	measures may represent that functional.  Conversely, suppose that the
	right and left polynomial solutions of the recurrence for \(T\), normalized
	by \(\xi\) and \(\nu\), form a weakly normal mixed-type system for a
	functional \(\boldsymbol{\mathcal G}\).  Then
	\(\boldsymbol{\mathcal G}=\boldsymbol{\mathcal F}_{T;\xi,\nu}\).
\end{proposition}

\begin{proof}
	The monomials form a basis of \(\mathbb C[x]\), and
	\eqref{eq:recurrence-initial-data-moment-functional} specifies the value
	of every entry of the functional on each monomial.  For
	\(f(x)=\sum_{n=0}^N f_nx^n\), linearity therefore forces
	\[
		\bigl\langle\boldsymbol{\mathcal F}_{T;\xi,\nu},f\bigr\rangle
		=
		\sum_{n=0}^N f_n
		\xi^{-1}\mathsf E_{[q]}T^n
		\mathsf E_{[p]}^{\top}\nu^{-\top},
	\]
	which also proves existence.  For the converse, the spectral reconstruction
	identity for a weakly normal mixed-type recurrence gives
	\[
		\bigl\langle\boldsymbol{\mathcal G},x^n\bigr\rangle
		=
		\xi^{-1}\mathsf E_{[q]}T^n
		\mathsf E_{[p]}^{\top}\nu^{-\top},
		\qquad n\in\mathbb N_0;
	\]
	see \cite[Proposition~9.4 and Remark~9.5]
	{BranquinhoFoulquieManas2023Spectral} and its unbounded extension
	\cite{BranquinhoFoulquieManas2026Unbounded}.  The two functionals therefore
	agree on every monomial and hence on \(\mathbb C[x]\).
\end{proof}

For \(q=1\) with monic normalization, \(\xi=1\), and the first identity in
\eqref{eq:mixed-block-output-identities} is exactly Sokal's Hessenberg
production/output identity \cite{Sokal2024}: \(T\) is the production matrix,
\(\mathcal O_1(T;1)=\mathscr L^{-1}\) is its output matrix, and
\(\mathscr L\) is the coefficient matrix of the monic type-II polynomial
sequence.  Thus it is \(T\), not either triangular matrix, that plays the role
of production matrix.  In the mixed-type problem, applying the construction
to \(T\), grouped into blocks of size \(q\), and to \(T^{\top}\), grouped into
blocks of size \(p\), gives the two output matrices associated with the two
polynomial families.

\begin{remark}[Two finite approximations of the Weyl matrix]
	The spectral construction in
	\cite[Proposition~4.1]{BranquinhoFoulquieManas2023Spectral} uses polynomial
	left and right eigenvectors of a principal truncation \(T^{[K]}\).  Its
	determinantal kernel reduces to a quotient involving two consecutive
	characteristic polynomials and controls the Christoffel numbers of that
	truncation, and hence gives a finite spectral approximation of the Weyl
	matrix.

	The continued-fraction convergents studied here provide a second finite
	approximation of the same Weyl matrix.  They retain the height blocks
	\(0,\ldots,N\) of the sparse matrix \(\mathscr H\), rather than a
	principal truncation of \(T\).  Thus both procedures have the same infinite
	spectral target and use different finite matrices.
\end{remark}

For \(a\in\{1,\ldots,M\}\) and \(r\in\mathbb N\), let \(\mathsf Q_r^{(a)}\) denote
the coordinate projection onto the states
\[
(a,0),(a,1),\ldots,(a,r-1)
\]
in the coordinate ordering used for \(\mathscr H\), and write
\(\mathsf Q_r\coloneq \mathsf Q_r^{(1)}\).  These coordinates need not be
consecutive; \(\mathsf Q_r^{(a)}\) selects them through their labels, and
\((\mathsf Q_r^{(a)})^{\top}\) embeds vectors into those coordinates.
Whenever a finite compression contains all the selected states, the same
symbol denotes the restriction of this projection to the retained coordinate
space, and its transpose denotes the corresponding finite-dimensional
embedding.
Proposition~\ref{prop:direct-cyclic-agreement} identifies \(\mathscr H\) with
the matrix of \(\mathscr C\) in this ordering.  Hence the projection
\(\Pi_a\) onto the \(a\)-th copy followed by \(\mathsf E_{[r]}\) becomes
\(\mathsf Q_r^{(a)}\).

\begin{proposition}[Rectangular projection of the matrix continued fraction]
	\label{prop:general-rectangular-projection}
	Let \(t=\tau^M\). Then
	\begin{equation}
		\label{eq:general-rectangular-projection}
		\mathsf G_{p,q}(t)
		=
		\mathsf Q_q(I-\tau\mathscr H)^{-1}\mathsf Q_p^{\top}.
	\end{equation}
	Thus the rectangular moment series of \(T\) is a rectangular block of
	\((I-\tau\mathscr H)^{-1}\), whose tail blocks give the matrix continued fraction
	in Theorem~\ref{thm:matrix-CF-height-blocks}.
\end{proposition}

\begin{proof}
	By Lemma~\ref{lem:general-cyclic-linearization},
	\(\Pi_1(I-\tau\mathscr C)^{-1}\Pi_1^{\top}=(I-\tau^MT)^{-1}\).
	Taking the first \(q\) rows and the first \(p\) columns, and reading the
	coordinate projections in the height-ordered basis, gives
	\eqref{eq:general-rectangular-projection} with \(t=\tau^M\).
\end{proof}

The first-copy projection recovers the moment series of the original
recurrence matrix.  Observing another position in the factor cycle gives the
moment series of the corresponding cyclic Darboux transform, with the
partial product between the two factor positions appearing explicitly in the
rectangular projected resolvent blocks.

\begin{corollary}[Rectangular moment series for the cyclic Darboux transforms]
	\label{cor:Darboux-shifted-rectangular-series}
	Let \(a,b\in\{1,\ldots,M\}\), and let \(s(a,b)\) and \(F_{a\to b}\) be as
	in Proposition~\ref{prop:Darboux-shifted-resolvent-blocks}. Then, with
	\(t=\tau^M\), for \(r,s\in\mathbb N\),
	\begin{align*}
		\mathsf Q_r^{(a)}(I-\tau\mathscr H)^{-1}
		\bigl(\mathsf Q_s^{(b)}\bigr)^{\top}
		&=
		\tau^{s(a,b)}
		\mathsf E_{[r]}F_{a\to b}
		(I-tT_{(b)})^{-1}\mathsf E_{[s]}^{\top}
		\\
		&=
		\tau^{s(a,b)}
		\mathsf E_{[r]}(I-tT_{(a)})^{-1}
		F_{a\to b}\mathsf E_{[s]}^{\top}.
	\end{align*}
	In particular, for \(a=b\), \(r=q\), \(s=p\),
	\begin{equation*}
		\mathsf Q_q^{(a)}(I-\tau\mathscr H)^{-1}
		\bigl(\mathsf Q_p^{(a)}\bigr)^{\top}
		=
		\mathsf G_{p,q}^{(a)}(t).
	\end{equation*}
	Thus the cyclic resolvent \((I-\tau\mathscr H)^{-1}\) contains
	the rectangular moment
	series of all cyclic Darboux transforms
	\(T_{(1)},\ldots,T_{(M)}\), together with the partial products
	\(F_{a\to b}\) connecting two factor positions.
\end{corollary}

\begin{proof}
	Apply Proposition~\ref{prop:Darboux-shifted-resolvent-blocks} and take the
	first \(r\) rows and first \(s\) columns; the two forms correspond to
	\eqref{eq:Darboux-shifted-resolvent-right} and
	\eqref{eq:Darboux-shifted-resolvent-left}.
\end{proof}

\begin{remark}[Coupling the tail to the initial height blocks]
	\label{rem:finite-window-extraction-Darboux}
	The coordinate projections above generally involve initial states
	distributed among several height blocks.  Use the state set
	\(\mathscr V_{[0,N]}=\mathscr V_0\cup\cdots\cup\mathscr V_N\), with
	\(N\) large enough so that the states selected by \(\mathsf Q_r^{(a)}\) and
	\(\mathsf Q_s^{(b)}\) lie in \(\mathscr V_{[0,N]}\). For \(K>N\), apply finite
	Schur complements to \(I-\tau\mathscr H^{[0,K]}\), eliminating the blocks
	\(\mathscr V_K,\ldots,\mathscr V_{N+1}\).  The corresponding projected block
	stabilizes coefficientwise as \(K\to\infty\), and its limit is a submatrix of
	\[
	\left(
	I_{|\mathscr V_{[0,N]}|}
	-
	\tau[\mathscr H]_{\mathscr V_{[0,N]}\times\mathscr V_{[0,N]}}
	-
	\tau^2
	[\mathscr H]_{\mathscr V_{[0,N]}\times\mathscr V_{N+1}}
	\Phi_{N+1}(\tau)
	[\mathscr H]_{\mathscr V_{N+1}\times\mathscr V_{[0,N]}}
	\right)^{-1}.
	\]
	Here \(\Phi_{N+1}\) is the coefficientwise limit of the finite terminal
	blocks \(\Phi_{N+1}^{[K]}\).  Thus the blocks beyond
	\(\mathscr V_{[0,N]}\) are represented by \(\Phi_{N+1}\), while the initial blocks determine which cyclic
	Darboux transform and which initial coordinates are observed. In
	particular, when the
	\(q\)-state tails beginning at \(m\ge m_*\) from
	Subsection~\ref{subsec:q-channel-branched-refinement} are used, these
	initial blocks connect the resolvent blocks \(\mathbf G_m\) to the rectangular projection at
	the origin.
\end{remark}

\subsection{Spectral representation and convergence of the rectangular approximants}
\label{subsec:general-Weyl-spectral-approximation}

We first treat the bounded case, in which the operator resolvent can be
identified directly with a Stieltjes transform.  By
Lemma~\ref{lem:height-operator-norm-bound}, boundedness of
\(\mathscr H\) is equivalent to uniform boundedness of the scalar
bidiagonal coefficients.  For a nonnegative factorization, it is also
equivalent to boundedness of \(T\).  The unbounded case, where the Weyl
matrix is defined from the Favard measure rather than from an operator
resolvent, is treated at the end of the subsection.

Let \(\xi\in\mathbb C^{q\times q}\) and
\(\nu\in\mathbb C^{p\times p}\) be nonsingular.  Whenever \(T\) is a
bounded operator, define its normalized rectangular Weyl matrix by
\begin{equation}
	\label{eq:general-operator-Weyl-definition}
	\mathsf S(z)
	\coloneq
	\xi^{-1}\mathsf E_{[q]}(zI-T)^{-1}
	\mathsf E_{[p]}^{\top}\nu^{-\top},
	\qquad z\in\rho(T).
\end{equation}

\begin{theorem}[Factor-resolved formula for the rectangular Weyl matrix]
	\label{thm:general-Weyl-CF}
	Assume that \(\mathscr H\) is bounded on
	\(\ell^2(\mathscr V_{p,q})\).  If \(z^{-1}=\tau^M\) and
	\(\lvert\tau\rvert<\|\mathscr H\|_2^{-1}\), then
	\begin{equation}
		\label{eq:general-factor-resolved-Weyl}
		\mathsf S(z)
		=
		\frac{1}{z}\,
		\xi^{-1}
		\mathsf Q_q(I-\tau\mathscr H)^{-1}\mathsf Q_p^{\top}
		\nu^{-\top}.
	\end{equation}
	The right-hand side is independent of the choice of the
	\(M\)-th root \(\tau\).
\end{theorem}

\begin{proof}
	Lemma~\ref{lem:height-operator-norm-bound} implies that \(T\) is bounded,
	and \(T\) is a coordinate compression of \(\mathscr H^M\).  Hence
	\[
		\|T\|_2\le \|\mathscr H\|_2^M<|z|,
	\]
	so \(z\in\rho(T)\) and
	\eqref{eq:general-operator-Weyl-definition} is defined.  The Neumann
	series
	for the resolvent of \(\mathscr H\) converges in operator norm.
	Proposition~\ref{prop:general-rectangular-projection} shows that the
	projected powers vanish unless their degree is divisible by \(M\), and
	that the surviving powers satisfy
	\[
		\mathsf Q_q\mathscr H^{Mn}\mathsf Q_p^{\top}
		=
		\mathsf E_{[q]}T^n\mathsf E_{[p]}^{\top}.
	\]
	Thus Proposition~\ref{prop:general-rectangular-projection}, with
	\(z^{-1}=\tau^M\), gives
	\[
		\mathsf E_{[q]}(zI-T)^{-1}\mathsf E_{[p]}^{\top}
		=
		\frac1z\,
		\mathsf Q_q(I-\tau\mathscr H)^{-1}\mathsf Q_p^{\top}.
	\]
	Multiplication by \(\xi^{-1}\) and \(\nu^{-\top}\) proves
	\eqref{eq:general-factor-resolved-Weyl}.
	The projected Neumann series contains only powers \(\tau^{Mn}\), so its
	value does not depend on the chosen \(M\)-th root of \(z^{-1}\).
\end{proof}

By Theorem~\ref{thm:matrix-CF-height-blocks}, the projected resolvent in
\eqref{eq:general-factor-resolved-Weyl} is evaluated by the height-block
matrix continued fraction.  Hence that continued fraction represents the
rectangular Weyl matrix \(\mathsf S\).

We next identify \(\mathsf S\) with a matrix Stieltjes transform.  This step
uses only the following moment identity from mixed-type spectral theory.
Suppose that a compactly supported \(q\times p\) matrix of finite measures
\(\mathrm d\Psi\) satisfies
\begin{equation}
	\label{eq:general-spectral-moments}
	\xi^{-1}\mathsf E_{[q]}T^n\mathsf E_{[p]}^{\top}\nu^{-\top}
	=
	\int x^n\,\mathrm d\Psi(x),
	\qquad
	n\in\mathbb N_0.
\end{equation}
Sufficient Favard hypotheses for the existence of \(\mathrm d\Psi\),
\(\xi\), and \(\nu\) are given in
\cite[Proposition~9.4 and Remark~9.5]{BranquinhoFoulquieManas2023Spectral}.

\begin{proposition}[Stieltjes representation of the Weyl matrix]
	\label{prop:general-Weyl-Stieltjes}
	Assume that \(T\) is bounded and that
	\(\mathrm d\Psi\) satisfies \eqref{eq:general-spectral-moments}.  Then
	\begin{equation}
		\mathsf S(z)
		=
		\int\frac{\mathrm d\Psi(x)}{z-x}
	\end{equation}
	entry by entry on the unbounded connected component of the common
	analytic domain.
\end{proposition}

\begin{proof}
	Choose \(R_\Psi\) such that every entry of \(\mathrm d\Psi\) is supported
	in \(\{x:\lvert x\rvert\le R_\Psi\}\).  For
	\(\lvert z\rvert>\max\{\|T\|_2,R_\Psi\}\), the Neumann series and
	\eqref{eq:general-spectral-moments} give, entry by entry,
	\begin{align*}
		\mathsf S(z)
		&=
		\sum_{n=0}^{\infty}
		\frac{
		\xi^{-1}\mathsf E_{[q]}T^n
		\mathsf E_{[p]}^{\top}\nu^{-\top}
		}{z^{n+1}}=
		\int\sum_{n=0}^{\infty}\frac{x^n}{z^{n+1}}\,
		\mathrm d\Psi(x)
		=
		\int\frac{\mathrm d\Psi(x)}{z-x}.
	\end{align*}
	Compact support justifies the interchange of sum and integral.  Analytic
	continuation proves the stated identity on the unbounded component.
\end{proof}

\begin{remark}[Domains of analytic convergence]
	\label{rem:general-CF-resolvent-domain}
	The norm disc in Theorem~\ref{thm:general-Weyl-CF} is only a sufficient
	domain.  Theorem~\ref{thm:finite-CF-convergence} gives locally uniform
	convergence of the finite tails there, with the estimate
	\eqref{eq:finite-CF-geometric-error}.
	Theorem~\ref{thm:finite-CF-stable-domain} extends convergence to compact
	subsets of the resolvent set on which the finite-truncation inverses are
	uniformly bounded, and
	Corollary~\ref{cor:finite-CF-numerical-range} gives a domain in terms of
	the numerical range.  When \(\mathscr H\) is nonnegative,
	Theorem~\ref{thm:positive-CF-convergents} gives entrywise monotone
	convergence on
	\[
		0\le\tau<\rho(\mathscr H)^{-1}.
	\]
\end{remark}

The preceding Weyl matrix is naturally the radial form of a Cauchy transform
on a star.  The cyclic linearization gives this representation directly from
the radial moment series.  The average over
the rays below is the rotationally invariant discrete Fourier projection described
after Proposition~\ref{prop:cyclic-star-symmetry}.  Recall that
\(\zeta_M=\exp(2\pi\mathrm i/M)\), and put
\[
\Sigma_M
\coloneq
\bigcup_{j=0}^{M-1}\zeta_M^j[0,\infty).
\]

When the diagonal normalization of
Definition~\ref{def:stochastic-height-normalization} is used, let
\(\mathscr D_h=\mathscr P_h\mathscr D\mathscr P_h^{\top}\) and define the
finite diagonal matrices on the selected coordinates by
\[
	\mathcal D_r^{(a)}
	\coloneq
	\mathsf Q_r^{(a)}\mathscr D_h
	\bigl(\mathsf Q_r^{(a)}\bigr)^{\top}.
\]
For every finite truncation containing these coordinates, diagonal
conjugation gives the exact identity
\begin{equation}
	\label{eq:projected-diagonal-normalization}
	\mathsf Q_q^{(a)}(I-\tau\mathscr H^{[0,N]})^{-1}
	\bigl(\mathsf Q_p^{(a)}\bigr)^{\top}
	=
	\mathcal D_q^{(a)}
	\mathsf Q_q^{(a)}(I-\tau\widehat{\mathscr H}^{[0,N]})^{-1}
	\bigl(\mathsf Q_p^{(a)}\bigr)^{\top}
	(\mathcal D_p^{(a)})^{-1}.
\end{equation}
Thus diagonal normalization changes a selected resolvent block only by the
two explicitly displayed coordinate factors.

\begin{definition}[Rotational lift of the matrix of measures]
	\label{def:rotational-lift-matrix-measure}
	Assume that \(\mathrm d\Psi\) is supported on \([0,\infty)\).  Its
	rotational lift to \(\Sigma_M\) is the \(q\times p\) matrix of measures
	\(\mathrm d\widehat\Psi\) defined, entry by entry, by
	\begin{equation}
		\label{eq:rotational-lift-matrix-measure}
		\int_{\Sigma_M}f(\lambda)\,
		\mathrm d\widehat\Psi(\lambda)
		\coloneq
		\frac1M\sum_{j=0}^{M-1}
		\int_{[0,\infty)}
		f\bigl(\zeta_M^jx^{1/M}\bigr)\,
		\mathrm d\Psi(x)
	\end{equation}
	for every continuous compactly supported scalar function \(f\) on
	\(\Sigma_M\).  Here \(x^{1/M}\) denotes the nonnegative root.
\end{definition}

For the measure \(\mathrm d\Psi\) and its rotational lift, write
\[
	\mathsf S_{\Psi}(z)
	\coloneq
	\int_{[0,\infty)}\frac{\mathrm d\Psi(x)}{z-x},
	\qquad
	\widehat{\mathsf S}(\lambda)
	\coloneq
	\int_{\Sigma_M}
	\frac{\mathrm d\widehat\Psi(\zeta)}{\lambda-\zeta}.
\]

\begin{theorem}[Rotational lift of the Weyl matrix]
	\label{thm:cyclic-star-Weyl}
	Let \(\mathrm d\Psi\) be a \(q\times p\) matrix of finite measures
	supported on \([0,\infty)\), with all moments finite, and let
	\(\mathrm d\widehat\Psi\) be its rotational lift.  Then
	\begin{equation}
		\label{eq:cyclic-star-Weyl}
		\widehat{\mathsf S}(\lambda)
		=
		\lambda^{M-1}\mathsf S_{\Psi}(\lambda^M),
		\qquad
		\lambda\in\mathbb C\setminus\Sigma_M.
	\end{equation}
	Moreover,
	\begin{equation}
		\label{eq:cyclic-star-moments}
		\int_{\Sigma_M}\zeta^n\,
		\mathrm d\widehat\Psi(\zeta)
		=
		\begin{cases}
		\displaystyle\int_{[0,\infty)}x^k\,\mathrm d\Psi(x),
		&n=Mk,\\[4pt]
		0,&M\nmid n.
		\end{cases}
	\end{equation}
\end{theorem}

\begin{proof}
	For \(r\ge0\) and \(\lambda^M\ne r^M\), the partial-fraction identity
	\[
	\frac1M\sum_{j=0}^{M-1}
	\frac1{\lambda-\zeta_M^jr}
	=
	\frac{\lambda^{M-1}}{\lambda^M-r^M}
	\]
	follows by taking the logarithmic derivative of
	\(\lambda^M-r^M\).  Apply it entry by entry in
	\eqref{eq:rotational-lift-matrix-measure}, with \(r=x^{1/M}\).  This proves
	\eqref{eq:cyclic-star-Weyl}.
	Finally,
	\[
	\frac1M\sum_{j=0}^{M-1}\zeta_M^{jn}
	=
	\begin{cases}1,&M\mid n,\\0,&M\nmid n,\end{cases}
	\]
	and substitution in the definition of the lifted measure proves
	\eqref{eq:cyclic-star-moments}.
\end{proof}

Theorem~\ref{thm:cyclic-star-Weyl} applies to finite measures with all
moments.  The following operator-resolvent representation adds the
boundedness hypothesis.

\begin{proposition}[Resolvent representation of the star Weyl matrix]
	\label{prop:cyclic-star-resolvent-representation}
	Assume that \(\mathscr H\) is bounded, that \(\mathrm d\Psi\) is compactly
	supported, and that \eqref{eq:general-spectral-moments} holds.  For all
	sufficiently large \(\lvert\lambda\rvert\),
	\begin{align}
		\widehat{\mathsf S}(\lambda)
		&=
		\xi^{-1}\mathsf E_{[q]}\Pi_1
		(\lambda I-\mathscr C)^{-1}
		\Pi_1^{\top}\mathsf E_{[p]}^{\top}\nu^{-\top}
		\label{eq:cyclic-star-Weyl-C}
		\\
		&=
		\frac1\lambda\,
		\xi^{-1}\mathsf Q_q
		(I-\lambda^{-1}\mathscr H)^{-1}
		\mathsf Q_p^{\top}\nu^{-\top}.
		\label{eq:cyclic-star-Weyl-H}
	\end{align}
\end{proposition}

\begin{proof}
	Proposition~\ref{prop:general-Weyl-Stieltjes} identifies
	\(\mathsf S_{\Psi}\) with \(\mathsf S\), and
	Theorem~\ref{thm:general-Weyl-CF} represents \(\mathsf S\) by the
	projected resolvent of \(\mathscr H\).  It remains to use
	\[
		\Pi_1(\lambda I-\mathscr C)^{-1}\Pi_1^{\top}
		=
		\lambda^{M-1}(\lambda^MI-T)^{-1},
	\]
	which follows from Lemma~\ref{lem:general-cyclic-linearization}.
	Taking the first \(q\) rows and \(p\) columns gives
	\eqref{eq:cyclic-star-Weyl-C}; conjugating by the height permutation gives
	\eqref{eq:cyclic-star-Weyl-H}.
\end{proof}

Since all three members of
\eqref{eq:cyclic-star-Weyl-C}--\eqref{eq:cyclic-star-Weyl-H} are analytic,
the equalities extend from large \(\lvert\lambda\rvert\) to the unbounded
connected component of their common domain.  In particular, the
factor-resolved continued fraction represents the rotationally symmetric
\(q\times p\) Weyl matrix without choosing an \(M\)-th root in the
\(\lambda\)-variable.

\begin{remark}[Diagonal normalization and spectral identification]
	For every finite truncation,
	\eqref{eq:projected-diagonal-normalization} gives the exact rescaling of
	the projected resolvent.  If the normalized finite resolvents converge,
	the same identity gives a limit for the original projected entries.
	A specified operator realization identifies this limit with its projected
	resolvent.  A representing measure identifies the same limit with its
	Cauchy transform only when the corresponding
	resolvent--Cauchy-transform equality has been established separately.
	Stochastic normalization supplies convergence under its stated hypotheses.
\end{remark}

\section{Moment systems for the recurrence
\texorpdfstring{\(\mathscr H_{(a)}\)}{H(a)} and cyclic Christoffel families}
\label{sec:balanced-star-Christoffel}

Fix \(a\in\mathcal A_{\mathrm{cyc}}\).  The terminal lower--upper PBF of
\(T_{(a)}\) supplies its Favard data; when the prescribed factors are
Christoffel connection matrices, these are the transformed data of
Theorem~\ref{thm:cyclic-Favard-transport}, including
\(T_{(p+1)}=BA\) with measure \(x\,\mathrm d\Psi\).  The minimum-height
factorization of the same product constructs \(\mathscr H_{(a)}\).  This
section builds its moment system, relates its polynomial solutions to the
cyclic-product recurrences, and derives its radial measures and finite star
spectra.  The common product and projected resolvent link the two
factorizations; the moment functional of \(\mathscr H_{(a)}\) is constructed
from its powers.

The moments of the system associated with \(\mathscr H_{(a)}\) are defined directly from powers of
the matrix \(\mathscr H_{(a)}\).  Its rotational symmetry separates them into congruence
classes modulo \(M=p+q\).  The resulting moment matrix has an explicit
triangular factorization and the matrix \(\mathscr H_{(a)}\) is exactly its recurrence
matrix.  The nonzero moments in each congruence class form a positive
Stieltjes moment sequence on the radial variable.  For fixed compactly
supported radial representing measures, the atom-at-zero and negative
fractional-moment conditions characterize, entry by entry, representation
of the corresponding Fourier components by finite complex measures supported
on the corresponding compact \(M\)-ray star.  In parallel, the cyclic
products are identified with
the recurrence matrices along their cyclic orbit.  The Christoffel
description is developed at three levels.  First, the PBF itself determines
matrix Christoffel transformations adapted to its individual factors and to
the initial conditions chosen at consecutive cyclic products.  Second, the
fixed multipliers \(\mathsf K_q(x)\) and
\(\mathsf K_p(x)^{\top}\) give a more explicit rectangle when its two
one-sided boundary chains are weakly normal.  Third, total positivity of the
scalar moment matrix makes this boundary normality automatic.

Subsection~\ref{subsec:balanced-star-normality-spectra} proves weak normality at
every scalar order in the coprime case.  In general, grouping the
\(d=\gcd(p,q)\) states of each height gives a block moment matrix with
invertible triangular factors whose recurrence is exactly
\(\mathscr H_{(a)}\).  It also locates the
spectra of the finite truncations of \(\mathscr H_{(a)}\) on the \(p+q\) rays.
For fixed compactly supported radial representing measures, it records the
criterion for lifting the algebraic moment functionals to finite complex
measures supported on the corresponding compact star.
Strong normality additionally requires every nontrivial component to attain
its prescribed maximal degree and is established under the separate
maximal-degree hypotheses stated in the individual results and applications.
Subsection~\ref{subsec:PBF-adapted-Christoffel} constructs the factor-adapted
Christoffel transformations along the cyclic PBF orbit and explains the role
of triangular changes of initial conditions.  It proves that a complete
factor circuit gives multiplication by \(x\), up to a constant triangular
normalization at the endpoint, and records the direct transport of the
polynomial recurrence solutions through the bidiagonal factors.
Subsection~\ref{subsec:complete-Christoffel-grid} then studies the special
rectangle obtained by multiplying the matrix of measures with the fixed
matrices \(\mathsf K_q(x)\) and \(\mathsf K_p(x)^{\top}\).  Under weak
normality of the two boundary chains, the bidiagonal connections propagate
weak normality throughout the rectangle and identify the measure, moment
matrix, recurrence matrix, and both polynomial families at every vertex.
If the scalar moment matrix is totally positive, all the shifted moment
matrices are totally positive submatrices, so the boundary hypothesis is
automatic.  Strong normality still requires the separate maximal-degree
conditions stated in the individual results and applications.
Subsection~\ref{subsec:Weyl-approximation-line-star} returns to rational
approximation and derives finite approximants, convergence estimates, and
Pad\'e-type orders for the Weyl matrices on the line and on the star.
Subsection~\ref{subsec:unbounded-PBF-Weyl} states separately the Poincar\'e
moment approximation for a measure-defined Weyl matrix and the conditional
operator-resolvent convergence of stable finite sections.

\subsection{Moment functionals for
\texorpdfstring{\(\mathscr H_{(a)}\)}{H(a)}, normality, and finite spectra}
\label{subsec:balanced-star-normality-spectra}

Theorem~\ref{thm:cyclic-star-Weyl} averages equally over the \(M\) rays and
therefore retains only the Fourier mode that is invariant under rotation.
The matrix \(\mathscr H_{(a)}\) contains more information: an entry joining factor
positions \(a\) and \(b\) can be nonzero only in powers congruent to
\(b-a\pmod M\).  This congruence selects a Fourier character on the star.
We first encode these moments as linear functionals on the polynomial space
\(\mathbb C[\lambda]\).  For a fixed compactly supported radial representing
measure, representation of the corresponding Fourier component by a finite
complex measure supported on the corresponding compact star is characterized
by the subsequent lifting criterion.

\begin{definition}[Moment functionals of the matrix \(\mathscr H_{(a)}\)]
	\label{def:balanced-character-star-functional}
	Order the coordinate states of \(\mathscr H_{(a)}\) by increasing
	height in the factor order used to construct it, using the convention of
	Definition~\ref{def:height-ordered-cyclic-operator}.  Let
	\(b_0,\ldots,b_{q-1}\) be its first \(q\) coordinate states and
	\(c_0,\ldots,c_{p-1}\) its first \(p\) coordinate states.
	Let \(\mathsf E_{[q]}^{(a)}\) and
	\(\mathsf E_{[p]}^{(a)}\) be the coordinate projections onto
	\((b_0,\ldots,b_{q-1})\) and \((c_0,\ldots,c_{p-1})\), respectively.
	For \(j\in\{0,\ldots,q-1\}\) and
	\(i\in\{0,\ldots,p-1\}\), define the linear functional
	\[
		\widehat\Psi_{j,i}:\mathbb C[\lambda]\longrightarrow\mathbb C
	\]
	by prescribing its moments
	\begin{equation}
		\bigl\langle\widehat\Psi_{j,i},\lambda^n\bigr\rangle
		\coloneq
		\left(
		\mathsf E_{[q]}^{(a)}
		(\mathscr H_{(a)})^n
		(\mathsf E_{[p]}^{(a)})^{\top}
		\right)_{j,i},
		\qquad n\in\mathbb N_0.
	\end{equation}
	Thus, if \(f(\lambda)=\sum_{n=0}^{N}f_n\lambda^n\), then
	\[
		\bigl\langle\widehat\Psi_{j,i},f\bigr\rangle
		=
		\sum_{n=0}^{N}f_n
		\bigl\langle\widehat\Psi_{j,i},\lambda^n\bigr\rangle.
	\]
	The \(q\times p\) matrix formed by these functionals is denoted by
	\(\widehat{\boldsymbol\Psi}\).
\end{definition}

If \(a(b_j)\) and \(a(c_i)\) are the factor-position labels of the two
states, put
\[
	s_{j,i}\equiv a(c_i)-a(b_j)\pmod M,
	\qquad 0\le s_{j,i}<M.
\]
Every path from \(b_j\) to \(c_i\) has length congruent to
\(s_{j,i}\) modulo \(M\).  Consequently,
\[
	\bigl\langle\widehat\Psi_{j,i},\lambda^n\bigr\rangle=0
	\qquad\text{if}\qquad n\not\equiv s_{j,i}\pmod M.
\]
Thus each entry retains one Fourier character of the rotational symmetry.
For later use, retain its nonzero moments in the radial functional
\(\mathcal R_{j,i}\) defined by
\begin{equation}
	\mathcal R_{j,i}(x^k)
	\coloneq
	\left(
	\mathsf E_{[q]}^{(a)}
	(\mathscr H_{(a)})^{s_{j,i}+Mk}
	(\mathsf E_{[p]}^{(a)})^{\top}
	\right)_{j,i},
	\qquad k\in\mathbb N_0.
\end{equation}

Before lifting these moments to the star, we determine their positivity on
the radial variable \(x=\lambda^M\).  For factor-copy labels
\(a,b\in\{1,\ldots,M\}\), let
\[
	s(a,b)\in\{0,\ldots,M-1\},
	\qquad
	b\equiv a+s(a,b)\pmod M,
\]
and denote by
\(\widetilde F_{a\to b}\) the product of the \(s(a,b)\) consecutive
refactorized factors beginning with \(\widetilde F_a\); when \(a=b\), this is
the identity.  The intertwining identity
\[
	\widetilde T_{(a)}\widetilde F_{a\to b}
	=
	\widetilde F_{a\to b}\widetilde T_{(b)}
\]
shows that the two expressions
\begin{equation}
	m_k^{a,b}
	\coloneq
	\left(
	\widetilde F_{a\to b}
	(\widetilde T_{(b)})^k
	\right)_{0,0}
	=
	\left(
	(\widetilde T_{(a)})^k
	\widetilde F_{a\to b}
	\right)_{0,0}
\end{equation}
coincide.

\begin{theorem}[Positive measures for all radial moment sequences]
	\label{thm:balanced-radial-Hankel-positivity}
	For every pair \(a,b\in\{1,\ldots,M\}\), every finite minor of the Hankel
	matrix
	\[
		\left[m_{r+c}^{a,b}\right]_{r,c\ge0}
	\]
	is a polynomial with nonnegative integer coefficients in the nontrivial
	entries of the refactorized bidiagonal factors.  For strictly positive
	numerical factors, the leading Hankel determinants and their one-step
	shifts satisfy
	\begin{equation}
		\label{eq:balanced-strict-radial-Hankel-determinants}
		\det\left[m_{r+c}^{a,b}\right]_{r,c=0}^{k-1}>0,
		\qquad
		\det\left[m_{r+c+1}^{a,b}\right]_{r,c=0}^{k-1}>0
	\end{equation}
	for every \(k\ge1\).  Consequently, there is a positive measure
	\(\mathrm d\rho^{a,b}\) on \([0,\infty)\), with infinite support, such
	that
	\[
		m_k^{a,b}
		=
		\int_{[0,\infty)}x^k\,\mathrm d\rho^{a,b}(x),
		\qquad k\in\mathbb N_0.
	\]
	If the refactorized factor coefficients are uniformly bounded, every
	representing measure is compactly supported and the measure is unique.
\end{theorem}

\begin{proof}
	Construct a periodic directed network with vertices
	\[
		(t,n)\in\mathbb Z\times\mathbb N_0.
	\]
	The layer between times \(t\) and \(t+1\) carries the refactorized factor
	whose copy label is congruent to \(a+t\) modulo \(M\).  If that factor is
	\(\widetilde F_c\), draw an edge from \((t,n)\) to \((t+1,m)\), with
	weight \((\widetilde F_c)_{n,m}\), whenever this matrix entry is nonzero.
	A lower-factor layer has the changes \(m-n=0,-1\), and an upper-factor
	layer has the changes \(m-n=0,+1\).  The network is therefore acyclic
	and planar in the upper half-plane.  Matrix multiplication identifies an
	entry of a consecutive factor product with the total weight of the
	corresponding paths.  Although the periodic network is infinite, fixed
	sources and targets involve only the finite subnetwork between their time
	coordinates and below the largest height reachable in that interval.
	Thus the finite-network form of the
	Lindstr\"om--Gessel--Viennot lemma applies to every minor considered below.
	The nested family used to prove strict positivity is shown schematically in
	Figure~\ref{fig:LGV-nested-paths} for three paths.

	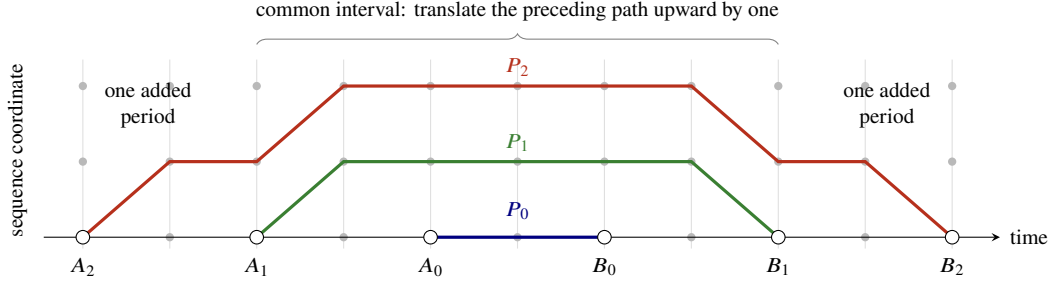
\begin{figure}[t]
		\centering
		\begin{tikzpicture}[
			x=1.15cm,
			y=1.0cm,
			>=stealth,
			vertex/.style={circle,fill=Gray!55,inner sep=1.1pt},
			endpoint/.style={circle,draw=black,fill=white,inner sep=1.8pt},
			pathzero/.style={very thick,NavyBlue},
			pathone/.style={very thick,OliveGreen},
			pathtwo/.style={very thick,BrickRed},
			lab/.style={font=\scriptsize}
		]
			\foreach \x in {-5,-4,...,5}{
				\draw[Gray!25] (\x,-0.15) -- (\x,2.35);
				\foreach \y in {0,1,2}{\node[vertex] at (\x,\y) {};}
			}
			\draw[->] (-5.45,0) -- (5.55,0) node[right,lab] {time};
			\node[lab,rotate=90] at (-5.75,1.15) {sequence coordinate};

			\draw[pathzero] (-1,0) -- (1,0);
			\draw[pathone] (-3,0) -- (-2,1) -- (2,1) -- (3,0);
			\draw[pathtwo]
				(-5,0) -- (-4,1) -- (-3,1) -- (-2,2) --
				(2,2) -- (3,1) -- (4,1) -- (5,0);

			\foreach \x/\name in {-5/A_2,-3/A_1,-1/A_0,1/B_0,3/B_1,5/B_2}{
				\node[endpoint] at (\x,0) {};
				\node[lab,below=4pt] at (\x,0) {\(\name\)};
			}

			\node[lab,text=NavyBlue] at (0,0.35) {\(P_0\)};
			\node[lab,text=OliveGreen] at (0,1.25) {\(P_1\)};
			\node[lab,text=BrickRed] at (0,2.25) {\(P_2\)};
			\draw[decorate,decoration={brace,amplitude=4pt},Gray]
				(-3,2.55) -- (3,2.55)
				node[midway,above=5pt,lab,text=black]
				{common interval: translate the preceding path upward by one};
			\node[lab,align=center] at (-4.25,1.75)
				{one added\\period};
			\node[lab,align=center] at (4.25,1.75)
				{one added\\period};
		\end{tikzpicture}
		\caption{The nested vertex-disjoint paths used in the LGV argument,
			shown for \(k=3\).  On the interval where \(P_r\) and
			\(P_{r-1}\) coexist, \(P_r\) is the vertical translate of
			\(P_{r-1}\) by one coordinate.  One period is added at each end,
			where the path rises from and returns to the lower boundary.  Thus
			\(P_r\) is disjoint from every earlier path and supplies a strictly
			positive monomial in the corresponding Hankel determinant.}
		\label{fig:LGV-nested-paths}
	\end{figure}

	Fix strictly increasing index sets
	\[
		I=\{i_1<\cdots<i_k\},
		\qquad
		J=\{j_1<\cdots<j_k\},
	\]
	and put \(s=s(a,b)\).  Choose sources and targets
	\[
		A_r=(-Mi_r,0),
		\qquad
		B_c=(s+Mj_c,0).
	\]
	Periodicity of the factors shows that the total weight of paths from
	\(A_r\) to \(B_c\) is \(m_{i_r+j_c}^{a,b}\).  Along the lower boundary
	the vertices occur in the order
	\[
		A_k,\ldots,A_1,B_1,\ldots,B_k.
	\]
	Thus a vertex-disjoint path family can only connect \(A_r\) to \(B_r\).
	The Lindstr\"om--Gessel--Viennot lemma
	\cite{Lindstrom1973,GesselViennot1985} gives
	\[
		\det\left[m_{i_r+j_c}^{a,b}\right]_{r,c=1}^{k}
		=
		\sum_{\mathcal P}\operatorname{wt}(\mathcal P),
	\]
	where the sum runs over these vertex-disjoint families.  Every summand is
	a monomial with a nonnegative integer coefficient.  Since \(I\) and
	\(J\) were arbitrary, this proves the first assertion.  The determinant
	may vanish when no such family exists, so strict positivity is established
	separately for the two families of determinants needed below.

	We next exhibit a positive family for each determinant in
	\eqref{eq:balanced-strict-radial-Hankel-determinants}.  For the first
	determinant, take the sources \((-Mr,0)\) and targets
	\((s+Mr,0)\), with \(r\in\{0,\ldots,k-1\}\).  Denote the path with these
	endpoints by \(P_r\).  The path \(P_0\) stays at coordinate \(0\) and uses
	only diagonal entries.

	Suppose that \(P_0,\ldots,P_{r-1}\) have been constructed.  On the time
	interval common to \(P_r\) and \(P_{r-1}\), define \(P_r\) to be the
	vertical translate of \(P_{r-1}\) by one coordinate.  Prefix this
	translated path by one additional period: use diagonal edges at coordinate
	\(0\) until an upper-factor layer, take its superdiagonal edge from
	coordinate \(0\) to coordinate \(1\), and then use diagonal edges at
	coordinate \(1\) until the common interval begins.  Suffix it by one
	additional period in the analogous way, using diagonal edges at coordinate
	\(1\) until a lower-factor layer, its subdiagonal edge from coordinate
	\(1\) to coordinate \(0\), and then diagonal edges at coordinate \(0\).
	Every period contains both an upper and a lower factor, so these choices
	are always available and have strictly positive weight.

	On the common interval, \(P_r\) lies one coordinate above \(P_{r-1}\),
	whereas all previously constructed paths lie on or below \(P_{r-1}\).
	Outside that interval none of the earlier paths is present.  Thus \(P_r\)
	is disjoint from every one of \(P_0,\ldots,P_{r-1}\).  Induction gives the
	required \(k\) vertex-disjoint paths
	of strictly positive weight.

	For the shifted determinant, use the sources \((-M(r+1),0)\) and the
	targets \((s+Mr,0)\).  The innermost path stays at coordinate \(0\)
	from time \(-M\) to time \(s\).  Each subsequent path is obtained from
	the preceding one by exactly the same three operations: add a left period,
	translate the whole common part upward by one coordinate, and add a right
	period.  The preceding argument again proves disjointness from all earlier
	paths.  Both determinants in
	\eqref{eq:balanced-strict-radial-Hankel-determinants} are therefore
	strictly positive.

	The Hankel matrix and its one-step shift are positive definite.  The
	Stieltjes moment criterion \cite{Akhiezer1965} gives a positive representing measure on
	\([0,\infty)\).  Strict positivity at every order forces infinite
	support.  If the factors are uniformly bounded, the partial product
	\(\widetilde F_{a\to b}\) and the cyclic product
	\(\widetilde T_{(b)}\) are bounded.  Hence there are constants
	\(C,R>0\) such that \(m_k^{a,b}\le CR^k\) for every \(k\).  Fix
	\(R_1>R\).  If a representing measure assigned positive mass to
	\([R_1,\infty)\), then
	\[
		m_k^{a,b}
		\ge
		R_1^k\rho^{a,b}([R_1,\infty)),
	\]
	which contradicts \(m_k^{a,b}\le CR^k\) for large \(k\).  Thus every
	representing measure is supported in \([0,R]\).  Polynomial density on
	this compact interval then gives uniqueness.
\end{proof}

\begin{corollary}[Stieltjes moment sequences throughout the prescribed cyclic orbit]
	\label{cor:prescribed-cyclic-radial-Hankel-positivity}
	Let \(F_1,\ldots,F_M\) be the prescribed positive bidiagonal factors and,
	for \(a,b\in\{1,\ldots,M\}\), define
	\[
		\mu_k^{a,b}
		\coloneq
		\left(F_{a\to b}(T_{(b)})^k\right)_{0,0},
		\qquad k\in\mathbb N_0.
	\]
	Every finite minor of the Hankel matrix
	\(\left[\mu_{r+c}^{a,b}\right]_{r,c\ge0}\) is a polynomial with
	nonnegative integer coefficients in the nontrivial entries of the
	prescribed factors.  Its leading Hankel determinants and those of its
	one-step shift are strictly positive.  Consequently, every sequence
	\((\mu_k^{a,b})_{k\ge0}\) is represented by a positive Stieltjes measure
	with infinite support.  In particular, this holds for the moments
	\[
		(T^k)_{0,0}
		\quad\text{and}\quad
		\bigl((T_{(a)})^k\bigr)_{0,0}
	\]
	of the original product and all its cyclic products.
\end{corollary}

\begin{proof}
	The proof of
	Theorem~\ref{thm:balanced-radial-Hankel-positivity} uses periodicity,
	positivity, and the lower or upper bidiagonal form of each factor.  The same
	periodic planar network argument applies to the prescribed factors
	\(F_1,\ldots,F_M\).  Taking
	\(a=b=1\) gives \(F_{1\to1}=I\) and \(T_{(1)}=T\); taking \(a=b\) gives
	the stated moments of every cyclic product.
\end{proof}

\begin{remark}
	The same argument applies to any periodic order of positive lower and
	upper bidiagonal factors.  For \(q=1\), it reduces to the coefficientwise
	Hankel total positivity underlying the branched Stieltjes fractions
	\cite{PetreolleSokalZhu2023}.  It extends to arbitrary \(p\) and \(q\)
	because the planar network retains the elementary factor steps instead of
	contracting them into the entries of the banded product.
\end{remark}

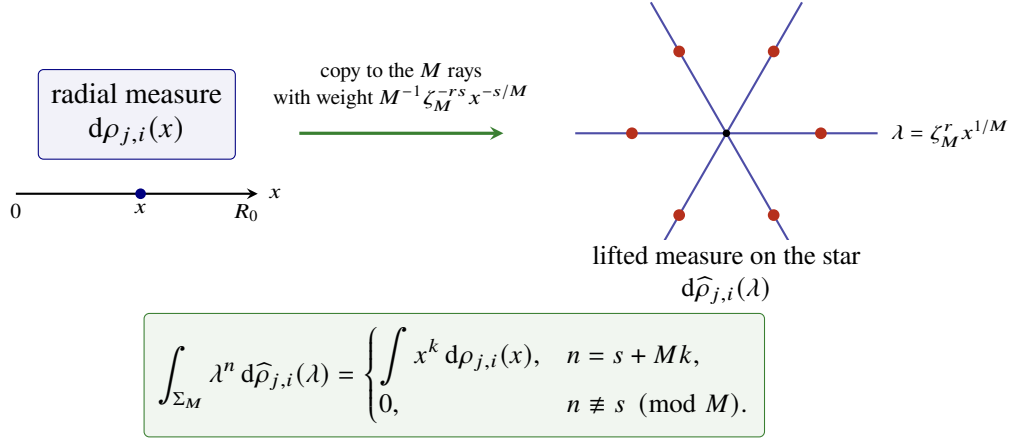
\begin{figure}[H]
	\centering
	\begin{tikzpicture}[
		x=1cm,
		y=1cm,
		>=stealth,
		lab/.style={font=\small},
		smalllab/.style={font=\scriptsize},
		box/.style={rounded corners=2pt,draw=NavyBlue,fill=NavyBlue!5,
			inner sep=5pt,align=center}
	]
		\node[box] (radial) at (-5.4,0.4)
			{radial measure\\\(\mathrm d\rho_{j,i}(x)\)};
		\draw[->,thick] (-7.0,-0.65) -- (-3.8,-0.65)
			node[right,smalllab] {\(x\)};
		\node[smalllab,below] at (-7.0,-0.65) {\(0\)};
		\node[smalllab,below] at (-3.95,-0.65) {\(R_0\)};
		\fill[NavyBlue] (-5.35,-0.65) circle (2pt);
		\node[smalllab,below] at (-5.35,-0.65) {\(x\)};

		\draw[->,very thick,OliveGreen]
			(-3.25,0.15) -- (-0.55,0.15)
			node[midway,above=4pt,smalllab,align=center,text=black]
			{copy to the \(M\) rays\\with weight
			\(M^{-1}\zeta_M^{-rs}x^{-s/M}\)};

		\coordinate (O) at (2.4,0.15);
		\foreach \r in {0,...,5}{
			\draw[thick,NavyBlue!70] (O) -- ++({60*\r}:2.0);
			\fill[BrickRed] ($(O)+({60*\r}:1.25)$) circle (2.2pt);
		}
		\fill (O) circle (1.4pt);
		\node[smalllab,anchor=west] at (4.45,0.15)
			{\(\lambda=\zeta_M^r x^{1/M}\)};
		\node[lab,align=center,fill=white,inner sep=2pt] at (2.4,-1.72)
			{lifted measure on the star\\
			\(\mathrm d\widehat\rho_{j,i}(\lambda)\)};

		\node[draw=OliveGreen,fill=OliveGreen!5,rounded corners=2pt,
			inner sep=5pt,align=center,lab] at (-1.2,-3.05)
			{\(\displaystyle
			\int_{\Sigma_M}\lambda^n\,\mathrm d\widehat\rho_{j,i}(\lambda)
			=
			\begin{cases}
			\displaystyle\int x^k\,\mathrm d\rho_{j,i}(x),&n=s+Mk,\\[2pt]
			0,&n\not\equiv s\pmod M.
			\end{cases}\)};
	\end{tikzpicture}
	\caption{Fourier lift of a radial measure to the \(M\)-ray star.  A
		point \(x\) is sent to its \(M\) roots
		\(\zeta_M^r x^{1/M}\), weighted by the \(s\)-th character of the
		cyclic rotation.  The root-of-unity average retains exactly the moments
		whose degrees are congruent to \(s\) modulo \(M\).  Six rays are drawn
		for illustration; the construction uses \(M=p+q\) rays.}
	\label{fig:Fourier-lift-to-star}
\end{figure}

For the boundary states in
Definition~\ref{def:balanced-character-star-functional}, one has
\[
	\mathcal R_{j,i}(x^k)
	=
	m_k^{a(b_j),a(c_i)}.
\]
Indeed, the first \(q\) and \(p\) states have heights below \(Q\) and
\(P\), respectively, and therefore all have sequence coordinate zero.
Theorem~\ref{thm:balanced-radial-Hankel-positivity} therefore supplies a
positive radial measure for every entry of the moment matrix of
\(\mathscr H_{(a)}\).  If the refactorized coefficients are uniformly
bounded, that measure is unique and compactly supported.  The next result
characterizes the additional integrability condition after a compactly
supported radial representative has been fixed.
Figure~\ref{fig:Fourier-lift-to-star} summarizes the lift and its selection of
one congruence class of moments.

\begin{proposition}[A measure on one ray variable and its Fourier lift to the star]
	\label{prop:balanced-character-measure-lift}
	Fix \((j,i)\) and put \(s=s_{j,i}\).
	Suppose that \(\mathcal R_{j,i}\) is represented by a finite complex
	measure \(\mathrm d\rho_{j,i}\) supported on \([0,R_0]\).  The linear
	functional \(\widehat\Psi_{j,i}\) is represented by a finite
	complex measure on \(\Sigma_M^{[R_0]}\) if and only if, when \(s>0\),
	\(\rho_{j,i}(\{0\})=0\) and
	\begin{equation}
		\label{eq:balanced-character-negative-moment-condition}
		\int_{(0,R_0]}x^{-s/M}\,
		\mathrm d\lvert\rho_{j,i}\rvert(x)<\infty.
	\end{equation}
	When these conditions hold, set
	\(\Sigma_M^{[R_0]}=\bigcup_{r=0}^{M-1}
	\zeta_M^r[0,R_0^{1/M}]\) and define one such representing measure by
	\begin{equation}
		\label{eq:balanced-character-measure-lift}
		\int_{\Sigma_M^{[R_0]}}f(\lambda)\,
		\mathrm d\widehat\rho_{j,i}(\lambda)
		\coloneq
		\frac1M\sum_{r=0}^{M-1}\zeta_M^{-rs}
		\int_{(0,R_0]}
		f\bigl(\zeta_M^r x^{1/M}\bigr)x^{-s/M}\,
		\mathrm d\rho_{j,i}(x).
	\end{equation}
	When \(s=0\), an atom of \(\rho_{j,i}\) at the origin is included in
	the right-hand side in the evident way.  Then
	\(\mathrm d\widehat\rho_{j,i}\) is a finite complex measure and
	\begin{equation}
		\label{eq:balanced-character-measure-moments}
		\int_{\Sigma_M^{[R_0]}}\lambda^n\,
		\mathrm d\widehat\rho_{j,i}(\lambda)
		=
		\begin{cases}
		\mathcal R_{j,i}(x^k),&n=s+Mk,\\
		0,&n\not\equiv s\pmod M.
		\end{cases}
	\end{equation}
	Its Cauchy transform is
	\begin{equation}
		\label{eq:balanced-character-measure-Cauchy-transform}
		\int_{\Sigma_M^{[R_0]}}
		\frac{\mathrm d\widehat\rho_{j,i}(\zeta)}{\lambda-\zeta}
		=
		\lambda^{M-s-1}
		\int_{[0,R_0]}\frac{\mathrm d\rho_{j,i}(x)}{\lambda^M-x},
		\qquad
		\lambda\in\mathbb C\setminus\Sigma_M^{[R_0]}.
	\end{equation}
\end{proposition}

\begin{proof}
	We first prove sufficiency.  Condition
	\eqref{eq:balanced-character-negative-moment-condition} shows that every
	term in \eqref{eq:balanced-character-measure-lift} has finite total
	variation.  Thus the formula defines a finite complex measure.  When
	\(s>0\), the hypothesis \(\rho_{j,i}(\{0\})=0\) ensures that no mass is
	placed where the factor \(x^{-s/M}\) is singular.
	Moreover, for every \(n\ge0\),
	\[
		x^{(n-s)/M}=x^{n/M}x^{-s/M}
	\]
	is integrable with respect to \(\lvert\rho_{j,i}\rvert\), because
	\(x^{n/M}\) is bounded on \([0,R_0]\) and the negative fractional moment
	is finite.  When \(s=0\), the same integrability follows directly from the
	finiteness and compact support of \(\rho_{j,i}\).
	For \(n\ge0\), the root-of-unity average is
	\[
		\frac1M\sum_{r=0}^{M-1}\zeta_M^{r(n-s)}
		=
		\begin{cases}1,&n\equiv s\pmod M,\\0,&n\not\equiv s\pmod M.
		\end{cases}
	\]
	If \(n=s+Mk\), the remaining power of \(x\) in
	\eqref{eq:balanced-character-measure-lift} is exactly \(x^k\).  This
	proves \eqref{eq:balanced-character-measure-moments}.  Finally, for
	\(r=x^{1/M}\) and \(\lvert\lambda\rvert>R_0^{1/M}\), expansion in powers
	of \(r/\lambda\) gives
	\[
		\frac1M\sum_{h=0}^{M-1}
		\frac{\zeta_M^{-hs}}{\lambda-\zeta_M^h r}
		=
		\frac{r^s\lambda^{M-s-1}}{\lambda^M-r^M}.
	\]
	Multiplication by \(x^{-s/M}=r^{-s}\) and integration prove
	\eqref{eq:balanced-character-measure-Cauchy-transform}; analytic
	continuation gives the identity throughout
	\(\mathbb C\setminus\Sigma_M^{[R_0]}\).

	For necessity, let \(\mathrm d\eta\) be any finite complex measure on
	\(\Sigma_M^{[R_0]}\) representing \(\widehat\Psi_{j,i}\).  Replace it by its
	\(s\)-th Fourier component under rotations:
	\[
		\int f(\lambda)\,\mathrm d\eta_s(\lambda)
		\coloneq
		\frac1M\sum_{r=0}^{M-1}\zeta_M^{-rs}
		\int f(\zeta_M^r\lambda)\,\mathrm d\eta(\lambda).
	\]
	This is again a finite measure and has the same moments as \(\eta\), because
	the prescribed moments vanish outside degrees congruent to \(s\) modulo
	\(M\).  Push the measure \(\lambda^s\,\mathrm d\eta_s(\lambda)\) forward
	under \(\lambda\mapsto\lambda^M\), and call the resulting measure
	\(\mathrm d\rho'\).  For every \(k\ge0\),
	\[
		\int_{[0,R_0]}x^k\,\mathrm d\rho'(x)
		=
		\int_{\Sigma_M^{[R_0]}}\lambda^{s+Mk}\,
		\mathrm d\eta_s(\lambda)
		=
		\mathcal R_{j,i}(x^k).
	\]
	Finite complex measures on the compact interval \([0,R_0]\) are determined
	by their moments.  Hence \(\rho'=\rho_{j,i}\).  If \(s>0\), multiplication
	by \(\lambda^s\) removes any atom at the origin, so
	\(\rho_{j,i}(\{0\})=0\).  Moreover, the total-variation inequality for a
	pushforward gives
	\[
		\int_{(0,R_0]}x^{-s/M}\,
		\mathrm d\lvert\rho_{j,i}\rvert(x)
		\le
		\int_{\Sigma_M^{[R_0]}\setminus\{0\}}
		\lvert\lambda\rvert^{-s}\lvert\lambda\rvert^s\,
		\mathrm d\lvert\eta_s\rvert(\lambda)
		\le \lVert\eta_s\rVert<\infty.
	\]
	This proves necessity.  For \(s=0\), the same construction requires no
	negative moment and retains a possible atom at the origin.
\end{proof}

\begin{proposition}[Obstruction to positive scalar representations]
	\label{prop:balanced-nonzero-Fourier-mode-not-positive}
	Let \(s\in\{1,\ldots,M-1\}\).  A nonzero functional whose moments
	vanish outside degrees congruent to \(s\) modulo \(M\) cannot be
	represented by a finite positive scalar measure on the star.
\end{proposition}

\begin{proof}
	The moment of degree zero vanishes.  If a positive measure represented
	the functional, that moment would be its total mass.  The measure would
	therefore be zero, contrary to the assumption that the functional is
	nonzero.
\end{proof}

The support of every lifted entry lies on the same \(M\) rays as the
rotationally invariant lift.  Since
\(s_{j,i}\equiv a(c_i)-a(b_j)\pmod M\), the phases on each ray separate into a
left and a right diagonal unitary factor.  The additional factor
\(x^{-s_{j,i}/M}\) in
\eqref{eq:balanced-character-measure-lift} shows that support reaching the
origin adds a negative fractional-moment condition to the full moment
problem on the star.  Positivity belongs to the radial measures
on \([0,\infty)\); the nonzero Fourier components on the star are
necessarily complex.

For nonsingular triangular matrices
\(\xi\in\mathbb C^{q\times q}\) and
\(\nu\in\mathbb C^{p\times p}\), normalize
\[
	\widehat{\boldsymbol\Psi}^{\xi,\nu}
	\coloneq
	\xi^{-1}\widehat{\boldsymbol\Psi}\nu^{-\top}.
\]
Form its scalar step-line moment matrix by
\[
	\mathscr M^{\xi,\nu}_{qu+j,pv+i}
	\coloneq
	\left(
		\bigl\langle
		\widehat{\boldsymbol\Psi}^{\xi,\nu},\lambda^{u+v}
		\bigr\rangle
	\right)_{j,i}.
\]
Write \(\mathscr M\coloneq\mathscr M^{I_q,I_p}\) for the corresponding
unnormalized scalar moment matrix.

Recall the normality terminology fixed in
Definition~\ref{def:weak-strong-step-line-normality}.  For the scalar moment
matrix above, the row positions are ordered as \(qu+j\), with
\(j\in\{0,\ldots,q-1\}\), and the column positions as \(pv+i\), with
\(i\in\{0,\ldots,p-1\}\).  Hence its leading \(N\times N\) block determines
the row and column multi-indices at the \(N\)-th position of the prescribed
step-line.

For the scalar step-line moment matrix used here, nonvanishing of the leading
minors of orders \(1,\ldots,N\) is equivalent to weak normality through order
\(N\).  The maximal-degree condition in strong normality requires additional
information.  When \(d>1\), the argument below first proves weak normality at
orders \(d,2d,3d,\ldots\), because those are precisely the cuts that retain
all \(d\) coordinate states having the last included height.

Terminology in the literature is not uniform.  The word \emph{normality} in
\cite{BranquinhoFoulquieManas2024BTP} denotes the maximal-degree property and
therefore corresponds to strong normality in
Definition~\ref{def:weak-strong-step-line-normality}.  In the
hypergeometric setting of \cite{Manas2026Hypergeometric}, normality denotes
existence and uniqueness, while the AT property established there also gives
the maximal degrees and hence strong normality.

The factorization used in both normality statements can be written before
any normalization is chosen.  Define two scalar path arrays by
\[
	\mathcal A_{qu+j,k}
	\coloneq
	\bigl((\mathscr H_{(a)})^u\bigr)_{j,k},
	\qquad
	\mathcal B_{pv+i,k}
	\coloneq
	\bigl((\mathscr H_{(a)})^v\bigr)_{k,i},
\]
where \(j\in\{0,\ldots,q-1\}\) and \(i\in\{0,\ldots,p-1\}\) refer to
the first \(q\) and \(p\) scalar states in exact-height order.  Then the
unnormalized scalar moment matrix satisfies
\begin{equation}
	\label{eq:balanced-path-array-moment-factorization}
	\mathscr M=\mathcal A\mathcal B^{\top}.
\end{equation}
Indeed, its \((qu+j,pv+i)\)-entry is
\[
	\sum_{k\ge0}
	\bigl((\mathscr H_{(a)})^u\bigr)_{j,k}
	\bigl((\mathscr H_{(a)})^v\bigr)_{k,i}
	=
	\bigl((\mathscr H_{(a)})^{u+v}\bigr)_{j,i}.
\]
For fixed \(u\) and \(v\), the sum is finite because
\(\mathscr H_{(a)}\) is banded.

\begin{proposition}[Weak normality after retaining complete height groups]
	\label{prop:balanced-complete-fiber-normality}
	For every \(K\ge1\), the leading principal minor of order \(dK\) of the
	unnormalized scalar moment matrix is nonzero.  Equivalently, weak normality
	holds after every truncation that retains all the coordinate states of
	heights \(0,\ldots,K-1\).
\end{proposition}

\begin{proof}
	Use the factorization
	\eqref{eq:balanced-path-array-moment-factorization}.  The first
	\(q=dQ\) row-boundary states consist of the \(d\) states at each of the
	heights
	\[
	0,\ldots,Q-1.
	\]
	After \(u\) successive upward steps, their endpoints are precisely the
	\(d\) states at each of the heights
	\[
	uQ,\ldots,(u+1)Q-1.
	\]
	No other path of \(u\) steps can end at a greater height.  Therefore,
	when the indices are grouped according to their common height,
	\(\mathcal A\) is block lower triangular, with blocks of order \(d\).
	Each diagonal block is a permutation matrix: between a boundary state
	and its maximal endpoint there is a unique path consisting entirely of
	upward steps, and every upward edge has weight \(1\).
	
	Similarly, the first \(p=dP\) column-boundary states consist of the
	\(d\) states at each of the heights
	\[
	0,\ldots,P-1.
	\]
	A path of \(v\) steps ending at one of these states can start no higher
	than the initial point of the unique path consisting entirely of
	downward steps.  These initial points are precisely the states at the
	heights
	\[
	vP,\ldots,(v+1)P-1.
	\]
	Hence \(\mathcal B\) is also block lower triangular.  Its diagonal
	blocks are positive monomial matrices, because they are supplied by the
	unique all-down paths.  In particular, every diagonal block of
	\(\mathcal A\) and \(\mathcal B\) is invertible.
	
	Fix \(K\ge1\).  The first \(dK\) indices comprise all
	the states of heights \(0,\ldots,K-1\).  Block lower triangularity shows
	that the first \(dK\) rows of both path arrays vanish in every column
	with index at least \(dK\).  Therefore the leading \(dK\times dK\) block of
	\[
	\mathscr M=\mathcal A\mathcal B^{\top}
	\]
	is the product of the corresponding leading blocks:
	\[
	\mathscr M^{[dK-1]}
	=
	\mathcal A^{[dK-1]}
	\bigl(\mathcal B^{[dK-1]}\bigr)^{\top}.
	\]
	Both factors are block triangular with invertible diagonal blocks.
	Consequently,
	\[
	\det\mathscr M^{[dK-1]}
	=
	\det\mathcal A^{[dK-1]}\det\mathcal B^{[dK-1]}
	\ne0.
	\]
\end{proof}

\begin{theorem}[Weak normality at every scalar order in the coprime case]
	\label{thm:balanced-scalar-normality}
	Assume that \(d=\gcd(p,q)=1\), and let \(\xi\) and \(\nu\) be lower
	unitriangular.  Every leading principal minor of
	\(\mathscr M^{\xi,\nu}\) is nonzero.  Hence the prescribed step-line is
	weakly normal at every scalar order.
\end{theorem}

\begin{proof}
	Start from the path-array factorization
	\eqref{eq:balanced-path-array-moment-factorization}.

	Since \(d=1\), every height contains a single coordinate state.  In
	increasing exact-height order, the row of \(\mathcal A\) numbered
	\(qu+j\) records paths of \(u\) steps beginning at height \(j\).
	Such a path can end no higher than \(j+uq=qu+j\), and equality is attained
	only by the all-up path, whose weight is one.  Hence \(\mathcal A\) is
	scalar lower unitriangular.
	Applying \(\xi^{-1}\) to each block of \(q\) consecutive rows preserves
	this scalar lower-triangular structure and its nonzero diagonal.  Denote
	the resulting array by \(\mathcal A^{\xi}\).

	The array \(\mathcal B\) is also scalar lower triangular.  Its row numbered
	\(pv+i\) records, in transposed path orientation, paths of \(v\) steps
	ending at the boundary state of height \(i\).  Since one step can decrease
	height by at most \(p\), a path beginning at height \(k\) can reach that
	state only if
	\[
		k\le i+vp=pv+i.
	\]
	Thus the entries in columns numbered above the row index vanish.  Equality
	requires every step to be downward, so the diagonal entry is the strictly
	positive weight of the unique all-down path.

	Right normalization of the moment matrix by \(\nu^{-\top}\) replaces
	\(\mathcal B\) by an array \(\mathcal B^{\nu}\) obtained by applying
	\(\nu^{-1}\) to each block of \(p\) consecutive rows.  Because \(\nu^{-1}\)
	is lower triangular, \(\mathcal B^{\nu}\) remains scalar lower triangular.
	The positive diagonal blocks described above and the unit diagonal of
	\(\nu^{-1}\) show that every diagonal entry of \(\mathcal B^{\nu}\) is
	nonzero.

	The normalized moment factorization is
	\[
		\mathscr M^{\xi,\nu}
		=
		\mathcal A^{\xi}(\mathcal B^{\nu})^{\top}.
	\]
	Because \(\mathcal A^{\xi}\) is scalar lower triangular, its first \(N\)
	rows vanish in every column numbered at least \(N\).  Therefore the
	leading \(N\times N\) block of the displayed product is the product of
	the leading blocks of \(\mathcal A^{\xi}\) and
	\((\mathcal B^{\nu})^{\top}\).  Both determinants are nonzero.  By
	Definition~\ref{def:weak-strong-step-line-normality}, this proves weak
	normality at every scalar order.
\end{proof}

\begin{remark}[Weak and strong normality for the sparse recurrence]
	\label{rem:balanced-weak-not-strong-normality}
	The preceding theorem establishes weak normality.  Strong normality may
	fail even if \(\xi\) and \(\nu\) are chosen to be Pascal-type totally
	positive matrices.  For
	example, take \(p=2\) and \(q=3\).  With the usual lower triangular
	normalization, the first right polynomial vector is
	\(\boldsymbol B_0=(1,0,0)\).  The row of the two-diagonal recurrence at
	height \(0\) gives
	\[
		\boldsymbol B_3(x)=x\boldsymbol B_0(x)=(x,0,0).
	\]
	At this step-line position the second and third components have nontrivial
	constant polynomial spaces, but they vanish.  Thus this index pair is
	weakly normal and not strongly normal.

	The BTP criterion of \cite{BranquinhoFoulquieManas2024BTP} applies to the
	cyclic products that have a PBF.  The sparse matrix \(\mathscr H_{(a)}\)
	has a zero main diagonal and zero diagonals between its two extreme nonzero
	diagonals, so it lies outside the BTP class.
\end{remark}

Group the scalar coordinate states of
\(\mathscr H_{(a)}\) according to their common height.  Each group
has dimension \(d\), and the path arrays \(\mathcal A\) and \(\mathcal B\)
become semi-infinite lower triangular block matrices with blocks of order
\(d\).  Set
\[
	\mathscr L_{\star}\coloneq\mathcal A^{-1},
	\qquad
	\mathscr U_{\star}\coloneq\mathcal B^{-\top}.
\]
Let \(\mathsf\Lambda_{[Q]}^{(d)}\) be the block shift whose nonzero blocks
are \(d\times d\) identity matrices in block positions \((k,k+Q)\), and
define \(\mathsf\Lambda_{[P]}^{(d)}\) analogously.

\begin{theorem}[The matrix \(\mathscr H_{(a)}\) as a block mixed-type recurrence]
	\label{thm:balanced-star-recurrence-matrix}
	With the preceding notation,
	\[
		\mathscr M
		=
		\mathscr L_{\star}^{-1}\mathscr U_{\star}^{-1},
		\qquad
		\mathscr L_{\star}^{-1}=\mathcal A,
		\qquad
		\mathscr U_{\star}^{-\top}=\mathcal B
	\]
	is a triangular factorization of the moment matrix by invertible lower and
	upper block matrices.  Equivalently, after a diagonal-block normalization,
	it is a block Gauss--Borel factorization with the same moment matrix.
	Moreover,
	\begin{equation}
		\label{eq:balanced-star-recurrence-matrix}
		\mathscr H_{(a)}
		=
		\mathscr L_{\star}\mathsf\Lambda_{[Q]}^{(d)}
		\mathscr L_{\star}^{-1}
		=
		\mathscr U_{\star}^{-1}
		\bigl(\mathsf\Lambda_{[P]}^{(d)}\bigr)^{\top}
		\mathscr U_{\star}.
	\end{equation}
\end{theorem}

\begin{proof}
	The proof of
	Proposition~\ref{prop:balanced-complete-fiber-normality} shows that
	\(\mathcal A\) and \(\mathcal B\), after grouping equal-height states,
	are lower triangular block matrices with invertible \(d\times d\)
	diagonal blocks.  They are therefore invertible as formal triangular
	block matrices.  The path-array identity gives
	\[
		\mathscr M
		=
		\mathcal A\mathcal B^{\top}
		=
		\mathscr L_{\star}^{-1}\mathscr U_{\star}^{-1}.
	\]

	The block row of \(\mathcal A\) numbered \(Qu+\alpha\), with
	\(\alpha\in\{0,\ldots,Q-1\}\), is the block row of
	\((\mathscr H_{(a)})^u\) beginning at height \(\alpha\).
	Deleting the first \(Q\) block rows therefore advances the power by one:
	\[
		\mathsf\Lambda_{[Q]}^{(d)}\mathcal A
		=
		\mathcal A\mathscr H_{(a)}.
	\]
	Substitution of \(\mathcal A=\mathscr L_{\star}^{-1}\) gives the
	first representation in
	\eqref{eq:balanced-star-recurrence-matrix}.

	Similarly, the block row of \(\mathcal B\) numbered \(Pv+\beta\), with
	\(\beta\in\{0,\ldots,P-1\}\), is the corresponding block row of
	\(((\mathscr H_{(a)})^{\top})^v\).  Hence
	\[
		\mathsf\Lambda_{[P]}^{(d)}\mathcal B
		=
		\mathcal B(\mathscr H_{(a)})^{\top}.
	\]
	Transposition and substitution of
	\(\mathcal B=\mathscr U_{\star}^{-\top}\) give the second
	representation.  These are the right and left recurrence identities of
	the polynomial systems associated with the block Gauss--Borel factors.
	Flattening the blocks recovers the scalar moment functionals already
	defined.  A representing matrix of measures, when it exists, changes the
	realization of those functionals but not these recurrence identities.
\end{proof}

The preceding factorization also makes precise the relation between the
rotational lift of the Favard measure of \(T_{(a)}\) and the matrix of star
functionals \(\widehat{\boldsymbol\Psi}\).  Define the right and left
polynomial arrays of the star recurrence by
\[
	\boldsymbol B_{\star}(\lambda)
	\coloneq
	\mathscr L_{\star}\boldsymbol\chi_{[q]}(\lambda),
	\qquad
	\boldsymbol A_{\star}(\lambda)
	\coloneq
	\boldsymbol\chi_{[p]}(\lambda)^{\top}\mathscr U_{\star}.
\]
Let \(\mathsf R_a\) be the factor-copy projection used in
Corollary~\ref{cor:admissible-cyclic-star-CF-family}, and put
\[
	\mathsf Q_{r,\star}^{(a)}
	\coloneq
	\mathsf E_{[r]}\mathsf R_a,
	\qquad
	\boldsymbol B_{a}(\lambda)
	\coloneq
	\mathsf Q_{q,\star}^{(a)}\boldsymbol B_{\star}(\lambda),
	\qquad
	\boldsymbol A_{a}(\lambda)
	\coloneq
	\boldsymbol A_{\star}(\lambda)
	\bigl(\mathsf Q_{p,\star}^{(a)}\bigr)^{\top}.
\]
Thus \(\boldsymbol B_a\) and \(\boldsymbol A_a\) are the polynomial
matrices attached to the first \(q\) and \(p\) coordinates, respectively,
of the factor copy whose complete circuit gives \(T_{(a)}\).

\begin{proposition}[The rotational lift inside the star moment system]
	\label{prop:rotational-lift-inside-star-system}
	Suppose that the matrix of functionals
	\(\widehat{\boldsymbol\Psi}\) is represented on the star by a
	\(q\times p\) matrix of finite measures
	\(\mathrm d\widehat{\boldsymbol\Psi}_{\star}\).  Suppose also that
	\(T_{(a)}\) has Favard data
	\(\mathrm d\Psi^{(a)},\xi^{(a)},\nu^{(a)}\), with
	\(\mathrm d\Psi^{(a)}\) supported on \([0,\infty)\), and let
	\(\mathrm d\widehat\Psi^{(a)}\) be the rotational lift of
	\(\mathrm d\Psi^{(a)}\) given by
	Definition~\ref{def:rotational-lift-matrix-measure}.  Then, for every
	\(n\in\mathbb N_0\),
	\begin{equation}
		\label{eq:rotational-lift-inside-star-system-moments}
		\int_{\Sigma_M}\lambda^n\,
		\mathrm d\widehat\Psi^{(a)}(\lambda)
		=
		(\xi^{(a)})^{-1}
		\int_{\Sigma_M}
		\boldsymbol B_a(\lambda)\lambda^n\,
		\mathrm d\widehat{\boldsymbol\Psi}_{\star}(\lambda)
		\boldsymbol A_a(\lambda)
		(\nu^{(a)})^{-\top}.
	\end{equation}
	Consequently, the rotational lift is the matrix spectral measure of the
	star system observed from the coordinates belonging to that factor copy.
	If the measures are compactly supported, the moment identity determines
	them and one has the measure identity
	\begin{equation}
		\label{eq:rotational-lift-inside-star-system-measures}
		\mathrm d\widehat\Psi^{(a)}(\lambda)
		=
		(\xi^{(a)})^{-1}
		\boldsymbol B_a(\lambda)\,
		\mathrm d\widehat{\boldsymbol\Psi}_{\star}(\lambda)
		\boldsymbol A_a(\lambda)
		(\nu^{(a)})^{-\top}.
	\end{equation}
\end{proposition}

\begin{proof}
	The moment factorization in
	Theorem~\ref{thm:balanced-star-recurrence-matrix} gives
	\[
		\int_{\Sigma_M}
		\boldsymbol B_{\star}(\lambda)\,
		\mathrm d\widehat{\boldsymbol\Psi}_{\star}(\lambda)
		\boldsymbol A_{\star}(\lambda)
		=I.
	\]
	The two recurrence identities in
	\eqref{eq:balanced-star-recurrence-matrix} imply
	\[
		\mathscr H_{(a)}\boldsymbol B_{\star}(\lambda)
		=\lambda\boldsymbol B_{\star}(\lambda),
		\qquad
		\boldsymbol A_{\star}(\lambda)\mathscr H_{(a)}
		=\lambda\boldsymbol A_{\star}(\lambda).
	\]
	It follows that
	\[
		\int_{\Sigma_M}
		\boldsymbol B_{\star}(\lambda)\lambda^n\,
		\mathrm d\widehat{\boldsymbol\Psi}_{\star}(\lambda)
		\boldsymbol A_{\star}(\lambda)
		=
		(\mathscr H_{(a)})^n.
	\]
	Selecting the indicated factor-copy coordinates therefore gives
	\begin{align*}
		&\int_{\Sigma_M}
		\boldsymbol B_a(\lambda)\lambda^n\,
		\mathrm d\widehat{\boldsymbol\Psi}_{\star}(\lambda)
		\boldsymbol A_a(\lambda)
		\\
		&\hspace{25mm}=
		\mathsf Q_{q,\star}^{(a)}
		(\mathscr H_{(a)})^n
		\bigl(\mathsf Q_{p,\star}^{(a)}\bigr)^{\top}.
	\end{align*}
	By \eqref{eq:admissible-cyclic-star-resolvent}, the last member is zero
	when \(M\nmid n\), while for \(n=Mk\) it equals
	\[
		\mathsf E_{[q]}(T_{(a)})^k\mathsf E_{[p]}^{\top}.
	\]
	After multiplication by the normalization matrices, these are precisely
	the moments of the rotational lift of \(\mathrm d\Psi^{(a)}\).  This
	proves \eqref{eq:rotational-lift-inside-star-system-moments}.  For compact
	support, polynomial moments determine finite measures on the finite union
	of rays \(\Sigma_M\), so the measure identity follows entry by entry.
\end{proof}

Thus \(\mathscr H_{(a)}\) is the recurrence matrix of the block
moment system constructed above, with \(Q\) right and \(P\) left boundary
components of size \(d\).  When \(d=1\), this is an ordinary weakly normal
mixed-type system along the full scalar step-line.  When \(d>1\), the result
gives block normality after complete groups of \(d\) consecutive states.
Proposition~\ref{prop:balanced-character-measure-lift} characterizes, under
its compact-support hypothesis, which of these functionals admit finite
representing measures supported on the corresponding compact star.

\begin{remark}[Linear functionals and measures on the star]
	\label{rem:balanced-star-functional-versus-measure}
	For \(\gcd(p,q)=1\),
	Theorem~\ref{thm:balanced-star-recurrence-matrix} gives a recurrence
	realization for the matrix of linear functionals
	\(\widehat{\boldsymbol\Psi}\).
	Theorem~\ref{thm:balanced-radial-Hankel-positivity} proves directly from
	the PBF that every radial moment sequence occurring in this matrix is a
	Stieltjes moment sequence, and therefore supplies a positive radial measure
	for every boundary pair.  When the chosen radial measure is compactly
	supported, Proposition~\ref{prop:balanced-character-measure-lift} gives a
	finite complex measure supported on the corresponding compact star for the
	Fourier component if and only if, in nonzero mode, the radial measure has no
	atom at the origin and satisfies the negative fractional-moment condition.
	No such equivalence
	for an arbitrary noncompact representative is asserted here; in an
	indeterminate radial moment problem the moments alone do not select that
	representative.  The lifted measures are complex in every nonzero Fourier mode, as
	Proposition~\ref{prop:balanced-nonzero-Fourier-mode-not-positive} shows.
	Verification of the compactly supported lifting condition is
	family-dependent.  The weak normality result in the coprime case, the block
	recurrence \eqref{eq:balanced-star-recurrence-matrix}, and the finite
	spectral results below are established at the level of the moment
	functionals.  When \(\gcd(p,q)>1\), the uniform conclusion is the block
	Gauss--Borel realization in
	Theorem~\ref{thm:balanced-star-recurrence-matrix}, with cuts taken after
	complete groups of equal-height states.  Strong normality requires
	additional information, as explained in
	Remark~\ref{rem:balanced-weak-not-strong-normality}.
\end{remark}

\begin{theorem}[Spectra of finite truncations of \(\mathscr H_{(a)}\)]
	\label{thm:balanced-finite-star-spectra}
	Let \(K_N\) be the finite principal matrix obtained from
	\(\mathscr H_{(a)}\) by retaining all coordinate states of heights
	\(0,\ldots,N\) in this order.  Thus \(K_N\) has order
	\(d(N+1)\).  There
	are an integer \(r_N\in\{0,\ldots,M-1\}\) and a polynomial
	\(\mathcal P_N\), with \(\mathcal P_N(0)\ne0\), such that
	\[
		\det(\lambda I-K_N)
		=
		\lambda^{r_N}\mathcal P_N(\lambda^M),
		\qquad
		r_N=d(N+1)-M\left\lfloor\frac{d(N+1)}M\right\rfloor.
	\]
	Moreover,
	\[
		\deg\mathcal P_N
		=
		\left\lfloor\frac{d(N+1)}M\right\rfloor
		=
		\left\lfloor\frac{N+1}{P+Q}\right\rfloor.
	\]
	Every zero of \(\mathcal P_N\) is positive.  Hence every nonzero pole of
	the finite resolvent of \(\mathscr H_{(a)}\) lies on one of the \(M=p+q\) rays through
	the \(M\)-th roots of unity.  More precisely, order the retained states
	by factor copy and let \(C_{a,N}\) be the diagonal block of
	\(K_N^M\) on the \(a\)-th copy.  If
	that copy contains \(n_{a,N}\) states and
	\(m_N=\lfloor d(N+1)/M\rfloor\), then
	\[
		\det(yI-C_{a,N})
		=
		y^{n_{a,N}-m_N}\mathcal P_N(y).
	\]
	Thus \(\mathcal P_N\) is the characteristic polynomial of every diagonal
	block of order \(m_N\) in the \(M\)-th power, and its zeros are the
	nonzero radial spectral parameters.
\end{theorem}

\begin{proof}
	Put
	\[
		D_N\coloneq d(N+1)=\dim K_N.
	\]
	The proof has three steps.

	\emph{Rotational form.}
	Every nonzero entry advances the factor-copy label by one modulo \(M\).
	If \(D\) multiplies the coordinates in copy \(a\) by \(\zeta_M^a\),
	then
	\[
		D^{-1}K_ND
		=
		\zeta_MK_N.
	\]
	Consequently, the nonzero terms in the characteristic polynomial have
	degrees congruent to \(D_N\) modulo \(M\).  Let \(r_N\) be the remainder
	of \(D_N\) on division by \(M\).  At this stage the symmetry gives
	\[
		\det(\lambda I-K_N)=\lambda^{r_N}\mathcal R_N(\lambda^M)
	\]
	for some polynomial \(\mathcal R_N\); it does not yet exclude an
	additional factor of \(\lambda^M\) in \(\mathcal R_N(\lambda^M)\).

	\emph{The \(M\)-th power.}
	Order the coordinates by factor copy.  The finite truncation of \(\mathscr H_{(a)}\) is then
	a cyclic block matrix whose block from one copy to the next is a leading
	rectangular submatrix of one of the positive bidiagonal factors.  Every
	such rectangular block is totally nonnegative.  Hence the \(M\)-th power
	is block diagonal, and each square diagonal block is a product of
	compatible rectangular totally nonnegative matrices.  Repeated
	Cauchy--Binet shows that these cyclic products are totally nonnegative, so
	their eigenvalues are real and nonnegative
	\cite{FallatJohnson2011,Pinkus2010}.

	Put \(m_N=\lfloor D_N/M\rfloor\) and write
	\(N=sR+t\), where \(s\in\mathbb N_0\) and
	\(t\in\{0,\ldots,R-1\}\).  The states in factor copy \(a\) have heights
	\[
		\widetilde\omega_a,
		\widetilde\omega_a+R,
		\widetilde\omega_a+2R,\ldots .
	\]
	Consequently, the number \(n_{a,N}\) of retained states in that copy is
	\(s+1\) when \(\widetilde\omega_a\le t\), and is \(s\) otherwise.
	Because every residue in \(\{0,\ldots,R-1\}\) occurs exactly \(d\) times
	among \(\widetilde\omega_1,\ldots,\widetilde\omega_M\), precisely
	\(d(t+1)\) copies have size \(s+1\).

	If \(t<R-1\), then \(m_N=s\) and
	\(r_N=d(t+1)\).  If \(t=R-1\), then every copy has size
	\(s+1=m_N\) and \(r_N=0\).  Thus every factor copy contains either
	\(m_N\) or \(m_N+1\) retained states, exactly \(r_N\) copies have the
	larger size, and at least one copy has size \(m_N\).  Choose
	one such copy.  In
	the repeated Cauchy--Binet expansion of the determinant of its
	\(m_N\times m_N\) cyclic product, select the first \(m_N\) indices at
	every intermediate copy.  This term is the product of leading square
	minors of positive bidiagonal factors and is strictly positive; all other
	terms are nonnegative.  This diagonal block of \(K_N^M\) is therefore
	nonsingular.

	\emph{Common nonzero spectrum.}
	The nonsingular block just constructed is totally nonnegative.  Its
	\(m_N\) eigenvalues are therefore positive.  Repeated use of the
	finite-dimensional \(AB/BA\) identity shows that all \(M\) cyclic blocks
	of \(K_N^M\) have these same \(m_N\) nonzero eigenvalues, counting
	algebraic multiplicity; any additional eigenvalues of a larger block are
	zero.  Hence \(K_N^M\) has exactly \(Mm_N\) nonzero eigenvalues.

	By the spectral mapping theorem, the eigenvalues of \(K_N^M\), with
	algebraic multiplicity, are the \(M\)-th powers of the eigenvalues of
	\(K_N\).  Therefore \(K_N\) itself has exactly \(Mm_N\) nonzero
	eigenvalues.  Its zero eigenvalue has algebraic
	multiplicity
	\[
		D_N-Mm_N=r_N.
	\]
	Consequently \(\mathcal R_N(0)\ne0\).  Set
	\(\mathcal P_N=\mathcal R_N\); this proves both the asserted factor
	\(\lambda^{r_N}\) and the formula
	\(\deg\mathcal P_N=m_N\).

	The block \(C_{a,N}\) has order
	\(n_{a,N}\in\{m_N,m_N+1\}\); by the preceding argument its nonzero
	eigenvalues are the same \(m_N\) positive numbers and its remaining
	eigenvalues are zero.  Since
	both characteristic polynomials are monic, this proves
	\[
		\det(yI-C_{a,N})
		=
		y^{n_{a,N}-m_N}\mathcal P_N(y).
	\]
	The block \(C_{a,N}\) is totally nonnegative, so all its eigenvalues are
	real and nonnegative.  The displayed characteristic-polynomial identity
	shows that the zeros of \(\mathcal P_N\) are exactly its nonzero
	eigenvalues.  They are therefore positive.  Finally, the spectral mapping
	theorem implies that every nonzero eigenvalue \(\lambda\) of \(K_N\)
	satisfies \(\lambda^M>0\), which proves the claimed ray localization.
\end{proof}

Proposition~\ref{prop:cyclic-star-resolvent-representation} identifies the
Weyl matrix observed from the first factor copy in the bounded case.  The
next corollary shows that a global PBF also supplies positive Favard data at
every factor position, whether or not the factors are bounded and whether or
not their origin is a Christoffel transformation.

\begin{lemma}[Compatible truncations of a global PBF]
	\label{lem:global-PBF-compatible-truncations}
	Let
	\[
		S=K_1\cdots K_pV_q\cdots V_1
	\]
	be a semi-infinite product in which the \(K_a\) are lower bidiagonal and
	the \(V_b\) are upper bidiagonal.  If \(S^{[N]}\), \(K_a^{[N]}\), and
	\(V_b^{[N]}\) denote their leading principal submatrices of order
	\(N+1\), then
	\begin{equation}
		\label{eq:global-PBF-compatible-truncations}
		S^{[N]}
		=
		K_1^{[N]}\cdots K_p^{[N]}
		V_q^{[N]}\cdots V_1^{[N]}.
	\end{equation}
	Consequently, if the semi-infinite factorization is positive, all its
	leading principal submatrices have PBFs with zero scalar shift.  These PBFs
	are coherent in the sense of the unbounded Favard theorem: the finite
	boundary matrices are independent of \(N\) once their fixed number of
	coordinates is present.
\end{lemma}

\begin{proof}
	An entry of the product is a sum over intermediate coordinate indices.
	Starting from a row index at most \(N\), the lower factors can only decrease
	the current index or leave it fixed, so the index is still at most \(N\)
	when the upper factors begin.  The upper factors can only increase the
	index or leave it fixed.  If the index then became larger than \(N\), it
	could not return to a final column index at most \(N\).  Hence every
	contributing intermediate index also lies in
	\(\{0,\ldots,N\}\), which proves
	\eqref{eq:global-PBF-compatible-truncations}.

	If the factors are positive, their leading submatrices give the asserted
	finite PBFs.  The normalization convention used in
	\cite{BranquinhoFoulquieManas2026Unbounded} places one positive diagonal
	matrix between the lower and upper factors and makes every upper factor
	unitriangular.  It is equivalent to the convention used here.  Indeed,
		put \(D_0=I\) and choose positive diagonal matrices
		\(D_1,\ldots,D_q\) recursively so that
		\[
			\overline U_b
			\coloneq
			D_b^{-1}U_bD_{b-1}
		\]
		has unit diagonal.  Equivalently, the diagonal of \(D_b\) is the
		entrywise product of the diagonals of \(U_b\) and \(D_{b-1}\).  Since
		\(U_bD_{b-1}=D_b\overline U_b\), successive cancellation gives
		\[
			U_q\cdots U_1
			=
			D_q\overline U_q\cdots\overline U_1.
		\]
		Thus \(D_q\) is precisely the positive diagonal factor required by that
		normalization, and all the new superdiagonal entries are positive.

	The boundary matrices of
	\cite[Definition~4.1]{BranquinhoFoulquieManas2026Unbounded} use only fixed
	upper-left portions of these factors.  Because the finite factors are the
	leading submatrices of one fixed semi-infinite factorization, those portions
	do not depend on \(N\) once \(N\ge\max\{p,q\}-1\).  Thus the coherence
	condition in that definition is satisfied: the same initial-condition
	matrices can be used for every \(N\).
\end{proof}

\begin{corollary}[Favard data for the complete cyclic orbit]
	\label{cor:positive-Favard-data-cyclic-orbit}
	Assume that the factorization
	\eqref{eq:general-two-sided-product} is positive.  For every
	\(a\in\{1,\ldots,M\}\), there exist nonsingular initial-condition
	matrices \(\xi_{\mathrm F}^{(a)}\) and
	\(\nu_{\mathrm F}^{(a)}\), two mixed-type polynomial systems, and an
	entrywise positive \(q\times p\) matrix of finite measures
	\(\mathrm d\Psi_{\mathrm F}^{(a)}\) such that
	\begin{equation}
		\label{eq:positive-Favard-data-cyclic-orbit}
		(\xi_{\mathrm F}^{(a)})^{-1}
		\mathsf E_{[q]}(T_{(a)})^n\mathsf E_{[p]}^{\top}
		(\nu_{\mathrm F}^{(a)})^{-\top}
		=
		\int x^n\,\mathrm d\Psi_{\mathrm F}^{(a)}(x),
		\qquad n\in\mathbb N_0.
	\end{equation}
	For every \(a\) for which \(T_{(a)}\) is bounded on
	\(\ell^2(\mathbb N_0)\), the measures may be
	chosen compactly supported and the corresponding diagonal projection of
	the common cyclic resolvent has this rectangular Stieltjes representation.
	If \(T_{(a)}\) is unbounded on \(\ell^2(\mathbb N_0)\), the moment identity
	remains valid.  An operator-resolvent representation then requires the
	additional analytic hypotheses stated below.
\end{corollary}

\begin{proof}
	Theorem~\ref{thm:cyclic-PBF-closure} gives a normalized positive
	bidiagonal factorization in standard lower--upper order for every
	\(T_{(a)}\).  Lemma~\ref{lem:global-PBF-compatible-truncations} shows that
	its leading principal submatrices have coherent PBFs with zero shift and
	the same initial-condition matrices.  Fix \(a\).  If \(T_{(a)}\) is
	bounded on \(\ell^2(\mathbb N_0)\), the last part of
	Lemma~\ref{lem:height-operator-norm-bound},
	applied to its nonnegative factorization, shows that the refactorized
	bidiagonal coefficients are uniformly bounded.  The bounded Favard theorem
	\cite{BranquinhoFoulquieManas2023Spectral} then gives an entrywise positive
	matrix of compactly supported measures.  If \(T_{(a)}\) is unbounded on
	\(\ell^2(\mathbb N_0)\), the
	unbounded Favard theorem
	\cite{BranquinhoFoulquieManas2026Unbounded} applies to the coherent family
	just described and gives an entrywise positive matrix of finite measures.
	In both cases one obtains the moment identity
	\eqref{eq:positive-Favard-data-cyclic-orbit}.  In the bounded case,
	Proposition~\ref{prop:Darboux-shifted-resolvent-blocks} also identifies the
	corresponding diagonal block of the cyclic resolvent with the resolvent of
	\(T_{(a)}\).
\end{proof}

\begin{remark}[Normality of the cyclic Favard systems]
	\label{rem:unbounded-cyclic-Favard-data}
	By Theorem~\ref{thm:cyclic-PBF-closure}, every cyclic product has a PBF and
	is therefore BTP.  The maximal-degree criterion of
	\cite{BranquinhoFoulquieManas2024BTP} then gives strong normality when
	\(\nu^{-\top}\), the inverse of the left initial block \(\nu^{\top}\),
	is upper \(\Delta\mathrm{TP}\), and \(\xi^{-1}\), the inverse of the
	right initial block \(\xi\), is lower \(\Delta\mathrm{TP}\).  Thus the
	Favard theorem supplies the measures
	in both the bounded and unbounded cases, while the BTP theorem supplies the
	maximal-degree condition and hence strong normality, in the terminology of
	Definition~\ref{def:weak-strong-step-line-normality}.
\end{remark}

\subsection{Christoffel transformations adapted to the cyclic PBF}
\label{subsec:PBF-adapted-Christoffel}

The Favard data supplied at a cyclic product depend on the two triangular
initial-condition matrices.  This freedom is part of the Christoffel
transformation: after a factor is moved cyclically, its natural matrix
polynomial need not be the fixed shift matrix \(\mathsf K_q(x)\) or
\(\mathsf K_p(x)^{\top}\).  The next proposition records exactly how a change
of initial conditions modifies the moments while leaving the recurrence
matrix unchanged.

\begin{proposition}[Moment functionals for different initial conditions]
	\label{prop:fixed-recurrence-initial-condition-change}
	Let \(\mathcal T\) be one of the semi-infinite \((p,q)\)-banded
	recurrence matrices considered above, admitting a positive bidiagonal
	factorization.  Let two pairs of nonsingular lower triangular
	initial-condition matrices reconstruct from the same matrix \(\mathcal T\)
	the block moments
	\[
		W_n
		=
		\xi^{-1}\mathsf E_{[q]}\mathcal T^n
		\mathsf E_{[p]}^{\top}\nu^{-\top},
		\qquad
		\widetilde W_n
		=
		\widetilde\xi^{-1}\mathsf E_{[q]}\mathcal T^n
		\mathsf E_{[p]}^{\top}\widetilde\nu^{-\top}.
	\]
	Then
	\begin{equation}
		\label{eq:fixed-recurrence-initial-condition-change}
		\widetilde W_n=G W_n H,
		\qquad
		G\coloneq\widetilde\xi^{-1}\xi,
		\qquad
		H\coloneq\nu^{\top}\widetilde\nu^{-\top}.
	\end{equation}
	The matrix \(G\) is nonsingular lower triangular and \(H\) is nonsingular
	upper triangular.

	Let \(\mathscr M\) and \(\widetilde{\mathscr M}\) be the corresponding
	scalar step-line moment matrices.  If
	\(\mathscr G_{\mathrm L}=\operatorname{diag}(G,G,\ldots)\) and
	\(\mathscr G_{\mathrm R}=\operatorname{diag}(H,H,\ldots)\), then
	\[
		\widetilde{\mathscr M}
		=
		\mathscr G_{\mathrm L}\mathscr M\mathscr G_{\mathrm R}.
	\]
	For every \(N\in\mathbb N_0\),
	\begin{equation}
		\label{eq:fixed-initial-condition-leading-minors}
		\det\widetilde{\mathscr M}^{[N]}
		=
		\det\mathscr G_{\mathrm L}^{[N]}\,
		\det\mathscr M^{[N]}\,
		\det\mathscr G_{\mathrm R}^{[N]}.
	\end{equation}
	Consequently, \(\mathscr M\) is weakly normal at every scalar order if
	and only if \(\widetilde{\mathscr M}\) is.  This triangular gauge
	statement concerns weak normality only; it does not in general preserve
	the maximal component degrees required for strong normality.
\end{proposition}

\begin{proof}
	Equation~\eqref{eq:fixed-recurrence-initial-condition-change} follows by
	inserting the definitions and cancelling
	\(\xi\xi^{-1}\) and
	\(\nu^{-\top}\nu^{\top}\).  Repeating this identity in every row and
	column block gives
	\(\widetilde{\mathscr M}
	=\mathscr G_{\mathrm L}\mathscr M\mathscr G_{\mathrm R}\).
	Products and inverses of nonsingular lower triangular matrices are lower
	triangular, so \(G\) is lower triangular; transposition gives that \(H\)
	is upper triangular.
	Because \(\mathscr G_{\mathrm L}\) is lower triangular and
	\(\mathscr G_{\mathrm R}\) is upper triangular, their first \(N+1\)
	rows and columns do not receive contributions from indices outside the
	leading \((N+1)\times(N+1)\) blocks.
	Taking determinants therefore proves
	\eqref{eq:fixed-initial-condition-leading-minors}.  The two triangular
	determinants are nonzero.
\end{proof}

For the remainder of this subsection, the initial conditions are part of the
polynomial data.  More precisely, a right polynomial array
\(\boldsymbol B\) for a recurrence matrix \(\mathcal T\), with right
initial-condition matrix \(\xi\), is the solution of
\(\mathcal T\boldsymbol B(x)=x\boldsymbol B(x)\) satisfying
\(\mathsf E_{[q]}\boldsymbol B(x)=\xi\).  Similarly, a left polynomial
array \(\boldsymbol A\), with left initial-condition matrix \(\nu\), satisfies
\(\boldsymbol A(x)\mathcal T=x\boldsymbol A(x)\) and
\(\boldsymbol A(x)\mathsf E_{[p]}^{\top}=\nu^{\top}\).  Once the
nonsingular matrices \(\xi\) and \(\nu\) have been fixed, the two recurrence
equations determine their polynomial arrays uniquely.  Strong normality means
in addition that these arrays attain the prescribed maximal step-line
degrees.

Consider a PBF
\(T_0=F_1\cdots F_{M-1}U\), where the last factor \(U\) is normalized
upper bidiagonal, and set
\(D\coloneq F_1\cdots F_{M-1}\) and \(T_1\coloneq UD\).  Hence
\(T_0=DU\), while \(T_1\) is the cyclic product obtained by moving \(U\)
to the left end.

Fix weakly normal biorthogonal data
\((\boldsymbol B_0,\boldsymbol A_0,\mathrm d\Psi_0)\) for \(T_0\), and
assume that the right array \(\boldsymbol B_0\) attains its prescribed
maximal step-line degrees.  Let \(\xi_0,\nu_0\) be its nonsingular lower
triangular initial-condition matrices, so that
\(\mathsf E_{[q]}\boldsymbol B_0=\xi_0\) and
\(\boldsymbol A_0\mathsf E_{[p]}^{\top}=\nu_0^{\top}\).  By
Theorem~\ref{thm:cyclic-PBF-closure}, \(T_1\) also has a PBF and is therefore
BTP.  Hence the maximal-degree criterion recalled in
Remark~\ref{rem:unbounded-cyclic-Favard-data} guarantees the existence of
admissible right initial conditions; in particular, it is enough to choose
\(\xi_1\) so that \(\xi_1^{-1}\) is lower \(\Delta\mathrm{TP}\).
Fix such a matrix \(\xi_1\), and let \(\boldsymbol B_1\) be the unique
solution of
\[
	T_1\boldsymbol B_1(x)=x\boldsymbol B_1(x),
	\qquad
	\mathsf E_{[q]}\boldsymbol B_1(x)=\xi_1,
\]
which then has the prescribed maximal step-line degrees.

\begin{theorem}[Christoffel transformation associated with an upper cyclic factor]
	\label{thm:PBF-adapted-cyclic-Christoffel}
	Under the preceding assumptions, the following statements hold.
	\begin{enumerate}[label=\textup{(\roman*)}]
		\item There is a unique \(q\times q\) matrix polynomial \(S(x)\)
		such that
		\begin{equation}
			\label{eq:PBF-adapted-upper-Christoffel-connection}
			U\boldsymbol B_0(x)=\boldsymbol B_1(x)S(x).
		\end{equation}
		It is given by
		\[
			S(x)=\xi_1^{-1}\mathsf E_{[q]}U\boldsymbol B_0(x).
		\]
		Moreover, \(S\) is regular and has degree one, its leading coefficient
		has rank one, and \(\det S(x)=cx\) for some \(c\ne0\).

		\item Set
		\begin{equation}
			\label{eq:PBF-adapted-upper-Christoffel-data}
			\boldsymbol A_1\coloneq\boldsymbol A_0U^{-1},
			\qquad
			\nu_1^{\top}\coloneq
			\boldsymbol A_1\mathsf E_{[p]}^{\top},
			\qquad
			\mathrm d\Psi_1(x)\coloneq
			S(x)\,\mathrm d\Psi_0(x).
		\end{equation}
		Then \(\nu_1\) is nonsingular lower triangular,
		\(\boldsymbol A_1(x)T_1=x\boldsymbol A_1(x)\), and
		\[
			\int\boldsymbol B_1(x)\,\mathrm d\Psi_1(x)\,
			\boldsymbol A_1(x)=I.
		\]
		The array \(\boldsymbol A_1\) belongs componentwise to the prescribed
		left step-line polynomial spaces, while \(\boldsymbol B_1\) attains
		its prescribed maximal degrees.  Consequently,
		\((\boldsymbol B_1,\boldsymbol A_1,\mathrm d\Psi_1)\) is a weakly
		normal biorthogonal system with recurrence matrix \(T_1\) and
		initial-condition matrices \(\xi_1,\nu_1\).

		If, in addition, \(\boldsymbol A_1\) attains its prescribed maximal
		step-line degrees, then the transformed system is strongly normal.
		In particular, this additional conclusion holds if
		\(\nu_0^{-\top}\) is upper \(\Delta\mathrm{TP}\).

		\item For any other pair of nonsingular lower triangular
		initial-condition matrices
		\((\widetilde\xi_1,\widetilde\nu_1)\), the moment functional
		reconstructed from the same recurrence matrix \(T_1\) is represented,
		whenever \(\mathrm d\Psi_1\) is a representing matrix of measures, by
		\[
			\mathrm d\widetilde\Psi_1
			=
			\widetilde\xi_1^{-1}\xi_1\,
			\mathrm d\Psi_1\,
			\nu_1^{\top}\widetilde\nu_1^{-\top}.
		\]
		Such a triangular change preserves weak normality.  Strong normality
		of the regauged polynomial arrays must be imposed or verified
		separately.
	\end{enumerate}
\end{theorem}

\begin{proof}
	The eigenvalue equation
	\(T_0\boldsymbol B_0(x)=x\boldsymbol B_0(x)\), together with
	\(T_0=DU\) and \(T_1=UD\), gives
	\(T_1(U\boldsymbol B_0)=UDU\boldsymbol B_0
	=UT_0\boldsymbol B_0=xU\boldsymbol B_0\).
	Thus \(U\boldsymbol B_0\) is a right polynomial solution of the
	recurrence associated with \(T_1\).
	
	A right polynomial solution is determined by its first \(q\) rows.
	Applying \(\mathsf E_{[q]}\) to
	\eqref{eq:PBF-adapted-upper-Christoffel-connection} and using
	\(\mathsf E_{[q]}\boldsymbol B_1=\xi_1\) gives necessarily
	\[
		S(x)=\xi_1^{-1}\mathsf E_{[q]}U\boldsymbol B_0(x).
	\]
	With this choice, the two sides of
	\eqref{eq:PBF-adapted-upper-Christoffel-connection} have the same first
	\(q\) rows and satisfy the same recurrence for \(T_1\).  Uniqueness of
	the polynomial solution therefore proves the identity in every row.
	
	We next prove that \(S\) is regular.  If \(\det S\) were identically
	zero, there would be a nonzero vector
	\(\boldsymbol v(x)\in\mathbb C(x)^q\) such that
	\(S(x)\boldsymbol v(x)=0\).  Equation
	\eqref{eq:PBF-adapted-upper-Christoffel-connection} would then imply
	\(U\boldsymbol B_0(x)\boldsymbol v(x)=0\).  Applying \(D\) and using
	\(T_0=DU\), one would obtain
	\(x\boldsymbol B_0(x)\boldsymbol v(x)=0\), and hence
	\(\boldsymbol B_0(x)\boldsymbol v(x)=0\).  Its first \(q\) rows form
	the nonsingular matrix \(\xi_0\), so
	\(\boldsymbol v(x)=0\), a contradiction.  Therefore
	\(\det S\not\equiv0\).
	
	Write \(\boldsymbol B_{0,n}(x)\) for the \(n\)-th row of
\(\boldsymbol B_0(x)\), and denote its \(s\)-th component by
\(B_{0,n}^{(s)}(x)\), where \(s\in\{0,\ldots,q-1\}\).
If \(n=qm+j\), with \(j\in\{0,\ldots,q-1\}\), the maximal
step-line degrees are
\[
\deg B_{0,n}^{(s)}
=
\begin{cases}
	m,   & 0\le s\le j,\\
	m-1, & j<s\le q-1,
\end{cases}
\]
with the usual convention that a negative upper bound means that the
corresponding polynomial vanishes.  In particular,
\(\boldsymbol B_{0,0},\ldots,\boldsymbol B_{0,q-1}\) are constant,
whereas \(\boldsymbol B_{0,q}\) has degree one and its coefficient of
\(x\) is nonzero only in its first component.

Since \(U\) is normalized upper bidiagonal,
\[
\bigl(U\boldsymbol B_0\bigr)_n
=
u_n\boldsymbol B_{0,n}
+
\boldsymbol B_{0,n+1}.
\]
Therefore the rows
\(\bigl(U\boldsymbol B_0\bigr)_0,\ldots,
\bigl(U\boldsymbol B_0\bigr)_{q-2}\) are constant, while
\(\bigl(U\boldsymbol B_0\bigr)_{q-1}\) has degree one.  Its
	coefficient of \(x\) is nonzero by the assumed maximal-degree property of
	\(\boldsymbol B_0\).  Consequently,
\[
\mathsf E_{[q]}U\boldsymbol B_0(x)=C_0+xC_1,
\]
where \(C_1\) has exactly one nonzero row and hence rank one.
Finally,
\[
S(x)
=
\xi_1^{-1}\mathsf E_{[q]}U\boldsymbol B_0(x)
=
S_0+xS_1.
\]
	Because \(\xi_1^{-1}\) is nonsingular,
	\(\operatorname{rank}S_1=\operatorname{rank}C_1=1\).  By
	multilinearity of the determinant,
\(\deg\det S\le1\).  
	
	It remains to locate the zero of \(\det S\).  The matrix
	\(D=F_1\cdots F_{M-1}\) is the product of the remaining bidiagonal
	factors.  Since \(U\) is the last upper factor of the PBF, the first row
	of \(D\) is supported in the first \(q\) coordinates and its restriction
	to those coordinates is nonzero.  From
	\(D\boldsymbol B_1(x)S(x)=DU\boldsymbol B_0(x)
	=x\boldsymbol B_0(x)\), evaluation at \(x=0\) gives
	\(D\boldsymbol B_1(0)S(0)=0\).  Since the first \(q\) rows of
	\(\boldsymbol B_1(0)\) form the nonsingular matrix \(\xi_1\), the first row of
	\(D\boldsymbol B_1(0)\) is nonzero and provides a nonzero left null
	vector of \(S(0)\).  Thus \(\det S(0)=0\).  Since \(S\) is regular and
	\(\deg\det S\le1\), it follows that \(\det S(x)=cx\) with \(c\ne0\).
	
	For the left polynomial solution,
	\(\boldsymbol A_0T_0=x\boldsymbol A_0\) and \(T_0=DU\) imply
	\(\boldsymbol A_0D=x\boldsymbol A_0U^{-1}\).  Therefore
	\(\boldsymbol A_1T_1
	=\boldsymbol A_0U^{-1}UD
	=\boldsymbol A_0D
	=x\boldsymbol A_0U^{-1}
	=x\boldsymbol A_1\).
	
	Finally, using
	\(\boldsymbol B_1S=U\boldsymbol B_0\) and
	\(\boldsymbol A_1=\boldsymbol A_0U^{-1}\), one obtains
	\begin{align*}
		\int
		\boldsymbol B_1\,
		\mathrm d\Psi_1\,
		\boldsymbol A_1
		&=
		\int
		\boldsymbol B_1S\,
		\mathrm d\Psi_0\,
		\boldsymbol A_0U^{-1} =
		U
		\left(
		\int
		\boldsymbol B_0\,
		\mathrm d\Psi_0\,
		\boldsymbol A_0
		\right)
		U^{-1}
		=I.
	\end{align*}
	It remains to record exactly what the triangular transformation implies
	for the left polynomial array.  Put
	\[
		U_{[p]}
		\coloneq
		\mathsf E_{[p]}U\mathsf E_{[p]}^{\top}.
	\]
	Since \(U^{-1}\) is upper triangular, its first \(p\) columns are
	supported in the first \(p\) rows, and hence
	\[
		U^{-1}\mathsf E_{[p]}^{\top}
		=
		\mathsf E_{[p]}^{\top}U_{[p]}^{-1}.
	\]
	Therefore
	\[
		\nu_1^{\top}
		=
		\boldsymbol A_0U^{-1}\mathsf E_{[p]}^{\top}
		=
		\nu_0^{\top}U_{[p]}^{-1}.
	\]
	Both factors on the right are nonsingular upper triangular, so
	\(\nu_1\) is nonsingular lower triangular.

	Moreover, for every \(n\ge0\),
	\[
		\boldsymbol A_{1,n}
		=
		\sum_{k=0}^{n}
		\boldsymbol A_{0,k}(U^{-1})_{k,n}.
	\]
	For each component, the maximal degree allowed at position \(k\) is no
	larger than the maximal degree allowed at position \(n\) whenever
	\(k\le n\).  Thus \(\boldsymbol A_{1,n}\) belongs to the prescribed
	left polynomial space.  This triangular combination need not, however,
	preserve a nonzero maximal-degree coefficient.

	The biorthonormality identity proved above makes every finite initial
	segment of each transformed family linearly independent.  Since the
	polynomials lie in the corresponding step-line flags, they form dual
	bases of those flags.  Hence the transformed system is weakly normal.
	The right array \(\boldsymbol B_1\) has maximal degrees by the choice of
	\(\xi_1\).

	Finally, if \(\nu_0^{-\top}\) is upper \(\Delta\mathrm{TP}\), then
	\[
		\nu_1^{-\top}=U_{[p]}\nu_0^{-\top}.
	\]
	The matrix \(U_{[p]}\) is upper bidiagonal and totally nonnegative, with
	positive diagonal.  For every upper-triangularly admissible minor with
	row and column sets \(I,J\), the \(K=I\) term in the Cauchy--Binet sum is
	\(\det U_{[p]}[I,I]\det\nu_0^{-\top}[I,J]>0\), while all other terms are
	nonnegative.  Thus \(\nu_1^{-\top}\) is upper
	\(\Delta\mathrm{TP}\).  Since
	\(T_1\) is BTP and \(\xi_1^{-1}\) is lower
	\(\Delta\mathrm{TP}\), the maximal-degree criterion recalled in
	Remark~\ref{rem:unbounded-cyclic-Favard-data} gives strong normality.
	The final gauge statement is
	Proposition~\ref{prop:fixed-recurrence-initial-condition-change} applied
	to the two pairs of lower triangular initial conditions.
	This completes the proof.
\end{proof}

\begin{corollary}[Christoffel transformation associated with a lower cyclic factor]
	\label{cor:PBF-adapted-lower-cyclic-Christoffel}
	Let \(T_0=LF_2\cdots F_M\), where \(L\) is normalized lower
	bidiagonal, and put \(D\coloneq F_2\cdots F_M\) and \(T_1\coloneq DL\).
	Fix weakly normal biorthogonal data
	\((\boldsymbol B_0,\boldsymbol A_0,\mathrm d\Psi_0)\) for \(T_0\),
	with nonsingular lower triangular initial-condition matrices
	\(\xi_0,\nu_0\),
	and assume that the left array \(\boldsymbol A_0\) attains its prescribed
	maximal step-line degrees.  Choose a nonsingular lower triangular left
	initial-condition matrix \(\nu_1\) for which the corresponding left
	polynomial array \(\boldsymbol A_1\) for \(T_1\) has its prescribed
	maximal step-line degrees.

	Then there is a regular \(p\times p\) matrix polynomial \(R(x)\) of
	degree one, with rank-one leading coefficient and
	\(\det R(x)=c'x\), \(c'\ne0\), such that
	\[
		\boldsymbol A_0(x)L=R(x)\boldsymbol A_1(x),
		\qquad
		\mathrm d\Psi_1(x)=\mathrm d\Psi_0(x)R(x).
	\]
	Define
	\[
		\boldsymbol B_1\coloneq L^{-1}\boldsymbol B_0,
		\qquad
		\xi_1\coloneq\mathsf E_{[q]}\boldsymbol B_1.
	\]
	Then \(\xi_1\) is nonsingular lower triangular, and
	\((\boldsymbol B_1,\boldsymbol A_1,\mathrm d\Psi_1)\) is a weakly
	normal biorthogonal system with recurrence matrix \(T_1\).  The left
	array \(\boldsymbol A_1\) has maximal step-line degrees by construction.

	If, in addition, \(\boldsymbol B_1\) attains its prescribed maximal
	step-line degrees, then the transformed system is strongly normal.  In
	particular, this additional conclusion holds if \(\xi_0^{-1}\) is lower
	\(\Delta\mathrm{TP}\).  Different nonsingular lower triangular initial
	conditions change this functional by the constant transformations of
	Proposition~\ref{prop:fixed-recurrence-initial-condition-change} and
	preserve weak normality; strong normality must be verified separately
	after such a change.
\end{corollary}

\begin{proof}
	Transpose the argument of
	Theorem~\ref{thm:PBF-adapted-cyclic-Christoffel}.  The connection
	polynomial now acts on the right of the matrix of measures.  Under
	transposition, the nonzero subdiagonal entries of \(L\) become a nonunit
	superdiagonal; the proof of regularity, the rank-one leading coefficient,
	and the simple zero of the determinant uses only their nonvanishing, so
	the same argument applies with these scalar factors retained.  The
	identities
	\[
		T_1L^{-1}=L^{-1}T_0,
		\qquad
		\boldsymbol A_0L=R\boldsymbol A_1
	\]
	give the transformed recurrence equations, and
	\[
		\int
		\boldsymbol B_1\,\mathrm d\Psi_1\,\boldsymbol A_1
		=
		L^{-1}
		\left(
		\int
		\boldsymbol B_0\,\mathrm d\Psi_0\,\boldsymbol A_0
		\right)L
		=I.
	\]
	Since \(L^{-1}\) is lower triangular,
	\(\boldsymbol B_{1,n}\) is a finite linear combination of
	\(\boldsymbol B_{0,0},\ldots,\boldsymbol B_{0,n}\); hence it belongs
	to its prescribed right polynomial space.  Its leading \(q\times q\)
	block is the product of two nonsingular lower triangular matrices, so
	\(\xi_1\) is nonsingular lower triangular.  Biorthonormality then gives
	weak normality.  Triangularity alone does not prevent cancellation of a
	maximal right degree, so strong normality requires the additional
	hypothesis stated in the corollary.

	For the sufficient condition, put
	\(L_{[q]}\coloneq
	\mathsf E_{[q]}L\mathsf E_{[q]}^{\top}\).  The initial blocks satisfy
	\(\xi_1=L_{[q]}^{-1}\xi_0\), and hence
	\[
		\xi_1^{-1}=\xi_0^{-1}L_{[q]}.
	\]
	If \(\xi_0^{-1}\) is lower \(\Delta\mathrm{TP}\), then for every
	lower-triangularly admissible minor with row and column sets \(I,J\), the
	\(K=J\) term in the Cauchy--Binet sum for
	\(\xi_0^{-1}L_{[q]}\) is strictly positive and all other terms are
	nonnegative.  Thus \(\xi_1^{-1}\) is lower
	\(\Delta\mathrm{TP}\).  The BTP criterion
	therefore supplies the missing maximal right degrees.
\end{proof}

\begin{corollary}[A complete factor circuit]
	\label{cor:PBF-adapted-complete-circuit}
	Recall that \(T=AB\), where \(A=L_1\cdots L_p\) and
	\(B=U_q\cdots U_1\).  Set \(a_0\coloneq1\) and
	\(a_j\coloneq M-j+1\) for \(j\in\{1,\ldots,q\}\).
	Moving successively the rightmost upper factor to the left carries
	\(T_{(a_{j-1})}\) into \(T_{(a_j)}\).  In particular,
	\(T_{(a_0)}=T\) and \(T_{(a_q)}=T_{(p+1)}=BA\).
	
	Fix a weakly normal biorthogonal system
	\((\boldsymbol B^{(a_0)},\boldsymbol A^{(a_0)},\mathrm d\Psi)\) at
	\(T_{(a_0)}=T\), and assume that its right array attains the prescribed
	maximal step-line degrees.  For each \(j\in\{1,\ldots,q\}\), choose a
	nonsingular lower triangular right initial-condition matrix
	\(\xi^{(a_j)}\) for which the corresponding right polynomial array
	\(\boldsymbol B^{(a_j)}\) has the prescribed maximal step-line degrees.
	Transport the left arrays and moment functionals successively by the
	upper-factor transformations of
	Theorem~\ref{thm:PBF-adapted-cyclic-Christoffel}.  Let \(S_j(x)\) be the
	connection matrix associated by that theorem with the transfer from
	\(T_{(a_{j-1})}\) to \(T_{(a_j)}\).  Then the following statements hold.
	\begin{enumerate}[label=\textup{(\roman*)}]
		\item There is a nonsingular constant lower triangular matrix \(H\)
		such that
		\begin{equation}
			\label{eq:PBF-adapted-complete-circuit}
			S_q(x)\cdots S_1(x)=xH.
		\end{equation}
		The matrix \(H\) is the lower triangular endpoint-normalization factor.
		
		\item The composition of the \(q\) Christoffel transformations gives
		the representative
		\(xH\,\mathrm d\Psi\) at \(T_{(a_q)}=BA\).
		Replacing the endpoint right initial condition by
		\(\widehat\xi^{(a_q)}=\xi^{(a_q)}H\) applies the constant left gauge
		\(H^{-1}\) to the functional, and the endpoint representative becomes
		\(x\,\mathrm d\Psi\).
		
		\item The circuit obtained by moving successively the lower factors
		has the transposed form: one chooses a left-maximal array at every
		vertex and transports the right arrays and functionals.  Its connection
		matrices act on the right and their product is \(xK\), where \(K\) is
		a nonsingular constant upper triangular matrix.
	\end{enumerate}
\end{corollary}

\begin{proof}
	Put \(P(x)\coloneq S_q(x)\cdots S_1(x)\).  At the \(j\)-th transfer,
	Theorem~\ref{thm:PBF-adapted-cyclic-Christoffel} relates the right
	polynomial arrays at \(T_{(a_{j-1})}\) and \(T_{(a_j)}\).
	Its one-sided conclusion is sufficient for iteration: the transported
	system remains weakly normal, while the right array at the next vertex is
	maximal by construction.
	Composing these \(q\) relations gives
	\(B\boldsymbol B^{(a_0)}(x)
	=\boldsymbol B^{(a_q)}(x)P(x)\).
	
	The first \(q\) rows of \(B\boldsymbol B^{(a_0)}(x)\) have degree at
	most one.  Indeed, \(B\) is a product of \(q\) upper bidiagonal matrices,
	so these rows involve only the first \(2q\) rows of
	\(\boldsymbol B^{(a_0)}\), whose maximal step-line degrees are at most
	one.  The first \(q\) rows of \(\boldsymbol B^{(a_q)}\) form the
	nonsingular constant matrix \(\xi^{(a_q)}\).  It follows that
	\(\deg P\le1\).
	
	The eigenvalue equation for \(\boldsymbol B^{(a_0)}\) gives
	\(A(B\boldsymbol B^{(a_0)}(x))
	=T\boldsymbol B^{(a_0)}(x)
	=x\boldsymbol B^{(a_0)}(x)\).
	Evaluating at \(x=0\) and using the invertibility of the lower triangular
	matrix \(A\) yields \(B\boldsymbol B^{(a_0)}(0)=0\).  The endpoint
	connection identity then gives
	\(\xi^{(a_q)}P(0)=0\).  Since \(\xi^{(a_q)}\) is nonsingular,
	\(P(0)=0\).  Therefore \(P(x)=xH\) for some constant matrix \(H\).
	
	For every \(j\in\{1,\ldots,q\}\),
	Theorem~\ref{thm:PBF-adapted-cyclic-Christoffel} gives
	\(\det S_j(x)=c_jx\), with \(c_j\ne0\).  Hence
	\(\det P(x)=(\prod_{j=1}^{q}c_j)x^q\).  On the other hand,
	\(P(x)=xH\) gives \(\det P(x)=x^q\det H\).  Thus
	\(\det H=\prod_{j=1}^{q}c_j\ne0\), and \(H\) is nonsingular.
	To see the triangular form explicitly, put
	\(A_{[q]}\coloneq
	\mathsf E_{[q]}A\mathsf E_{[q]}^{\top}\).  Since
	\(A(B\boldsymbol B^{(a_0)})=x\boldsymbol B^{(a_0)}\), the endpoint
	connection identity and \(P=xH\) give, after cancelling \(x\),
	\[
		\boldsymbol B^{(a_q)}H=A^{-1}\boldsymbol B^{(a_0)}.
	\]
	Taking the first \(q\) rows yields
	\[
		H
		=
		\bigl(\xi^{(a_q)}\bigr)^{-1}
		A_{[q]}^{-1}\xi^{(a_0)}.
	\]
	All three factors are lower triangular, so \(H\) is lower triangular.
	This proves \textup{(i)}.
	
	At each transfer the matrix of measures is multiplied on the left by
	the corresponding connection matrix.  After the \(q\) transfers it is
	therefore \(P(x)\,\mathrm d\Psi(x)=xH\,\mathrm d\Psi(x)\).
	The formula above gives
	\(\widehat\xi^{(a_q)}=\xi^{(a_q)}H
	=A_{[q]}^{-1}\xi^{(a_0)}\), which is nonsingular lower triangular.
	Moreover,
	\[
		\bigl(\widehat\xi^{(a_q)}\bigr)^{-1}\xi^{(a_q)}=H^{-1}.
	\]
	Proposition~\ref{prop:fixed-recurrence-initial-condition-change} therefore
	applies the constant left gauge \(H^{-1}\) to \(xH\,\mathrm d\Psi\)
	and produces the representative \(x\,\mathrm d\Psi\).  This proves
	\textup{(ii)}.
	
	Applying the same argument to the left polynomial arrays and the lower
	bidiagonal factors proves \textup{(iii)}.
\end{proof}

The connection matrices depend covariantly on the chosen initial conditions.
Indeed, if \(\boldsymbol B_j\) is replaced by
\(\boldsymbol B_j^{\sharp}=\boldsymbol B_jC_j\), with \(C_j\) a constant
nonsingular lower triangular matrix, then
\[
	S_j^{\sharp}(x)=C_j^{-1}S_j(x)C_{j-1}.
\]
Thus all intermediate gauges cancel in a complete circuit; only the two
endpoint gauges change the constant matrix \(H\).

\begin{corollary}[Comparison of endpoint functionals]
	\label{cor:PBF-adapted-path-independence}
	Consider two factor-transfer constructions of the preceding type, starting
	from the same strongly normal biorthogonal system with recurrence matrix
	\(T_{\mathrm{in}}\) and moment functional \(\Psi_{\mathrm{in}}\).  At
	every upper, respectively lower, step choose the right-, respectively
	left-maximal data required by
	Theorem~\ref{thm:PBF-adapted-cyclic-Christoffel} and
	Corollary~\ref{cor:PBF-adapted-lower-cyclic-Christoffel}.  The
	intermediate recurrence matrices and their lower triangular
	initial-condition matrices may be different.  Suppose that both
	constructions end at the same recurrence matrix
	\(T_{\mathrm{out}}\).
	
	Let \(\Psi_1\) and \(\Psi_2\) be the two endpoint moment functionals, and
	let \((\xi_1,\nu_1)\) and \((\xi_2,\nu_2)\) be their respective endpoint
	initial-condition matrices.  Then
	\[
	\Psi_2
	=
	G\Psi_1H,
	\qquad
	G=\xi_2^{-1}\xi_1,
	\qquad
	H=\nu_1^{\top}\nu_2^{-\top}.
	\]
	The matrix \(G\) is nonsingular lower triangular and \(H\) is
	nonsingular upper triangular.  If the two constructions use the same
	initial-condition matrices at \(T_{\mathrm{out}}\), then
	\(G=I_q\), \(H=I_p\), and \(\Psi_1=\Psi_2\).
\end{corollary}

\begin{proof}
	For \(i\in\{1,2\}\), reconstruction from the endpoint recurrence and
	initial conditions gives
	\[
		\bigl\langle\Psi_i,x^n\bigr\rangle
		=
		\xi_i^{-1}\mathsf E_{[q]}T_{\mathrm{out}}^n
		\mathsf E_{[p]}^{\top}\nu_i^{-\top},
		\qquad n\in\mathbb N_0.
	\]
	Substitution shows on every monomial that
	\(\Psi_2=(\xi_2^{-1}\xi_1)\Psi_1
	(\nu_1^{\top}\nu_2^{-\top})\), and hence proves the displayed
	functional identity.  The standard lower triangular orientation of the
	initial-condition matrices makes \(G\) lower triangular and \(H\) upper
	triangular.  If the endpoint matrices coincide, both constant
	transformations are the identity matrices.
\end{proof}

The recurrence solutions can also be transported directly through the
bidiagonal factors.  Under the Christoffel hypotheses introduced in the next
subsection, the transported solutions belong to the standard step-line
polynomial spaces and carry the normalization of the transformed mixed-type
orthogonal system.  For the first
cyclic step, write
\[
T=F_1G_1,
\qquad
G_1=F_2\cdots F_M,
\qquad
T_{(2)}=G_1F_1.
\]
If \(\boldsymbol B(x)\) and \(\boldsymbol A(x)\) are respectively right and
left polynomial solutions of the recurrence equations for \(T\), then
\begin{equation}
	\label{eq:first-cyclic-polynomial-transport}
	\widetilde{\boldsymbol B}^{(2)}(x)
	=F_1^{-1}\boldsymbol B(x)
	=\frac1xG_1\boldsymbol B(x),
	\qquad
	\widetilde{\boldsymbol A}^{(2)}(x)
	=\boldsymbol A(x)F_1.
\end{equation}
Indeed, \(TF_1=F_1T_{(2)}\).  Therefore
\[
	T_{(2)}F_1^{-1}\boldsymbol B
	=F_1^{-1}T\boldsymbol B
	=xF_1^{-1}\boldsymbol B,
\]
and the analogous calculation on the left proves the two recurrence
equations.  The identity
\(G_1\boldsymbol B=xF_1^{-1}\boldsymbol B\) follows from
\(T=F_1G_1\).  Since \(F_1=L_1\) is lower unitriangular,
\(F_1^{-1}\) acts on each component by a finite triangular sum, so the first
vector in \eqref{eq:first-cyclic-polynomial-transport} is polynomial.

At a general factor position, define
\[
	G_a\coloneq F_{a+1}\cdots F_MF_1\cdots F_{a-1},
	\qquad
	T_{(a)}=F_aG_a,
	\qquad
	T_{(a+1)}=G_aF_a,
\]
with cyclic indices.  The same intertwining calculation shows that
\(G_a\boldsymbol B^{(a)}\) and
\(\boldsymbol A^{(a)}F_a\) satisfy the right and left recurrence equations
for \(T_{(a+1)}\).  When the factors form the weakly normal Christoffel chain
of Theorem~\ref{thm:cyclic-Favard-transport}, the nonvanishing Christoffel
determinants ensure existence and divisibility, while the initial blocks
\(\xi^{(a+1)}\) and \(\nu^{(a+1)}\) select the standard normalization and
identify the transformed measures.  For a general factorization, the same
calculation gives the transport of recurrence solutions.

This observation has different content in the two applications below.  For
the Pi\~neiro system, each Christoffel step only shifts a power parameter, so
\eqref{eq:first-cyclic-polynomial-transport} becomes a contiguous relation
between already known hypergeometric formulas.  In the Jacobi-like family,
the lower steps are explicit but the upper steps require the upper
bidiagonal factors, whose entries are ratios of moment minors.  Once those
minors have been evaluated, the complete cyclic chain of polynomial systems
becomes explicit.

\begin{remark}
	This factor-adapted construction requires no separate nonvanishing condition
for the fixed moment shifts: its matrix polynomial is determined from the
two cyclic recurrence systems.  The fixed shift matrices considered next
give a more explicit statement, but their one-sided transforms must exist.
\end{remark}

\subsection{The two-sided Christoffel rectangle}
\label{subsec:complete-Christoffel-grid}

This subsection answers the following question.  Starting from a matrix of
measures and its bidiagonally factorized recurrence matrix, can the elementary
row and column Christoffel transformations be performed simultaneously while
retaining positive bidiagonal recurrence factorizations?

\begin{figure}[H]
	\centering
	\begin{tikzpicture}[
		scale=0.82,
		transform shape,
		>=stealth,
		vertex/.style={draw=Gray,fill=Gray!4,rounded corners=2pt,
			inner sep=4pt,font=\small},
		columnvertex/.style={vertex,draw=NavyBlue,fill=NavyBlue!7},
		rowvertex/.style={vertex,draw=BrickRed,fill=BrickRed!7},
		bothvertex/.style={vertex,draw=OliveGreen,fill=OliveGreen!9},
		harrow/.style={->,thick,NavyBlue},
		varrow/.style={->,thick,BrickRed}
	]
		\node[vertex] (p00) at (0,3) {
			\(\mathrm d\Psi_{(0,0)}=\mathrm d\Psi\)
		};
		\node[columnvertex] (p01) at (3.2,3) {
			\(\mathrm d\Psi_{(0,1)}\)
		};
		\node[text=NavyBlue] (pdot1) at (5.3,3) {\(\cdots\)};
		\node[columnvertex] (p0p) at (7.6,3) {
			\(\mathrm d\Psi_{(0,p)}=x\,\mathrm d\Psi\)
		};

		\node[text=BrickRed] (vdot1) at (0,1.7) {\(\vdots\)};
		\node[text=BrickRed] (vdot2) at (3.2,1.7) {\(\vdots\)};
		\node[text=BrickRed] (vdot3) at (7.6,1.7) {\(\vdots\)};

		\node[rowvertex] (pq0) at (0,0.4) {
			\(\mathrm d\Psi_{(q,0)}=x\,\mathrm d\Psi\)
		};
		\node[bothvertex] (pq1) at (3.2,0.4) {
			\(\mathrm d\Psi_{(q,1)}\)
		};
		\node[text=OliveGreen] (pdot2) at (5.3,0.4) {\(\cdots\)};
		\node[bothvertex] (pqp) at (7.6,0.4) {
			\(\mathrm d\Psi_{(q,p)}=x^2\,\mathrm d\Psi\)
		};

		\draw[harrow] (p00) -- node[above,text=NavyBlue] {
			\(\cdot\,\mathsf K_p^{\top}\)
		} (p01);
		\draw[harrow] (p01) -- (pdot1);
		\draw[harrow] (pdot1) -- (p0p);
		\draw[harrow] (pq0) -- (pq1);
		\draw[harrow] (pq1) -- (pdot2);
		\draw[harrow] (pdot2) -- (pqp);
		\draw[varrow] (p00) -- node[left,text=BrickRed] {
			\(\mathsf K_q\,\cdot\)
		} (vdot1);
		\draw[varrow] (vdot1) -- (pq0);
		\draw[varrow] (p01) -- (vdot2);
		\draw[varrow] (vdot2) -- (pq1);
		\draw[varrow] (p0p) -- (vdot3);
		\draw[varrow] (vdot3) -- (pqp);

		\node[vertex] (s00) at (10.2,2.8) {
			\(\mathrm d\Psi_{(b-1,a)}\)
		};
		\node[columnvertex] (s01) at (13.8,2.8) {
			\(\mathrm d\Psi_{(b-1,a+1)}\)
		};
		\node[rowvertex] (s10) at (10.2,0.6) {
			\(\mathrm d\Psi_{(b,a)}\)
		};
		\node[bothvertex] (s11) at (13.8,0.6) {
			\(\mathrm d\Psi_{(b,a+1)}\)
		};

		\draw[harrow] (s00) -- node[above,text=NavyBlue] {
			\(\mathcal L_{b-1,a}\)
		} (s01);
		\draw[harrow] (s10) -- node[below,text=NavyBlue] {
			\(\mathcal L_{b,a}\)
		} (s11);
		\draw[varrow] (s00) -- node[left,text=BrickRed] {
			\(\mathcal U_{b,a}\)
		} (s10);
		\draw[varrow] (s01) -- node[right,text=BrickRed] {
			\(\mathcal U_{b,a+1}\)
		} (s11);

		\node[rounded corners=2pt,draw=OliveGreen,fill=OliveGreen!7,
			inner sep=3pt] at (12,-0.25) {
			\(
			\mathcal U_{b,a}\mathcal L_{b-1,a}
			=
			\mathcal L_{b,a}\mathcal U_{b,a+1}
			\)
		};
	\end{tikzpicture}
	\caption{The two-sided Christoffel rectangle and its elementary commuting
	square.  Horizontal arrows are column-side Christoffel transformations and
	downward arrows are row-side Christoffel transformations; they are shown in
	blue and red, respectively, while green marks vertices reached in both
	directions.  A factorization
	path from \((q,0)\) to \((0,p)\) traverses horizontal arrows forward and
	vertical arrows backward.}
	\label{fig:two-sided-Christoffel-rectangle}
\end{figure}
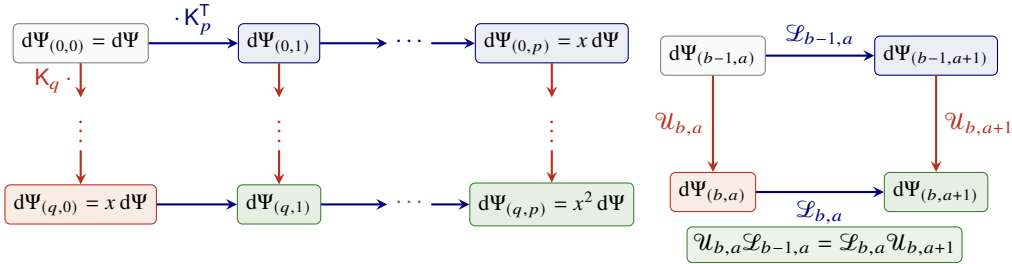

For \(b\in\{0,\ldots,q\}\) and \(a\in\{0,\ldots,p\}\), the vertices of the
rectangle will be represented by
\[
\mathrm d\Psi_{(b,a)}(x)
=
\mathsf K_q(x)^b\,\mathrm d\Psi(x)\,
\bigl(\mathsf K_p(x)^a\bigr)^{\top}.
\]
Its scalar moment matrix is obtained from the original moment matrix
\(\mathscr M\) by shifting the row index by \(b\) and the column index by
\(a\).  The main result shows that weak normality on the two boundary chains
propagates to every vertex.  At each vertex it identifies the transformed
measure, the shifted moment matrix, the recurrence matrix, its positive
bidiagonal factors, and the two mixed-type polynomial families.  The proof
reduces to the positive refactorization attached to one elementary square.
Figure~\ref{fig:two-sided-Christoffel-rectangle} displays both the full
rectangle and that elementary square.

\subsubsection{Boundary chains and shifted moment matrices}

Corollary~\ref{cor:positive-Favard-data-cyclic-orbit} gives positive spectral
data for every cyclic product.  We now describe the more specific situation
in which the prescribed bidiagonal factors are the connection matrices of
two one-sided Christoffel chains.  For \(r\ge1\), define
\[
\mathsf K_r(x)
\coloneq
\begin{bNiceMatrix}[margin=2pt]
	0_{(r-1)\times1}&I_{r-1}\\
	x&0_{1\times(r-1)}
\end{bNiceMatrix},
\qquad
\mathsf K_r(x)^r=xI_r.
\]
Starting from a \(q\times p\) matrix of measures \(\mathrm d\Psi\), define
the successive right Christoffel transforms by
\(\mathrm d\Psi_{(0,a)}(x)
\coloneq\mathrm d\Psi(x)(\mathsf K_p(x)^a)^{\top}\), for
\(a\in\{0,\ldots,p\}\), and the successive left Christoffel transforms by
\(\mathrm d\Psi_{(b,0)}(x)
\coloneq\mathsf K_q(x)^b\,\mathrm d\Psi(x)\), for
\(b\in\{0,\ldots,q\}\).
Thus \(\mathrm d\Psi_{(0,0)}=\mathrm d\Psi\), and the two chains have the
common endpoint
\(\mathrm d\Psi_{(0,p)}=\mathrm d\Psi_{(q,0)}=x\,\mathrm d\Psi\).

The chains are weakly normal when all their step-line orthogonality problems
exist and are unique up to nonzero constants.  At each degree, the
Christoffel formula is obtained from a finite linear system expressing the
required cancellation at \(x=0\).  The determinants of these systems are the
Christoffel \(\tau\)-determinants.  The Christoffel factorization theorem
\cite{BranquinhoFoulquieManas2026Factorization} states that weak normality of
the two chains is equivalent to the nonvanishing of the corresponding
\(\tau\)-determinants at every degree.

\begin{definition}[One-step Christoffel connection matrices]
	\label{def:one-step-Christoffel-connection-matrices}
	Let \(\mathrm d\Psi_{0,0}\) be a matrix of measures with a weakly normal
	mixed-type system \((\boldsymbol A_{0,0},\boldsymbol B_{0,0})\) and fixed
	initial conditions.  Suppose that the two elementary Christoffel
	transforms
	\[
		\mathrm d\Psi_{0,1}
		=\mathrm d\Psi_{0,0}\mathsf K_p^{\top},
		\qquad
		\mathrm d\Psi_{1,0}
		=\mathsf K_q\mathrm d\Psi_{0,0}
	\]
	are weakly normal, and equip their polynomial arrays with the initial
	conditions obtained from the corresponding one-step Christoffel formulas.
	The normalized lower bidiagonal matrix \(L\) characterized by
	\[
		\boldsymbol A_{0,0}L
		=\mathsf K_p^{\top}\boldsymbol A_{0,1},
		\qquad
		L\boldsymbol B_{0,1}=\boldsymbol B_{0,0}
	\]
	is called the one-step column-side Christoffel connection matrix.  The
	normalized upper bidiagonal matrix \(U\) characterized by
	\[
		U\boldsymbol B_{0,0}
		=\boldsymbol B_{1,0}\mathsf K_q,
		\qquad
		\boldsymbol A_{1,0}U=\boldsymbol A_{0,0}
	\]
	is called the one-step row-side Christoffel connection matrix.  The same
	terminology is used at every successive vertex of either one-sided chain.
\end{definition}

Denote by \(T_{(0,a)}\) and \(T_{(b,0)}\) the recurrence matrices of the
systems represented by \(\mathrm d\Psi_{(0,a)}\) and
\(\mathrm d\Psi_{(b,0)}\), respectively.  Their initial-condition matrices
are obtained successively from those at \((0,0)\) by the one-step
Christoffel formulas.

For \(a\in\{1,\ldots,M\}\), set
\[
\mathrm d\Psi^{(a)}
\coloneq
\begin{cases}
\mathrm d\Psi_{(0,a-1)},&1\le a\le p+1,\\
\mathrm d\Psi_{(M+1-a,0)},&p+2\le a\le M.
\end{cases}
\]
Let \(\mathscr M^{(a)}\) be its moment matrix and let
\(\xi^{(a)}\), \(\nu^{(a)}\) be the initial-condition matrices transported
by the corresponding one-step Christoffel transformations.

\begin{theorem}[Recurrence matrices under cyclic Christoffel shifts]
	\label{thm:cyclic-Favard-transport}
	Let \(\mathrm d\Psi\) be a \(q\times p\) matrix of measures with all
	moments finite, and suppose that its moment matrix is weakly normal along
	the step-line and has
	recurrence matrix \(T\) with initial-condition matrices \(\xi,\nu\).
	Suppose that both one-sided chains defined above are weakly normal.  Assume,
	in the sense of
	Definition~\ref{def:one-step-Christoffel-connection-matrices}, that
	\(L_a\) is the normalized column-side connection matrix from
	\(\mathrm d\Psi_{(0,a-1)}\) to \(\mathrm d\Psi_{(0,a)}\), for every
	\(a\in\{1,\ldots,p\}\), and that \(U_b\) is the normalized row-side
	connection matrix from \(\mathrm d\Psi_{(b-1,0)}\) to
	\(\mathrm d\Psi_{(b,0)}\), for every \(b\in\{1,\ldots,q\}\).
	Then, for every \(a\in\{1,\ldots,M\}\), \(T_{(a)}\) is the recurrence
	matrix of the transformed system with data
	\((\mathrm d\Psi^{(a)},\xi^{(a)},\nu^{(a)})\).  At the first position,
	\[
	\mathrm d\Psi^{(1)}=\mathrm d\Psi,
	\qquad
	\xi^{(1)}=\xi,
	\qquad
	\nu^{(1)}=\nu.
	\]
	At the common endpoint of the row and column Christoffel chains,
	\[
		T_{(p+1)}=BA,
		\qquad
		\mathrm d\Psi^{(p+1)}=x\,\mathrm d\Psi.
	\]
	For every \(a\in\{1,\ldots,M\}\), these data satisfy
	\begin{equation}
		\label{eq:cyclic-Favard-transport-moments}
		(\xi^{(a)})^{-1}\mathsf E_{[q]}(T_{(a)})^n
		\mathsf E_{[p]}^{\top}(\nu^{(a)})^{-\top}
		=
		\int x^n\,\mathrm d\Psi^{(a)}(x),
		\qquad n\in\mathbb N_0.
	\end{equation}
	The right and left polynomial systems at that position are the Christoffel
	transforms determined by the transported initial-condition matrices
	\(\xi^{(a)}\) and \(\nu^{(a)}\).
	These conclusions hold for both bounded and unbounded factors.
\end{theorem}

\begin{proof}
	The weak normality hypothesis ensures that the two one-sided transformed
	systems and their transported initial conditions are defined at every
	step-line degree.  The one-step Christoffel connection
	formulas give the intertwining identities
	\[
	T_{(0,a-1)}L_a=L_aT_{(0,a)},
	\qquad
	U_bT_{(b-1,0)}=T_{(b,0)}U_b.
	\]
	Starting from \(T=L_1\cdots L_pU_q\cdots U_1\), induction therefore yields
	\begin{align*}
	T_{(0,a)}
	&=L_{a+1}\cdots L_pU_q\cdots U_1L_1\cdots L_a,
	&&0\le a\le p,\\
	T_{(b,0)}
	&=U_b\cdots U_1L_1\cdots L_pU_q\cdots U_{b+1},
	&&1\le b\le q.
	\end{align*}
	Since \(\mathrm d\Psi_{(0,p)}=\mathrm d\Psi_{(q,0)}
	=x\,\mathrm d\Psi\), these two lists join at their common endpoint.
	Comparison with
	the definition of \(T_{(a)}\) now gives
	\[
	T_{(a+1)}=T_{(0,a)}\quad(0\le a\le p),
	\qquad
	T_{(M+1-b)}=T_{(b,0)}\quad(1\le b\le q).
	\]
	Thus the cyclic products are exactly the recurrence matrices of the
	transformed systems listed above.  At each one-step transformation the
	Christoffel formula transports the recurrence solutions, their initial
	conditions, and the represented functional.  Starting from
	\((T,\xi,\nu,\mathrm d\Psi)\) and iterating around the two boundary chains
	therefore gives
	\eqref{eq:cyclic-Favard-transport-moments}.  Equivalently,
	Proposition~\ref{prop:recurrence-initial-data-moment-functional} identifies
	\(\mathrm d\Psi^{(a)}\) with the unique functional determined by
	\(T_{(a)},\xi^{(a)},\nu^{(a)}\).  The spectral representations are valid
	in both the bounded and unbounded settings under the hypotheses of
	\cite{BranquinhoFoulquieManas2023Spectral,
	BranquinhoFoulquieManas2026Unbounded}.
\end{proof}

The initial conditions throughout the rectangle are fixed by transport.  Starting
from those at the vertex \((0,0)\), the one-sided Christoffel formulas
determine the initial conditions at the boundary vertices and then at every
interior vertex.  Thus the polynomial arrays at different vertices are not
independently regauged.  In the statements below, a normalized connection
matrix always refers to this transported choice of initial conditions,
together with the bidiagonal normalization fixed in
Section~\ref{sec:positive-refactorization-sparse-matrices}.

For a recurrence matrix \(S\) with initial-condition matrices
\(\xi_S\) and \(\nu_S\), let
\(\mathscr R(S;\xi_S,\nu_S)\) denote the scalar moment matrix reconstructed
from the functional in
Proposition~\ref{prop:recurrence-initial-data-moment-functional}.  Explicitly,
for \(u,v\in\mathbb N_0\), \(j\in\{0,\ldots,q-1\}\), and
\(i\in\{0,\ldots,p-1\}\),
\begin{equation}
	\label{eq:reconstructed-scalar-moment-matrix}
	\mathscr R(S;\xi_S,\nu_S)_{qu+j,pv+i}
	=
	\left[
	\xi_S^{-1}\mathsf E_{[q]}S^{u+v}
	\mathsf E_{[p]}^{\top}\nu_S^{-\top}
	\right]_{j+1,i+1}.
\end{equation}

\begin{lemma}[Measures representing row and column moment shifts]
	\label{lem:balanced-Christoffel-shifted-moments}
	For \(b\in\{0,\ldots,q\}\) and
	\(a\in\{0,\ldots,p\}\), define the shifted scalar matrix
	\[
		\mathscr M^{\langle b,a\rangle}
		=
		\bigl[\mathscr M_{r+b,s+a}\bigr]_{r,s\ge0}.
	\]
	If \(\mathscr M\) is represented by a \(q\times p\) matrix of measures
	\(\mathrm d\Psi\), then
	\(\mathscr M^{\langle b,a\rangle}\) is represented by
	\begin{equation}
		\label{eq:shifted-moment-representing-measure}
		\mathsf K_q(x)^b\,\mathrm d\Psi(x)\,
		\bigl(\mathsf K_p(x)^a\bigr)^{\top}.
	\end{equation}
\end{lemma}

\begin{proof}
	Left multiplication by \(\mathsf K_q(x)\) advances the
	scalar step-line row index by one, including the factor \(x\) when the
	component index passes from \(q\) back to \(1\).  Likewise, right
	multiplication by \(\mathsf K_p(x)^{\top}\) advances the scalar column
	index by one.  Iterating these two identities shows, entry by entry, that
	the moments of \eqref{eq:shifted-moment-representing-measure} are
	\(\mathscr M_{r+b,s+a}\), as claimed.
\end{proof}

\subsubsection{The elementary square and construction of the rectangle}

\begin{lemma}[Permutability of one row and one column Christoffel step]
	\label{lem:one-step-functional-transport}
	Under the hypotheses and with the notation of
	Definition~\ref{def:one-step-Christoffel-connection-matrices}, let
	\begin{equation}
		\label{eq:one-square-bidiagonal-refactorization}
		UL=L'U'
	\end{equation}
	be its normalized positive lower--upper refactorization.  Then
	\[
		\mathrm d\Psi_{1,1}
		=\mathsf K_q\mathrm d\Psi_{0,0}\mathsf K_p^{\top}
	\]
	has a weakly normal mixed-type system.  Its polynomial arrays are
	\begin{equation}
		\label{eq:outgoing-Christoffel-polynomials}
		\boldsymbol B_{1,1}=(L')^{-1}\boldsymbol B_{1,0},
		\qquad
		\boldsymbol A_{1,1}=\boldsymbol A_{0,1}(U')^{-1},
	\end{equation}
	and the outgoing connection relations are
	\begin{equation}
		\label{eq:outgoing-Christoffel-connections}
		\begin{aligned}
			\boldsymbol A_{1,0}L'
			&=\mathsf K_p^{\top}\boldsymbol A_{1,1},
		&\qquad
			L'\boldsymbol B_{1,1}&=\boldsymbol B_{1,0},\\
			U'\boldsymbol B_{0,1}
			&=\boldsymbol B_{1,1}\mathsf K_q,
		&
			\boldsymbol A_{1,1}U'&=\boldsymbol A_{0,1}.
		\end{aligned}
	\end{equation}
\end{lemma}

\begin{proof}
	First we verify the degrees.  Since \(L'\) is normalized lower
	bidiagonal, every finite leading block of \((L')^{-1}\) is lower
	unitriangular.  Hence, for every \(n\),
	\[
		\boldsymbol B_{1,1,n}
		=\boldsymbol B_{1,0,n}
		+\sum_{k<n}c_{n,k}\boldsymbol B_{1,0,k}
	\]
	for suitable scalars \(c_{n,k}\).  The column-side shift from
	\((1,0)\) to \((1,1)\) does not change the prescribed right step-line
	polynomial space.  Since the degree bounds in that flag are nondecreasing
	with \(n\), the displayed triangular combination shows that
	\(\boldsymbol B_{1,1,n}\) belongs to the required space.  No
	maximal-degree conclusion is needed here.

	Likewise, \(U'\) is upper bidiagonal with a nonzero diagonal.  If its
	\(n\)-th diagonal entry is \(u'_n\), then the \(n\)-th column of the second
	formula in \eqref{eq:outgoing-Christoffel-polynomials} gives
	\[
		\boldsymbol A_{1,1,n}
		=\frac{1}{u'_n}\boldsymbol A_{0,1,n}
		+\sum_{k<n}d_{k,n}\boldsymbol A_{0,1,k}.
	\]
	The row-side shift from \((0,1)\) to \((1,1)\) does not change the
	prescribed left step-line polynomial space.  The same monotonicity of the
	degree bounds shows that \(\boldsymbol A_{1,1,n}\) belongs to that space.
	The triangular combination need not preserve a nonzero maximal-degree
	coefficient, and weak normality does not require it.

	The initial blocks remain nonsingular with the standard lower triangular
	orientation.  Indeed, if \(\xi_{1,0}\) and \(\nu_{0,1}\) denote the
	boundary initial-condition matrices and
	\[
		L'_{[q]}=\mathsf E_{[q]}L'\mathsf E_{[q]}^{\top},
		\qquad
		U'_{[p]}=\mathsf E_{[p]}U'\mathsf E_{[p]}^{\top},
	\]
	then triangularity gives
	\[
		\mathsf E_{[q]}\boldsymbol B_{1,1}
		=(L'_{[q]})^{-1}\xi_{1,0},
		\qquad
		\boldsymbol A_{1,1}\mathsf E_{[p]}^{\top}
		=\nu_{0,1}^{\top}(U'_{[p]})^{-1}.
	\]

	Two relations in \eqref{eq:outgoing-Christoffel-connections} follow
	immediately from \eqref{eq:outgoing-Christoffel-polynomials}.  For the
	other two, use the incoming relations:
	\[
		UL\boldsymbol B_{0,1}
		=U\boldsymbol B_{0,0}
		=\boldsymbol B_{1,0}\mathsf K_q.
	\]
	Substitution of \(UL=L'U'\), followed by multiplication on the left by
	\((L')^{-1}\), gives
	\(U'\boldsymbol B_{0,1}=\boldsymbol B_{1,1}\mathsf K_q\).  Similarly,
	\[
		\boldsymbol A_{1,0}UL
		=\boldsymbol A_{0,0}L
		=\mathsf K_p^{\top}\boldsymbol A_{0,1}.
	\]
	Substitution of \(UL=L'U'\), followed by multiplication on the right by
	\((U')^{-1}\), gives
	\(\boldsymbol A_{1,0}L'=\mathsf K_p^{\top}\boldsymbol A_{1,1}\).

	It remains to verify the orthogonality directly.  Put
	\[
		\mathcal I
		\coloneq
		\int
		\boldsymbol B_{1,1}(x)\,
		\mathrm d\Psi_{1,0}(x)\,
		\boldsymbol A_{1,0}(x).
	\]
	Since \(L'\boldsymbol B_{1,1}=\boldsymbol B_{1,0}\) and the system at
	\((1,0)\) is biorthonormal,
	\[
		L'\mathcal I
		=
		\int
		\boldsymbol B_{1,0}(x)\,
		\mathrm d\Psi_{1,0}(x)\,
		\boldsymbol A_{1,0}(x)
		=I.
	\]
	Thus \(\mathcal I=(L')^{-1}\), and
	\begin{align*}
		\int
		\boldsymbol B_{1,1}(x)\,
		\mathrm d\Psi_{1,1}(x)\,
		\boldsymbol A_{1,1}(x)
		&=
		\int
		\boldsymbol B_{1,1}(x)\,
		\mathrm d\Psi_{1,0}(x)\,
		\mathsf K_p(x)^{\top}\boldsymbol A_{1,1}(x)\\
		&=\mathcal I L'=I.
	\end{align*}
	The two arrays are therefore biorthonormal for
	\(\mathrm d\Psi_{1,1}\).  The degree statements proved at the beginning of
	the argument show that they belong to the prescribed step-line polynomial
	spaces.  Hence the system at \((1,1)\) is weakly normal.
\end{proof}

\begin{corollary}[The elementary \texorpdfstring{$\tau$}{tau}-determinants]
	\label{cor:elementary-grid-tau-determinants}
	In the setting of Lemma~\ref{lem:one-step-functional-transport}, the
	Christoffel \(\tau\)-determinants on the outgoing column and row sides are
	\[
		\tau_n^{\mathrm{col}}=A_{1,0,n}^{(1)}(0),
		\qquad
		\tau_n^{\mathrm{row}}=B_{0,1,n}^{(1)}(0),
	\]
	and are nonzero for every \(n\in\mathbb N_0\).  Equivalently, the first
	row of the upper Gauss--Borel coefficient matrix at \((1,0)\) and the
	first column of the lower Gauss--Borel coefficient matrix at \((0,1)\)
	have no zero entries.
	The superscripts \(\mathrm{col}\) and \(\mathrm{row}\) refer to the
	outgoing edges controlled by these determinants: the evaluation at
	\((1,0)\) controls the column step to \((1,1)\), whereas the evaluation at
	\((0,1)\) controls the row step to \((1,1)\).
\end{corollary}

\begin{proof}
	Both \(\mathsf K_p^{\top}\) and \(\mathsf K_q\) have the single simple
	determinantal zero \(x=0\).  The relevant left null vector of
	\(\mathsf K_p(0)^{\top}\), and the relevant right null vector of
	\(\mathsf K_q(0)\), are in both cases the first coordinate vector.  Hence the general
	determinants of \cite{ManasRojas2026Christoffel} reduce to the two scalar
	evaluations displayed above.

	Write \(L'_{n+1,n}=\ell'_n>0\), \(U'_{n,n}=u'_n>0\), and
	\(U'_{n,n+1}=1\).  Evaluating the first and third relations in
	\eqref{eq:outgoing-Christoffel-connections} at \(x=0\), and selecting the
	first polynomial component, gives
	\[
		A_{1,0,n}^{(1)}(0)
		+\ell'_nA_{1,0,n+1}^{(1)}(0)=0,
		\qquad
		u'_nB_{0,1,n}^{(1)}(0)
		+B_{0,1,n+1}^{(1)}(0)=0.
	\]
	The initial evaluations are nonzero because the initial-condition matrices
	are triangular and nonsingular.  Iterating the two preceding identities
	gives, for every \(n\ge1\),
	\begin{equation}
		\label{eq:elementary-grid-tau-products}
		A_{1,0,n}^{(1)}(0)
		=
		\frac{(-1)^nA_{1,0,0}^{(1)}(0)}
		{\displaystyle\prod_{k=0}^{n-1}\ell'_k},
		\qquad
		B_{0,1,n}^{(1)}(0)
		=
		(-1)^n
		\left(\prod_{k=0}^{n-1}u'_k\right)
		B_{0,1,0}^{(1)}(0).
	\end{equation}
	Thus the two sequences alternate in sign and have no zero entries.  Notice
	that this conclusion uses only the connection identities, the nonzero
	initial evaluations, and the strict positivity of the bidiagonal
	coefficients; no total-positivity assumption on the moment matrix is used.
	The same calculation applies to an edge on either boundary whenever its
	prescribed factor is already known to be the normalized Christoffel
	connection matrix.  Finally,
	\(\boldsymbol B(0)=\mathscr L\mathsf E_{[q]}^{\top}\) and
	\(\boldsymbol A(0)=\mathsf E_{[p]}\mathscr U\), so these evaluations are
	the entries of the stated column and row.
\end{proof}

\begin{proposition}[Positive bidiagonal factors on the moment-shift rectangle]
	\label{prop:two-sided-grid-edge-factors}
	There is a unique collection of matrices
	\(\mathcal U_{b,a}\) and \(\mathcal L_{b,a}\), indexed by the edges of
	the moment-shift rectangle, that satisfies the boundary conditions
	\begin{equation}
		\label{eq:two-sided-grid-boundary-edge-factors}
		\mathcal U_{b,0}=U_b,
		\qquad b\in\{1,\ldots,q\},
		\qquad
		\mathcal L_{0,a}=L_{a+1},
		\qquad a\in\{0,\ldots,p-1\}.
	\end{equation}
	At every elementary square these matrices satisfy
	\begin{equation}
		\label{eq:two-sided-grid-local-square}
		\mathcal U_{b,a}\mathcal L_{b-1,a}
		=
		\mathcal L_{b,a}\mathcal U_{b,a+1},
		\qquad
		b\in\{1,\ldots,q\},\quad
		a\in\{0,\ldots,p-1\}.
	\end{equation}
	Every \(\mathcal U_{b,a}\) is normalized upper bidiagonal and every
	\(\mathcal L_{b,a}\) is normalized lower bidiagonal.  All their
	nontrivial entries are strictly positive.
\end{proposition}

\begin{proof}
	Suppose that row \(b-1\) has already been constructed.  At the beginning
	of row \(b\), the factor \(\mathcal U_{b,0}=U_b\) is known.  For
	\(a=0\), the left-hand side of
	\eqref{eq:two-sided-grid-local-square} is therefore a known normalized
	positive upper--lower product.  Applying
	Lemma~\ref{lem:positive-adjacent-refactorization} determines a unique
	normalized positive lower--upper pair
	\(\mathcal L_{b,0}\mathcal U_{b,1}\).  The newly obtained upper factor
	is then used at \(a=1\), and the argument continues until the end of the
	row.  Repeating this procedure for
	\(b=1,\ldots,q\) constructs the entire rectangle.  The same lemma proves
	at every step the asserted normalization, positivity, and uniqueness.
\end{proof}

\subsubsection{Recurrences and polynomial systems at the vertices}

\begin{theorem}[Moment systems on the two-sided Christoffel rectangle]
	\label{thm:balanced-Christoffel-grid-path}
	Suppose that the prescribed PBF is the Christoffel factorization of a
	weakly normal system represented by \(\mathrm d\Psi\), with recurrence
	matrix \(T\) and initial-condition matrices \(\xi,\nu\).  Suppose also
	that the two one-sided Christoffel chains along the lower and left
	boundaries of the rectangle are weakly normal, or equivalently that their
	Christoffel \(\tau\)-determinants do not vanish.

	For every vertex \((b,a)\), the matrices
	\(\mathcal L_{b,a}\) and \(\mathcal U_{b,a}\) are the normalized
	column-side and row-side connection matrices on the corresponding edges.
	The right and left recurrence polynomials and their initial-condition
	matrices obtained by transporting \((T,\xi,\nu)\) along a path from
	\((0,0)\) to \((b,a)\) are independent of that path.  Denote them by
	\(\boldsymbol B_{(b,a)},\boldsymbol A_{(b,a)}\),
	\(\xi_{(b,a)},\nu_{(b,a)}\), and denote their recurrence matrix by
	\(T_{(b,a)}\).  Their reconstructed moment matrix is
	\begin{equation}
		\label{eq:two-sided-grid-reconstructed-moment-matrix}
		\mathscr R\bigl(T_{(b,a)};
		\xi_{(b,a)},\nu_{(b,a)}\bigr)
		=
		\mathscr M^{\langle b,a\rangle}
		=
		\bigl[\mathscr M_{r+b,s+a}\bigr]_{r,s\ge0}.
	\end{equation}
	The matrix of measures
	\begin{equation}
		\label{eq:two-sided-grid-measure-identification}
		\mathrm d\Psi_{(b,a)}(x)
		\coloneq
		\mathsf K_q(x)^b\,\mathrm d\Psi(x)\,
		\bigl(\mathsf K_p(x)^a\bigr)^{\top}
	\end{equation}
	represents the reconstructed moment matrix.  In particular, for every
	\(n\in\mathbb N_0\),
	\begin{equation}
		\label{eq:two-sided-grid-moment-identification}
		\xi_{(b,a)}^{-1}\mathsf E_{[q]}T_{(b,a)}^n
		\mathsf E_{[p]}^{\top}\nu_{(b,a)}^{-\top}
		=
		\int x^n\,\mathrm d\Psi_{(b,a)}(x).
	\end{equation}
	The right and left polynomial solutions of the recurrence
	for \(T_{(b,a)}\), normalized by \(\xi_{(b,a)}\) and
	\(\nu_{(b,a)}\), are the mixed-type multiple orthogonal polynomials for
	\(\mathrm d\Psi_{(b,a)}\).
\end{theorem}

\begin{proof}
	We first identify the systems at the vertices and the factors on the
	interior edges.  The systems on the two boundary chains and their
	connection factors are known by hypothesis.  We fill the rectangle one
	elementary square at a time, from left to right in each row and then from
	the bottom row to the top.  Suppose that the systems at the lower-left,
	lower-right, and upper-left vertices of a square have already been
	constructed.  In the notation of
	Lemma~\ref{lem:one-step-functional-transport}, their measures are
	\[
		\mathrm d\Psi_{0,0},
		\qquad
		\mathrm d\Psi_{0,1}
		=
		\mathrm d\Psi_{0,0}\mathsf K_p^{\top},
		\qquad
		\mathrm d\Psi_{1,0}
		=
		\mathsf K_q\mathrm d\Psi_{0,0}.
	\]
	The two known incoming edge factors form the upper--lower product on the
	left-hand side of \eqref{eq:two-sided-grid-local-square}.  Its unique
	positive lower--upper refactorization is the pair of outgoing edge factors
	constructed in
	Proposition~\ref{prop:two-sided-grid-edge-factors}.  Lemma~\ref{lem:one-step-functional-transport}
	then constructs the two polynomial arrays at the fourth vertex, proves
	their biorthogonality for
	\[
		\mathsf K_q\mathrm d\Psi_{0,0}\mathsf K_p^{\top},
	\]
	and proves all four connection identities on the square.  This also proves
	weak normality at the fourth vertex.  Induction therefore gives, at every
	vertex,
	\begin{equation}
		\label{eq:two-sided-grid-inductive-measure}
		\mathrm d\Psi_{(b,a)}(x)
		=
		\mathsf K_q(x)^b\,\mathrm d\Psi(x)\,
		\bigl(\mathsf K_p(x)^a\bigr)^{\top},
	\end{equation}
	and identifies every factor in
	Proposition~\ref{prop:two-sided-grid-edge-factors} as the corresponding
	normalized Christoffel connection matrix.  Corollary~\ref{cor:elementary-grid-tau-determinants}
	shows at the same time that all elementary \(\tau\)-determinants on the
	new edges are nonzero.

	Any two monotone paths from \((0,0)\) to a fixed vertex \((b,a)\) are
	connected by elementary-square moves.  The four connection identities in
	Lemma~\ref{lem:one-step-functional-transport} show that replacing either
	two-edge route across one square by the other produces the same normalized
	polynomial arrays and the same transported initial conditions.  Hence
	\(\boldsymbol B_{(b,a)},\boldsymbol A_{(b,a)},\xi_{(b,a)}\), and
	\(\nu_{(b,a)}\) are independent of the chosen path.  Their recurrence
	equations determine \(T_{(b,a)}\) uniquely.

	It remains to identify the scalar moment matrix reconstructed from these
	recurrence data.  By
	Lemma~\ref{lem:balanced-Christoffel-shifted-moments}, the moments of
	\eqref{eq:two-sided-grid-inductive-measure} form precisely
	\[
		\mathscr M^{\langle b,a\rangle}
		=
		\bigl[\mathscr M_{r+b,s+a}\bigr]_{r,s\ge0}.
	\]
	On the other hand,
	Proposition~\ref{prop:recurrence-initial-data-moment-functional}, applied
	to the normalized endpoint system, reconstructs those same moments from
	\(T_{(b,a)},\xi_{(b,a)},\nu_{(b,a)}\).  This proves
	\eqref{eq:two-sided-grid-reconstructed-moment-matrix}; inserting the
	definition \eqref{eq:reconstructed-scalar-moment-matrix} gives
	\eqref{eq:two-sided-grid-moment-identification}.  Hence
	\eqref{eq:two-sided-grid-measure-identification} is a representing matrix
	of measures for the reconstructed functional, and the endpoint arrays are
	its two mixed-type multiple orthogonal polynomial families.

\end{proof}

\begin{corollary}[Path factorizations of the Darboux product]
	\label{cor:two-sided-grid-path-factorizations}
	For a monotone path \(\gamma\) from \((q,0)\) to \((0,p)\), let
	\(\mathcal F(\gamma)\) be the ordered product obtained by assigning
	\(\mathcal U_{b,a}\) to each step
	\((b,a)\mapsto(b-1,a)\) and \(\mathcal L_{b,a}\) to every step
	\((b,a)\mapsto(b,a+1)\).  Then
	\[
		\mathcal F(\gamma)
		=
		BA
		=
		U_q\cdots U_1L_1\cdots L_p,
	\]
	so \(\mathcal F(\gamma)\) is independent of \(\gamma\).  If \(\gamma\)
	passes through
	\((b,a)\), reading its factors cyclically from that vertex gives the
	recurrence matrix \(T_{(b,a)}\) of
	Theorem~\ref{thm:balanced-Christoffel-grid-path}; in particular, the cut
	product is also independent of \(\gamma\).  The
	\(\binom{p+q}{p}\) monotone paths therefore give positive bidiagonal
	factorizations of \(BA\) with every interlacing of the lower and upper
	factor types that preserves their relative orders.
\end{corollary}

\begin{proof}
	For the boundary path consisting first of the \(q\) upward steps and then
	of the \(p\) rightward steps, the edge factors are
	\(U_q,\ldots,U_1,L_1,\ldots,L_p\), and their product is \(BA\).  Any two
	monotone paths are connected by elementary-square moves.  Each move replaces
	\(\mathcal U_{b,a}\mathcal L_{b-1,a}\) by
	\(\mathcal L_{b,a}\mathcal U_{b,a+1}\), or conversely, and therefore leaves
	the product unchanged by \eqref{eq:two-sided-grid-local-square}.

	Now cut a path through \((b,a)\) at that vertex and denote its cyclic
	product by \(T_{(b,a)}^\gamma\).  Composing the right connection identities
	around the path makes the horizontal connections telescope, while the
	\(q\) vertical connections contribute \(\mathsf K_q(x)^q=xI_q\).  Hence
	\[
		T_{(b,a)}^\gamma\boldsymbol B_{(b,a)}(x)
		=x\boldsymbol B_{(b,a)}(x).
	\]
	Composing the left connection identities in the opposite order gives
	\[
		\boldsymbol A_{(b,a)}(x)T_{(b,a)}^\gamma
		=x\boldsymbol A_{(b,a)}(x),
	\]
	because \(\mathsf K_p(x)^p=xI_p\).  Thus the cut product is the recurrence
	matrix of the normalized weakly normal system at \((b,a)\).  The uniqueness
	of that recurrence matrix gives \(T_{(b,a)}^\gamma=T_{(b,a)}\), independently
	of \(\gamma\).
\end{proof}

\begin{corollary}[Minimum-height factorizations represented in the rectangle]
	\label{cor:two-sided-grid-star-refactorizations}
	Fix \(a\in\mathcal A_{\mathrm{cyc}}=\{a_-,\ldots,a_+\}\).  The positive
	refactorization used to construct \(\mathscr H_{(a)}\) is obtained by
	cutting one of the paths in
	Corollary~\ref{cor:two-sided-grid-path-factorizations} at the vertex
	immediately preceding the first factor of the cyclic product \(T_{(a)}\).
	The cyclic products of the resulting factorization are the recurrence
	matrices \(T_{(b,a)}\) of the simultaneously shifted systems encountered at
	the successive vertices of that path.
\end{corollary}

\begin{proof}
	Start with the boundary path and mark the vertex immediately preceding the
	first factor of \(T_{(a)}\).  Perform the adjacent exchanges used in
	Theorem~\ref{thm:positive-balanced-refactorization}, without moving a factor
	across the marked cut.  Each exchange replaces the two edges of one route
	through an elementary square by the other two edges and is therefore exactly
	\eqref{eq:two-sided-grid-local-square}.  Uniqueness in
	Lemma~\ref{lem:positive-adjacent-refactorization} identifies the resulting
	edge factors with those used to construct \(\mathscr H_{(a)}\).  The claim
	about their cyclic products now follows from
	Corollary~\ref{cor:two-sided-grid-path-factorizations}.
\end{proof}

\subsubsection{Total positivity and the Pi\~neiro and Jacobi-like systems}

\begin{corollary}[The totally positive moment case]
	\label{cor:TP-moment-shift-rectangle}
	In the setting of
	Theorem~\ref{thm:balanced-Christoffel-grid-path}, assume that the scalar
	step-line moment matrix \(\mathscr M\) is totally positive in the
	TN--TP terminology fixed above.  Then, for every
	\(b\in\{0,\ldots,q\}\) and \(a\in\{0,\ldots,p\}\), the shifted matrix
	\[
		\mathscr M^{\langle b,a\rangle}
		=
		\bigl[\mathscr M_{r+b,s+a}\bigr]_{r,s\ge0}
	\]
	is totally positive.  Consequently, the two boundary chains and every
	interior vertex are weakly normal, all their Christoffel
	\(\tau\)-determinants are nonzero, and every conclusion of
	Theorem~\ref{thm:balanced-Christoffel-grid-path} holds without separately
	assuming weak normality of the boundary chains.
\end{corollary}

\begin{proof}
	Every finite minor of
	\(\mathscr M^{\langle b,a\rangle}\) is a finite minor of
	\(\mathscr M\), obtained by translating its row indices by \(b\) and its
	column indices by \(a\).  It is therefore strictly positive.  In
	particular, all leading principal minors of every shifted matrix are
	positive, which is equivalent to weak normality along its step-line.
	The one-step existence criterion of
	\cite{ManasRojas2026Christoffel} identifies this weak normality with the
	nonvanishing of the corresponding Christoffel \(\tau\)-determinants.
	The hypotheses of
	Theorem~\ref{thm:balanced-Christoffel-grid-path} are therefore satisfied.
\end{proof}

\begin{remark}[Explicit formulas on the boundary chains]
	The preceding theorem identifies the moment functional, recurrence matrix,
	and two mixed-type polynomial families at every vertex of the rectangle;
	in particular, it establishes the full bilateral Christoffel
	transformation.  Closed determinant formulas expressing the transformed
	polynomials in terms of the original ones are presently available
	for a polynomial perturbation on one side of the rectangular matrix of
	measures; in the present notation they cover the boundary chains
	\(\mathrm d\Psi_{(b,0)}\), with \(b\in\{0,\ldots,q\}\), and
	\(\mathrm d\Psi_{(0,a)}\), with \(a\in\{0,\ldots,p\}\); see
	\cite{ManasRojas2026Christoffel}.  At an interior vertex the polynomials
	are obtained here by transporting the recurrence and its initial
	conditions through the bidiagonal relations along a path.
	Finally, the moment identity determines the polynomial functional; in an
	indeterminate moment problem, several matrices of measures may represent
	the same functional.
\end{remark}

\begin{remark}[Direction of the transformations along a path]
	Each vertex \(\mathrm d\Psi_{(b,a)}\) is defined by applying the displayed
	row and column Christoffel multipliers to \(\mathrm d\Psi\).  A path from
	\((q,0)\) to \((0,p)\), however, traverses these transformations in two
	different directions.  A step \((b,a)\mapsto(b,a+1)\) is an elementary
	column-side Christoffel transformation, whereas
	\((b,a)\mapsto(b-1,a)\) reverses an elementary row-side Christoffel
	transformation and is therefore a Geronimus step.  The functionals and the
	initial conditions at both endpoints have already been fixed by
	\(\mathrm d\Psi_{(b,a)}\).  The uniqueness of the normalized recurrence
	solutions therefore determines the Geronimus parameter and the connecting
	factor.
\end{remark}

\begin{corollary}[The Pi\~neiro moment-shift rectangle]
	\label{cor:balanced-Pineiro-complete-grid}
	For the Pi\~neiro matrix
	\[
		\mathrm d\Psi_{j,i}(x)=x^{\alpha_i+\beta_j}\,\mathrm dx,
		\qquad 0<x<1,
	\]
	assume
	\[
		\alpha_1<\cdots<\alpha_p<\alpha_1+1,
		\qquad
		\beta_1<\cdots<\beta_q<\beta_1+1,
		\qquad
		\alpha_1+\beta_1+1>0.
	\]
	Then the scalar moment matrix is totally positive.  Every shifted
	system in the rectangle is therefore weakly normal, and the cyclic product
	obtained by cutting any path at a vertex has precisely the corresponding
	Pi\~neiro measure and its transported initial-condition matrices.  Moreover, the row and column
	power systems at every vertex are AT systems; hence all these transformed
	systems are strongly normal.
\end{corollary}

\begin{proof}
	With
	\[
		x_{qu+j}=u+\beta_{j+1},
		\qquad
		y_{pv+i}=v+\alpha_{i+1}+1,
	\]
	the moment matrix is the Cauchy matrix
	\(\mathscr M_{r,s}=1/(x_r+y_s)\).  The parameter inequalities make both
	sequences strictly increasing and every sum positive, so every minor is a
	positive Cauchy determinant.  Weak normality and the identification of the
	shifted measures follow from
	Corollary~\ref{cor:TP-moment-shift-rectangle}, and the cut-product
	identification follows from
	Corollary~\ref{cor:two-sided-grid-path-factorizations}.  The left and right
	Christoffel multiplications replace the exponent lists by the affine
	cyclic shifts described in
	\cite{BranquinhoFoulquieManas2026Pineiro} and used independently in
	\cite[Section~7 and Appendix~A]{ManasMarkov2026}.  The displayed strict
	inequalities are preserved by those shifts.  Hence the two transformed
	power vectors remain AT systems.  The mixed AT theorem then gives the
	maximal component degrees and strong normality.
\end{proof}

\begin{corollary}[The Jacobi-like moment-shift rectangle]
	\label{cor:balanced-Jacobi-complete-grid}
	In the positive factorization chamber of
	\cite[Theorem~9.4]{ManasMarkov2026}, the scalar Jacobi-like moment matrix
	is totally positive.  Hence every simultaneous row and column
	moment shift in the rectangle is weakly normal, and the cyclic product
	obtained by cutting any path at a vertex is the recurrence matrix of the
	corresponding shifted Jacobi-like measure.
	A transformed system is strongly normal whenever its two transformed
	weight vectors satisfy the AT hypotheses of
	\cite{Manas2026Hypergeometric}.  Thus the AT range gives strong normality,
	while the larger positive-factorization chamber gives the PBF and weak
	normality conclusions.
\end{corollary}

\begin{proof}
	The proof of the cited theorem factors the scalar moment matrix as
	\(\mathscr M^{\mathrm J}=\mathcal C\mathcal G\), where \(\mathcal C\)
	is block diagonal with nonsingular totally nonnegative finite blocks, while
	\(\mathcal G\) is totally positive in the TN--TP terminology fixed above.
	Fix finite row and column
	sets \(I,J\) of the same cardinality.  The rows of \(\mathcal C\) indexed
	by \(I\) are linearly independent, so some column set \(K\) satisfies
	\(\det\mathcal C_{I,K}>0\).  Since
	\(\det\mathcal G_{K,J}>0\), this is a strictly positive term in the
	Cauchy--Binet expansion of
	\(\det\mathscr M^{\mathrm J}_{I,J}\).  Only finitely many column sets
	occur because \(\mathcal C\) is block diagonal with finite blocks, and
	every other term is nonnegative.
	Thus all minors of \(\mathscr M^{\mathrm J}\) are positive.  The
	moment-shift and recurrence identifications follow from
	Corollaries~\ref{cor:TP-moment-shift-rectangle}
	and~\ref{cor:two-sided-grid-path-factorizations}.
\end{proof}

\subsection{Weyl approximation on the line and on the star}
\label{subsec:Weyl-approximation-line-star}

The preceding results identify both the radial Weyl matrix
\(\mathsf S(z)\) and its star lift \(\widehat{\mathsf S}(\lambda)\) in the
bounded case.
Choose \(N_0\) so that all states selected by \(\mathsf Q_q\) and
\(\mathsf Q_p\) belong to \(\mathscr V_{[0,N_0]}\).  For \(N\ge N_0\), define
\begin{equation}
	\label{eq:finite-CF-Weyl-convergent}
	\mathsf S^{[N]}(z)
	\coloneq
	\frac1z\,\xi^{-1}
	\mathsf Q_q
	\left(I-\tau\mathscr H^{[0,N]}\right)^{-1}
	\mathsf Q_p^{\top}\nu^{-\top},
	\qquad \tau^M=z^{-1}.
\end{equation}
Here \(\mathscr H^{[0,N]}\) is the finite principal submatrix from
Definition~\ref{def:principal-compression}, obtained by retaining the height
blocks \(0,\ldots,N\).  Thus the superscript \([N]\) on
\(\mathsf S^{[N]}\) records the terminal height block.  This definition uses
only a finite matrix and is therefore meaningful even if \(\mathscr H\) is
unbounded.  Under the hypotheses of Theorem~\ref{thm:general-Weyl-CF}, the
convergence results of Subsection~\ref{subsec:finite-convergents} apply to
these rectangular approximants.

The corresponding approximation on the star is obtained without choosing a
root:
\begin{equation}
	\widehat{\mathsf S}^{[N]}(\lambda)
	\coloneq
	\frac1\lambda\,\xi^{-1}\mathsf Q_q
	\left(I-\lambda^{-1}\mathscr H^{[0,N]}\right)^{-1}
	\mathsf Q_p^{\top}\nu^{-\top}
	=
	\lambda^{M-1}\mathsf S^{[N]}(\lambda^M).
\end{equation}
Thus the same backward Schur recursion computes rational approximants either
radially, in the variable \(z=\lambda^M\), or directly on the star, in the
variable \(\lambda\).  In the latter form their poles are eigenvalues of the
finite two-diagonal height matrix itself, before taking \(M\)-th powers.

In the next corollary, \(\|\cdot\|_2\) denotes the subordinate norm between
the indicated Euclidean coordinate spaces, also for rectangular matrices.

\begin{corollary}[Convergence to the rectangular Weyl matrix]
	\label{cor:finite-CF-Weyl-convergence}
	Assume the hypotheses of Theorem~\ref{thm:general-Weyl-CF}.  With the
	notation above, \(\mathsf S^{[N]}(z)\) is independent of the
	chosen \(M\)-th root \(\tau\).
	If
	\[
		\theta_z
		\coloneq
		\|\mathscr H\|_2\,|z|^{-1/M}
		<1,
	\]
	then
	\begin{equation}
		\label{eq:finite-CF-Weyl-error}
		\left\|
		\mathsf S(z)-\mathsf S^{[N]}(z)
		\right\|_2
		\le
		\frac{2}{|z|}
		\|\xi^{-1}\|_2\,\|\nu^{-\top}\|_2\,
		\frac{\theta_z^{2(N-N_0+1)}}{1-\theta_z}.
	\end{equation}
	Consequently, the rectangular approximants
	\(\mathsf S^{[N]}(z)\) converge locally uniformly to
	\(\mathsf S(z)\) in
	\(|z|>\|\mathscr H\|_2^M\).

	If \(\mathscr H\) is nonnegative, set
	\(\varrho\coloneq\rho_{\ell^2}(\mathscr H)\).  For every
	\(\gamma>\varrho\) satisfying \(\gamma|z|^{-1/M}<1\), define
	\[
		C_{\gamma,2}
		\coloneq
		\sup_{k\ge0}\frac{\|\mathscr H^k\|_2}{\gamma^k}<\infty.
	\]
	Then
	\begin{equation}
		\label{eq:finite-CF-Weyl-spectral-radius-error}
		\left\|
		\mathsf S(z)-\mathsf S^{[N]}(z)
		\right\|_2
		\le
		\frac{C_{\gamma,2}}{|z|}
		\|\xi^{-1}\|_2\,\|\nu^{-\top}\|_2\,
		\frac{
		\left(\gamma|z|^{-1/M}\right)^{2(N-N_0+1)}
		}{
		1-\gamma|z|^{-1/M}
		}.
	\end{equation}
	Consequently, convergence holds throughout
	\(|z|>\varrho^M\).
\end{corollary}

\begin{proof}
	The factor position advances by one at every multiplication by
	\(\mathscr H\), and the principal height compression only deletes
	edges; it does not alter their factor-position labels.  Hence every surviving
	path
	from the first position back to the first position has length divisible by
	\(M\).  The rectangular block in
	\eqref{eq:finite-CF-Weyl-convergent} therefore depends on \(\tau\) only
	through \(\tau^M=z^{-1}\).

	A product contributing to the difference between the infinite and finite
	rectangular resolvents must leave \(\mathscr V_{[0,N]}\).  Both selected
	endpoints lie in blocks with index at most \(N_0\).  Such a product must
	make at least \(N-N_0+1\) outward block transitions and the same number
	of inward transitions.  Thus the coefficients of the two rectangular
	resolvents agree below degree \(2(N-N_0+1)\).  Repeating the Neumann-series
	estimate in the proof of Theorem~\ref{thm:finite-CF-convergence} yields
	\[
		\left\|
		\mathsf Q_q
		\left[
		(I-\tau\mathscr H)^{-1}
		-
		(I-\tau\mathscr H^{[0,N]})^{-1}
		\right]
		\mathsf Q_p^{\top}
		\right\|_2
		\le
		\frac{2\theta_z^{2(N-N_0+1)}}{1-\theta_z}.
	\]
	Multiplication by \(z^{-1}\xi^{-1}\) and \(\nu^{-\top}\) gives
	\eqref{eq:finite-CF-Weyl-error}.

	For nonnegative \(\mathscr H\), replace the norm estimate by the
	pathwise spectral-radius estimate
	\eqref{eq:positive-CF-spectral-radius-error}.  The same endpoint argument
	starts the omitted paths at degree \(2(N-N_0+1)\), and gives
	\eqref{eq:finite-CF-Weyl-spectral-radius-error}.
\end{proof}

Besides analytic convergence, the approximant \(\mathsf S^{[N]}\) reproduces
an exact range of moment coefficients by the path argument in
Proposition~\ref{prop:finite-CF-leading-error}.  Expressed in
the variable \(z\), this becomes a rectangular Pad\'e-type approximation at
infinity.  Put
\[
\operatorname{ord}_N
\coloneq
\left\lceil
\frac{2(N-N_0+1)}{M}
\right\rceil.
\]

	\begin{corollary}[Rectangular Pad\'e-type approximation at infinity]
		\label{cor:finite-CF-Weyl-Pade}
		Assume the hypotheses of
		Theorem~\ref{thm:general-Weyl-CF}.  Then the matrix
		\(\mathsf S^{[N]}(z)\) is a rational \(q\times p\) matrix, independent
	of the chosen \(M\)-th root of \(z^{-1}\), and
	\begin{equation}
		\label{eq:finite-CF-Weyl-Pade-order}
		\mathsf S(z)-\mathsf S^{[N]}(z)
		=
		\mathrm{O}\left(z^{-\operatorname{ord}_N-1}\right),
		\qquad z\longrightarrow\infty.
	\end{equation}
	Thus \(\mathsf S^{[N]}\) is a rectangular Pad\'e-type approximant of
	moment order \(\operatorname{ord}_N\): it reproduces
	\(\int x^n\,\mathrm d\Psi(x)\) for
	\(0\le n<\operatorname{ord}_N\).
	Its poles are contained among the \(M\)-th powers of the nonzero
	eigenvalues of \(\mathscr H^{[0,N]}\), after cancellation.
\end{corollary}

\begin{proof}
	The proof of Corollary~\ref{cor:finite-CF-Weyl-convergence} shows that the
	projected resolvents agree in every power \(\tau^k\) with
	\(k<2(N-N_0+1)\).  Only powers divisible by \(M\) survive the projection
	onto one copy of the cyclic linearization.
	Consequently, the coefficients of
	\(\tau^{Mn}=z^{-n}\) agree for
	\[
	Mn<2(N-N_0+1),
	\]
	which is equivalent to \(0\le n<\operatorname{ord}_N\).  The additional factor
	\(z^{-1}\) in \eqref{eq:finite-CF-Weyl-convergent} gives
	\eqref{eq:finite-CF-Weyl-Pade-order}.

	Proposition~\ref{prop:finite-CF-determinants} places the poles of the
	finite resolvent at reciprocals of the nonzero eigenvalues of the height
	compression in the variable \(\tau\).  Since the projected rational matrix
	is independent of the chosen \(M\)-th root \(\tau\) of \(z^{-1}\), these
	poles occur in \(M\)-cycles and map to the \(M\)-th powers of
	the corresponding eigenvalues in the variable \(z=\tau^{-M}\).
\end{proof}

\subsection{Unbounded positive bidiagonal factorizations}
\label{subsec:unbounded-PBF-Weyl}

The preceding norm estimates concern bounded realizations.  In the
nonnegative PBF setting, Lemma~\ref{lem:height-operator-norm-bound} identifies
this condition with boundedness of the original matrix \(T\).  For unbounded
\(T\), the formal continued fraction, every finite convergent, and the moment
identities remain well defined.  Analytic finite-section convergence requires
a specified closed realization with the stated core and stability properties,
and its target is the projected resolvent of that realization.  The spectral
Weyl matrix is defined separately from a chosen Favard measure; equality of
these two analytic objects is an additional identification.

\begin{definition}[Spectral Weyl matrix]
	\label{def:unbounded-spectral-Weyl}
	Let \(\mathrm d\Psi\) be a \(q\times p\) matrix of entrywise nonnegative
	finite measures
	supported on \(\mathbb R\), and assume that all its moments are finite.
	Its spectral Weyl matrix is
	\begin{equation*}
	\mathsf S_{\Psi}(z)
	\coloneq
	\int_{\mathbb R}\frac{\mathrm d\Psi(x)}{z-x},
	\qquad z\in\mathbb C\setminus\mathbb R,
	\end{equation*}
	where the integral is taken entry by entry.  If \(\mathrm d\Psi\) is
	compactly supported, satisfies
	\eqref{eq:general-spectral-moments}, and \(\mathscr H\) is bounded, then
	Proposition~\ref{prop:general-Weyl-Stieltjes} gives
	\(\mathsf S_{\Psi}=\mathsf S\) on their common domain.
	In an indeterminate moment problem, \(\mathsf S_{\Psi}\) depends on the
	chosen representing measure; the moment identity alone does not identify
	it with a projected operator resolvent.
\end{definition}

Let \(\mathcal H_0\) be a closed, densely defined realization of
\(\mathscr H\) on \(X_0=\ell^2(\mathscr V_{\ge0})\), assume that
\(c_{00}(\mathscr V_{\ge0})\) is a core, and let
\(K\Subset\Omega_0\setminus\{0\}\) be a compact set on which the finite
sections satisfy \eqref{eq:finite-CF-uniform-stability}.  Define
\begin{equation}
	\label{eq:unbounded-projected-resolvent}
	\mathsf R_{\mathcal H_0}(\tau)
	\coloneq
	\tau^M\xi^{-1}\mathsf Q_q
	\left(I-\tau\mathcal H_0\right)^{-1}
	\mathsf Q_p^{\top}\nu^{-\top}.
\end{equation}
The finite-window assertion
\eqref{eq:finite-CF-stable-window-convergence}, with \(m=0\), gives
\begin{equation}
	\label{eq:unbounded-rectangular-stability}
	\sup_{\tau\in K}
	\left\|
	\mathsf S^{[N]}(\tau^{-M})-
	\mathsf R_{\mathcal H_0}(\tau)
	\right\|_2
	\longrightarrow0.
\end{equation}
If
\[
	K^{-M}\coloneq\{\tau^{-M}:\tau\in K\}
	\subset\mathbb C\setminus\mathbb R
\]
and the chosen Favard
measure additionally satisfies
\begin{equation}
	\label{eq:unbounded-resolvent-Cauchy-identification}
	\mathsf S_{\Psi}(\tau^{-M})
	=
	\mathsf R_{\mathcal H_0}(\tau),
	\qquad \tau\in K,
\end{equation}
then the same convergence is to the measure-defined spectral Weyl matrix.

For measures with unbounded support, approximation at infinity is formulated
in the Poincar\'e sense.  The following elementary remainder formula gives
the corresponding Pad\'e-type statement.

\begin{proposition}[Moment expansion for an unbounded spectral measure]
	\label{prop:unbounded-Weyl-Poincare}
	Let \(\mathrm d\Psi\) satisfy
	Definition~\ref{def:unbounded-spectral-Weyl}, and put
	\(\boldsymbol\mu_n\coloneq\int x^n\,\mathrm d\Psi(x)\).  For every
	integer \(r\ge1\),
	\begin{equation}
	\label{eq:unbounded-Weyl-Poincare}
	\mathsf S_{\Psi}(z)
	=
	\sum_{n=0}^{r-1}\frac{\boldsymbol\mu_n}{z^{n+1}}
	+
	\mathrm{O}\left(|z|^{-r-1}\right),
	\qquad z\longrightarrow\infty,
	\end{equation}
	entry by entry, uniformly in every closed region satisfying
	\begin{equation*}
	\operatorname{dist}(z,\mathbb R)\ge\delta|z|
	\qquad(\delta>0).
	\end{equation*}
	Consequently, if the moment identity in
	\eqref{eq:general-spectral-moments} holds without the compact-support
	assumption, then, for every fixed \(N\), the finite height approximant
	\(\mathsf S^{[N]}\) in
	\eqref{eq:finite-CF-Weyl-convergent} still satisfies
	\begin{equation*}
	\mathsf S_{\Psi}(z)-\mathsf S^{[N]}(z)
	=
	\mathrm{O}
	\left(|z|^{-\operatorname{ord}_N-1}\right),
	\qquad z\longrightarrow\infty,
	\end{equation*}
	in the same regions.  Both Landau
	estimates are entrywise.  In the first, the implied constant may depend on
	\(\delta\), \(r\), and \(\mathrm d\Psi\); in the second, it may depend on
	\(\delta\), \(N\), and the fixed spectral data, but in neither case does it
	depend on \(z\).
\end{proposition}

\begin{proof}
	The scalar identity
	\begin{equation*}
	\frac1{z-x}
	=
	\sum_{n=0}^{r-1}\frac{x^n}{z^{n+1}}
	+
	\frac{x^r}{z^r(z-x)}
	\end{equation*}
	gives, entry by entry,
	\begin{equation*}
	\mathsf S_{\Psi}(z)
	-
	\sum_{n=0}^{r-1}\frac{\boldsymbol\mu_n}{z^{n+1}}
	=
	\frac1{z^r}
	\int_{\mathbb R}\frac{x^r}{z-x}\,\mathrm d\Psi(x).
	\end{equation*}
	In the stated region, \(|z-x|\ge\delta|z|\) for \(x\in\mathbb R\).  Since
	the measures are positive, finiteness of the even moments implies
	finiteness of every absolute moment.  The absolute value of each entry of
	the remainder is therefore bounded by its finite absolute \(r\)-th moment
	divided by \(\delta|z|^{r+1}\).  This proves
	\eqref{eq:unbounded-Weyl-Poincare}.  The path argument in the proof of
	Corollary~\ref{cor:finite-CF-Weyl-Pade} shows algebraically, without a
	boundedness assumption, that \(\mathsf S^{[N]}\) has the same first
	\(\operatorname{ord}_N\) moment coefficients.  The matrix
	\(\mathsf S^{[N]}\) is rational and analytic at infinity.  Consequently,
	for fixed \(N\), its Laurent expansion may be truncated after
	\(z^{-\operatorname{ord}_N}\), with a rational remainder that is
	\(\mathrm O(|z|^{-\operatorname{ord}_N-1})\) for all sufficiently large
	\(|z|\).  Subtracting this truncated expansion from
	\eqref{eq:unbounded-Weyl-Poincare}, with
	\(r=\operatorname{ord}_N\), proves the last claim.
\end{proof}

\medskip
The mixed-type Favard theorem, first established in the bounded setting and
then extended to unbounded matrices, supplies a positive spectral measure
when the finite shifted PBFs form a coherent family with fixed
initial-condition matrices
\cite{BranquinhoFoulquieManas2023Spectral,
BranquinhoFoulquieManas2026Unbounded}; \(\Delta\mathrm{TP}\) initial
conditions give the strong-normality conclusion described in
\cite{BranquinhoFoulquieManas2024BTP}.  These results justify the existence
of the measure-defined spectral Weyl matrix used above.  Independently,
\eqref{eq:unbounded-rectangular-stability} shows that uniformly stable height
sections converge to the projected resolvent \(\mathsf R_{\mathcal H_0}\) of
a specified closed realization.  This limit is the spectral Weyl matrix
\(\mathsf S_{\Psi}\) only under the additional identification
\eqref{eq:unbounded-resolvent-Cauchy-identification}.  The formal
moment-matching statement holds independently of these analytic hypotheses.

\begin{remark}[Weyl matrices of the cyclic Darboux transforms]
	\label{rem:Darboux-shifted-Weyl-matrices}
When the PBF is the mixed-type Christoffel factorization and its
\(\tau\)-determinants do not vanish,
Theorem~\ref{thm:cyclic-Favard-transport} supplies the transformed measure,
spectral data, and normalization matrices of every cyclic Darboux transform.
Its spectral Weyl matrix is defined from that measure by
	\[
	\mathsf S^{(a)}(z)
	\coloneq
	\int_{\mathbb R}
	\frac{\mathrm d\Psi^{(a)}(x)}{z-x},
	\qquad z\in\mathbb C\setminus\mathbb R.
	\]
	If \(T_{(a)}\) has a specified closed operator realization and its
	normalized projected operator resolvent equals this Cauchy transform---in
	particular under the bounded hypotheses used above---then
	Corollary~\ref{cor:Darboux-shifted-rectangular-series} gives, for
	\(z^{-1}=\tau^M\),
	\[
	\mathsf S^{(a)}(z)
	=
	\frac{1}{z}
	(\xi^{(a)})^{-1}
	\mathsf Q_q^{(a)}(I-\tau\mathscr H)^{-1}
	\bigl(\mathsf Q_p^{(a)}\bigr)^{\top}
	(\nu^{(a)})^{-\top}.
	\]
	Thus, after replacing \(\mathsf Q_q,\mathsf Q_p,\xi,\nu\) by
	\(\mathsf Q_q^{(a)},\mathsf Q_p^{(a)},\xi^{(a)},\nu^{(a)}\), the continued-fraction
	representation derived from \((I-\tau\mathscr H)^{-1}\) gives the Weyl
	matrices of all cyclic Darboux transforms.  For an unbounded PBF, the
	Favard measures and their spectral Weyl matrices are defined by the measure
	integral.  The displayed resolvent formula additionally requires the stated
	operator realization and its identification with that integral.  For a
	general bidiagonal factorization, the projected-resolvent identities retain
	their formal power-series meaning.
\end{remark}

\section{The mixed-type Pi\~neiro Weyl fraction}
\label{sec:case-study-mixed-Pineiro}

The general construction is now specialized to the mixed-type Pi\~neiro
system of~\cite{BranquinhoFoulquieManas2026Pineiro,ManasMarkov2026}.  The
cited works provide the polynomial forms and bidiagonal coefficients used
here; \cite[Section~7 and Appendix~A]{ManasMarkov2026} gives an independent,
self-contained derivation of the factor entries and the exact PBF chamber.
The rational
example below evaluates these formulas completely: the measure, every
required factor entry, the product matrix, and the finite resolvents are
displayed and checked directly.  The positive endpoint vector then
turns the resulting fraction into convergent Pad\'e-type approximants to the
explicit Pi\~neiro Weyl matrix throughout \(|z|>1\).
The canonical Pi\~neiro bidiagonal factors may be unbounded on the standard
sequence space.  Their Favard representation and every finite continued
fraction remain explicit, while analytic convergence follows from the
stochastic diagonal normalization used below.
Subsection~\ref{subsec:mixed-Pineiro-specialization-CF} starts from its matrix of
power weights and their moments, recalls the positive-factor chamber, and
identifies the factorization used as input.
Subsection~\ref{subsec:Pineiro-height-block-CF} gives the resulting height-block
continued fraction and tests backward evaluation, convergence, and residual
error in a rational \((3,2)\) example.  It also gives a two-state recursion
with two different return depths, genuinely distinct from the height-block
recursion.
Subsection~\ref{subsec:Pineiro-branched-qstate-tails} specializes the \(q\)-state
first-return recursion, and
Subsection~\ref{subsec:Pineiro-rectangular-Weyl} identifies the rectangular Weyl matrix
extracted from \((I-\tau\mathscr H_{p,q}^{\mathrm P})^{-1}\) and proves the
stochastic convergence statement just described.

\subsection{Measures, explicit bidiagonal factorization, and positivity}
\label{subsec:mixed-Pineiro-specialization-CF}

Consider a concrete matrix of measures for which the bidiagonal coefficients entering
the preceding constructions can be written explicitly.  This example is
the mixed-type Pi\~neiro specialization of the Jacobi-like systems
\cite{Manas2026Hypergeometric}, and it will be used
throughout this section.  It is useful for two reasons.  First, it
gives a nontrivial family where the continued fractions in the individual
bidiagonal coefficients contain
no unevaluated determinants or moments.  Second, under the explicit parameter
conditions given below, all bidiagonal coefficients are positive.  The
resulting continued fraction therefore has positive coefficients and
represents the corresponding mixed-type Weyl matrix whenever
Theorem~\ref{thm:general-Weyl-CF} and
Proposition~\ref{prop:general-Weyl-Stieltjes} apply.

The cyclic Christoffel transforms remain within the Pi\~neiro
special-function family.  Multiplying a row or a column of the matrix weight
by \(x\) increments the corresponding power parameter, and the
polynomial transport in
\eqref{eq:first-cyclic-polynomial-transport} reduces to the associated
contiguous hypergeometric identities.  The role of the Pi\~neiro section is
instead to provide a determinant-free test of the cyclic matrix, the two
continued-fraction recursions, the Weyl projection, and the positivity and
convergence results.

Let
\[
\boldsymbol\alpha=(\alpha_1,\ldots,\alpha_p),
\qquad
\boldsymbol\beta=(\beta_1,\ldots,\beta_q),
\]
and consider the mixed-type Pi\~neiro matrix of measures on \((0,1)\)
\begin{equation}
	\label{eq:mixed-Pineiro-matrix-CF-paper}
	\mathrm d\boldsymbol\mu^{\mathrm P}(x)
	=
	\left[
		x^{\alpha_i+\beta_j}
	\right]_{
		\substack{j\in\{1,\ldots,q\}\\ i\in\{1,\ldots,p\}}
	}
	\dx.
\end{equation}
The moment entries are therefore
\begin{equation}
	\label{eq:mixed-Pineiro-moments-CF-paper}
	\int_0^1 x^n\,\mathrm d\mu^{\mathrm P}_{j,i}(x)
	=
	\frac{1}{n+\alpha_i+\beta_j+1},
\end{equation}
whenever the exponents are integrable at the origin.  This is the degeneration
of the Jacobi-like matrix in which both sides of the matrix of measures are
power vectors.

Throughout this section, assume that \(\alpha_i+\beta_j>-1\) for every
\(i\in\{1,\ldots,p\}\) and \(j\in\{1,\ldots,q\}\).  Assume also that the
step-line multi-indices of the original system and of all the cyclic
Christoffel transforms used below are weakly normal.  Equivalently, the
corresponding leading principal minors of the moment matrices are nonzero.
This hypothesis guarantees the existence of the
Gauss--Borel factorizations, the recurrence matrix, and the Christoffel
connections.  It also excludes the zeros of the denominators in the explicit
formulas below.  The parameter conditions
\eqref{eq:Pineiro-positive-chamber-CF-paper} give a concrete sufficient
condition for weak normality.  In that chamber the two power systems are AT
systems, so the corresponding step-line problems are in fact strongly normal
in the sense of Definition~\ref{def:weak-strong-step-line-normality}.

The step-line recurrence matrix has the bidiagonal factorization
\begin{equation}
	\label{eq:Pineiro-PBF-CF-paper}
	T^{\mathrm P}
	=
	L_1^{\mathrm P}\cdots L_p^{\mathrm P}
	U_q^{\mathrm P}\cdots U_1^{\mathrm P}.
\end{equation}
The determinant-free Gamma--Pochhammer formulas for all its nontrivial
entries are proved in
\cite[Corollary~7.5 and equations~(133)--(136)]{Manas2026Hypergeometric};
see also the Christoffel-factorization theorem in
\cite{BranquinhoFoulquieManas2026Factorization}.  The example below evaluates
them at rational parameters and displays the factors, the resulting
recurrence matrix, and a finite continued fraction.  Write
\[
\ell_N^{\mathrm P,(a)}\coloneq(L_a^{\mathrm P})_{N+1,N},
\qquad
u_N^{\mathrm P,(b)}\coloneq(U_b^{\mathrm P})_{N,N}
\]
for the nontrivial entries used in the computation below.

The factorization is positive in the chamber
\begin{equation}
	\label{eq:Pineiro-positive-chamber-CF-paper}
	\alpha_1<\cdots<\alpha_p<\alpha_1+1,
	\qquad
	\beta_1<\cdots<\beta_q<\beta_1+1,
	\qquad
	\alpha_1+\beta_1+1>0,
\end{equation}
as follows from the endpoint and factor-positivity results in
\cite{ManasMarkov2026}.  Those results use the auxiliary normalization
\(\alpha_1,\beta_1>-1\), which causes no loss here.  Indeed, the opposite
translation
\[
\alpha_i\longmapsto\alpha_i+c,
\qquad
\beta_j\longmapsto\beta_j-c
\]
leaves the measures and the factor quotients unchanged.  Since
\(\alpha_1+\beta_1+1>0\), the interval
\((-1-\alpha_1,\beta_1+1)\) is nonempty, and any \(c\) in this interval
makes both translated first parameters greater than \(-1\).  Therefore, in
\eqref{eq:Pineiro-positive-chamber-CF-paper}, all entries
\[
\ell_N^{\mathrm P,(k)}>0,
\qquad
u_N^{\mathrm P,(b)}>0,
\qquad
N\in\Nzero,
\]
are positive.  If \(B_N^{\mathrm P}\) is the \(N\)-th right step-line form,
then
\begin{equation*}
v_N^{\mathrm P}
\coloneq
B_N^{\mathrm P}(1)
=
\sum_{j=1}^{q}B_N^{(j),\mathrm P}(1),
\qquad N\in\mathbb N_0.
\end{equation*}
The endpoint theorem gives \(v_N^{\mathrm P}>0\), and evaluation of the
recurrence at \(x=1\) gives
\[
T^{\mathrm P}\boldsymbol v^{\mathrm P}=\boldsymbol v^{\mathrm P}.
\]
The terminating component formulas in
\cite{BranquinhoFoulquieManas2026Pineiro,ManasMarkov2026} express every
\(v_N^{\mathrm P}\) as a finite Gamma--Pochhammer sum.  This formula supplies
the positive eigenvector used in the stochastic normalization below.  Since
each factor is nonnegative and has a fixed positive
normalizing entry in every row, all intermediate factor images of
\(\boldsymbol v^{\mathrm P}\) are strictly positive.

\subsection{The height-block fraction and an exact rational example}
\label{subsec:Pineiro-height-block-CF}

The general construction of
Sections~\ref{sec:positive-refactorization-sparse-matrices}--\ref{sec:branched-q-state-CF} applies directly to
\eqref{eq:Pineiro-PBF-CF-paper}.  The rational example follows the chain
\[
\text{bidiagonal factors}
\longrightarrow T^{\mathrm P}
\longrightarrow \text{step-line recurrence}
\longrightarrow \text{finite matrix continued fraction}.
\]
Positivity and the endpoint vector then yield numerical convergence.
For a semi-infinite matrix \(F\), write
\[
\operatorname{pr}_5(F)\coloneq[F_{n,k}]_{n,k=0}^{4}.
\]
In the prescribed factor order all lower factors precede all upper factors.
Consequently, a path contributing to this leading principal block cannot
leave the indices \(0,1,2,3,4\) and later return, and hence
\[
\operatorname{pr}_5(T)
=\operatorname{pr}_5(L_1)\cdots\operatorname{pr}_5(L_p)
 \operatorname{pr}_5(U_q)\cdots\operatorname{pr}_5(U_1).
\]
This permits the finite factorization and its coefficient sequences to be
displayed together.

\begin{example}[Backward evaluation for a rational \texorpdfstring{$(3,2)$}{(3,2)}
	Pi\~neiro system]
	\label{ex:Pineiro-CF-numerical}
	Take
	\[
	\boldsymbol\alpha=
	\begin{bNiceMatrix}-\frac13&0&\frac13\end{bNiceMatrix},
	\qquad
	\boldsymbol\beta=
	\begin{bNiceMatrix}0&\frac13\end{bNiceMatrix}.
	\]
	These parameters satisfy
	\eqref{eq:Pineiro-positive-chamber-CF-paper}.  In this case the matrix of
	measures is
	\begin{equation*}
		\mathrm d\boldsymbol\mu^{\mathrm P}(x)
		=
		\begin{bNiceMatrix}[margin=2pt,cell-space-limits=2pt]
			x^{-1/3}&1&x^{1/3}\\
			1&x^{1/3}&x^{2/3}
		\end{bNiceMatrix}\mathrm dx,
		\qquad 0<x<1.
	\end{equation*}
	Thus every moment, factor coefficient, and finite continued-fraction block
	is rational.  Evaluation of the explicit factor formulas cited above gives
	\begin{align*}
		\operatorname{pr}_5(L_1^{\mathrm P})
		&=\begin{bNiceMatrix}[margin=2pt,cell-space-limits=2pt]
		1&0&0&0&0\\ \frac1{12}&1&0&0&0\\ 0&\frac8{35}&1&0&0\\
		0&0&\frac1{12}&1&0\\ 0&0&0&\frac7{22}&1
		\end{bNiceMatrix},
		&
		\operatorname{pr}_5(L_2^{\mathrm P})
		&=\begin{bNiceMatrix}[margin=2pt,cell-space-limits=2pt]
		1&0&0&0&0\\ \frac1{20}&1&0&0&0\\ 0&\frac5{28}&1&0&0\\
		0&0&\frac1{15}&1&0\\ 0&0&0&\frac{40}{143}&1
		\end{bNiceMatrix},
		\\
		\operatorname{pr}_5(L_3^{\mathrm P})
		&=\begin{bNiceMatrix}[margin=2pt,cell-space-limits=2pt]
		1&0&0&0&0\\ \frac1{30}&1&0&0&0\\ 0&\frac17&1&0&0\\
		0&0&\frac3{55}&1&0\\ 0&0&0&\frac{45}{182}&1
		\end{bNiceMatrix},
		&
		\operatorname{pr}_5(U_2^{\mathrm P})
		&=\begin{bNiceMatrix}[margin=2pt,cell-space-limits=2pt]
		\frac35&1&0&0&0\\ 0&\frac5{21}&1&0&0\\ 0&0&\frac7{12}&1&0\\
		0&0&0&\frac8{55}&1\\ 0&0&0&0&\frac{15}{26}
		\end{bNiceMatrix},
		\\
		\operatorname{pr}_5(U_1^{\mathrm P})
		&=\begin{bNiceMatrix}[margin=2pt,cell-space-limits=2pt]
		\frac23&1&0&0&0\\ 0&\frac35&1&0&0\\ 0&0&\frac5{28}&1&0\\
		0&0&0&\frac{35}{66}&1\\ 0&0&0&0&\frac8{65}
		\end{bNiceMatrix}.
	\end{align*}
	Every entry of the five finite bidiagonal factors is now explicit.  Direct
	multiplication gives
	\begin{equation*}
		\operatorname{pr}_5(T^{\mathrm P})
		=
		\operatorname{pr}_5(L_1^{\mathrm P})
		\operatorname{pr}_5(L_2^{\mathrm P})
		\operatorname{pr}_5(L_3^{\mathrm P})
		\operatorname{pr}_5(U_2^{\mathrm P})
		\operatorname{pr}_5(U_1^{\mathrm P})
		=
		\begin{bNiceMatrix}[margin=2pt,cell-space-limits=2pt]
			\frac25&\frac65&1&0&0\\
			\frac1{15}&\frac{12}{35}&\frac7{12}&1&0\\
			\frac1{100}&\frac{19}{175}&\frac{43}{120}&\frac{183}{110}&1\\
			\frac1{5040}&\frac{17}{2940}&\frac{13}{352}&\frac{1051}{3080}&\frac{123}{260}\\
			0&\frac1{2310}&\frac{389}{56628}&\frac{1008}{7865}&\frac{5282}{15015}
		\end{bNiceMatrix}.
	\end{equation*}
	The zero pattern displays the complete \((p,q)=(3,2)\) band:
	\(T^{\mathrm P}_{n,k}=0\) outside \(n-3\le k\le n+2\).  In particular,
	the first three entries on the outer superdiagonal are \(1\).  The lower
	left corner also shows the band boundary: \(T^{\mathrm P}_{3,0}=1/5040>0\)
	but \(T^{\mathrm P}_{4,0}=0\), while
	\(T^{\mathrm P}_{4,1}=1/2310>0\).
	For example, the first row is the concrete recurrence
	\[
	xB_0^{\mathrm P}(x)=\frac25B_0^{\mathrm P}(x)
	+\frac65B_1^{\mathrm P}(x)+B_2^{\mathrm P}(x).
	\]
	Thus the upper and lower recurrence coefficients are recovered from the
	bidiagonal factors rather than inserted as independent data.
	Let \(v_N\) be the value at \(x=1\) of the \(N\)-th right step-line
	form, and put \(\boldsymbol v=(v_N)_{N\ge0}\).  Its first values are
	\[
		v_0=1,\qquad v_1=\frac13,\qquad v_2=\frac15,\qquad
		v_3=\frac1{28},\qquad v_4=\frac1{44},\qquad
		v_5=\frac7{2145}.
	\]
	The recurrence at \(x=1\) gives
	\(T^{\mathrm P}\boldsymbol v=\boldsymbol v\), and the ordered
	positive factorization gives positive intermediate vectors.  The diagonal
	normalization of Definition~\ref{def:stochastic-height-normalization}, with
	the properties proved in Proposition~\ref{prop:stochastic-height-gauge},
	therefore converts all five factors into row-stochastic bidiagonal
	matrices.  Hence \(\|\widehat{\mathscr H}_{3,2}^{\mathrm P}\|_\infty=1\).

	For \((p,q)=(3,2)\), the reduced height parameters are
	\(P=3\), \(Q=2\), \(R=5\), and the block width is \(3\).
	With the factor order
	\((L_1,L_2,L_3,U_2,U_1)\), the first height blocks are
	\[
		\begin{aligned}
			\mathscr V_0&=((1,0),(2,0)),&
			\mathscr V_1&=((5,0),(3,0),(1,1)),\\
			\mathscr V_2&=((4,0),(2,1),(5,1)),&
			\mathscr V_3&=((3,1),(1,2),(4,1)).
		\end{aligned}
	\]
	Before stochastic normalization, the cut after \(\mathscr V_3\) is the
	following explicit block-tridiagonal matrix.  The horizontal and vertical
	lines separate the coordinates belonging to
	\(\mathscr V_0,\mathscr V_1,\mathscr V_2,\mathscr V_3\), in that order:
	\begin{equation*}
		(\mathscr H_{3,2}^{\mathrm P})^{[0,3]}
		=
		\begin{bmatrix}[margin=4pt,cell-space-limits=2pt]
			\mathsf A_0^{\mathrm P}&\mathsf B_0^{\mathrm P}&0&0\\
			\mathsf C_0^{\mathrm P}&\mathsf A_1^{\mathrm P}&\mathsf B_1^{\mathrm P}&0\\
			0&\mathsf C_1^{\mathrm P}&\mathsf A_2^{\mathrm P}&\mathsf B_2^{\mathrm P}\\
			0&0&\mathsf C_2^{\mathrm P}&\mathsf A_3^{\mathrm P}
		\end{bmatrix}
		=
		\begin{bNiceArray}{cc|ccc|ccc|ccc}[margin=4pt,cell-space-limits=2pt]
			0&1&0&0&0&0&0&0&0&0&0\\
			0&0&0&1&0&0&0&0&0&0&0\\
			\hline
			\frac23&0&0&0&1&0&0&0&0&0&0\\
			0&0&0&0&0&1&0&0&0&0&0\\
			0&\frac1{12}&0&0&0&0&1&0&0&0&0\\
			\hline
			0&0&\frac35&0&0&0&0&1&0&0&0\\
			0&0&0&\frac1{20}&0&0&0&0&1&0&0\\
			0&0&0&0&\frac35&0&0&0&0&1&0\\
			\hline
			0&0&0&0&0&\frac1{30}&0&0&0&0&1\\
			0&0&0&0&0&0&\frac8{35}&0&0&0&0\\
			0&0&0&0&0&0&0&\frac5{21}&0&0&0
		\end{bNiceArray}.
	\end{equation*}
	Thus the individual factor coefficients occur as the nonzero transitions
	of the height matrix.  Define the Schur denominators of the displayed finite
	continued fraction backwards by
	\begin{align*}
	\mathsf D_3(\tau)
	&=I_3-\tau
	\begin{bNiceMatrix}0&0&1\\0&0&0\\0&0&0\end{bNiceMatrix},\\
	\mathsf D_2(\tau)
	&=I_3-\tau
	\begin{bNiceMatrix}0&0&1\\0&0&0\\0&0&0\end{bNiceMatrix}
	-\tau^2
	\begin{bNiceMatrix}0&0&0\\1&0&0\\0&1&0\end{bNiceMatrix}
	\mathsf D_3(\tau)^{-1}
	\begin{bNiceMatrix}\frac1{30}&0&0\\0&\frac8{35}&0\\0&0&\frac5{21}\end{bNiceMatrix},\\
	\mathsf D_1(\tau)
	&=I_3-\tau
	\begin{bNiceMatrix}0&0&1\\0&0&0\\0&0&0\end{bNiceMatrix}
	-\tau^2
	\begin{bNiceMatrix}0&0&0\\1&0&0\\0&1&0\end{bNiceMatrix}
	\mathsf D_2(\tau)^{-1}
	\begin{bNiceMatrix}\frac35&0&0\\0&\frac1{20}&0\\0&0&\frac35\end{bNiceMatrix},\\
	\mathsf D_0(\tau)
	&=I_2-\tau
	\begin{bNiceMatrix}0&1\\0&0\end{bNiceMatrix}
	-\tau^2
	\begin{bNiceMatrix}0&0&0\\0&1&0\end{bNiceMatrix}
	\mathsf D_1(\tau)^{-1}
	\begin{bNiceMatrix}\frac23&0\\0&0\\0&\frac1{12}\end{bNiceMatrix}.
	\end{align*}
	The concrete finite fraction is
	\begin{equation}
		\Phi_0^{\mathrm P,[3]}(\tau)=\mathsf D_0(\tau)^{-1}.
	\end{equation}
	Thus every block entering \(\Phi_0^{\mathrm P,[3]}\) is visible.  Direct
	inversion of the displayed \(11\times11\) truncation and backward
	evaluation of the continued fraction
	give the same exact value:
	\begin{equation*}
		\left[
		\left(I_{11}-\frac12
		(\mathscr H_{3,2}^{\mathrm P})^{[0,3]}\right)^{-1}
		\right]_{\mathscr V_0\times\mathscr V_0}
		=\Phi_0^{\mathrm P,[3]}\!\left(\frac12\right)
		=
		\begin{bNiceMatrix}[margin=2pt,cell-space-limits=2pt]
			\frac{14857}{14670}&\frac{414}{815}\\
			\frac{187}{7335}&\frac{828}{815}
		\end{bNiceMatrix}.
	\end{equation*}
	This identity is the finite calculation behind the general Schur-complement
	recursion: the left-hand side uses the complete truncated height matrix,
	whereas the middle expression evaluates only its successive diagonal
	blocks.
	This small cut displays the block structure and the finite fraction
	explicitly.  In this example the states selected by \(\mathsf Q_2\) and
	\(\mathsf Q_3\) lie in the height blocks \(0\), \(1\), and \(3\), so the
	minimal index in the general Weyl approximation is \(N_0=3\).  Consequently,
	the cut \(N=3\) has
	\(\operatorname{ord}_3=\lceil 2/5\rceil=1\) and reproduces only the
	zeroth moment.  The deeper cuts in
	Table~\ref{tab:Pineiro-32-convergence} illustrate the increasing contact
	order together with analytic convergence.

	The same four-block display in the stochastic basis is
	\begin{equation*}
	(\widehat{\mathscr H}_{3,2}^{\mathrm P})^{[0,3]}
	=
	\begin{bNiceArray}{cc|ccc|ccc|ccc}[margin=4pt,cell-space-limits=2pt]
	0&1&0&0&0&0&0&0&0&0&0\\
	0&0&0&1&0&0&0&0&0&0&0\\
	\hline
	\frac23&0&0&0&\frac13&0&0&0&0&0&0\\
	0&0&0&0&0&1&0&0&0&0&0\\
	0&\frac14&0&0&0&0&\frac34&0&0&0&0\\
	\hline
	0&0&\frac35&0&0&0&0&\frac25&0&0&0\\
	0&0&0&\frac15&0&0&0&0&\frac45&0&0\\
	0&0&0&0&\frac12&0&0&0&0&\frac12&0\\
	\hline
	0&0&0&0&0&\frac16&0&0&0&0&\frac56\\
	0&0&0&0&0&0&\frac27&0&0&0&0\\
	0&0&0&0&0&0&0&\frac47&0&0&0
	\end{bNiceArray}.
	\end{equation*}
	The semi-infinite matrix \(\widehat{\mathscr H}_{3,2}^{\mathrm P}\) is
	row-stochastic, but the displayed principal truncation is only
	substochastic at its terminal boundary.  Its last two rows sum to
	\(2/7\) and \(4/7\), because their outgoing normalizing transitions enter
	\(\mathscr V_4\) and are deleted by the cut.  The full operator includes
	these transitions and satisfies
	\(\|\widehat{\mathscr H}_{3,2}^{\mathrm P}\|_\infty=1\).
	Its finite continued fraction is equally explicit.  With backward
	denominators
	\begin{align*}
	\widehat{\mathsf D}_3(\tau)
	&=I_3-\tau
	\begin{bNiceMatrix}0&0&\frac56\\0&0&0\\0&0&0\end{bNiceMatrix},\\
	\widehat{\mathsf D}_2(\tau)
	&=I_3-\tau
	\begin{bNiceMatrix}0&0&\frac25\\0&0&0\\0&0&0\end{bNiceMatrix}
	-\tau^2
	\begin{bNiceMatrix}0&0&0\\\frac45&0&0\\0&\frac12&0\end{bNiceMatrix}
	\widehat{\mathsf D}_3(\tau)^{-1}
	\begin{bNiceMatrix}\frac16&0&0\\0&\frac27&0\\0&0&\frac47\end{bNiceMatrix},\\
	\widehat{\mathsf D}_1(\tau)
	&=I_3-\tau
	\begin{bNiceMatrix}0&0&\frac13\\0&0&0\\0&0&0\end{bNiceMatrix}
	-\tau^2
	\begin{bNiceMatrix}0&0&0\\1&0&0\\0&\frac34&0\end{bNiceMatrix}
	\widehat{\mathsf D}_2(\tau)^{-1}
	\begin{bNiceMatrix}\frac35&0&0\\0&\frac15&0\\0&0&\frac12\end{bNiceMatrix},\\
	\widehat{\mathsf D}_0(\tau)
	&=I_2-\tau
	\begin{bNiceMatrix}0&1\\0&0\end{bNiceMatrix}
	-\tau^2
	\begin{bNiceMatrix}0&0&0\\0&1&0\end{bNiceMatrix}
	\widehat{\mathsf D}_1(\tau)^{-1}
	\begin{bNiceMatrix}\frac23&0\\0&0\\0&\frac14\end{bNiceMatrix},
	\end{align*}
	one has
	\[
	\widehat\Phi_0^{[3]}(\tau)=\widehat{\mathsf D}_0(\tau)^{-1}.
	\]
	At the same test point,
	\[
	\left[
	\left(I_{11}-\frac12
	(\widehat{\mathscr H}_{3,2}^{\mathrm P})^{[0,3]}\right)^{-1}
	\right]_{\mathscr V_0\times\mathscr V_0}
	=\widehat\Phi_0^{[3]}\!\left(\frac12\right)
	=
	\begin{bNiceMatrix}[margin=4pt,cell-space-limits=2pt]
	\frac{14857}{14670}&\frac{414}{815}\\
	\frac{187}{7335}&\frac{828}{815}
	\end{bNiceMatrix}.
	\]
	The equality of this leading resolvent block with the unnormalized one is
	exact: the zeroth coordinate of every intermediate vector
	\(\boldsymbol v^{[r]}\) is \(1\), and both states of
	\(\mathscr V_0\) have sequence index \(0\); hence the height-ordered
	diagonal normalization restricts to \(I_2\) on \(\mathscr V_0\).  Finite
	diagonal conjugation therefore forces the two leading blocks to agree.  The
	stochastic calculation instead displays the normalized transition
	probabilities and verifies the corresponding backward denominators, with
	every block substituted.
	To make the notation
	in the numerical comparison explicit, set
	\begin{align*}
		\widehat\Phi_m(\tau)
		&\coloneq
		\left[
		\left(I-\tau
		(\widehat{\mathscr H}_{3,2}^{\mathrm P})^{\langle m\rangle}
		\right)^{-1}
		\right]_{\mathscr V_m\times\mathscr V_m},
		\\
		\widehat\Phi_m^{[N]}(\tau)
		&\coloneq
		\left[
		\left(I-\tau
		(\widehat{\mathscr H}_{3,2}^{\mathrm P})^{[m,N]}
		\right)^{-1}
		\right]_{\mathscr V_m\times\mathscr V_m}.
	\end{align*}
	For the infinite tail, Theorem~\ref{thm:matrix-CF-height-blocks} gives, at the first level,
	\[
		\widehat\Phi_0(\tau)
		=
		\left(
		I_2-\tau\widehat{\mathsf A}_0^{\mathrm P}
		-\tau^2\widehat{\mathsf B}_0^{\mathrm P}
		\widehat\Phi_1(\tau)\widehat{\mathsf C}_0^{\mathrm P}
		\right)^{-1}.
	\]
	The subsequent levels are obtained by continuing the same explicit
	backward recursion.  Moreover,
	\(\|\widehat{\mathscr H}_{3,2}^{\mathrm P}\|_\infty=1\).
	All resolvents, errors, and bounds in the remainder of this example are
	computed directly for the normalized operator
	\(\widehat{\mathscr H}_{3,2}^{\mathrm P}\); finite truncations are related
	to the unnormalized ones by the diagonal conjugations above.

	At \(\tau=0.75\), cutting after block \(220\) gives the reference
	\[
	\widehat\Phi_0^{[220]}(0.75)
	=
	\begin{bNiceMatrix}
	1.111994065342&0.856167208327\\
	0.149325420456&1.141556277770
	\end{bNiceMatrix}.
	\]
	Table~\ref{tab:Pineiro-32-convergence} compares shallower cuts with this
	reference.
	The observed error is measured in the matrix infinity norm.  The a priori column uses
	\eqref{eq:finite-CF-geometric-error} with \(\theta=0.75\), while the residual
	column uses \eqref{eq:finite-CF-a-posteriori-bound}.

	The factors were generated exactly and converted to double precision for
	the linear algebra.  Direct inversion and backward block evaluation were
	used independently to check every displayed value.

	\begin{table}[H]
	\centering
	\caption{Convergence and certified error bounds for the rational
		\((3,2)\) Pi\~neiro system at \(\tau=0.75\).  The columns give the
		truncation block \(N\), the finite-section dimension, the observed
		absolute matrix infinity-norm difference from the truncation after block
		\(220\), the a priori bound, and the a posteriori residual bound.}
	\label{tab:Pineiro-32-convergence}
	\(
	\begin{NiceArray}{r r c c c}[cell-space-limits=3pt]
	N&\dim&
	\textnormal{observed error}&
	\textnormal{a priori}&
	\textnormal{residual}\\
	\hline
	2&8&8.09\mathbin{\times}10^{-3}&1.42&2.51\mathbin{\times}10^{-1}\\
	4&14&1.91\mathbin{\times}10^{-4}&4.51\mathbin{\times}10^{-1}&2.49\mathbin{\times}10^{-2}\\
	8&26&6.78\mathbin{\times}10^{-8}&4.51\mathbin{\times}10^{-2}&4.92\mathbin{\times}10^{-4}\\
	12&38&2.86\mathbin{\times}10^{-11}&4.52\mathbin{\times}10^{-3}&1.07\mathbin{\times}10^{-5}\\
	32&98&\textnormal{roundoff}&4.54\mathbin{\times}10^{-8}&6.10\mathbin{\times}10^{-14}\\
	64&194&\textnormal{roundoff}&4.58\mathbin{\times}10^{-16}&5.80\mathbin{\times}10^{-27}\\
	192&578&\textnormal{roundoff}&4.75\mathbin{\times}10^{-48}&1.16\mathbin{\times}10^{-78}
	\end{NiceArray}
	\)
	\end{table}
	The direct inversion and backward evaluation agree to double-precision
	roundoff in every row.  The last printed digits depend on the order of the
	floating-point operations and are not used as reproducibility targets.  At
	\(N=32\) and beyond in the table, the observed difference from the reference
	is at the level of floating-point roundoff;
	the theoretical and residual bounds continue to decrease, certifying
	convergence beyond the precision visible in the observed-error column.
\end{example}

\subsection{Branched and \texorpdfstring{$q$}{q}-state tails}
\label{subsec:Pineiro-branched-qstate-tails}

The second grouping retains the \(q\) labeled states at the lower boundary of
each tail.  These \emph{boundary states} form the space \(\mathcal U_m\).
This grouping is also explicit for the mixed-type
Pi\~neiro family with arbitrary \(p\) and \(q\).  Let
\[
P=\kappa Q+\sigma,
\qquad
\kappa\in\mathbb N_0,
\qquad
\sigma\in\{0,\ldots,Q-1\},
\]
and, for \(m\ge m_*=p\), form the \(q\times q\) matrices
\[
\mathbf A_m^{\mathrm P},
\qquad
\mathbf B_m^{\mathrm P},
\qquad
\mathbf D_m^{\mathrm P,(\eta)},
\qquad
\eta\in\{\kappa,\kappa+1\},
\]
from \(\mathscr H_{p,q}^{\mathrm P}\) by the projections in
\eqref{eq:q-channel-A-B-D-definitions}.  The nonzero entries of
\(\mathbf B_m^{\mathrm P}\) are equal to \(1\).  This statement refers to
the unnormalized Pi\~neiro factors.  After a diagonal normalization, the
corresponding rise entries remain positive but need not equal \(1\), as the
	matrix in Example~\ref{ex:Pineiro-CF-numerical} shows.  Every nonzero coefficient entry
of \(\mathbf A_m^{\mathrm P}\) or
\(\mathbf D_m^{\mathrm P,(\eta)}\) is one of the factor coefficients
\[
\ell_N^{\mathrm P,(a)}
\qquad\textnormal{or}\qquad
u_N^{\mathrm P,(b)}
\]
of \eqref{eq:Pineiro-PBF-CF-paper}.
As in Subsection~\ref{subsec:q-channel-branched-refinement}, the superscript
\(\langle m\rangle_Q\) denotes principal compression to the tail formed by the
spaces \(\mathcal U_m,\mathcal U_{m+1},\ldots\), each formed from \(Q\)
consecutive heights.  For \(m\ge p\), define the Pi\~neiro \(q\)-state tail by
\[
\mathbf G_m^{\mathrm P}(\tau)
\coloneq
E_m^{\top}
\left(I-\tau(\mathscr H_{p,q}^{\mathrm P})^{\langle m\rangle_Q}\right)^{-1}
E_m.
\]

Theorem~\ref{thm:q-channel-matrix-BCF}, applied to the Pi\~neiro
factorization, now gives the \(q\)-state recursion without repeating its
general formula.  Its possible return depths are \(\kappa\) and
\(\kappa+1\), as determined by the Euclidean division of \(P\) by \(Q\).
If \(\sigma=0\), the second depth is absent; if \(\kappa=0\), the
depth-zero falls are contained in \(\mathbf A_m^{\mathrm P}\).  The weights
retain the complete dependence on \(\boldsymbol\alpha\) and
\(\boldsymbol\beta\) through the explicit factor sequences.  They are
strictly positive under
\eqref{eq:Pineiro-positive-chamber-CF-paper}.

\medskip

\begin{example}[The two return depths in the rational
\texorpdfstring{$(3,2)$}{(3,2)} system]
	\label{ex:Pineiro-32-q-state-depths}
	Return to the rational system in
	Example~\ref{ex:Pineiro-CF-numerical}.  Here
	\[
		d=1,
		\qquad
		P=3,
		\qquad
		Q=2,
		\qquad
		P=1\cdot Q+1.
	\]
	Thus \(\kappa=1\), \(\sigma=1\), and both possible return depths occur.
	The first complete boundary-state space, ordered by increasing height, is
	\[
		\mathcal U_3=((4,0),(2,1)).
	\]
	The pairs record the factor label and the nonnegative sequence index.  In
	this basis, and in the corresponding increasing-height bases of
	\(\mathcal U_4\) and \(\mathcal U_5\), direct substitution of the
	factor data used in Example~\ref{ex:Pineiro-CF-numerical} gives
	\[
		\mathbf A_3^{\mathrm P}=0,
		\qquad
		\mathbf B_3^{\mathrm P}=\mathbf B_4^{\mathrm P}=I_2,
		\qquad
		\mathbf D_3^{\mathrm P,(1)}
		=
		\begin{bNiceMatrix}0&0\\[1mm]\frac1{30}&0\end{bNiceMatrix},
		\qquad
		\mathbf D_3^{\mathrm P,(2)}
		=
		\begin{bNiceMatrix}0&\frac8{35}\\[1mm]0&0\end{bNiceMatrix}.
	\]
	The next two pairs of nonzero return coefficients are
	\[
		\begin{NiceArray}{c|ccc}[margin=4pt,cell-space-limits=2pt]
		m&3&4&5\\ \hline
		(\mathbf D_m^{\mathrm P,(1)})_{2,1}
		&\frac1{30}&\frac5{21}&\frac5{28}\\[1mm]
		(\mathbf D_m^{\mathrm P,(2)})_{1,2}
		&\frac8{35}&\frac5{28}&\frac17
		\end{NiceArray}.
	\]
	In particular, the first tail satisfies
	\[
		\mathbf G_3^{\mathrm P}(\tau)
		=
		\left(
		I_2
		-\tau^2\mathbf G_4^{\mathrm P}(\tau)
		\mathbf D_3^{\mathrm P,(1)}
		-\tau^3\mathbf G_4^{\mathrm P}(\tau)
		\mathbf G_5^{\mathrm P}(\tau)
		\mathbf D_3^{\mathrm P,(2)}
		\right)^{-1}.
	\]
	Since \(Q<P\), each \(\mathcal U_m\) groups two consecutive heights,
	whereas the block-tridiagonal construction groups three.  Thus the two
	recursions use genuinely different groupings in this example.
\end{example}

For \(q=1\), Corollary~\ref{cor:q-channel-recovers-scalar-BCF} reduces
this specialization of Theorem~\ref{thm:q-channel-matrix-BCF} to the
Jacobi--Pi\~neiro
\(p\)-branched Stieltjes fraction, with scalar weights given by the same
coefficients \(u_N^{\mathrm P,(1)}\) and
\(\ell_N^{\mathrm P,(a)}\).

\subsection{The rectangular Weyl matrix}
\label{subsec:Pineiro-rectangular-Weyl}

It remains to identify the rectangular Weyl matrix represented by the
general mixed-type Pi\~neiro continued fraction.  The moments of
\eqref{eq:mixed-Pineiro-matrix-CF-paper} are
\[
\int_0^1 x^n\,\mathrm d\boldsymbol\mu^{\mathrm P}(x)
=
\left[
\frac{1}{n+\alpha_i+\beta_j+1}
\right]_{
\substack{j\in\{1,\ldots,q\}\\i\in\{1,\ldots,p\}}},
\qquad
n\in\Nzero.
\]
For \(|z|>1\), the corresponding Stieltjes matrix is
\begin{equation}
	\label{eq:Pineiro-general-Weyl-entries-CF-paper}
	\mathsf S_{p,q}^{\mathrm P}(z)
	=
	\int_0^1
	\frac{\mathrm d\boldsymbol\mu^{\mathrm P}(x)}{z-x}
=
	\left[
	\frac{1}{z(\alpha_i+\beta_j+1)}
		\pFq{2}{1}
		{1,\alpha_i+\beta_j+1}
		{\alpha_i+\beta_j+2}
		{\dfrac1z}
	\right]_{
	\substack{j\in\{1,\ldots,q\}\\i\in\{1,\ldots,p\}}}.
\end{equation}

The weak-normality hypothesis gives the step-line Gauss--Borel factorization
\[
	\mathscr M^{\mathrm P}
	=
	(\mathscr L^{\mathrm P})^{-1}
	(\mathscr U^{\mathrm P})^{-1},
\]
and hence the two finite normalizations
\[
	\xi^{\mathrm P}
	\coloneq
	\mathsf E_{[q]}\mathscr L^{\mathrm P}\mathsf E_{[q]}^{\top},
	\qquad
	(\nu^{\mathrm P})^{\top}
	\coloneq
	\mathsf E_{[p]}\mathscr U^{\mathrm P}\mathsf E_{[p]}^{\top}.
\]

\begin{proposition}[Pi\~neiro moment normalization]
	\label{prop:Pineiro-moment-normalization-CF-paper}
	For every \(n\in\mathbb N_0\), the matrices just defined satisfy
	\[
		(\xi^{\mathrm P})^{-1}
		\mathsf E_{[q]}(T^{\mathrm P})^n\mathsf E_{[p]}^{\top}
		(\nu^{\mathrm P})^{-\top}
		=
		\int_0^1 x^n\,\mathrm d\boldsymbol\mu^{\mathrm P}(x).
	\]
\end{proposition}

\begin{proof}
	This is Theorem~\ref{thm:mixed-block-output-matrices} applied to the
	Pi\~neiro moment matrix, with the normalization determined by its leading
	Gauss--Borel blocks.
\end{proof}

Proposition~\ref{prop:Pineiro-moment-normalization-CF-paper} verifies the
moment hypothesis of Proposition~\ref{prop:general-Weyl-Stieltjes} with the concrete
matrix of measures \eqref{eq:mixed-Pineiro-matrix-CF-paper}.  Its support is contained in
\([0,1]\).

The preceding normalization also gives an explicit rotational lift of the
Pi\~neiro matrix of measures and identifies its Cauchy transform with a
projected resolvent of the sparse matrix associated with the prescribed
factor order.
For the two compactly supported applications below, write
\(\Sigma_M^{[1]}\coloneq
\bigcup_{k=0}^{M-1}\zeta_M^k[0,1]\).

\begin{proposition}[Rotational lift of the Pi\~neiro measure and its Weyl matrix]
\label{prop:Pineiro-star-system}
The rotational lift of \(\mathrm d\boldsymbol\mu^{\mathrm P}\) to
\(\Sigma_M^{[1]}\) is the positive
\(q\times p\) matrix of measures determined entrywise by
\begin{equation}
\label{eq:Pineiro-star-measure}
\int_{\Sigma_M^{[1]}}f(\lambda)\,
\mathrm d\widehat\mu_{j,i}^{\mathrm P}(\lambda)
=
\frac1M\sum_{k=0}^{M-1}\int_0^1
f\bigl(\zeta_M^k x^{1/M}\bigr)x^{\alpha_i+\beta_j}\,\mathrm dx.
\end{equation}
Its moments and Weyl matrix are
\begin{align}
\int_{\Sigma_M^{[1]}}\lambda^n\,
\mathrm d\widehat\mu_{j,i}^{\mathrm P}(\lambda)
&=
\begin{cases}
\displaystyle\frac1{r+\alpha_i+\beta_j+1},&n=Mr,\\[5pt]
0,&M\nmid n,
\end{cases}
\label{eq:Pineiro-star-moments}
\\
\widehat{\mathsf S}_{p,q}^{\mathrm P}(\lambda)
&=
\lambda^{M-1}\mathsf S_{p,q}^{\mathrm P}(\lambda^M),
\qquad
\lambda\in\mathbb C\setminus\Sigma_M^{[1]}.
\label{eq:Pineiro-star-Weyl}
\end{align}
This rotational lift is the rotationally invariant component of the cyclic
spectral system.  The factor-resolved representation and convergence of its
Weyl matrix are stated after the endpoint normalization in
Corollary~\ref{cor:Pineiro-stochastic-Weyl-convergence}.  If
	\(\gcd(p,q)=1\), the recurrence for the complete matrix of moment
	functionals is instead \(\mathscr H_{(a)}\) of
Theorem~\ref{thm:balanced-star-recurrence-matrix}; in exact-height order,
only its \(q\)-th superdiagonal and \(p\)-th subdiagonal are nonzero.  If
\(\gcd(p,q)>1\), every height contains \(d=\gcd(p,q)\) coordinate states
and the matrix \(\mathscr H_{(a)}\) has the two monomial block diagonals described in
Theorem~\ref{thm:balanced-two-extreme-diagonals}.  Grouping those states gives
the block mixed-type recurrence of
Theorem~\ref{thm:balanced-star-recurrence-matrix}, with normality at cuts
formed by complete height groups.

At every cyclic factor position, the corresponding Favard system is the
usual Pi\~neiro Christoffel transform: one of the power parameters is
increased by one according to its position in the right or left Christoffel
chain.  Thus all the transformed matrix measures and polynomial systems are
again explicit Pi\~neiro systems, and their rotational lifts are obtained
from \eqref{eq:Pineiro-star-measure} with the shifted parameter.
\end{proposition}

\begin{proof}
Equation~\eqref{eq:Pineiro-star-measure} is the specialization of
Definition~\ref{def:rotational-lift-matrix-measure}.  Summing the \(M\)-th
roots of unity gives \eqref{eq:Pineiro-star-moments}, and
Theorem~\ref{thm:cyclic-star-Weyl} gives
\eqref{eq:Pineiro-star-Weyl}.  The assertions about the matrix \(\mathscr H_{(a)}\)
follow from Theorems~\ref{thm:balanced-two-extreme-diagonals}
and~\ref{thm:balanced-star-recurrence-matrix};
Corollary~\ref{cor:balanced-coprime-double-band} gives the scalar form when
\(d=1\).  Finally,
Theorem~\ref{thm:cyclic-Favard-transport} identifies the other cyclic
positions with the successive Christoffel transforms.  For the power matrix
\eqref{eq:mixed-Pineiro-matrix-CF-paper}, these transforms increment the
corresponding exponent, so the family is preserved.
\end{proof}

Proposition~\ref{prop:Pineiro-star-system} produces an explicit rotationally
invariant positive matrix of measures on the star and its Weyl matrix from
the original interval system.  This positive lift represents the invariant
Fourier component; the remaining character components are the generally
complex measures described above.  In the rational \((3,2)\) example,
\(\gcd(p,q)=1\).  Its prescribed-order exact-height matrix is scalar
two-diagonal after a finite initial part, whereas \(\mathscr H_{(a)}\) has
that scalar two-diagonal form from height \(0\).

In the positive chamber, the endpoint vector and stochastic normalization
yield the following global convergence statement.

\begin{corollary}[Pi\~neiro Weyl convergence in the positive chamber]
\label{cor:Pineiro-stochastic-Weyl-convergence}
Assume \eqref{eq:Pineiro-positive-chamber-CF-paper}.  Apply
Definition~\ref{def:stochastic-height-normalization} to
\(\boldsymbol v^{\mathrm P}\), and denote the resulting height-ordered
diagonal matrix and stochastic height matrix by
\(\mathscr D_h^{\mathrm P}\) and
\(\widehat{\mathscr H}_{p,q}^{\mathrm P}\), respectively.  For
\(z^{-1}=\tau^{p+q}\) and \(|\tau|<1\),
\begin{equation}
\label{eq:Pineiro-stochastic-Weyl-factor-resolved}
\mathsf S_{p,q}^{\mathrm P}(z)
=
\frac1z
(\xi^{\mathrm P})^{-1}
\mathsf Q_q\mathscr D_h^{\mathrm P}
\left(I-\tau\widehat{\mathscr H}_{p,q}^{\mathrm P}\right)^{-1}
(\mathscr D_h^{\mathrm P})^{-1}\mathsf Q_p^{\top}
(\nu^{\mathrm P})^{-\top},
\qquad
z^{-1}=\tau^{p+q},
\quad |\tau|<1.
\end{equation}
The right-hand side is independent of the chosen \((p+q)\)-th root \(\tau\)
of \(z^{-1}\), because the two coordinate projections select paths that
complete whole factor cycles.
The corresponding finite rectangular projections converge locally uniformly
to \(\mathsf S_{p,q}^{\mathrm P}(z)\) for \(|z|>1\), with the error estimate
of \eqref{eq:stochastic-normalized-CF-error} before projection.
For real \(0\le\tau<1\), equivalently \(z=\tau^{-(p+q)}>1\), the normalized
tail blocks increase entrywise and their Schur denominators are nonsingular
\(M\)-matrices.
\end{corollary}

\begin{proof}
The endpoint vector verifies the hypotheses of
Proposition~\ref{prop:stochastic-height-gauge}, so
\(\|\widehat{\mathscr H}_{p,q}^{\mathrm P}\|_\infty=1\) and
Corollary~\ref{cor:stochastic-normalized-CF-convergence} applies.  The finite
diagonal conjugation of
Remark~\ref{rem:finite-diagonal-normalization-resolvents}, followed by the
coordinate projections, gives the finite version of
\eqref{eq:Pineiro-stochastic-Weyl-factor-resolved}; only finitely many entries
of \(\mathscr D_h^{\mathrm P}\) occur beside those projections.  Passing to
the limit and using
Proposition~\ref{prop:Pineiro-moment-normalization-CF-paper} identifies the
result with the Stieltjes matrix
\eqref{eq:Pineiro-general-Weyl-entries-CF-paper}.
\end{proof}

Thus the family realizes the complete general construction: its explicit
Pochhammer factors determine the finite fractions, and their rectangular
projections converge to the explicit \(q\times p\) Weyl matrix.  The cyclic
polynomial systems are the Pi\~neiro formulas with shifted power parameters
and give a direct verification of the Darboux transport.

\section{Jacobi-like extension of the Pi\~neiro Weyl fraction}
\label{sec:Jacobi-like-CF}

The mixed-type Pi\~neiro system of
\cite{BranquinhoFoulquieManas2026Pineiro,ManasMarkov2026} is the special case of the
broader Jacobi-like family of~\cite{Manas2026Hypergeometric} in which all
Gamma factors cancel.
In the complete-cancellation specialization, this relation holds for the
recurrence matrix, the rectangular Weyl matrix, and its continued-fraction
approximants.  The weight and factor formulas of
\cite{Manas2026Hypergeometric} provide the input for the construction.  This
section identifies the beta-convolution measure, computes its hypergeometric
Weyl matrix, and gives its rotational lift to the star.  When the parameters
additionally belong to the positive-factorization chamber stated below, the
canonical PBF yields the factor-resolved continued fraction and its finite
Weyl approximants.

In the Jacobi-like family the lower factors have closed Gamma--Pochhammer
entries, whereas the upper factors depend on moment minors.  Consequently,
the transformed systems can leave the original special-function family,
whereas the Pi\~neiro cyclic chain is described by direct parameter shifts.
Within that positive-factorization chamber, the essential input from
\cite[Theorem~9.4]{ManasMarkov2026} is stronger
than positivity of the original recurrence: it proves weak normality throughout
the cyclic Christoffel chain.  The present paper turns that simultaneous
weak normality into positive Favard systems for the complete Darboux orbit and
into their rotationally symmetric systems on \(p+q\) rays.
The canonical Jacobi-like factors can be unbounded on the standard sequence
space in part of the positive region.  The unbounded Favard theorem still
provides the matrix of measures and the weakly normal Christoffel chain;
boundedness enters only in the unnormalized operator-resolvent
representation.  The subsection below
derives the beta-type convolution weights, Gamma-quotient moments, and
generalized hypergeometric expression for the Weyl matrix and its rotational
lift.  Cancellation of all Gamma quotients recovers the Pi\~neiro measure
and Weyl matrix.

The standing condition \(a_1>-1\) is a convenient sufficient condition for
the beta-convolution construction.  In the complete-cancellation case it is
slightly stronger than the full Pi\~neiro integrability condition
\(\alpha_1+\beta_1+1>0\).

\subsection{Beta-convolution weights and the Weyl matrix}
\label{subsec:Jacobi-like-Weyl-CF}

We use the vector notation
\[
	\Gamma(\boldsymbol c)
	\coloneq
	\prod_{h=1}^{r}\Gamma(c_h),
	\qquad
	(\boldsymbol c)_n
	\coloneq
	\prod_{h=1}^{r}(c_h)_n,
	\qquad
	(c)_n
	\coloneq
	\frac{\Gamma(c+n)}{\Gamma(c)},
	\qquad
	\boldsymbol c=(c_1,\ldots,c_r).
\]
The generalized hypergeometric function is normalized by
\[
	{}_rF_s
	\left(
	\begin{matrix}
		a_1,\ldots,a_r\\
		b_1,\ldots,b_s
	\end{matrix}
	;z
	\right)
	\coloneq
	\sum_{n=0}^{\infty}
	\frac{(a_1)_n\cdots(a_r)_n}
	     {(b_1)_n\cdots(b_s)_n}
	\frac{z^n}{n!},
\]
in its disk of convergence, with the usual terminating or analytically
continued interpretation when applicable.  The Meijer function appearing
below is taken in the Mellin--Barnes normalization
\[
	G_{q,q}^{q,0}
	\left(
	x\,\middle|\,
	\begin{matrix}
		\boldsymbol b+\boldsymbol e_j^{(q)}\\
		\boldsymbol a
	\end{matrix}
	\right)
	\coloneq
	\frac{1}{2\pi\mathrm i}
	\int_{\mathcal L}
	\frac{\Gamma(s\one_q+\boldsymbol a)}
	     {\Gamma(s\one_q+\boldsymbol b+\boldsymbol e_j^{(q)})}
	x^{-s}\,\mathrm ds,
\]
where \(\mathcal L\) is an admissible Mellin--Barnes contour.  Equivalently,
this normalization is characterized by the Mellin transform
\[
	\int_0^1 x^{s-1}
	G_{q,q}^{q,0}
	\left(
	x\,\middle|\,
	\begin{matrix}
		\boldsymbol b+\boldsymbol e_j^{(q)}\\
		\boldsymbol a
	\end{matrix}
	\right)\dx
	=
	\frac{\Gamma(s\one_q+\boldsymbol a)}
	     {\Gamma(s\one_q+\boldsymbol b+\boldsymbol e_j^{(q)})}
\]
throughout the common fundamental strip.

Let
\[
\boldsymbol\alpha=(\alpha_1,\ldots,\alpha_p),\qquad
\boldsymbol a=(a_1,\ldots,a_q),\qquad
\boldsymbol b=(b_1,\ldots,b_q),
\]
and assume
\begin{equation}
\label{eq:Jacobi-like-CF-basic-order}
\alpha_1<\cdots<\alpha_p<\alpha_1+1,
\qquad
a_1<\cdots<a_q,
\qquad
b_1<\cdots<b_q,
\qquad
a_h\le b_h.
\end{equation}
Assume also that \(a_1>-1\) and
\(\alpha_i+a_h>-1\) for all \(i,h\).  Define the row weights by
their Mellin transforms
\begin{equation}
\label{eq:Jacobi-like-CF-weight-Mellin}
\int_0^1x^{s-1}w_j^{\mathrm J}(x)\dx
=
\frac{1}{s+b_j}
\frac{\Gamma(s\one_q+\boldsymbol a)}
{\Gamma(s\one_q+\boldsymbol b)},
\qquad 1\le j\le q.
\end{equation}
When all inequalities \(a_h<b_h\) are strict, these are positive Mellin
convolutions of beta-type kernels.  In the non-strict case, common Gamma
factors can be cancelled.  Indeed, \(a_i=b_j\) implies \(j\le i\), since
\(j>i\) would give \(a_i=b_j>b_i\).  Deleting this common value preserves
the componentwise inequalities between the remaining ordered lists: before
position \(j\) nothing changes, between \(j\) and \(i\) each surviving
\(a\)-parameter is paired with the next larger \(b\)-parameter, and after
position \(i\) both indices shift by one.  Repeating the cancellation leaves
strict inequalities among all surviving pairs.  Hence the remaining Gamma
quotient is again the Mellin transform of a positive beta-type convolution.
If all factors cancel, then \(w_j^{\mathrm J}(x)=x^{b_j}\).  Thus the row
weights remain positive.  In standard Meijer notation,
\begin{equation*}
w_j^{\mathrm J}(x)
=
G_{q,q}^{q,0}
\left(
x\,\middle|\,
\begin{matrix}
\boldsymbol b+\boldsymbol e_j^{(q)}\\
\boldsymbol a
\end{matrix}
\right),
\qquad 0<x<1.
\end{equation*}
Here \(\boldsymbol e_j^{(q)}\in\mathbb Z^q\) has a \(1\) in position \(j\)
and \(0\) in every other position.

The Jacobi-like matrix of measures is
\begin{equation}
\label{eq:Jacobi-like-CF-matrix-measure}
\mathrm d\boldsymbol\mu^{\mathrm J}(x)
=
\left[
x^{\alpha_i}w_j^{\mathrm J}(x)
\right]_{
\substack{1\le j\le q\\1\le i\le p}}
\dx.
\end{equation}
Its moments are explicit:
\begin{equation}
\label{eq:Jacobi-like-CF-moments}
\int_0^1x^n\,\mathrm d\mu_{j,i}^{\mathrm J}(x)
=
\frac{1}{n+\alpha_i+b_j+1}
\frac{
\Gamma((n+\alpha_i+1)\one_q+\boldsymbol a)
}{
\Gamma((n+\alpha_i+1)\one_q+\boldsymbol b)
},
\qquad n\in\Nzero.
\end{equation}
Consequently the Weyl matrix can be written without an unevaluated
integral.  Put
\[
\widetilde{\boldsymbol a}_i\coloneq(\alpha_i+1)\one_q+\boldsymbol a
=(\widetilde a_{i,1},\ldots,\widetilde a_{i,q}),
\qquad
\widetilde{\boldsymbol b}_i\coloneq(\alpha_i+1)\one_q+\boldsymbol b
=(\widetilde b_{i,1},\ldots,\widetilde b_{i,q}),
\qquad
	c_{j,i}^{\mathrm J}
	\coloneq
\frac{1}{\widetilde b_{i,j}}
\frac{\Gamma(\widetilde{\boldsymbol a}_i)}
{\Gamma(\widetilde{\boldsymbol b}_i)}.
\]
For \(|z|>1\), expansion of \((z-x)^{-1}\) and
\eqref{eq:Jacobi-like-CF-moments} give
\begin{equation}
\label{eq:Jacobi-like-CF-Weyl-hypergeometric}
\left(\mathsf S_{p,q}^{\mathrm J}(z)\right)_{j,i}
\coloneq
\int_0^1
\frac{x^{\alpha_i}w_j^{\mathrm J}(x)}{z-x}\dx
=
	\frac{c_{j,i}^{\mathrm J}}{z}
\pFq{q+1}{q}
{1,\widetilde{\boldsymbol a}_i}
{\widetilde{\boldsymbol b}_i+\boldsymbol e_j^{(q)}}
{\dfrac1z}.
\end{equation}
Indeed,
\[
\frac{1}{n+\widetilde b_{i,j}}
=
\frac{1}{\widetilde b_{i,j}}
\frac{(\widetilde b_{i,j})_n}
{(\widetilde b_{i,j}+1)_n},
\]
and the remaining Gamma quotients in
\eqref{eq:Jacobi-like-CF-moments} become the Pochhammer quotients in
\eqref{eq:Jacobi-like-CF-Weyl-hypergeometric}.  The identity
\((1)_n/n!=1\) supplies the first upper hypergeometric parameter.

The Jacobi-like data likewise determine both a concrete rotational lift of
the matrix of measures and the corresponding change of variable in its Weyl
matrix.

\begin{proposition}[Rotational lift of the Jacobi-like measure and its Weyl matrix]
\label{prop:Jacobi-like-star-system}
The positive \(q\times p\) matrix of measures on
\(\Sigma_M^{[1]}\) is defined by
\begin{equation}
\label{eq:Jacobi-like-star-measure}
\int_{\Sigma_M^{[1]}}f(\lambda)\,
\mathrm d\widehat\mu_{j,i}^{\mathrm J}(\lambda)
=
\frac1M\sum_{k=0}^{M-1}\int_0^1
f\bigl(\zeta_M^k x^{1/M}\bigr)
x^{\alpha_i}w_j^{\mathrm J}(x)\,\mathrm dx.
\end{equation}
Its moments vanish outside multiples of \(M\); for \(n=Mr\),
\begin{equation}
\label{eq:Jacobi-like-star-moments}
\int_{\Sigma_M^{[1]}}\lambda^{Mr}\,
\mathrm d\widehat\mu_{j,i}^{\mathrm J}(\lambda)
=
\frac{1}{r+\alpha_i+b_j+1}
\frac{
\Gamma((r+\alpha_i+1)\one_q+\boldsymbol a)
}{
\Gamma((r+\alpha_i+1)\one_q+\boldsymbol b)
}.
\end{equation}
Its Weyl matrix is the explicit hypergeometric matrix
\begin{equation}
\label{eq:Jacobi-like-star-Weyl}
\widehat{\mathsf S}_{p,q}^{\mathrm J}(\lambda)
=
\lambda^{M-1}\mathsf S_{p,q}^{\mathrm J}(\lambda^M),
\qquad
\lambda\in\mathbb C\setminus\Sigma_M^{[1]}.
\end{equation}
For \(\lvert\lambda\rvert>1\), the entries on the right are given by the
convergent hypergeometric series in
\eqref{eq:Jacobi-like-CF-Weyl-hypergeometric}.  Elsewhere in
\(\mathbb C\setminus\Sigma_M^{[1]}\), the right-hand side denotes the
single-valued analytic continuation determined by the Stieltjes integral
from infinity.  This rotational lift is the rotationally invariant component
of the cyclic spectral system.
\end{proposition}

\begin{proof}
Positivity of the radial weights was proved above.
Equation~\eqref{eq:Jacobi-like-star-measure} is therefore a positive
specialization of Definition~\ref{def:rotational-lift-matrix-measure}; the
root-of-unity sum and \eqref{eq:Jacobi-like-CF-moments} give
\eqref{eq:Jacobi-like-star-moments}.  Finally,
Theorem~\ref{thm:cyclic-star-Weyl} and
\eqref{eq:Jacobi-like-CF-Weyl-hypergeometric} give
\eqref{eq:Jacobi-like-star-Weyl} on the stated domain.
\end{proof}

We now impose the additional hypotheses needed for the PBF conclusions and
assume in addition that \(q\ge2\).
Let \(r_{\mathrm J}\in\{0,\ldots,q\}\) be the largest integer such that
every member of the initial string \(b_1,\ldots,b_{r_{\mathrm J}}\) occurs
in the \(a\)-list.  Delete those values from both lists and denote the
remaining increasing \(a\)-list by
\(\widehat a_1<\cdots<\widehat a_m\), where
\(m=q-r_{\mathrm J}\).  In addition to
\eqref{eq:Jacobi-like-CF-basic-order} and the preceding integrability
conditions, assume
\begin{equation}
\label{eq:Jacobi-like-CF-positive-chamber}
a_q<b_1+1,
\qquad
\widehat a_{m-1}\le b_{r_{\mathrm J}+1}
\quad\text{when }m\ge2.
\end{equation}

\begin{corollary}[Jacobi-like PBF realization of the Weyl matrices]
\label{cor:Jacobi-like-PBF-Weyl-realization}
Under \eqref{eq:Jacobi-like-CF-positive-chamber}, form the scalar step-line
moment matrix
\[
\mathscr M^{\mathrm J}_{qu+j,pv+i}
\coloneq
\int_0^1x^{u+v}\,\mathrm d\mu_{j+1,i+1}^{\mathrm J}(x),
\]
where \(u,v\in\mathbb N_0\), \(0\le j<q\), and \(0\le i<p\).
It has a Gauss--Borel factorization
\[
\mathscr M^{\mathrm J}
=
(\mathscr L^{\mathrm J})^{-1}(\mathscr U^{\mathrm J})^{-1},
\]
and its recurrence matrix has the canonical positive bidiagonal
factorization
\begin{equation}
\label{eq:Jacobi-like-CF-PBF}
T^{\mathrm J}
=
L_1^{\mathrm J}\cdots L_p^{\mathrm J}
U_q^{\mathrm J}\cdots U_1^{\mathrm J},
\qquad
T^{\mathrm J}
=
\mathscr L^{\mathrm J}\mathsf\Lambda_{[q]}
(\mathscr L^{\mathrm J})^{-1}
=
(\mathscr U^{\mathrm J})^{-1}
\mathsf\Lambda_{[p]}^{\top}\mathscr U^{\mathrm J}.
\end{equation}
With
\[
\xi^{\mathrm J}
\coloneq
\mathsf E_{[q]}\mathscr L^{\mathrm J}\mathsf E_{[q]}^{\top},
\qquad
(\nu^{\mathrm J})^{\top}
\coloneq
\mathsf E_{[p]}\mathscr U^{\mathrm J}\mathsf E_{[p]}^{\top},
\]
one has
\begin{equation}
\label{eq:Jacobi-like-moment-normalization}
(\xi^{\mathrm J})^{-1}\mathsf E_{[q]}(T^{\mathrm J})^n
\mathsf E_{[p]}^{\top}(\nu^{\mathrm J})^{-\top}
=
\int_0^1x^n\,\mathrm d\boldsymbol\mu^{\mathrm J}(x),
\qquad n\in\mathbb N_0.
\end{equation}

Let \(\mathscr H_{\mathrm{ord}}^{\mathrm J}\) be the prescribed-order
height matrix formed from \eqref{eq:Jacobi-like-CF-PBF}.  If it is bounded,
then, whenever \(z^{-1}=\tau^M\) and
\(\lvert\tau\rvert<\lVert\mathscr H_{\mathrm{ord}}^{\mathrm J}\rVert_2^{-1}\),
\[
\mathsf S_{p,q}^{\mathrm J}(z)
=
\frac1z(\xi^{\mathrm J})^{-1}\mathsf Q_q
\bigl(I-\tau\mathscr H_{\mathrm{ord}}^{\mathrm J}\bigr)^{-1}
\mathsf Q_p^{\top}(\nu^{\mathrm J})^{-\top}.
\]
The height-block Schur recursion therefore gives its factor-resolved matrix
continued fraction and finite Weyl approximants.  Without boundedness, the
same identity holds coefficientwise as a formal Laurent series at infinity.
Under the closed-realization, core, and uniform-stability hypotheses,
analytic finite-section convergence is to the corresponding projected
operator resolvent.  Identifying that limit with the Weyl matrix of a chosen
Favard measure requires the additional equality
\eqref{eq:unbounded-resolvent-Cauchy-identification}.

For each \(a\in\mathcal A_{\mathrm{cyc}}\), let
\(\mathscr H_{(a)}^{\mathrm J}\) be the minimum-height matrix obtained from
the positive refactorization of \(T_{(a)}^{\mathrm J}\).  After the
transported normalizations, its factor-copy moment functional belonging to
\(T_{(a)}^{\mathrm J}\) agrees, at the level of moments, with the rotational
lift of the corresponding shifted Favard measure.  If \(\gcd(p,q)=1\),
\(\mathscr H_{(a)}^{\mathrm J}\) is the scalar recurrence of the complete
matrix of star moment functionals.  If \(d=\gcd(p,q)>1\), it is instead the
block mixed-type recurrence, with normality asserted at cuts containing
complete groups of the \(d\) equal-height coordinates.
\end{corollary}

\begin{proof}
The positive-factorization theorem
\cite[Theorem~9.4]{ManasMarkov2026} and
Corollary~\ref{cor:balanced-Jacobi-complete-grid} give the Gauss--Borel and
positive bidiagonal factorizations.  Theorem~\ref{thm:mixed-block-output-matrices}
gives \eqref{eq:Jacobi-like-moment-normalization}, while
Proposition~\ref{prop:general-Weyl-Stieltjes} and
Theorem~\ref{thm:general-Weyl-CF} identify the explicit Weyl matrix and give
the factor-resolved formula.
The assertions for the selected products follow from
Corollary~\ref{cor:positive-Favard-data-cyclic-orbit},
Theorem~\ref{thm:cyclic-Favard-transport},
\eqref{eq:admissible-cyclic-star-resolvent}, and
Theorems~\ref{thm:balanced-two-extreme-diagonals}
and~\ref{thm:balanced-star-recurrence-matrix}.
\end{proof}

Thus Proposition~\ref{prop:Jacobi-like-star-system} supplies the rotational
measure and its Weyl matrix under the basic beta-convolution hypotheses,
whereas Corollary~\ref{cor:Jacobi-like-PBF-Weyl-realization} supplies the
recurrence and factor-resolved conclusions only under the additional
positive-factorization chamber.  The other cyclic products have their own
shifted positive Favard data, and their upper Christoffel transforms need
not remain Jacobi-like.

\section{Integrable deformations encoded by the PBF}
\label{sec:integrable-hierarchies-PBF}

The sparse recurrences constructed above do more than provide static spectral
models.  A scalar two-term difference operator with exponents \(q\) and
\(-p\) carries a hierarchy of commuting Lax flows for arbitrary positive
\(p,q\).  When \(p\) and \(q\) are coprime, \(\mathscr H_{(a)}\) is precisely
the semi-infinite matrix of such an operator: the PBF selects its initial
values, while the star moment functionals give its exact determinant solution.
The same construction evolves all cyclic products of the refactorized factors
with common Toda times.  Moreover, the continuous flows commute with the
discrete Christoffel transformations generated by the factors.

The case \(q=1\) includes the Narita--Itoh--Bogoyavlensky lattice.  Moment
solutions and Miura-type transformations for that one-sided setting have been
studied in \cite{Svinin2011Miura,Osipov2024Bogoyavlensky}.  After two general
algebraic results, the PBF realization below uses the two-sided factorization
with arbitrary coprime \(p\) and \(q\), and simultaneously retains its complete
cyclic family.

\subsection{The arithmetic-support hierarchy}
\label{subsec:arithmetic-support-PBF-hierarchy}

The next two algebraic results do not require \(p\) and \(q\) to be
coprime.  Let \(p,q\ge1\), put
\[
d\coloneq\gcd(p,q),
\qquad
M\coloneq p+q,
\qquad
R\coloneq\frac{M}{d},
\]
and let
\(\mathcal S\) be the forward shift on scalar sequences,
\((\mathcal Sf)_h=f_{h+1}\).  Let \(c=(c_h)_{h\in\mathbb Z}\) be a scalar
sequence, with \(c_h=0\) for \(h<p\), and define
\begin{equation}
	\label{eq:coprime-sparse-Lax-operator}
	\mathcal L
	\coloneq
	\mathcal S^q+c\mathcal S^{-p}
\end{equation}
Multiplication is determined by
\(\mathcal S^j c=c_{\,\cdot+j}\mathcal S^j\).  For a difference operator
\(X\), write \(X_+\) for the sum of its nonnegative powers of
\(\mathcal S\), and define explicitly
\(\operatorname{res}_{\mathcal S}X\coloneq[\mathcal S^0]X\), the
coefficient of \(\mathcal S^0\).

For \(d>1\), the minimal positive powers that generically preserve the
two-term support are \(\mathcal L^{R\ell}\): all exponents in
\(\mathcal L^n\) are congruent to \(-np\) modulo \(M\), and the required
tangency occurs when \(M\mid np\).  We deliberately retain below the
subsequence
\[
\mathcal L^{M\ell}=\mathcal L^{R(d\ell)},
\]
because its first residue expansion has \(M=p+q\) positions and therefore
matches the \(M\) bidiagonal factors.  Thus, when \(d>1\), the hierarchy
displayed below is the factor-compatible commuting subhierarchy; when
\(d=1\), one has \(R=M\) and no support-preserving positive flows are
omitted.

\begin{theorem}[The factor-compatible sparse arithmetic-support hierarchy]
	\label{thm:coprime-sparse-PBF-hierarchy}
	For every \(\ell\ge1\), the equation
	\begin{equation}
		\label{eq:coprime-sparse-PBF-Lax-flows}
		\frac{\partial\mathcal L}{\partial t_\ell}
		=
		\left[
		\bigl(\mathcal L^{M\ell}\bigr)_+,
		\mathcal L
		\right]
	\end{equation}
	defines a formal flow that preserves the two exponents \(q\) and \(-p\),
	and these formal flows commute.  On the Banach space
	\[
	\mathcal X_\infty
	\coloneq
	\left\{c\in\ell^\infty(\mathbb N_0;\mathbb C):
	c_h=0\ \text{for }0\le h<p\right\},
	\]
	with each sequence extended by zero to negative indices,
	the right-hand side is locally Lipschitz; hence, for bounded initial data,
	each equation defines a unique local \(C^1\) flow.  If
	\begin{equation}
		\label{eq:coprime-sparse-PBF-residue}
		\varrho_h^{(\ell)}
		\coloneq
		\left(
		\operatorname{res}_{\mathcal S}
		\mathcal L^{M\ell}
		\right)_h,
	\end{equation}
	then
	\begin{equation}
		\label{eq:coprime-sparse-PBF-coefficient-flow}
		\frac{\partial c_h}{\partial t_\ell}
		=
		c_h\left(
		\varrho_h^{(\ell)}-\varrho_{h-p}^{(\ell)}
		\right).
	\end{equation}
	Consequently, the boundary values \(c_h=0\), \(h<p\), remain zero.
	Along every real solution, strict positivity of the remaining coefficients
	is preserved throughout its interval of existence.
\end{theorem}

\begin{proof}
	Work first in the bilateral difference-operator algebra and put
	\[
	X_\ell\coloneq(\mathcal L^{M\ell})_+,
	\qquad
	Y_\ell\coloneq(\mathcal L^{M\ell})_-.
	\]
	A monomial in \(\mathcal L^{M\ell}\) containing \(k\) upward shifts has
	exponent
	\[
	kq-(M\ell-k)p=M(k-\ell p),
	\]
	so every exponent is a multiple of \(M\).  Since
	\([X_\ell+Y_\ell,\mathcal L]=0\),
	\[
	[X_\ell,\mathcal L]=-[Y_\ell,\mathcal L].
	\]
	The left-hand side has no exponent below \(-p\), whereas the right-hand
	side has no exponent above \(q-M=-p\).  Hence both sides are supported
	only at \(\mathcal S^{-p}\); in particular, the coefficient of
	\(\mathcal S^q\) vanishes.  Differentiating \(X_\ell\), taking
	nonnegative parts, and using
	\([\mathcal L^{Mk},\mathcal L^{M\ell}]=0\) gives the zero-curvature
	identity
	\[
	\partial_{t_k}X_\ell-
	\partial_{t_\ell}X_k-[X_k,X_\ell]=0,
	\]
	which proves that the formal flows commute.

	The identities are polynomial in finitely many shifts of the
	coefficients.  Their restriction to nonnegative indices produces no
	boundary term: a path that starts at a nonnegative height and visits a
	negative height must first use a downward edge from some height \(h<p\),
	whose weight is \(c_h=0\).

	Therefore, in the coefficient of \(\mathcal S^{-p}\) in the commutator,
	only the coefficient of \(\mathcal S^0\) in
	\((\mathcal L^{M\ell})_+\) contributes.  Since
	\[
	\varrho^{(\ell)}c\mathcal S^{-p}
	-c\mathcal S^{-p}\varrho^{(\ell)}
	=
	c\bigl(\varrho^{(\ell)}-
	\mathcal S^{-p}(\varrho^{(\ell)})\bigr)\mathcal S^{-p},
	\]
	equation \eqref{eq:coprime-sparse-PBF-coefficient-flow} follows.

	For fixed \(\ell\), each \(\varrho_h^{(\ell)}\) is a polynomial in a
	fixed finite set of shifts of \(c\), uniformly in \(h\).  Thus the
	right-hand side of \eqref{eq:coprime-sparse-PBF-coefficient-flow} defines
	a locally Lipschitz map \(\mathcal X_\infty\to\mathcal X_\infty\).
	The Banach-space Picard--Lindel\"of theorem gives the asserted local flow.
	Finally,
	\[
	c_h(t)=c_h(0)\exp\!\left(
	\int_0^t
	(\varrho_h^{(\ell)}-\varrho_{h-p}^{(\ell)})(s)\,\mathrm ds
	\right)
	\]
	proves preservation of the boundary values and of strict positivity along
	real solutions.
\end{proof}

General arithmetic-support hierarchies are classified in
\cite{Carpentier2026Arithmetic}; the preceding two-term proof is included to
keep the support-preservation argument used here self-contained.

For the first flow of this factor-compatible subhierarchy, the residue has a
direct path expansion.  Let
\(\mathcal P_{p,q}\) consist of the sequences
\(\boldsymbol\epsilon=(\epsilon_1,\ldots,\epsilon_M)\), with
\(\epsilon_j\in\{q,-p\}\), containing \(p\) occurrences of \(q\) and
\(q\) occurrences of \(-p\).  Put \(\xi_0=0\) and
\(\xi_j=\epsilon_1+\cdots+\epsilon_j\).  Identifying \(+q\) with \(L\)
and \(-p\) with \(U\) identifies \(\mathcal P_{p,q}\) with the set of all
factor-type words of Definition~\ref{def:factor-position-heights}.  If
\(\eta_j^{\mathrm{fac}}\) denotes the factor-position height there, then
\(\xi_j=d\eta_j^{\mathrm{fac}}\), where \(d=\gcd(p,q)\).  Thus the two
statistics coincide in the coprime case.  The sum below runs over all such
words, not only over \(\sigma_{\mathrm{min}}\) or its rotations.

\begin{proposition}[The first factor-compatible nonlinear flow]
	\label{prop:first-PBF-flow-path-polynomial}
	The polynomial in the first flow is
	\begin{equation}
		\label{eq:first-PBF-flow-path-polynomial}
		\varrho_h^{(1)}
		=
		\sum_{\boldsymbol\epsilon\in\mathcal P_{p,q}}
		\prod_{\substack{1\le j\le M\\ \epsilon_j=-p}}
		c_{h+\xi_{j-1}}.
	\end{equation}
	Before the boundary values are imposed, it is homogeneous of degree \(q\)
	and has positive integer coefficients.  If \(q=1\), then
	\begin{equation}
		\label{eq:Narita-Itoh-Bogoyavlensky-reduction}
		\frac{\partial c_h}{\partial t_1}
		=
		c_h\left(
		\sum_{j=1}^{p}c_{h+j}
		-
		\sum_{j=1}^{p}c_{h-j}
		\right),
	\end{equation}
	which is the Narita--Itoh--Bogoyavlensky lattice.
\end{proposition}

\begin{proof}
	A term of exponent zero in \(\mathcal L^M\) uses \(p\) shifts
	\(+q\) and \(q\) shifts \(-p\).  A downward shift in position \(j\)
	contributes the coefficient at height \(h+\xi_{j-1}\), which gives
	\eqref{eq:first-PBF-flow-path-polynomial}.  When \(q=1\), placing the
	unique downward shift in each of the \(p+1\) possible positions gives
	\(\varrho_h^{(1)}=\sum_{j=0}^{p}c_{h+j}\).  Substitution in
	\eqref{eq:coprime-sparse-PBF-coefficient-flow} proves
	\eqref{eq:Narita-Itoh-Bogoyavlensky-reduction}; see also
	\cite{Wang2012NaritaItohBogoyavlensky}.
\end{proof}

\subsection{Exact determinant solution from the star moments}
\label{subsec:PBF-arithmetic-tau-solution}

We now return to the PBF realization within this subsection.  Assume
\(\gcd(p,q)=1\), fix \(a\in\mathcal A_{\mathrm{cyc}}\), set
\(c_h=c_h^{(a)}\) for \(h\ge p\) and \(c_h=0\) for \(h<p\), and write
\[
\mathcal L_{(a)}(\boldsymbol0)
\coloneq
\mathcal S^q+c^{(a)}\mathcal S^{-p}.
\]
By Corollary~\ref{cor:balanced-coprime-double-band}, its semi-infinite
matrix is \(\mathscr H_{(a)}\).

Let \(\widehat{\boldsymbol\Psi}_{(a)}\) be the matrix of star moment
functionals whose recurrence is \(\mathscr H_{(a)}\), as constructed in
Section~\ref{sec:balanced-star-Christoffel}.  Put
\[
E_{\boldsymbol t}(\lambda)
\coloneq
\exp\!\left(\sum_{\ell\ge1}t_\ell\lambda^{M\ell}\right).
\]
All formal assertions below are understood coefficientwise after restriction
to an arbitrary finite set of active times.  Define the radial deformation
of the functionals coefficientwise by
\begin{equation}
	\label{eq:star-radial-Toda-deformation}
	\left\langle
	\widehat{\boldsymbol\Psi}_{(a)}(\boldsymbol t),f
	\right\rangle
	\coloneq
	\left\langle
	\widehat{\boldsymbol\Psi}_{(a)},
	f(\lambda)E_{\boldsymbol t}(\lambda)
	\right\rangle,
	\qquad f\in\mathbb C[\lambda].
\end{equation}
When representing measures exist and the displayed exponential moments are
finite, this is equivalently
\[
\mathrm d\widehat{\boldsymbol\Psi}_{(a)}
(\lambda;\boldsymbol t)
=
E_{\boldsymbol t}(\lambda)\,
\mathrm d\widehat{\boldsymbol\Psi}_{(a)}(\lambda).
\]
For an analytic interpretation, fix \(L<\infty\), set \(t_\ell=0\) for
\(\ell>L\), and fix a connected open neighbourhood
\(D\subset\mathbb C^L\) of \(\boldsymbol0\) on which all pairings in
\eqref{eq:star-radial-Toda-deformation} exist and depend holomorphically on
\(\boldsymbol t\).  This is automatic when the corresponding height
operator is bounded, and also when the relevant representing measures are
compactly supported.  The deformation is radial: it depends on \(\lambda\)
only through \(x=\lambda^M\) and hence preserves every Fourier character of
the star functional.

Let \(\widehat{\mathscr M}_{(a)}(\boldsymbol t)\) be the scalar step-line
moment matrix of \eqref{eq:star-radial-Toda-deformation}.  Thus, for
\(r=qu+j\), \(s=pv+i\),
\begin{equation}
	\label{eq:deformed-star-step-line-moments}
	\bigl(\widehat{\mathscr M}_{(a)}(\boldsymbol t)\bigr)_{r,s}
	=
	\left\langle
	\widehat\Psi_{(a),j,i}\,,\,
	\lambda^{u+v}
	\exp\!\left(\sum_{\ell\ge1}t_\ell\lambda^{M\ell}\right)
	\right\rangle.
\end{equation}
Set \(\tau_0^{(a)}(\boldsymbol t)=1\) and, for \(N\ge1\),
\begin{equation}
	\label{eq:PBF-arithmetic-tau-determinants}
	\tau_N^{(a)}(\boldsymbol t)
	\coloneq
	\det\left(
	\widehat{\mathscr M}_{(a)}(\boldsymbol t)^{[N-1]}
	\right).
\end{equation}

These are the integrable tau functions of the time-deformed star moment
matrix.  They should be distinguished from the elementary Christoffel
\(\tau\)-determinants of
Corollary~\ref{cor:elementary-grid-tau-determinants}, which test the
existence of the discrete one-step connections.  Both families arise from
moment determinants, but they play different roles.

\begin{theorem}[Determinant solution selected by the PBF]
	\label{thm:PBF-arithmetic-tau-solution}
	Assume \(\gcd(p,q)=1\).  In every fixed finite set of times, let
	\(\mathcal L_{(a)}(\boldsymbol t)\) be the unique coefficientwise formal
	solution of \eqref{eq:coprime-sparse-PBF-Lax-flows} with initial operator
	\(\mathcal L_{(a)}(\boldsymbol0)\) defined above, whose semi-infinite
	matrix is \(\mathscr H_{(a)}\).  This formal solution is reconstructed
	coefficientwise from \eqref{eq:deformed-star-step-line-moments}.  Its
	nontrivial coefficients are
	\begin{equation}
		\label{eq:PBF-arithmetic-tau-coefficients}
		c_h^{(a)}(\boldsymbol t)
		=
		\frac{
		\tau_{h-p}^{(a)}(\boldsymbol t)
		\tau_{h+1}^{(a)}(\boldsymbol t)
		}{
		\tau_{h-p+1}^{(a)}(\boldsymbol t)
		\tau_h^{(a)}(\boldsymbol t)
		},
		\qquad h\ge p.
	\end{equation}
	For the analytic interpretation, let \(D\) be a domain satisfying the
	hypotheses stated above.  For each \(h\ge p\), let \(D_h\) be the connected
	component containing \(\boldsymbol0\) of
	\[
	\left\{
	\boldsymbol t\in D:
	\tau_N^{(a)}(\boldsymbol t)\ne0,\
	1\le N\le h+1
	\right\}.
	\]
	Then the finite Gauss--Borel construction through row \(h\) exists
	holomorphically on \(D_h\), and
	\eqref{eq:PBF-arithmetic-tau-coefficients} holds there.  In particular,
	every fixed coefficient \(c_h^{(a)}\) is analytic in a neighbourhood of
	\(\boldsymbol0\), which may depend on \(h\).  If all leading principal
	minors \(\tau_N^{(a)}\), \(N\ge1\), remain nonzero throughout \(D\), then
	the full recurrence is coefficientwise holomorphic on \(D\).
\end{theorem}

\begin{proof}
	Theorem~\ref{thm:balanced-scalar-normality} gives
	\(\tau_N^{(a)}(\boldsymbol0)\ne0\) for every \(N\), so the formal
	Gauss--Borel factorization below exists and is unique.
	In the analytic setting, ordinary finite-dimensional Gaussian elimination
	shows that the factors through row \(h\) depend holomorphically on
	\(\boldsymbol t\) throughout \(D_h\), since precisely the required pivots
	remain nonzero.  No convergence of an infinite triangular factorization is
	needed for this coefficientwise assertion.
	Let \(\mathsf\Lambda_{[q]}\) and \(\mathsf\Lambda_{[p]}\) be the two
	step-line block shifts.  Equation
	\eqref{eq:deformed-star-step-line-moments} gives the equivalent moment
	evolutions
	\begin{align*}
	\widehat{\mathscr M}_{(a)}(\boldsymbol t)
	&=
	\exp\!\left(\sum_{\ell\ge1}
	t_\ell\mathsf\Lambda_{[q]}^{M\ell}\right)
	\widehat{\mathscr M}_{(a)}(\boldsymbol0),
	\\
	\widehat{\mathscr M}_{(a)}(\boldsymbol t)
	&=
	\widehat{\mathscr M}_{(a)}(\boldsymbol0)
	\exp\!\left(\sum_{\ell\ge1}
	t_\ell
	\bigl((\mathsf\Lambda_{[p]})^{\top}\bigr)^{M\ell}\right).
	\end{align*}
	Factor the moment matrix as
	\(\widehat{\mathscr M}_{(a)}=\mathscr L^{-1}\mathscr U^{-1}\), with
	\(\mathscr L\) unit lower triangular and \(\mathscr U\) upper
	triangular.  The corresponding recurrence is
	\[
	\mathcal L_{(a)}(\boldsymbol t)
	=
	\mathscr L\mathsf\Lambda_{[q]}\mathscr L^{-1}
	=
	\mathscr U^{-1}(\mathsf\Lambda_{[p]})^{\top}\mathscr U.
	\]
	Differentiating the first moment evolution and separating strictly lower
	and upper triangular parts gives
	\((\partial_{t_\ell}\mathscr L)\mathscr L^{-1}
	=-(\mathcal L_{(a)}^{M\ell})_-\).  Therefore
	\[
	\frac{\partial\mathcal L_{(a)}}{\partial t_\ell}
	=
	[(\mathcal L_{(a)}^{M\ell})_+,\mathcal L_{(a)}].
	\]
	The recurrence produced by the moment factorization belongs to the fixed
	band \([-p,q]\).  In that band the Lax right-hand side is a finite-range
	polynomial vector field, so its formal initial-value problem is unique
	coefficientwise.  Theorem~\ref{thm:coprime-sparse-PBF-hierarchy} shows
	that the two-term submanifold is invariant.  Since the recurrence from
	the moment factorization has initial operator
	\(\mathcal L_{(a)}(\boldsymbol0)\) (equivalently, initial semi-infinite
	matrix \(\mathscr H_{(a)}\)),
	uniqueness forces it to remain on that submanifold.  It is therefore the
	unique formal solution with the prescribed initial value.

	The diagonal pivots of the Gauss--Borel factorization satisfy
	\[
	\mathscr U_{n,n}
	=
	\frac{\tau_n^{(a)}}{\tau_{n+1}^{(a)}}.
	\]
	In the second recurrence representation, the entry in position
	\((h,h-p)\) has only the diagonal contribution
	\[
	c_h^{(a)}
	=
	\frac{\mathscr U_{h-p,h-p}}{\mathscr U_{h,h}}.
	\]
	Substituting the pivot formula proves
	\eqref{eq:PBF-arithmetic-tau-coefficients}.
\end{proof}

The theorem distinguishes the time-dependent solution from the static PBF
formulas.  The PBF coefficients displayed in the Pi\~neiro and Jacobi-like
sections determine, through the selected minimum-height refactorizations, the
initial coefficients \(c_h^{(a)}(\boldsymbol0)\).  Formula
\eqref{eq:PBF-arithmetic-tau-coefficients} supplies their canonical formal
evolution and, on the domains specified in the theorem, their coefficientwise
analytic evolution.

\subsection{Discrete Christoffel transformations and synchronized Toda flows}
\label{subsec:PBF-Baecklund-synchronized-Toda}

The continuous deformation above is compatible with the discrete
transformations generated by the bidiagonal factors.  The compatibility is
most transparent at the level of the polynomial arrays.

\begin{proposition}[A Christoffel transfer commutes with the Toda times]
	\label{prop:Christoffel-transfer-commutes-Toda}
	Let
	\[
	E_{\boldsymbol t}(x)
	\coloneq
	\exp\!\left(\sum_{k\ge1}t_kx^k\right),
	\qquad
	\mathrm d\Psi_j(x;\boldsymbol t)
	=
	E_{\boldsymbol t}(x)\,\mathrm d\Psi_j(x;\boldsymbol0),
	\quad j\in\{0,1\},
	\]
	where \(\mathrm d\Psi_0(x;\boldsymbol0)\) and
	\(\mathrm d\Psi_1(x;\boldsymbol0)\) are two consecutive functionals
	related by a fixed row- or column-side Christoffel multiplication.
	Assume that both deformed systems are weakly normal as formal series.  For
	the analytic statement, assume instead a common domain on which the two
	deformed moment matrices are entrywise holomorphic and all leading minors
	required for their Gauss--Borel factorizations remain nonzero.  Let
	\(T_j(\boldsymbol t)\) and
	\(\boldsymbol B_j(x;\boldsymbol t)\) be their recurrence matrices and
	right polynomial arrays in Sato-compatible Gauss--Borel normalizations.
	These normalizations are chosen by fixed constant gauges at
	\(\boldsymbol t=\boldsymbol0\); their numerical initial blocks are not
	assumed to remain constant in time.  Let \(F(\boldsymbol t)\) be the normalized
	bidiagonal connection matrix and \(\mathsf K(x)\) the fixed matrix
	polynomial in the corresponding right-array relation, after orienting the
	labels \(0,1\) so that
	\[
	F(\boldsymbol t)\boldsymbol B_0(x;\boldsymbol t)
	=
	\boldsymbol B_1(x;\boldsymbol t)\mathsf K(x).
	\]
	For a column-side step, this orientation takes system \(0\) to be the
	transformed system and system \(1\) the preceding one; then \(F=L\) and
	\(\mathsf K=I\).  A row-side step uses the forward orientation displayed
	above.
	For every \(\ell\ge1\),
	\begin{equation}
		\label{eq:Christoffel-Toda-auxiliary-equation}
		\frac{\partial F}{\partial t_\ell}
		=
		(T_1^\ell)_+F-F(T_0^\ell)_+.
	\end{equation}
	Consequently, the Sato-compatible normalized connections around every
	elementary square of the Christoffel rectangle satisfy the same qd
	refactorization throughout that common time domain.
\end{proposition}

\begin{proof}
	Since \(E_{\boldsymbol t}(x)\) is scalar, it commutes with every fixed
	row- or column-side Christoffel multiplier.  Hence the Christoffel
	transform of the deformed functional is the deformation of the
	Christoffel-transformed functional, and the matrix polynomial
	\(\mathsf K(x)\) is independent of \(\boldsymbol t\).  Weak normality and
	the chosen Sato-compatible Gauss--Borel normalizations give the unique
	connection
	\[
	F\boldsymbol B_0=\boldsymbol B_1\mathsf K
	\]
	throughout the common time domain.

	In this normalization the Sato equations are
	\[
	\partial_{t_\ell}\boldsymbol B_j
	=-(T_j^\ell)_-\boldsymbol B_j.
	\]
	Differentiating the connection gives
	\[
	\left(
	\frac{\partial F}{\partial t_\ell}
	-F(T_0^\ell)_-
	+(T_1^\ell)_-F
	\right)\boldsymbol B_0=0.
	\]
	Since \(\boldsymbol B_0\) is a polynomial basis, its coefficient matrix is
	invertible lower triangular and it can be cancelled.  The recurrence
	equations also give \(T_1^\ell F=FT_0^\ell\).  Splitting both sides into
	their upper and strictly lower triangular parts yields
	\[
	\frac{\partial F}{\partial t_\ell}
	=(T_1^\ell)_+F-F(T_0^\ell)_+.
	\]

	Around an elementary square the fixed row and column Christoffel
	multipliers act on opposite sides and hence commute; the common scalar
	deformation commutes with both.  The two routes therefore give the same
	deformed moment functional.  Uniqueness of the chosen Gauss--Borel
	factorizations and of the normalized connections gives the time-dependent
	qd identity.
\end{proof}

Recall that the cyclic products of the refactorized factors belonging to
\(\mathscr H_{(a)}\) are denoted by
\(\widetilde T_{(a;1)},\ldots,\widetilde T_{(a;M)}\).

\begin{theorem}[One sparse Lax operator for the cyclic Toda family]
	\label{thm:sparse-root-synchronized-cyclic-Toda}
	Assume \(\gcd(p,q)=1\), and let
	\(\mathcal L_{(a)}(\boldsymbol t)\) be the PBF-selected solution of
	Theorem~\ref{thm:PBF-arithmetic-tau-solution}.  Then every cyclic block satisfies
	\begin{equation}
		\label{eq:synchronized-cyclic-Toda-flows}
		\frac{\partial\widetilde T_{(a;r)}}{\partial t_\ell}
		=
		\left[
		(\widetilde T_{(a;r)}^{\,\ell})_+,
		\widetilde T_{(a;r)}
		\right],
		\qquad
		r\in\{1,\ldots,M\}.
	\end{equation}
	The factors connecting consecutive blocks satisfy the auxiliary equations
	\begin{equation}
		\label{eq:factor-copy-Toda-auxiliary-equations}
		\frac{\partial \widetilde F_r^{(a)}}{\partial t_\ell}
		=
		(\widetilde T_{(a;r)}^{\,\ell})_+\widetilde F_r^{(a)}
		-
		\widetilde F_r^{(a)}
		(\widetilde T_{(a;r+1)}^{\,\ell})_+,
	\end{equation}
	where the indices are read cyclically and the orientation is fixed by
	\(\widetilde T_{(a;r)}\widetilde F_r^{(a)}
	=\widetilde F_r^{(a)}\widetilde T_{(a;r+1)}\).
	Thus the semi-infinite matrix \(\mathscr H_{(a)}(\boldsymbol t)\) of
	\(\mathcal L_{(a)}(\boldsymbol t)\) is a common integrable \(M\)-th root of the cyclic
	Toda family.  If, in addition, the factor-copy polynomial arrays are taken
	in the Sato-compatible Gauss--Borel normalizations of
	Proposition~\ref{prop:Christoffel-transfer-commutes-Toda}, so that the
	off-diagonal blocks \(\widetilde F_r^{(a)}\) are the corresponding
	Christoffel connections, then the factor transfers are discrete B\"acklund
	directions and \eqref{eq:factor-copy-Toda-auxiliary-equations} is precisely
	\eqref{eq:Christoffel-Toda-auxiliary-equation}.  Under other time-dependent
	initial-condition gauges, the compatibility is gauge-covariant and contains
	the corresponding gauge terms.
\end{theorem}

\begin{proof}
	Let \(\mathscr P_a\) denote the fixed permutation from factor-copy order
	to height order, so that
	\[
	\mathscr H_{(a)}
	=
	\mathscr P_a\mathscr C_{(a)}\mathscr P_a^{-1}.
	\]
	In factor-copy order,
	\[
	\mathscr P_a^{-1}\mathscr H_{(a)}^M\mathscr P_a
	=
	\mathscr C_{(a)}^M
	=
	\operatorname{diag}
	\bigl(
	\widetilde T_{(a;1)},\ldots,\widetilde T_{(a;M)}
	\bigr).
	\]
	The coordinates in one factor copy occur at heights separated by \(M\).
	Hence nonnegative height displacement is equivalent to an upper triangular
	entry inside that copy, and
	\[
	\mathscr P_a^{-1}
	(\mathscr H_{(a)}^{M\ell})_+
	\mathscr P_a
	=
	\operatorname{diag}
	\bigl(
	(\widetilde T_{(a;1)}^{\,\ell})_+,\ldots,
	(\widetilde T_{(a;M)}^{\,\ell})_+
	\bigr).
	\]
	Conjugating the Lax equation for \(\mathscr H_{(a)}\) by
	\(\mathscr P_a^{-1}\) gives the corresponding equation for
	\(\mathscr C_{(a)}\).  Its cyclic off-diagonal blocks are the factors
	\(\widetilde F_r^{(a)}\) and give
	\eqref{eq:factor-copy-Toda-auxiliary-equations}.  Taking the \(M\)-th
	power and then the diagonal factor-copy blocks gives
	\eqref{eq:synchronized-cyclic-Toda-flows}.  Under the additional
	normalization hypothesis in the statement,
	Proposition~\ref{prop:Christoffel-transfer-commutes-Toda} identifies these
	auxiliary equations with the continuous--discrete Christoffel
	compatibility; the statement for other gauges is the usual change-of-gauge
	form of the same identity.
\end{proof}

The band of each cyclic product extends from the \(p\)-th subdiagonal to the
\(q\)-th superdiagonal, so its equations belong to the
\((q,p)\)-bigraded Toda setting \cite{Carlet2006BigradedToda}.  The new
feature here is the factor-resolved sparse root: one scalar sequence
\(c_h^{(a)}(\boldsymbol t)\) evolves the entire cyclic family and retains the
individual bidiagonal connections between its members.

\subsection{Hypergeometric determinant solutions}
\label{subsec:PBF-hypergeometric-integrable-solutions}

The determinant solution becomes explicit once the cyclic spectral data are
known.  Fix two factor copies \(r,s\), let
\(\widetilde F_{r\to s}\) be their consecutive partial product, and let
\((\boldsymbol B^{(s)},\boldsymbol A^{(s)},\mathrm d\Psi^{(s)})\) be
biorthogonal data for \(\widetilde T_{(a;s)}\).  The radial moments of
Section~\ref{sec:balanced-star-Christoffel} have the spectral representation
\begin{equation}
	\label{eq:radial-moments-cyclic-spectral-data}
	m_k^{r,s}
	=
	\int
	\bigl(e_0^{\top}\widetilde F_{r\to s}
	\boldsymbol B^{(s)}(x)\bigr)
	x^k\,\mathrm d\Psi^{(s)}(x)
	\bigl(\boldsymbol A^{(s)}(x)e_0\bigr).
\end{equation}

Denote the scalar measure occurring on the right-hand side by
\[
\mathrm d\sigma^{r,s}(x)
\coloneq
\bigl(e_0^{\top}\widetilde F_{r\to s}
\boldsymbol B^{(s)}(x)\bigr)\,
\mathrm d\Psi^{(s)}(x)\,
\bigl(\boldsymbol A^{(s)}(x)e_0\bigr).
\]
It is a finite signed or complex measure whenever the displayed spectral
integral is defined; it is not positive in general.

\begin{proposition}[Radial deformation in cyclic spectral coordinates]
	\label{prop:radial-deformation-cyclic-spectral-coordinates}
	Equation \eqref{eq:radial-moments-cyclic-spectral-data} says that
	\(\mathrm d\sigma^{r,s}\) and the positive radial measure
	\(\mathrm d\rho^{r,s}\) of
	Theorem~\ref{thm:balanced-radial-Hankel-positivity} represent the same
	moment functional.  If \(\mathrm d\sigma^{r,s}\) is positive on
	\([0,\infty)\) and the corresponding Stieltjes moment problem is
	determinate, then
	\(\mathrm d\sigma^{r,s}=\mathrm d\rho^{r,s}\).  The same equality follows
	without positivity whenever both measures are compactly supported on the
	real line.  Under the first radial time \(t=t_1\), the formal deformed
	moments are
	\begin{equation}
		\label{eq:deformed-radial-moments-cyclic-data}
		m_k^{r,s}(t)
		=
		\int
		\bigl(e_0^{\top}\widetilde F_{r\to s}
		\boldsymbol B^{(s)}(x)\bigr)
		x^k \e^{tx}\,\mathrm d\Psi^{(s)}(x)
		\bigl(\boldsymbol A^{(s)}(x)e_0\bigr).
	\end{equation}
	Thus every entry of the deformed star moment matrix is a finite
	polynomial dressing of the deformed moments of a cyclic Favard system.
	In the compactly supported case the identity is entire in
	\(t\in\mathbb C\).
\end{proposition}

\begin{proof}
	Biorthogonality and either recurrence equation give
	\[
	(\widetilde T_{(a;s)})^k
	=
	\int
	\boldsymbol B^{(s)}(x)x^k\,
	\mathrm d\Psi^{(s)}(x)\boldsymbol A^{(s)}(x).
	\]
	Multiplication on the left by
	\(e_0^{\top}\widetilde F_{r\to s}\) and on the right by \(e_0\)
	gives \eqref{eq:radial-moments-cyclic-spectral-data}.  Its moments are
	\((\widetilde F_{r\to s}(\widetilde T_{(a;s)})^k)_{0,0}\), so it is the
	moment functional of the radial measure constructed in
	Theorem~\ref{thm:balanced-radial-Hankel-positivity}.  This proves equality
	of the two moment functionals.  Stieltjes determinacy gives equality of the
	measures when the dressed measure is positive; compact support on the real
	line gives equality for finite signed or complex measures by polynomial
	density.  Multiplication by \(\e^{tx}\), interpreted coefficientwise in
	\(t\), gives
	\eqref{eq:deformed-radial-moments-cyclic-data}; compact support makes this
	an entire identity.
\end{proof}

For the Pi\~neiro system, the first-time deformation of the scalar step-line
moment matrix is completely explicit.  With
\(r=qu+j\), \(s=pv+i\), put
\(x_r=u+\beta_{j+1}\) and \(y_s=v+\alpha_{i+1}+1\).  Under
\eqref{eq:Pineiro-positive-chamber-CF-paper},
\[
x_r+y_s\ge\alpha_1+\beta_1+1>0.
\]
Thus the following integral converges for every \(r,s\) and defines an
entire function of \(t\):
\begin{equation}
	\label{eq:deformed-Pineiro-moment-Kummer}
	\mathscr M_{r,s}^{\mathrm P}(t)
	=
	\int_0^1x^{x_r+y_s-1}\e^{tx}\,\mathrm dx
	=
	\frac{1}{x_r+y_s}
	\pFq{1}{1}{x_r+y_s}{x_r+y_s+1}{t}.
\end{equation}

\begin{theorem}[Global positive coefficientwise solutions for the Pi\~neiro system]
	\label{thm:global-positive-Pineiro-integrable-solutions}
	Under \eqref{eq:Pineiro-positive-chamber-CF-paper}, the deformed Pi\~neiro
	moment matrix in \eqref{eq:deformed-Pineiro-moment-Kummer} is totally
	positive for every real \(t\).  All moment factorizations and normalized
	Christoffel connections in its finite rectangle therefore exist and remain
	positive for every real \(t\); these assertions hold for arbitrary
	\(p,q\).  If, in addition, \(\gcd(p,q)=1\), the selected star
	refactorizations give a global positive coefficientwise solution of the
	first scalar flow
	in \eqref{eq:coprime-sparse-PBF-Lax-flows}, and their cyclic products
	satisfy \eqref{eq:synchronized-cyclic-Toda-flows} with \(\ell=1\).
	More generally, replacing \(\e^{tx}\) by
	\(\exp(\sum_{\ell=1}^{L}t_\ell x^\ell)\), with finitely many real times,
	gives, in the coprime case, global positive coefficientwise solutions of
	the corresponding commuting scalar flows.

	When \(\gcd(p,q)=1\), every coefficient \(c_h^{(a)}(t)\) of the
	first-time deformation can be obtained by finitely many
	determinants and rational operations from the Kummer entries in
	\eqref{eq:deformed-Pineiro-moment-Kummer}.  The initial coefficients
	\(c_h^{(a)}(0)\) are those of the selected minimum-height refactorizations
	of the Pi\~neiro PBF.  For the rational system in
	Example~\ref{ex:Pineiro-CF-numerical}, they are rational.
\end{theorem}

\begin{proof}
	For finite increasing row and column sets
	\(I=\{r_1<\cdots<r_N\}\) and \(J=\{s_1<\cdots<s_N\}\), expand the two
	determinants below.  Although an individual function
	\(z^{y_{s_v}-1}\) need not be integrable at the origin, every monomial in
	the finite expansion has, in each variable, an exponent of the form
	\(x_{r_u}+y_{s_v}-1>-1\).  The expansion is therefore absolutely
	integrable, and Andr\'eief's identity gives
	\begin{equation*}
	\det\left[\mathscr M_{r_u,s_v}^{\mathrm P}(t)\right]_{u,v=1}^{N}
=
	\frac1{N!}\int_{(0,1)^N}
	\det\left[z_k^{x_{r_u}}\right]_{u,k=1}^{N}
	\det\left[z_k^{y_{s_v}-1}\right]_{v,k=1}^{N}
	\prod_{k=1}^{N}\e^{tz_k}\,\mathrm dz_k.
	\end{equation*}
	The inequalities in
	\eqref{eq:Pineiro-positive-chamber-CF-paper} imply that the exponent
	sequences are strictly increasing.  Indeed, successive differences in
	\((x_r)\) are either \(\beta_{j+2}-\beta_{j+1}>0\) or
	\(1+\beta_1-\beta_q>0\), and the corresponding differences in \((y_s)\)
	are either \(\alpha_{i+2}-\alpha_{i+1}>0\) or
	\(1+\alpha_1-\alpha_p>0\).
	On \(0<z_1<\cdots<z_N<1\), the two generalized Vandermonde
	determinants have the same strictly positive sign because both exponent
	sequences are strictly increasing.  The remaining factor is positive for
	every real \(t\), so the minor is strictly positive.  The same proof applies
	with \(\e^{tz_k}\) replaced by
	\(\exp(\sum_{\ell=1}^{L}t_\ell z_k^\ell)\), which proves the multi-time
	assertion.

	Total positivity gives the Gauss--Borel factorization at every vertex of
	the finite rectangle for every real \(t\).  Choose there the
	Sato-compatible Gauss--Borel normalizations of
	Proposition~\ref{prop:Christoffel-transfer-commutes-Toda}; the resulting
	Christoffel connection matrices have entries that are ratios of minors of
	\eqref{eq:deformed-Pineiro-moment-Kummer}.  Starting from the chosen cyclic
	product, the finite directed sequence of adjacent qd refactorizations then
	constructs the factors of \(\mathscr H_{(a)}(t)\); each of their entries is
	a rational expression in those Christoffel coefficients and remains
	strictly positive.

	The connections satisfy
	Proposition~\ref{prop:Christoffel-transfer-commutes-Toda}, so their cyclic
	matrix satisfies the factor-copy Lax equation.  If \(\gcd(p,q)=1\),
	conjugating by the fixed height permutation gives
	\eqref{eq:coprime-sparse-PBF-Lax-flows} for \(\mathscr H_{(a)}(t)\).
	In that coprime case, at \(t=0\) this is the prescribed sparse matrix, so
	uniqueness of the
	formal initial-value problem identifies it with the determinant solution of
	Theorem~\ref{thm:PBF-arithmetic-tau-solution}.  Since the construction
	exists and stays positive for every real \(t\), the coefficientwise solution
	is global on the real first-time axis.  The same finite sequence of
	determinants and
	rational operations proves the asserted Kummer description.
\end{proof}

In the coprime case, these solutions are not finite-soliton solutions in
general.  The Pi\~neiro
measure is absolutely continuous on \((0,1)\).  On the first-time axis its
\(\tau\)-functions are determinants of Kummer moments of this continuous
spectral measure; for several times the entries are the corresponding
exponential-polynomial moments.  When only finitely many times are nonzero,
these tau functions are entire functions of those complex times.  Total
positivity excludes their zeros for all real values of those times, so every
coefficient \(c_h^{(a)}(\boldsymbol t)\) is positive and real analytic there.
For complex times, each fixed coefficient is meromorphic, with possible poles
only at zeros of the two denominator tau functions in
\eqref{eq:PBF-arithmetic-tau-coefficients}.  Thus the Pi\~neiro construction
supplies global positive coefficientwise solutions, represented by
continuous-spectrum determinants and selected explicitly by the original
PBF.

The Jacobi-like system has an analogous explicit deformation.  If
\(s=n+\alpha_i+1\), then
\begin{equation}
	\label{eq:deformed-Jacobi-like-moment-hypergeometric}
	\int_0^1x^n \e^{tx}\,\mathrm d\mu_{j,i}^{\mathrm J}(x)
	=
	\left(
	\int_0^1x^n\,\mathrm d\mu_{j,i}^{\mathrm J}(x)
	\right)
	\pFq{q}{q}
	{s\one_q+\boldsymbol a}
	{s\one_q+\boldsymbol b+\boldsymbol e_j^{(q)}}
	{t}.
\end{equation}
Indeed, expanding \(\e^{tx}\) and using
\eqref{eq:Jacobi-like-CF-moments} gives, for \(k\ge0\),
\[
\frac{\displaystyle\int_0^1x^{n+k}\,
\mathrm d\mu_{j,i}^{\mathrm J}(x)}
{\displaystyle\int_0^1x^n\,
\mathrm d\mu_{j,i}^{\mathrm J}(x)}
=
\frac{(s+b_j)_k}{(s+b_j+1)_k}
\frac{(s\one_q+\boldsymbol a)_k}
{(s\one_q+\boldsymbol b)_k}.
\]
The factor \((s+b_j)_k\) cancels the \(j\)-th component of the second
denominator.  The remaining lower vector is
\(s\one_q+\boldsymbol b+\boldsymbol e_j^{(q)}\), which gives exactly the
\({}_qF_q\) in
\eqref{eq:deformed-Jacobi-like-moment-hypergeometric}.

At every vertex for which the two weight systems at \(t=0\) satisfy the AT
hypotheses, their deformations remain AT for every real \(t\).  Indeed,
multiplying all weights in a linear form by the common factor
\(\e^{tx}>0\) does not change its zeros or their multiplicities.  Assume now
that \(q\ge2\), that the parameters lie in the Jacobi-like
positive-factorization chamber
\eqref{eq:Jacobi-like-CF-positive-chamber} and that, at every vertex required
by the selected Christoffel rectangle, both transformed weight systems
satisfy the AT hypotheses.  When \(\gcd(p,q)=1\), the same argument gives
global positive coefficientwise determinant solutions, built from the
\({}_qF_q\) entries in
\eqref{eq:deformed-Jacobi-like-moment-hypergeometric} and their corresponding
row and column shifts.  Indeed, the AT hypotheses prevent the required
determinants from vanishing; continuity then preserves the positive signs of
the PBF coefficients from \(t=0\) along the real axis.

In the larger Jacobi-like PBF chamber with \(\gcd(p,q)=1\),
Theorem~\ref{thm:PBF-arithmetic-tau-solution} still gives the formal solution.
Since every deformed moment entry is entire, for each fixed \(h\) it also
gives a coefficient \(c_h^{(a)}\) analytic in a neighbourhood of
\(\boldsymbol0\) on which the finitely many minors through order \(h+1\) do
not vanish.  This neighbourhood may depend on \(h\); neither a common
neighbourhood for all coefficients nor global nonvanishing of all deformed
moment minors is asserted.

The radial deformation \eqref{eq:star-radial-Toda-deformation} is a
common-time reduction of the multi-component 2D Toda deformation associated
with mixed-type moment matrices \cite{AlvarezFidalgoManas2011}.  Its role here
is to evolve the PBF coefficients while synchronizing the complete cyclic
family through Theorem~\ref{thm:sparse-root-synchronized-cyclic-Toda}.  For
\(d=\gcd(p,q)>1\), the scalar two-term operator still carries the
factor-compatible support-preserving commuting subhierarchy above and
decomposes into \(d\) invariant
residue classes.  What fails is its identification with
\(\mathscr H_{(a)}\), which has \(d\) coordinate states at each height and
block displacements \(+Q\) and \(-P\).  A PBF-adapted matrix
tau-function description of this block problem, and its relation to the
independent times of the full multi-component \(2D\) Toda hierarchy, is left
open.

\section{Conclusions}
\label{sec:conclusions}

The article shows that a positive bidiagonal factorization is an organizing
datum in its own right, finer than the banded matrix obtained by multiplying
the factors.  Its principal conclusions fall into the same three groups as
the results stated in the Introduction.

\begin{enumerate}[label=\textnormal{(\arabic*)},leftmargin=*]
	\item Every cyclic permutation of the prescribed factors again admits a
	normalized PBF in lower--upper order, for bounded and unbounded factors
	alike.  With suitable one-sided maximal-degree data at consecutive cyclic
	products, each cyclic transfer determines a degree-one matrix Christoffel
	polynomial with determinant \(cx\),
	\(c\ne0\), and the transformed functional has the next cyclic product as
	recurrence matrix.  Different nonsingular lower triangular initial
	conditions produce constant lower/upper triangular changes of the
	functional and preserve weak normality; maximal component degrees require
	separate admissibility.  A complete factor circuit gives multiplication by
	\(x\), up to the terminal normalization.
	For the fixed multipliers \(\mathsf K_q(x)\) and
	\(\mathsf K_p(x)^{\top}\), weak normality of the two boundary chains gives
	the finite rectangle of simultaneous moment shifts.  Total positivity of
	the scalar moment matrix makes this additional boundary hypothesis
	automatic.

	\item Exactly the cyclic products
	\(T_{(a)}\), \(a\in\mathcal A_{\mathrm{cyc}}=\{a_-,\ldots,a_+\}\),
	can reach the minimum-height factor pattern through the directed positive
	adjacent exchanges used here.  Their associated sparse matrices
	\(\mathscr H_{(a)}\) have only the block displacements \(+Q\) and \(-P\),
	with \(d=\gcd(p,q)\) coordinates at each height; in the coprime case only
	the \(q\)-th superdiagonal and the \(p\)-th subdiagonal remain.  Their
	projected resolvents recover the resolvents of the cyclic products, their
	finite nonzero spectra lie on a \((p+q)\)-ray star, and their explicit block
	moment factorizations give positive radial Stieltjes moment sequences.
	For fixed compactly supported radial representing measures, the Fourier
	components admit finite complex measures supported on the corresponding
	compact star precisely when
	the corresponding atom-at-zero and negative fractional-moment conditions
	of Proposition~\ref{prop:balanced-character-measure-lift} hold entrywise.
	In the coprime
	case, each \(\mathscr H_{(a)}\) also carries a local arithmetic-support
	Lax hierarchy.  A radial exponential deformation of the star moments gives
	the evolving coefficients by an explicit quotient of four leading minors.
	The \((p+q)\)-th power of the sparse Lax matrix synchronizes the Toda flows
	of the cyclic products of
	the refactorized factors, and the Christoffel transfers commute with the
	continuous times.  On the real first-time axis, the Pi\~neiro chamber
	produces global positive coefficientwise solutions expressed by Kummer
	determinants; the intersection of the Jacobi-like positive-factorization
	chamber with the requisite AT region produces analogous global positive
	coefficientwise
	solutions expressed by \({}_qF_q\) determinants.

	\item The sparse matrix built from the prescribed factor order retains
	every original bidiagonal coefficient.  Grouping consecutive heights makes
	it block tridiagonal and gives a factor-resolved matrix continued fraction.
	Its finite height-block truncations have explicit Pad\'e-type contact,
	denominator factorization, backward evaluation, and residual error bounds.
	For bounded nonnegative factors these convergents are entrywise monotone,
	their denominators are nonsingular \(M\)-matrices, and their limit is the
	least nonnegative solution.  Stochastic normalization gives convergence
	throughout the unit disk.  A second grouping produces the \(q\)-state
	recursion, which for \(q=1\) is the classical \(p\)-branched Stieltjes
	recursion.  A rectangular projection of the bounded resolvent is the Weyl
	matrix of the mixed-type system, so all its Stieltjes transforms are
	approximated simultaneously.  For unbounded factorizations the
	coefficientwise approximation and finite algorithm remain valid.  A
	specified closed realization satisfying the core and uniform finite-section
	stability hypotheses gives convergence to its projected resolvent; this
	limit is a measure-defined Weyl matrix only under the additional spectral
	identification.
\end{enumerate}

Together these results extend the classical passage
\[
\text{Jacobi matrix}
\longrightarrow
\text{Stieltjes fraction}
\longrightarrow
\text{scalar Weyl function}
\]
to a two-sided mixed-type setting.  The role of the Jacobi coefficients is
played by the individual positive bidiagonal factors.  Their prescribed
order produces rational approximants to the rectangular Weyl matrix in the
bounded or spectrally identified setting, and formal moment approximants in
general; their minimum-height factorizations produce the two-diagonal
recurrences and their moment systems.
The Pi\~neiro example displays the full construction and its certified
numerical error.  In the Jacobi-like specialization, cyclic transforms may
leave the original special-function family even though the PBF and its
spectral consequences persist.

The remaining questions include
uniform finite-section stability outside the domains controlled by the norm,
	spectral radius, and numerical range; ratio asymptotics and convergence
	rates of the matrix convergents; comparison with finite spectral
	approximants obtained from principal truncations in the unbounded case; and
	an extension of the star-lifting criterion beyond the compactly supported,
	determinate radial case, including intrinsic conditions that select
	noncompact representatives with no atom at the origin in nonzero modes and
	the required negative fractional moments.
	On the integrable side, the principal open
problem is not the existence of scalar support-preserving flows for
noncoprime \(p,q\), but their PBF-relevant block realization when
\(d=\gcd(p,q)>1\), including a matrix determinant or quasideterminant
solution and its relation to the independent times of the full
multi-component \(2D\) Toda hierarchy.

The occurrence of Christoffel words suggests a further combinatorial and
spectral direction.  If coprime pairs \((P_n,Q_n)\) are chosen along the
continued-fraction convergents to an irrational slope, the corresponding
Christoffel patterns have a natural Sturmian limit.  It remains to determine whether the
associated sparse recurrences and matrix continued fractions admit a
corresponding aperiodic limit, and whether the classical connection between
Christoffel words, ordinary continued fractions, and Markoff theory has a
spectral counterpart for positive bidiagonal factorizations
\cite{BerstelLauveReutenauerSaliola2009}.

\section*{Acknowledgements}

The author thanks Am\'\i lcar Branquinho and Ana Foulqui\'e-Moreno for valuable
discussions on mixed-type orthogonality.

\section*{Declarations}

\noindent\textbf{Funding.}
The author was supported by the Spanish State Research Agency
(Agencia Estatal de Investigaci\'on) through research projects
PID2021-122154NB-I00,
\emph{Ortogonalidad y aproximaci\'on con aplicaciones en machine learning y
teor\'\i a de la probabilidad}, and PID2024-155133NB-I00,
\emph{Ortogonalidad, aproximaci\'on e integrabilidad: aplicaciones en procesos
estoc\'asticos cl\'asicos y cu\'anticos}.

\medskip
\noindent\textbf{Competing interests.}
The author declares no competing interests.

\medskip
\noindent\textbf{Data availability.}
No external data set is used.  The Python and Sage scripts and their generated
output are supplied as Supplementary Material with the manuscript.  The
numerical script implements the cited Pi\~neiro Gamma--Pochhammer formulas and
reproduces Table~\ref{tab:Pineiro-32-convergence}, including the 12-digit
reference resolvent block, the finite-section dimensions, and the a priori and
residual bounds.  The exact-arithmetic regression suite checks the cyclic
powers and Darboux intertwinings, the height multiplicities and block grouping,
the finite Schur recursion and determinant factorization, the first omitted
path, the \(q\)-state recursion, the rectangular moment projection,
representative Jacobi-like Weyl coefficients, finite monotonicity, and the
exact residual identity.  Two additional exact-arithmetic Maple scripts
reconstruct the determinant solution of the sparse Lax hierarchy.  One checks
the coefficient formula for several coprime pairs \((p,q)\); the other verifies
the first nonlinear flow coefficientwise for the rational Pi\~neiro
\((3,2)\) example.  A further exact Sage check verifies, for the same set of
coprime pairs, that restriction of the bilateral Lax identity to the
semi-infinite zero boundary introduces no correction.

\printbibliography

\end{document}